\documentclass[11pt, french, a4paper]{smfbook}

\usepackage[T1]{fontenc}
\usepackage[utf8]{inputenc}
\usepackage[french,english]{babel}
\usepackage{textgreek}
\usepackage{lmodern}
\usepackage{ucs} %

\let\printindex\undef

\usepackage{splitidx}
\makeindex

\usepackage{latexsym}
\usepackage{amsmath}
\usepackage{amsthm}
\usepackage{amssymb}
\usepackage{amsfonts}
\usepackage{esint}
\usepackage{graphicx}
\usepackage{accents}
\usepackage[all]{xy}

\let\ORIGxymatrix\xymatrix
\def\xymatrix#1#{\begingroup\shorthandoff{;}\INNERxymatrix{#1}}
\def\INNERxymatrix#1#2{\ORIGxymatrix#1{#2}\endgroup}

\usepackage{xpatch}
\usepackage[hidelinks,hyperindex,pdfencoding=unicode]{hyperref}
\xyoption{2cell}
\UseAllTwocells
\RequirePackage{tikz}						%
\usetikzlibrary{decorations.markings}

\usepackage{float}
\usepackage{todonotes}
\usepackage{mathtools}
\usepackage{faktor}
\usepackage{chngcntr}
\usepackage{enumitem}
\usepackage{manfnt}
\usepackage{calc}

\usepackage{ifmtarg}
\usepackage{bold-extra}
\usepackage[pass]{geometry}

\let\thmnewline\relax

\makeatletter
\patchcmd{\smf@captionsfrench}{Bibliographie}{Références}{}{}
\makeatother

\newcommand{\dhorline}[3][0]{%
    \tikz[baseline=-2pt]{\path[decoration={markings, 
      mark=between positions 0 and 1 step 2*#3
      with {\node[color=blueamu, fill, circle, minimum width=#3, inner sep=0pt, anchor=south west] {};}},postaction={decorate}]  (0,#1) -- ++(#2,0);}}
\newcommand{\dvertline}[3][0]{%
    \tikz[baseline=2em]{\path[decoration={markings,
      mark=between positions 0 and 1 step 2*#2
      with {\node[color=black!50, fill, circle, minimum width=#2, inner sep=0pt, anchor=south west] {};}},postaction={decorate}] (0, #1) -- ++(0,#3);}}

\makeatletter\newcommand\HUGE{\@setfontsize\Huge{28}{0}}\makeatother

\theoremstyle{plain}
\newtheorem{thm}{Th\'eor\`eme}[section]
\newtheorem*{thm*}{Th\'eor\`eme}
\newtheorem{prop}[thm]{Proposition}
\newtheorem*{prop*}{Proposition}
\newtheorem{lemme}[thm]{Lemme}
\newtheorem{coro}[thm]{Corollaire}
\theoremstyle{remark}

\theoremstyle{definition}
\newtheorem{defin}[thm]{Définition}
\newtheorem{paragr}[thm]{}
\newtheorem{paragrA}{}[chapter]
\newtheorem{propA}[paragrA]{Proposition}
\newtheorem{exampleA}[paragrA]{Exemple}
\newtheorem{definA}[paragrA]{Définition}
\newtheorem{lemmeA}[paragrA]{Lemme}
\newtheorem{remarkA}[paragrA]{Remarque}
\newtheorem{coroA}[paragrA]{Corollaire}
\newtheorem*{defin*}{Définition}

\theoremstyle{plain}

\numberwithin{equation}{thm}

\newcommand{\mypar}[1]{%
  \bigbreak
  \noindent\textbf{#1. ---}%
}
\makeatletter
\newif\ifsection
\preto\section{\sectiontrue}
\preto\subsection{\sectionfalse}
\xapptocmd\@sect{%
  \ifsection
    \numberwithin{thm}{section}
  \else
    \numberwithin{thm}{subsection}
  \fi
  \setcounter{thm}{0}\relax}
  {}{}
\makeatother

\newcommand\letenv[2]{%
\expandafter\expandafter\expandafter\let\expandafter\expandafter
\csname #1\endcsname\csname #2\endcsname
\expandafter\expandafter\expandafter\let\expandafter\expandafter
\csname end#1\endcsname\csname end#2\endcsname
}

\letenv{paragraph}{paragr}
\letenv{proposition}{prop}
\letenv{corollary}{coro}
\letenv{theorem}{thm}
\letenv{theorem*}{thm*}
\letenv{remark}{rem}
\letenv{example}{exem}

\newindex[Index des notations]{not}
\newindex[Index terminologique]{term}
\AtWriteToIndex{not}{\let\thepage\thenoentry}
\AtWriteToIndex{term}{\let\thepage\thenoentry}

\makeatletter

\newcommand\termindex[1]{%
  \Hy@raisedlink{\hypertarget{\thenoentry}{}}%
  \label{\thenoentry}%
  \expandafter\xdef\csname idx\thenoentry\endcsname{\thethm}%
  \sindex[term]{#1|linkindex}%
  \stepcounter{noentry}%
}

\newcounter{noentry}

\newcommand\notindex[1]{%
  \Hy@raisedlink{\hypertarget{\thenoentry}{}}%
  \label{\thenoentry}%
  \expandafter\xdef\csname idx\thenoentry\endcsname{\thethm}%
  \sindex[not]{\thenoentry @\detokenize{#1}|linkindex}%
  \stepcounter{noentry}%
}

\newcommand\linkindex[1]{%
  \ifnum\pdf@strcmp{\@nameuse{idx#1}}{0.0}=0%
    \hyperlink{#1}{\hbox{p.~\pageref*{#1}, \S~Notations et terminologie}}%
  \else
    \hyperlink{#1}{p.~\pageref*{#1},~\S~\@nameuse{idx#1}}%
  \fi
}
\newcommand\ndef[2][\@empty]{%
  \@ifmtarg{#1}{}{%
  \Hy@raisedlink{\hypertarget{\thethm}{}}%
  \ifx\@empty#1%
    \termindex{#2}%
  \else%
    \termindex{#1}%
  \fi%
  }%
  \emph{#2}%
}
\newcommand\nnot[2][\@empty]{%
  \@ifmtarg{#1}{}{%
  \Hy@raisedlink{\hypertarget{\thethm}{}}%
  \ifx\@empty#1%
    \notindex{#2}%
  \else%
    \notindex{#1}%
  \fi%
  }%
  #2%
}

\DeclareFontFamily{OT1}{pzc}{}
\DeclareFontShape{OT1}{pzc}{m}{it}%
             {<-> s * [1,150] pzcmi7t}{}
\DeclareMathAlphabet{\mathpzc}{OT1}{pzc}%
                                 {m}{it}

\renewcommand\leq\leqslant
\renewcommand\geq\geqslant

\newcommand\Z{\mathbb{Z}}
\newcommand\N{\mathbb{N}}
\newcommand\G{\mathbb{G}}
\newcommand\Grefl{\mathbb{G}_{\text{ref}}}
\newcommand{\Ens}{{\mathcal{E}\mspace{-2.mu}\it{ns}}}
\newcommand{\Top}{{\mathcal{T}\mspace{-2.mu}\it{op}}}
\newcommand{\Cat}{{\mathcal{C}\mspace{-2.mu}\it{at}}}
\newcommand{\Gpd}{{\mathcal{G}\mspace{-2.mu}\it{pd}}}
\newcommand{\CAT}{{\mathcal{C}\mspace{-2.mu}\it{AT}}}

\newcommand{\Ch}{{\mathcal{C}\mspace{-2.mu}\it{h}_{+}}}
\newcommand{\Chm}{{\mathcal{C}\mspace{-2.mu}\it{h}^{+}}}
\newcommand{\Chu}{{\mathcal{C}\mspace{-2.mu}\it{h}}}

\newcommand{\C}{\mathcal{C}}
\newcommand\D{\mathcal{D}}
\let\limind\varinjlim
\let\limproj\varprojlim
\newcommand{\Ob}{\operatorname{\mathsf{Ob}}}
\newcommand{\Fl}{\operatorname{\mathsf{Fl}}}
\newcommand{\Hom}{\operatorname{\mathsf{Hom}}}
\newcommand{\Homi}{\operatorname{\kern.5truept\underline{\kern-.5truept\mathsf{Hom}\kern-.5truept}\kern1truept}}
\newcommand{\id}{\mathrm{id}}

\newcommand{\adjpair}[4]{#1:#2 \rightleftarrows #3:#4}

\DeclareMathOperator{\op}{\mathsf{op}}
\newcommand{\pref}[1]{{\widehat{ #1 }}}
\newcommand{\tranche}[2]{#1/{#2}}
\newcommand{\cotranche}[2]{#1\backslash{#2}}
\DeclareMathOperator{\pr}{pr}

\DeclareMathOperator{\ev}{ev}

\newcommand{\W}{\mathcal{W}}
\newcommand{\M}{\mathcal{M}}
\renewcommand{\D}{\mathbb{D}}

\makeatletter
\newcommand{\hocolim@}[2]{%
  \vtop{\m@th\ialign{##\cr
    \hfil$#1\operator@font holim$\hfil\cr
    \noalign{\nointerlineskip\kern1.5\ex@}#2\cr
    \noalign{\nointerlineskip\kern-\ex@}\cr}}%
}
\newcommand{\hocolim}{%
  \mathop{\mathpalette\hocolim@{\rightarrowfill@\textstyle}}\nmlimits@\nolimits
}

\newcommand{\holim@}[2]{%
  \vtop{\m@th\ialign{##\cr
    \hfil$#1\operator@font holim$\hfil\cr
    \noalign{\nointerlineskip\kern1.5\ex@}#2\cr
    \noalign{\nointerlineskip\kern-\ex@}\cr}}%
}
\newcommand{\holim}{%
  \mathop{\mathpalette\holim@{\leftarrowfill@\textstyle}}\nmlimits@
}
\makeatother

\def\TO#1{\mathrel{\hbox to #1pt{\rightarrowfill}}}
\def\OT#1{\mathrel{\hbox to #1pt{\leftarrowfill}}}

\newcommand{\Hotab}{\operatorname{\mathsf{Hotab}}}

\newcommand{\Hot}{\operatorname{\mathsf{Hot}}}

\newcommand{\Der}{\operatorname{\mathbb{D}}}
\newcommand{\DerHot}{\operatorname{\mathpzc{Hot}}}
\newcommand{\DerHotab}{\operatorname{\mathpzc{Hotab}}}
\newcommand{\orient}{\mathcal{O}}

\newcommand{\const}{\operatorname{\mathsf{const}}}

\newcommand{\Ab}{{\mathsf{Ab}}}
\newcommand{\A}{{\mathcal{A}}}
\newcommand{\ab}{{\mathsf{ab}}}

\newcommand{\Tot}{\operatorname{\mathsf{Tot}}}
\newcommand{\diag}{\operatorname{\mathsf{diag}}}
\DeclareMathOperator{\Add}{\mathrm{Add}}
\newcommand{\Addinf}{\operatorname{\mathrm{Addinf}}}
\newcommand{\Homadd}{\operatorname{\kern.5truept\underline{\kern-.5truept\mathsf{Hom}\kern-.5truept}\kern1truept_{\mathrm{add}}}}
\newcommand{\Homaddinf}{\operatorname{\kern.5truept\underline{\kern-.5truept\mathsf{Hom}\kern-.5truept}\kern1truept_{\mathrm{addinf}}}}
\newcommand{\Wh}[2]{\operatorname{{\mathsf{Wh}}_{#1}}\!{\left(#2\right)}}
\newcommand{\Whf}[1]{\operatorname{\mathsf{Wh}_{#1}}}

\renewcommand{\H}[2]{\operatorname{{H}}(#1,#2)}
\newcommand{\Hder}[3]{\operatorname{{H}}_{\star}^{#1}(#2,#3)}
\newcommand{\coH}[2]{\operatorname{{H^*}}(#1,#2)}
\newcommand{\coHder}[3]{\operatorname{{H}}^{\star}_{#1}(#2,#3)}
\newcommand{\Hf}[1]{\operatorname{H_{#1}}}
\newcommand{\Lodot}{\overset{L}{\odot}}
\newcommand{\Lotimes}{\overset{L}{\otimes}}
\newcommand{\prefab}[1]{\widehat{ #1 }^{\mathsf{ab}}}
\newcommand{\sAb}{\prefab{\Delta}}
\newcommand{\Wab}{\W^{\ab}}

\DeclareMathOperator{\im}{Im}
\DeclareMathOperator{\U}{U}
\DeclareMathOperator{\dk}{\mathsf{C}_N}
\DeclareMathOperator{\Hsing}{\dk\nerf}
\newcommand{\Lbk}[1]{{\ell_{#1}}}
\newcommand{\bourne}[1]{\mathsf{B}(#1)}
\newcommand{\bournef}{\mathsf{B}}
\DeclareMathOperator{\wcat}{\omega\Cat}
\DeclareMathOperator{\wgpd}{\omega\Gpd}

\DeclareMathSymbol{\bang}{\mathclose}{operators}{"21}

\newcommand{\cub}{\square}
\newcommand{\cubc}{{\square^\mathsf{c}}}
\DeclareMathOperator{\dkcc}{\mathsf{C}_N^{\square^{\,\mathclap{\mathsf{c}}}}}
\DeclareMathOperator{\dkccg}{\mathsf{C}_{N,\gamma}^{\square^\mathsf{c}}}
\DeclareMathOperator{\dkc}{\mathsf{C}_N^{\square}}
\newcommand{\normcc}{N_{\cubc}}

\DeclareMathOperator{\ximuf}{\mu^{\prime}}
\newcommand{\xim}{m^{\prime}}
\newcommand{\ximu}[2]{\ximuf(#1,#2)}
\DeclareMathOperator{\xiwr}{\wr^{\prime}}

\DeclareMathOperator{\codim}{\mathop{codim}}
\DeclareMathOperator{\sg}{sg}

\newcommand\EnsSimp{\widehat{\Delta}}

\newcommand{\nerf}{\mathsf{N}}

\let\epsilon\varepsilon

\newcommand{\mdvirg}{\text{ ,}\hspace{10pt}}
\newcommand\zbox[1]{\makebox[0pt][l]{#1}}
\newcommand\pbox[1]{\zbox{\quad #1}} %
\newcommand\annot[1]{\zbox{$\,\scriptstyle #1$}}

\newcommand{\hdash}{\rotatebox[origin=c]{90}{\ensuremath\vdash}}

\title{Correspondances de Dold-Kan homotopiques}
\author{Léo Hubert}

\begin{document}
\frontmatter
\newgeometry{left=3em,right=2em,top=2em,bottom=2em} %

\vspace{1em}

\begin{center}
	\begin{minipage}[c]{.5\linewidth}
		\raggedright\includegraphics[height=5em]{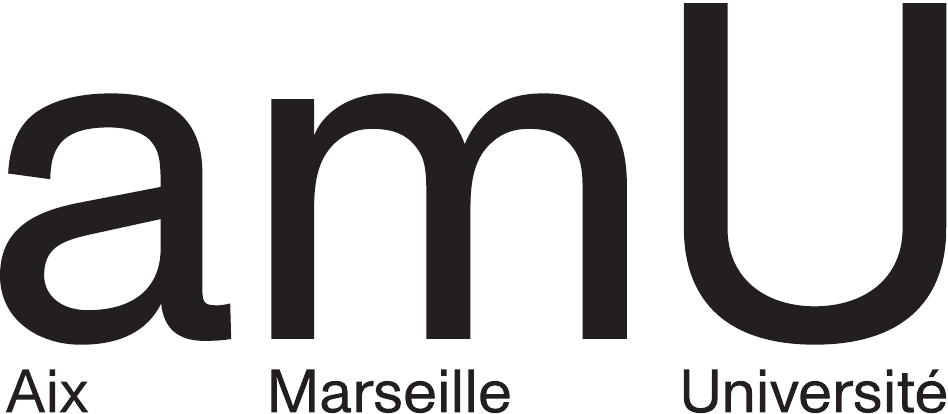}
	\end{minipage}\hfill
	\begin{minipage}[t]{.5\linewidth}
	\end{minipage}\hfill
\end{center}

\vspace{1em}

\begin{raggedleft}
	\begin{minipage}[c]{.81\linewidth}
		\textcolor{black}{\noindent\rule{\textwidth}{.7pt}}
	\end{minipage}
  \begin{minipage}[c]{.18\linewidth}
		\raggedright\SMALL\textsf{NNT : 2025AIXM0038}
  \end{minipage}
\end{raggedleft}

\vspace{1em}

\begin{raggedright}
    \Huge\textsc{Thèse de doctorat}\\
    \vspace{-.5em}
	  {\large\normalfont{Soutenue à AMU ― Aix-Marseille Université}}\\
	  \large\normalfont{le 13 février 2025 par}\\

\end{raggedright}
\vspace{3em}

\begin{center}
  	{\LARGE\textsc{Léo Hubert}}\\

    \vspace{3em}

    \begin{minipage}[c][][c]{0.77\linewidth}
		\textcolor{black}{\noindent\rule{\textwidth}{.7pt}}
    \vspace{1em}
\begin{center}
	  \LARGE \bfseries{CORRESPONDANCES} {DE} {DOLD-KAN} {HOMOTOPIQUES}\\ 
\end{center}

    \vspace{1.2em}
		\textcolor{black}{\noindent\rule{\textwidth}{.7pt}}

    \vspace{3em}

    \end{minipage}

\end{center}

\vspace{3em}

\begin{center}
	\begin{minipage}[t]{.25\linewidth}
\begin{raggedleft}
    	    \vspace{.5em}
        	\textbf{Discipline} \\

          Mathématiques
        	
    	    \vspace{2em}
        	\textbf{École doctorale}\\

          Mathématiques et informatique \\ (ED184)
        	
    	    \vspace{2em}
        	\textbf{Laboratoire} \\

           Institut de Mathématiques de Marseille \\ (I2M, UMR 7373)

\end{raggedleft}
\vspace{3em}
\begin{center}
\includegraphics[width=0.7\textwidth]{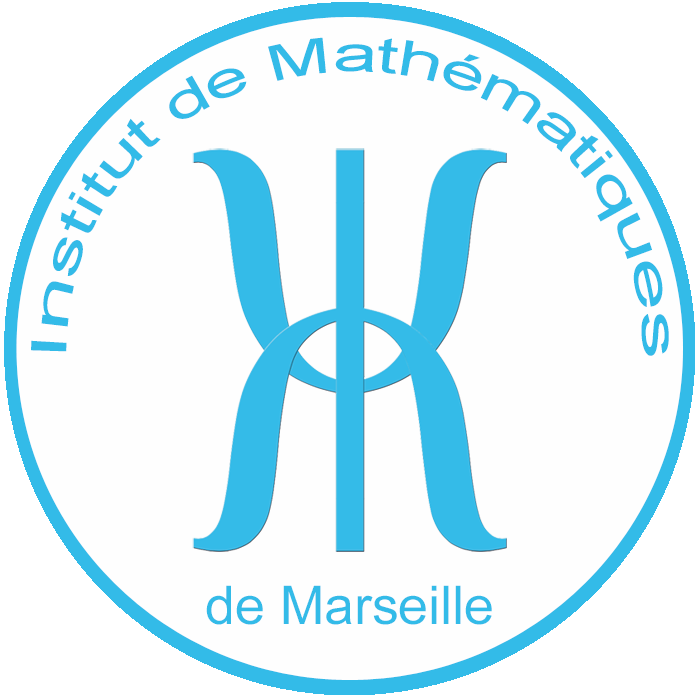}
\end{center}
	\end{minipage}\hspace{1em}\hfill
	\begin{minipage}[t]{.05\linewidth}
	    \rule[-290pt]{.7pt}{290pt}
	\end{minipage}\hfill
	\begin{minipage}[t]{.64\linewidth}
	    \vspace{.5em}
    	\textbf{Composition du jury}

	    \vspace{1em}
    	
        \begin{tabular}{p{15em} p{10em}}
        	\large{Clemens Berger} & Rapporteur \\
        	\small{PR, Université Côte d'Azur} & \small{Président du jury}\\
            \\
        	\large{Muriel Livernet} & Rapporteure \\
        	\small{PR, Université Paris Cité}\\
        	\\
        	\large{Georges Maltsiniotis} & Examinateur \\
        	\small{DR, CNRS, Université Paris Cité}\\
            \\
        	\large{Ieke Moerdijk} & Examinateur \\
        	\small{PR, Université d'Utrecht }\\
            \\
          \large{Christine Vespa} & Examinatrice \\
        	\small{PR, Aix-Marseille Université} \\
            \\
        	\large{Dimitri Ara} & Directeur de thèse \\
        	\small{MCF, Aix-Marseille Université} \\
           \\
        	\large{Yves Lafont} & Codirecteur de thèse \\
        	\small{PR, Aix-Marseille Université} \\
        \end{tabular}
        
	\end{minipage}\hfill
\end{center}       

\restoregeometry
\newpage\thispagestyle{empty}~

\newpage
\thispagestyle{empty}
\newgeometry{left=8em,right=8em,top=8em,bottom=4em} %

\mypar{\textit{Résumé}}
Ce travail trouve son origine dans les chapitres V et VII du manuscrit de
Grothendieck \emph{Pursuing Stacks}, qui contiennent une série de
questions ainsi qu'un formalisme jusqu'ici resté inexploré,
concernant l'interaction entre la notion de catégorie test et
l'homologie.

L'objectif principal de cette thèse est d'exhiber, dans le cadre des
catégories test, des
\emph{correspondances de Dold-Kan homotopiques}. Plus
précisément, on introduit, selon \hbox{Grothendieck,} un foncteur
généralisant le foncteur d'homologie simpliciale, associant aux
préfaisceaux en groupes abéliens sur une petite catégorie quelconque
un type d'homologie, c'est-à-dire un élément de la catégorie dérivée
des groupes abéliens en degré homologique positif. 
On cherche alors des conditions pour que ce foncteur induise une
équivalence de catégories, après localisation par la classe des morphismes dont
l'image dans la catégorie dérivée est un isomorphisme.

En général, il existe une seconde classe d'équivalences faibles, issue
de la théorie des catégories test, sur la
catégorie des préfaisceaux abéliens, et on appelle \emph{catégories de
\hbox{Whitehead}} les petites catégories pour lesquelles les deux
classes coïncident, généralisant ainsi le cas de $\Delta$. On
montre que des exemples importants de catégories test sont de
\hbox{Whitehead}, notamment la catégorie $\Theta$ de Joyal. On construit, pour
toute catégorie test locale de \hbox{Whitehead}, une structure de catégorie
de modèles sur sa catégorie des préfaisceaux abéliens
munie des équivalences faibles évoquées ci-dessus. On démontre alors
que pour toute catégorie test de \hbox{Whitehead}, le foncteur d'homologie
induit bien une équivalence entre les catégories localisées. On
obtient ainsi de nombreux exemples de correspondances de
Dold\nobreakdash-Kan homotopiques incluant, entre autres, la catégorie $\Theta$.

\mypar{\textit{Mots clés}} catégories test, théorème de Dold-Kan,
préfaisceaux en groupes abéliens, produit tensoriel de foncteurs,
homologie des préfaisceaux abéliens, types d'homologie.

\mypar{\textit{Abstract}}
This work originates from chapters V and VII of Grothendieck's manuscript
\emph{Pursuing Stacks}, which contains a series of questions, as well as a
previously unexplored formalism, concerning the interactions between
the notion of test categories and homology. 

The main objective of this thesis is to exhibit \emph{homotopical Dold-Kan
correspondences} in the context of test categories. More precisely, we
introduce, following Grothendieck, a functor generalizing simplicial
homology, from the category of abelian presheaves over any small
category to the derived category of abelian groups in non-negative homological
degree. We then look for conditions ensuring that this functor induces an
equivalence of categories, after localization by the class of morphisms whose image in the derived category is an isomorphism. 

Generally, there exists a second class of weak equivalences, arising
from the theory of test categories, on the category of abelian
presheaves, and we call \emph{Whitehead categories} those small
categories for which these two classes coincide, generalizing the case
of $\Delta$. We show that important examples of test categories are
Whitehead categories, notably Joyal's category $\Theta$. We construct,
for any Whitehead local test category, a model category structure
on its category of abelian presheaves with the weak
equivalences mentioned above. We then prove that for any Whitehead
test category, the homology functor does induce an equivalence
between the localized categories. We obtain this way many examples of
homotopical Dold-Kan correspondences, including, among others, the
category $\Theta$.

\mypar{\textit{Keywords}} 
abelian presheaves,
Dold-Kan theorem,
homology of abelian presheaves,
homology types,
tensor product of functors,
test categories.

\restoregeometry

\thispagestyle{empty}
\newgeometry{left=8em,right=8em,top=8em,bottom=8em} %
\chapter*{Remerciements}

Un travail de cette nature ne peut qu'être le fruit d'une heureuse
coïncidence de lieux et de rencontres, de
contextes que nous habitons avec attention et qui nous permettent, au
bout du compte, de les réaliser. J'ai conscience de la chance que j'ai
eue de pouvoir suivre le chemin (mouvementé) menant jusqu'à la
rédaction de ce manuscrit, et tenter d'exprimer ma reconnaissance à
chacune des belles circonstances qui l'ont construit est un exercice
difficile mais essentiel. 

Il m'est difficile de trouver des mots pour exprimer ma gratitude
envers Dimitri Ara, qui, bien qu'il fut un directeur de thèse disponible, rigoureux,
patient, perspicace et pédagogue, représente énormément
plus que tout ça. J'ai eu l'immense privilège d'être formé par un mathématicien extrêmement lucide et un être humain particulièrement sensible. Pour tout cela, je te remercie profondément.
Aussi, je remercie sincèrement Yves Lafont qui s'est toujours montré intéressé, accueillant et disponible pour m'aider. Merci pour tes remarques très pertinentes pendant la préparation de la soutenance. 

Je remercie chaleureusement Clemens Berger et Muriel Livernet, qui ont
accepté de rapporter cette thèse. Je tiens à remercier
particulièrement Muriel, dont les nombreux retours ont été
extrêmement précieux et ont permis de grandement améliorer ce manuscrit. 

Je remercie également Christine Vespa, Ieke Moerdijk et Georges
Maltsiniotis d'avoir bien voulu participer au jury de cette thèse.
Sans la tâche colossale que Georges s'est donnée d'éditer
\emph{Pursuing Stacks}, la question à l'origine de cette thèse
n'aurait peut être jamais émergée des quelques six cents pages de réflexions que contient ce manuscrit. Certains éléments de cette thèse, en particulier dans le
chapitre $2$, sont tirés des notes manuscrites de Georges résumant
avec beaucoup de clarté les propos de Grothendieck. Merci
également pour ta disponibilité.

Je tiens aussi à remercier François Métayer et Lionel Vaux Auclair
d'avoir participé, avec beaucoup d'attention et de bienveillance, à
mon comité de suvi de thèse.

J'ai baigné, depuis les fins timides du confinement jusqu'à aujourd'hui,
dans un environnement particulièrement doux et amical dans lequel
découvrir la pratique mathématique fut extrêmement agréable. Merci à Lionel pour ta
bienveillance et la chaleur que tu apportes à l'équipe, Laurent pour ta désinvolture et ton
énergie toujours sincère, Alexey pour les questions qui
énervent Laurent et les images de la mer Baltique. Merci Patrick pour
ton impertinence. Merci Olivier et Raphaël.
Merci Christine de t'être rendue disponible pour mes questions et de
m'avoir donné, pour la première fois, l'impression que ce que je
faisais pouvait intéresser quelqu'un d'autre que Dimitri. 
Merci Sary,
mes premières expériences d'enseignement auraient été catastrophiques
sans tes conseils et ton implication. Merci Myriam, Lionel et Étienne,
j'ai adoré enseigner à vos côtés. Merci Enea, Nathalie, et toute
l'équipe de l'IREM pour le si beau projet que vous portez.

Je tiens à remercier particulièrement celles et ceux qui, en plus de l'habiter,
construisent (ou construisaient) cette maison : Peter, Jean-Bruno, Pierre,
Corinne, David, Sonia, Jessica, Yosr, et tous ceux dont je n'ai pas réalisé
le rôle pourtant fondamental qu'ils jouaient.

Merci aux doctorant$\cdot$e$\cdot$s de l'I2M qui m'ont permis de les rencontrer malgré
ma timidité : Tommy-Lee, Juliette, Axel, Cédric, Adrien, ainsi que celles et ceux devant qui je n'ai pas réussi à enlever mon
masque de fantôme.

Merci à Rémy et Zoé pour toutes les discussions, et le beau cocon dont
nous nous sommes enveloppés. Je suis ému de devenir le \og petit
dernier\fg{}. 

Merci à vous qui, en dehors de Marseille, m'ont prouvé que les
mathématiques pouvaient aussi être un lieu d'amitiés et de poésie :
Saad, Jérémie, Léonard, Félix, Victor, Francesca.

À vous qui êtes le cœur de mon village : Samuel, Matthew, Gaya, Jérémi, Luke, Guillaume, Roxane, Clément, David, Domi, Bissera,
Antoine, Antonin, et toutes les cigales Marseillaises.

À Sarah-Charlotte, Nico, Baptiste, Tom, Simon, qui n'ont pas laissé la distance nous séparer. 

\smallskip
À Tim qui partout murmure la vie.
\smallskip

À mes parents et à ma famille qui m'ont toujours encouragé et
valorisé, que j'étudie les maths ou que je joue dans les marchés pour
quelques légumes ramollis de soleil. 

À Marseille, aux rochers de Malmousque et aux cafés qui m'ont
soufflé la majeure partie des réflexions présentes dans ce manuscrit.
Merci à la mer d'être une si patiente confidente et à la montagne d'être un si puissant refuge.
\smallskip

Merci Aurélie d'être la mer et la montagne, et celle qui les émerveille. 

\smallskip

\restoregeometry

\tableofcontents 
\mainmatter

\chapter*{Introduction}

\section*{Types d'homotopie et préfaisceaux}
En 1983, Grothendieck envoie une lettre à Quillen, qui donnera
naissance à \emph{Pursuing Stacks} \cite{pursuingstacks}. La plus
grande partie de ce manuscrit de près de six cents pages concerne la
\og poursuite des types d'homotopie \fg{}. Concrètement, on
note~$\Top$ la catégorie des espaces topologiques, et
\[
\W_\Top \subset \Fl(\Top)
\]
la classe des équivalences faibles d'homotopie, c'est-à-dire des
applications continues induisant un isomorphisme entre tous les
groupes d'homotopie. On peut alors définir une nouvelle catégorie
notée
\[
\Hot := \W_\Top^{-1} \Top
\]
obtenue en ajoutant formellement des inverses aux morphismes de
$\W_\Top$, qui deviennent donc des isomorphismes : on dit qu'on
\emph{localise} la catégorie $\Top$ par la classe $\W_\Top$. On
appelle alors \emph{types d'homotopie} les objets de la
catégorie~$\Hot$. 

La poursuite des \emph{modèles} des types d'homotopie consiste à
chercher d'autres catégories que celle des espaces topologiques
permettant, par le même procédé de localisation, d'obtenir une
catégorie équivalente à $\Hot$. On appelle \emph{modélisateur} la
donnée d'un couple $(\M,\W)$ où~$\M$ est une catégorie et $\W$ est une
classe de flèches de $\M$, telle que la catégorie~$\W^{-1}\M$ obtenue
en localisant $\M$ par~$\W$ est équivalente à la catégorie $\Hot$ des
types d'homotopie. Par exemple, un résultat classique de topologie
affirme que le couple formé des CW-complexes et des équivalences
(fortes) d'homotopie entre ceux-ci est un modélisateur. 

Un autre exemple est donné par la catégorie $\EnsSimp$ des
\emph{ensembles simpliciaux}. On dispose d'un foncteur de
\emph{réalisation géométrique}
\[
|-| : \EnsSimp \to \Top 
\]
grâce auquel on peut définir une classe d'\emph{équivalences faibles
simpliciales} notée 
\[
\W_\Delta := \lbrace f \in \Fl(\EnsSimp) \mid |f| \in \W_\Top\rbrace \mdvirg
\]
et on pose alors 
\[
\Hot_\Delta := \W_\Delta^{-1}\EnsSimp \pbox{.}
\]
Milnor prouve en 1957 dans \cite{milnor1957realization} que le couple
$(\pref{\Delta},\W_\Delta)$ est un modélisateur, et que le foncteur de
réalisation géométrique induit une équivalence de catégories 
\[
\Hot_\Delta \to \Hot \pbox{.}
\]
Les ensembles simpliciaux, munis des équivalences faibles
simpliciales, sont souvent considérés comme le modèle standard des
types d'homotopie. 

Pour Grothendieck, le modèle privilégié est donné par la catégorie
$\Cat$ des petites catégories : le foncteur \emph{nerf simplicial} 
\[
\nerf : \Cat \to \EnsSimp \pbox{,}
\]
introduit par Grothendieck dans \cite{grothendieck1961techniques},
associe à toute petite catégorie un ensemble simplicial, et donc un
type d'homotopie. On note alors
\[
  \W_\infty := \lbrace f \in \Fl(\Cat) \mid \nerf(f) \in \W_\Delta
  \rbrace
\]
la classe des foncteurs entre petites catégories dont l'image par le
foncteur nerf est une équivalence faible simpliciale, et on appelle
\emph{équivalences de Thomason} ses éléments. En 1972, Illusie
démontre dans sa thèse \cite{illusie1971cotangent}, tout en attribuant
le résultat à Quillen, que le foncteur nerf simplicial induit une
équivalence de catégories 
\[
\nerf : \W_\infty^{-1}\Cat \to \Hot_\Delta \pbox{.}
\]
Le couple $(\Cat,\W_\infty)$ est donc aussi un modélisateur, que
Grothendieck nomme \emph{modélisateur fondamental}. On note alors
simplement 
\[
\Hot=\W_\infty^{-1}\Cat 
\]
la catégorie localisée associée. Bien que le foncteur nerf admette un
adjoint à gauche $\mathsf{c} : \pref{\Delta}\to\Cat$ appelé foncteur
de \emph{réalisation catégorique}, ce dernier ne préserve pas les
équivalences faibles et n'induit donc pas de foncteur entre les
catégories localisées. Un quasi-inverse au foncteur nerf est donné par
le foncteur associant à tout ensemble simplicial sa \emph{catégorie
des éléments} 
\[
i_\Delta : \EnsSimp \to \Cat \mdvirg X 
\mapsto \tranche{\Delta}{X} \pbox{.}
\]
On peut toutefois noter le résultat de Fritsch et Latch
\cite{fritsch1981homotopyinverse} qui affirme qu'en précomposant deux
fois le foncteur de réalisation catégorique avec le foncteur de
subdivision, on obtient bien un foncteur entre les catégories
localisées, qui induit un autre quasi-inverse au foncteur nerf
simplicial.

Se donnant d'abord l'objectif de caractériser \emph{tous} les
modélisateurs, Grothendieck s'intéresse plus particulièrement aux
modélisateur $(\M,\W)$ où~$\M$ est une catégorie de préfaisceaux sur
une petite catégorie $A$, et où~$\W$ est définie en tirant en arrière
les équivalences de Thomason le long du foncteur catégorie des
éléments 
\[
i_A : \pref{A}\to\Cat \mdvirg X \mapsto \tranche{A}{X} \mdvirg 
\]
généralisant ainsi le cas de la catégorie $\Delta$. On note alors 
\[
  \W_A := \lbrace f \in \Fl(\pref{A}) \mid i_A(f) \in \W_\infty\rbrace
\]
la classe des \emph{équivalences test} de $\pref{A}$. 

Grothendieck appelle alors \emph{catégories pseudo-test} les
catégories $A$ telles que le couple~$(\pref{A},\W_A)$ est un
modélisateur. Plus précisément, en notant 
\[
\Hot_A := \W_A^{-1}\pref{A} 
\]
la catégorie localisée associée au couple $(\pref{A}, \W_A)$, la
petite catégorie $A$ est une catégorie pseudo-test si le foncteur 
\[
i_A : \Hot_A \to \Hot 
\]
induit par le foncteur catégorie des éléments est une équivalence de
catégories. Il introduit également les notions plus fortes de
catégories test faibles, test, et test strictes, que l'on détaillera
dans la section \ref{secCatTest}, et en trouve des caractérisations
particulièrement simples. Le slogan de Grothendieck est alors que
toute catégorie test est \og aussi bonne \fg{} que $\Delta$ pour faire
de la théorie de l'homotopie.

Parmi les exemples de catégories test (et donc de catégories
pseudo-test), on peut citer, en plus de la catégorie $\Delta$, les
catégories $\Delta^n$ pour $n \geq 1$, la catégorie cubique $\cub$
ainsi que la catégorie cubique avec connexion $\cubc$. Un exemple
important est celui de la catégorie $\Theta$ de Joyal
\cite{joyal1997disks} qui, en plus d'être une catégorie test, est
intimement liée à diverses notions de catégories supérieures. On
introduira aussi, pour tout entier $n \geq 1$, une sous-catégorie de
$\Theta_n$ notée~$\Xi_n$, définie comme l'image dans $\Theta_n$ de la
catégorie $\Delta^n$, et on montrera que celle-ci est une catégorie
test stricte.

\section*{Abélianisation des types d'homotopie}

Les chapitres V et VII de \emph{Pursuing Stacks} sont consacrés à la
question de l'abélianisation des types d'homotopie. On note
\[
\Hotab = \W_{qis}^{-1}\Ch(\Ab) 
\]
la catégorie des \emph{types d'homologie}, obtenue en localisant la
catégorie des complexes de chaînes de groupes abéliens concentrés en
degré positif par les quasi-isomorphismes.

Grothendieck pose la question de l'existence d'un foncteur canonique 
\[
\Hot \to \Hotab 
\]
associant à tout type d'homotopie son type d'homologie, sans utiliser
l'intermédiaire d'une catégorie test spécifique.

La catégorie $\Delta$ permet en effet de définir un tel foncteur, de
la manière suivante. On sait que si~$X$ est un groupe abélien
simplicial (c'est-à-dire, un préfaisceau en groupes abéliens sur
$\Delta$), on peut définir ses groupes d'homologie comme les groupes
d'homologie de son image par le foncteur \emph{complexe normalisé} 
\[
\dk : \prefab{\Delta} \to \Ch(\Ab) 
\]
où $\prefab{\Delta}$ désigne la catégorie des préfaisceaux en groupes
abéliens sur $\Delta$. En notant $\U : \prefab{\Delta} \to
\pref{\Delta}$ le foncteur d'oubli de la structure de groupe abélien,
on peut montrer que si $f : X \to Y$ est un morphisme de
$\prefab{\Delta}$ tel que $\U(f)$ est une équivalence faible
simpliciale, alors le morphisme de complexes $\dk(f)$ est un
quasi-isomorphisme. En d'autres termes, en notant 
\[
\Hotab^{\U}_\Delta = (\U^{-1}\W_\Delta)^{-1}\prefab{\Delta} 
\]
la catégorie obtenue en localisant la catégorie $\prefab{\Delta}$ par les
morphismes dont le morphisme d'ensembles simpliciaux
sous-jacent est une équivalence faible simpliciale, le foncteur
complexe normalisé induit un foncteur encore noté 
\[
\dk : \Hotab^{\U}_\Delta \to \Hotab \pbox{.}
\]
De plus, un résultat classique affirme que les groupes d'homologie d'un groupe abélien
simplicial correspondent aux groupes d'homotopie de son ensemble
simplicial sous-jacent. On peut alors montrer que 
le foncteur d'abélianisation
\[
\pref{\Delta}\to\prefab{\Delta} \mdvirg X \mapsto \Z^{(X)}
\]
envoie les éléments de $\W_\Delta$ sur des éléments de
$\U^{-1}\W_\Delta$, et induit donc un foncteur que Grothendieck
appelle \emph{foncteur de Whitehead} 
\[
\Whf{\Delta}:\Hot_\Delta \to \Hotab_{\Delta}^{\U} \pbox{.}
\]

Pour toute petite catégorie $A$, on obtient alors un diagramme commutatif à isomorphisme naturel près
\[
\xymatrix{
\pref{A} \ar[r]^{i_A} \ar[d]^{} 
& \Cat \ar[r]^{\nerf} \ar[d]^{} 
& \EnsSimp \ar[d]^{} \ar[r]^{\Whf{\Delta}} 
& \prefab{\Delta} \ar[d]^{} \ar[r]^{\dk} 
& \Ch(\Ab) \ar[d]^{} 
\\
\Hot_A \ar[r]_{i_A} 
& \Hot \ar[r]_{\nerf} & \Hot_\Delta \ar[r]_{\Whf{\Delta}}
&\Hotab^{\U}_\Delta \ar[r]_{\dk}
& \Hotab
} 
\]
où les flèches verticales
correspondent aux projections canoniques,
décrivant le processus \emph{d'abélianisation des types
d'homotopie}.

En revanche, le complexe associé au nerf de la catégorie des éléments
d'un préfaisceau est assez compliqué à manipuler, et il n'est pas
raisonnable d'espérer faire des calculs d'homologie explicites avec
cette description. 

\medskip

Après avoir introduit et étudié les catégories test, Grothendieck formule alors une série de questions
concernant l'abélianisation des types d'homotopie :
\begin{enumerate}
\item Étant donnée une petite catégorie quelconque $A$, peut-on
toujours définir un foncteur 
\[
\Hf{A} : \prefab{A}\to \Hotab 
\]
de sorte que le diagramme 
\[
\xymatrix{
\pref{A} \ar[d]_{\Whf{A}}  \ar[r]^{can} 
& \Hot_A \ar[r]^{i_A}
& \Hot \ar[d]^{\Whf{}} \\
\prefab{A} \ar@{-->}[rr]_{\Hf{A}}
&& \Hotab \pbox{,}
} 
\]
où la flèche de droite désigne le foncteur d'abélianisation
utilisant l'intermédiaire de la catégorie $\Delta$ défini ci-dessus,
soit commutatif à isomorphisme naturel près ?

\item Si un tel foncteur existe, est-il toujours possible de le
factoriser par la catégorie 
\[
\Hotab_A^{\U} := (\U^{-1}\W_A)^{-1}\prefab{A}
\]
de sorte que le diagramme 
\[
\xymatrix{
\prefab{A} \ar[rr]^-{\Hf{A}}="a" 
\ar[rd]_{can} 
&& \Hotab\\
& \Hotab_A^{\U} \ar@{-->}[ru]_(.6){\widetilde{\Hf{A}}}
} 
\]
soit commutatif ? En d'autres termes, est-ce que le foncteur $\Hf{A}$
envoie les éléments de $\U^{-1}\W_A$ sur des isomorphismes de $\Hotab$
?

\item Dans le cas où la réponse à la question précédente est
affirmative, on dispose alors d'un diagramme commutatif à isomorphisme
près
\[
\xymatrix{
{ \Hot_A } \ar[r]^{ i_A } \ar[d]_{ \Whf{A} } & { \Hot } \ar[d]^{
\Whf{} } \\
{ \Hotab_A^{\U} } \ar[r]_{ \Hf{A} } & { \Hotab }} 
\]
et la théorie des catégories test donne des caractérisations des
petites catégories~$A$ telles que la flèche du haut est une
équivalence de catégories. Peut-on alors trouver des conditions pour
que la flèche du bas dans ce diagramme soit une équivalence de
catégories ? 
\end{enumerate}

\section*{Deux classes d'équivalences faibles abéliennes}

La stratégie que nous allons explorer dans cette thèse est la
suivante. On va définir, selon Grothendieck,
un foncteur canonique 
\[
\Hf{A} : \prefab{A}\to\Hotab 
\]
pour toute petite catégorie $A$, généralisant l'homologie simpliciale
dans le cas où $A$ est la catégorie~$\Delta$, ainsi que l'homologie
des groupes dans le cas où
$A$ est le groupoïde associé à un groupe. Pour tout
préfaisceau en groupes abéliens~$X$ sur $A$, on dira que~$\Hf{A}(X)$ est
\emph{l'homologie de $A$ à coefficients dans~$X$}. Ce foncteur, introduit
par Grothendieck dans \emph{Pursuing Stacks}, correspond au foncteur
dérivé à gauche du foncteur de limite inductive 
\[
\limind\nolimits_{{A}^{\op}}^\Ab : \prefab{A} \to \Ab \mdvirg
\]
que l'on peut aussi voir comme le foncteur associant à tout
préfaisceau en groupes abéliens la limite inductive homotopique du
diagramme de complexes de chaînes concentrés en degré $0$ associé. De
ce point de vue, le foncteur $\Hf{A}$ est un équivalent abélien au
foncteur $i_A$, puisque ce dernier associe à tout préfaisceau
d'ensembles la limite inductive homotopique du diagramme d'ensembles
simpliciaux discrets associé. Il est intéressant de noter
que le foncteur~$\Hf{A}$ peut aussi être obtenu à partir du
\emph{produit tensoriel de foncteurs} dérivé utilisé dans la
communauté de l'homologie des foncteurs (voir par exemple
\cite{pirashvili2003introduction,djament2010homology}),
qui est exprimé ici dans le
langage des préfaisceaux, et donc des types d'homotopie, plutôt que
dans celui des modules à droite sur une catégorie additive.

Ce foncteur permet alors d'introduire, pour toute petite catégorie
$A$, une classe d'\emph{équivalences faibles abéliennes} sur la
catégorie $\prefab{A}$ en tirant en arrière les isomorphismes de
$\Hotab$ par le foncteur $\Hf{A}$. On note alors
\[
  \W_A^\ab \subset \Fl(\prefab{A})
\]
la classe des équivalences faibles abéliennes de $\prefab{A}$,
analogue abélien de la classe $\W_A$ des équivalences test sur
$\pref{A}$.
On peut alors introduire la catégorie 
\[
\Hotab_A := {\W_A^\ab}^{-1}\prefab{A} 
\]
obtenue en localisant la catégorie $\prefab{A}$ par la classe
$\W_A^\ab$. 

On montrera alors dans la proposition
\ref{homologiePrefaisceauxLibres} qu'on obtient bien de cette manière,
pour toute petite catégorie $A$, un diagramme commutatif à
isomorphisme naturel près 
\[
\tag{1}
\label{introHomologieCatElements}
\xymatrix{
\pref{A} \ar[d]_{\Whf{A}}  \ar[r]^-{can} 
& \Hot_A \ar[r]^-{i_A} \ar[d]^{\Whf{A}} 
& \Hot \ar[d]^{\Whf{}} \\
\prefab{A} \ar[r]_-{can}
&\Hotab_A \ar[r]_{\Hf{A}}
& \Hotab 
} 
\]
où la flèche de droite désigne à nouveau le foncteur d'abélianisation
des types d'homotopie utilisant l'intermédiaire de la catégorie
$\Delta$.

Étant donnée une petite catégorie $A$, on dispose donc \emph{a priori}
de deux catégories différentes 
\[
\Hotab_A \mdvirg \Hotab_A^{\U} 
\]
obtenues respectivement en localisant la catégorie $\prefab{A}$ par
les classes de flèches 
\[
\W_A^\ab \mdvirg \U^{-1}\W_A \mdvirg
\]
et on va s'intéresser aux conditions sous lesquelles ces deux classes
de flèches coïncident. On dira que $A$ est une \emph{catégorie de
Whitehead} lorsque c'est le cas. C'est par exemple le cas de la
catégorie $\Delta$. On montre de plus le résultat suivant.
\begin{prop*}
Soit $u : A \to B$ un foncteur
entre petites catégories tel que pour tout objet $b$ de $B$, la
catégorie~$\tranche{A}{b}$ a le type d'homotopie du point (on dit que
$u$ est \emph{asphérique}). Alors si $A$ est une catégorie de
Whitehead, il en est de même pour $B$.
\end{prop*} 
Puisqu'il existe un foncteur asphérique $\Delta \to \Theta_n$ pour
tout entier $n\geq1$, la proposition ci-dessus implique en particulier que
$\Theta_n$ est une catégorie de Whitehead. En revanche, on ne dispose
pas d'un foncteur asphérique $\Delta \to \Theta$. En utilisant une
variation de la notion de foncteur asphérique pour les foncteurs de
but une catégorie de préfaisceaux, on montre qu'il existe un
zig-zag de foncteurs asphériques 
\[
\Theta \leftarrow C \rightarrow \Delta 
\]
où $C$ est une catégorie de Whitehead, et la proposition précédente
permet alors d'affirmer que
$\Theta$ est une catégorie de Whitehead. 

En revanche, il existe des petites catégories qui ne sont pas des catégories
de Whitehead : on montre par exemple de la catégorie $\Grefl$ des
globes réflexifs n'est pas une catégorie de Whitehead.

\medskip

Les catégories de Whitehead présentent de très bonnes propriétés
dans le cadre de notre étude. Par exemple, Cisinski a montré que pour
toute catégorie test locale $A$, la catégorie des préfaisceaux
d'ensembles sur $A$ peut être munie d'une structure de catégorie de
modèles dont les équivalences faibles sont les équivalences test, et
dont les cofibrations sont des monomorphismes. Si $A$ est également
une catégorie de Whitehead, on montre qu'on peut transférer cette
structure de catégorie de modèles en utilisant le foncteur d'oubli $\U
: \prefab{A}\to\pref{A}$, et on obtient le résultat suivant.
\begin{theorem*}
Soit $A$ une catégorie test locale de Whitehead. Alors la catégorie
$\prefab{A}$ peut être munie d'une structure de catégorie de modèles
dont les équivalences faibles sont les éléments de
$\Wab_A=\U^{-1}\W_A$.
\end{theorem*}

Pour cette structure de catégorie de modèles, le foncteur $\Hf{A}$ est
alors le foncteur dérivé à gauche d'un foncteur de Quillen à gauche.

\section*{Types d'homologie et préfaisceaux}
Revenons maintenant à la troisième question de Grothendieck, en la
reformulant légèrement : on dispose, pour toute petite catégorie $A$,
d'un diagramme commutatif 
\[
\xymatrix{
\Hot_A \ar[r]^{i_A} \ar[d]_{\Whf{A}} & \Hot \ar[d]^{\Whf{}} \\
\Hotab_A \ar[r]_{\Hf{A}} & \Hotab \pbox{,}
} 
\]
où on souligne qu'on a remplacé, par rapport à la question originale,
la catégorie $\Hotab_A^{\U}$ par la catégorie $\Hotab_A$. 
\begin{defin*}
On dit que $A$ est une \emph{catégorie pseudo-test homologique} si le
foncteur $\Hf{A}$ ci-dessus est une équivalence de catégories.
\end{defin*}

Par exemple, dans le cas où $A$ est la catégorie $\Delta$, le célèbre théorème de Dold-Kan
\cite{kan1958functors,dold1958homology} affirme que le foncteur
complexe normalisé 
\[
\dk : \prefab{\Delta}\to\Ch(\Ab) 
\]
est une équivalence de catégories. Puisque le foncteur $\dk$ coïncide,
après localisation, avec le foncteur $\Hf{\Delta}$, cela implique
que le foncteur 
\[
\Hf{\Delta} : \Hotab_\Delta \to \Hotab 
\]
est aussi une équivalence de catégories. La catégorie $\Delta$, en
plus d'être une catégorie test, est aussi une catégorie pseudo-test
homologique. Dans le cas où~$A$ est la
catégorie cubique avec connexions $\cubc$ ou la catégorie $\Grefl$ des
globes réflexifs, on dispose aussi d'un théorème de Dold-Kan \og
strict \fg{}, et on verra que les foncteurs en jeu coïncident avec le
foncteur $\Hf{A}$, ce qui implique que ces deux catégories sont aussi
des catégories pseudo-test homologiques.

L'objectif principal de ce travail est d'exhiber une large classe de
catégories pseudo-test homologiques : c'est en ce sens que l'on parle
\emph{correspondances de Dold-Kan homotopiques}. 

Nous sommes particulièrement intéressés au cas des catégories test.
Notre résultat principal dans cette direction est le suivant, utilisant
la structure de catégorie de modèles évoquée précédemment.

\begin{theorem*}
Si $A$ est une catégorie test de Whitehead, alors $A$ est une
catégorie pseudo-test homologique.
\end{theorem*}

Ce théorème fournit des exemples importants de catégories pseudo-test
homologiques, comprenant : 
\begin{itemize}
\item les catégories
$\Delta^n$ pour $n>1$;
\item la catégorie $\cubc$ des cubes avec
connexions;
\item la catégorie~$\Theta$ de Joyal ainsi que les catégories
$\Theta_n$ pour $n\geq 1$;
\item les sous-catégories notées~$\Xi_n$
de~$\Theta_n$ définies
comme l'image de la catégorie~$\Delta^n$ dans $\Theta_n$.
\end{itemize}

En revanche, ceci ne suffit pas à caractériser toutes les catégories
pseudo-test homologiques : il existe des catégories
$A$ qui ne sont pas des catégories de Whitehead, et telles que le
foncteur $\Hf{A}$ est une équivalence de catégories. C'est par exemple
le cas de la catégorie $\Grefl$ des globes réflexifs. Même en limitant
la recherche aux catégories test, on peut trouver des catégories test
qui sont des catégories pseudo-test homologiques, mais qui ne sont pas
des catégories de Whitehead. Un tel exemple est donné par la catégorie
$\tranche{\Delta}{\nerf(\Grefl)}$. 

\bigbreak 

Il est important de souligner un point de confusion possible :
\emph{il ne s'agit pas} de chercher à caractériser les catégories $A$
dont la catégorie des préfaisceaux \emph{d'ensembles} modélise les types
d'homologie des espaces topologiques. Une réponse partielle, mais
peut-être suffisante, à cette question est déjà connue. On note
\[
  \W_\infty^\ab \subset \Fl(\Cat)
\]
la classe des foncteurs entre petites catégories induisant un
isomorphisme en homologie, c'est-à-dire les foncteurs $A \to B$
telles que le morphisme de complexes de chaînes $\dk \nerf A \to \dk
\nerf B$ est un quasi-isomorphisme. On montre dans l'annexe
\ref{annexeDerivateurs} que $\W_\infty^\ab$ est un \emph{localisateur
fondamental}, ce qui signifie entre autres que $\W_\infty^\ab$ vérifie
un analogue du théorème $A$ de Quillen \cite{quillenktheory}. Cela
implique en particulier, grâce au théorème de Cisinski-Grothendieck
\cite{cisinski2004localisateur} affirmant que tout localisateur
fondamental contient la classe $\W_\infty$, que toute catégorie test
est également une \og $\W_\infty^\ab$-catégorie test \fg{}. En d'autres termes,
la catégorie obtenue en localisant la catégorie $\pref{A}$ par les
morphismes dont l'image par le foncteur $i_A$ est dans $\W_\infty^\ab$
est équivalente à la catégorie~${\W_\infty^\ab}^{-1}\Cat$, qui à son
tour est équivalente à la catégorie des types d'homologie des espaces
topologiques. 

Toutefois, la classe~$\W_\infty^\ab$ joue un rôle
fondamental dans ce travail. En particulier, on peut grâce à elle
définir la notion de foncteur asphérique en homologie que l'on va
maintenant détailler.

\section*{Cofinalité en homologie}
Une des différences fondamentales entre le cadre des catégories test
(et donc des préfaisceaux d'ensembles) et le nôtre est qu'on ne
dispose pas d'un foncteur canonique 
\[
\prefab{A} \to \Ch(\Ab) 
\]
factorisant, pour toute petite catégorie $A$, le foncteur
$\Hf{A}: \prefab{A}\to\Hotab$,
c'est-à-dire d'un foncteur pouvant jouer le rôle du foncteur 
\[
  i_A : \pref{A}\to\Cat
\]
pour la limite inductive homotopique des préfaisceaux d'ensembles. 

On peut trouver des expressions simples de~$\Hf{A}$ en construisant
des \emph{intégrateurs} sur $A$, c'est-à-dire des foncteurs $L_A : A
\to \Ch(\Ab)$ qui, vus comme des complexes de foncteurs $A\to\Ab$,
forment une résolution projective du foncteur constant de valeur $\Z$
sur $A$. Étant donné un intégrateur $L_A$, on montre que son extension
de Kan à gauche
\[
{L_A}_! : \prefab{A} \to \Ch(\Ab) 
\]
induit, après localisation, un foncteur isomorphe au foncteur
$\Hf{A}$. Une expression particulièrement simple est possible dans le
cas où, pour tout entier $n\geq 0$, le foncteur 
\[
(L_A)_n : A \to \Ab 
\]
est une somme directe de foncteurs représentables. L'extension de Kan
à gauche de $L_A$ est alors donnée, en tout degré, comme une somme de
foncteurs d'évaluation.
On dit alors que $L_A$ est un \emph{intégrateur libre} (Grothendieck
utilise le terme d'intégrateur \og quasi-spécial \fg{}). 

En particulier, on montre que pour une classe de petites catégories
admettant une notion de dimension, on peut chercher
à construire un intégrateur libre en donnant un signe à tous les
monomorphismes de codimension~$1$. On peut alors exprimer l'homologie
de tout préfaisceau en groupes abéliens~$X$ sur~$A$ comme l'homologie
d'un complexe de la forme 
\[
\cdots \leftarrow \bigoplus_{\dim a = n} Xa \leftarrow
\bigoplus_{\dim a = n+1} Xa \leftarrow \cdots \pbox{,}
\]
où la différentielle est donnée par la somme des images par
$X$ des monomorphismes de codimension $1$, comptés avec leur signe. On
retrouve de cette manière le complexe non normalisé associé à un
groupe abélien simplicial, et on donnera un tel système de signes pour
la catégorie $\Xi_n$ dans la section~\ref{secXi}.

\medskip

Par contre, on est toujours confronté à la question du choix d'un tel
intégrateur. Bien qu'on puisse toujours choisir de travailler avec
l'intégrateur que l'on appelle \emph{intégrateur de Bousfield-Kan},
défini pour toute petite catégorie $A$ par la composition des
foncteurs 
\[
A \hookrightarrow \pref{A} \xrightarrow{i_A} \Cat \xrightarrow{\nerf}
\EnsSimp \xrightarrow{\dk} \Ch(\Ab) \pbox{,}
\]
le complexe obtenu de cette manière est difficile à manipuler, comme
on l'a déjà dit, et on
préfère chercher des manières intrinsèques d'exprimer l'homologie des
préfaisceaux abéliens sur $\prefab{A}$.

Une des difficultés venant avec ce choix est la suivante. Dans le cas
des préfaisceaux d'ensembles, on dispose, pour tout foncteur
$u : A \to B$ entre petites catégories, d'une
transformation naturelle s'insérant dans le diagramme
\[
\xymatrix{
{\pref{B}} 
\ar[rr]^{u^*} \ar[rd]_{{i_B}}^{}="c"
&& {\pref{A}} 
\ar@{}"c"|(.4){}="d"|(.9){}="e"
\ar@2"d";"e"_{\lambda_u}
\ar[ld]^{{i_A}} \\
& {\Cat} & \pbox{,} 
} 
\]
où $u^*$ désigne le foncteur de précomposition associant à un
préfaisceau $X$ sur~$B$ le préfaisceau $X\circ u$ sur $A$.
On peut alors montrer que $\lambda_u$ est argument par argument dans $\W_\infty$
(autrement dit, le foncteur~$u^*$ commute avec le foncteur de limite
inductive homotopique) \emph{si
et seulement si}~$u$ est un foncteur asphérique, c'est-à-dire si pour
tout objet $b$ de $B$, la catégorie~$\tranche{A}{b}$ a le type
d'homotopie du point. 

Dans le cas des préfaisceaux en groupes abéliens, on introduit la
notion de \emph{foncteur asphérique en homologie} : on dit qu'un
foncteur entre petites catégories $u : A \to B$ est asphérique en
homologie si pour tout objet $b$ de~$B$, la catégorie~$\tranche{A}{b}$
a le type d'homologie du point. On montre, en utilisant des notions
élémentaires de théorie des dérivateurs, que si $u$ est asphérique en
homologie, alors on obtient un diagramme commutatif à isomorphisme
naturel près 
\[
\xymatrix{
{\Hotab_B} \ar[rr]^{u^*}="a" \ar[rd]_{\Hf{B}} \ar@{}"a";[rd]|{\simeq} && {\Hotab_A}
\ar[ld]^{\Hf{A}} \\
& {\Hotab} & \pbox{.}
}  
\]
Mais, sans relations de compatibilité entre le choix des intégrateurs
sur $A$ et sur $B$, on n'obtient évidemment pas de transformation
naturelle explicite s'insérant dans le diagramme 
\[\label{diag:embarrasduchoix}
\tag{2}
\xymatrix{
{\prefab{B}} 
\ar[rr]^{u^*} \ar[rd]_{{L_B}_!}^{}="c"
&& {\prefab{A}} 
\ar@{}"c"|(.4){}="d"|(.9){}="e"
\ar@{==>}"d";"e"_{\lambda_u}
\ar[ld]^{{L_A}_!} \\
& {\Ch(\Ab)}
} 
\]
et à laquelle on pourrait demander d'être un quasi-isomorphisme argument
par argument.

Une stratégie que l'on explore pour exhiber une telle transformation
naturelle est la suivante. Si $u : A \to B$ est un foncteur, on
dispose d'un couple de foncteurs adjoints 
\[
\xymatrix{
\Hom(A,\Ab) 
\ar@<1ex>[r]^{({u}^{\op})_!^\ab}="a"
& 
\Hom(B,\Ab)
\ar@<1ex>[l]^{(u^{\op})^*}="b"
\ar@{}"a";"b"|-{{\hdash}}
}
\]
et on
peut alors chercher des conditions pour que ces foncteurs préservent les
intégrateurs (c'est-à-dire, les résolutions projectives du
foncteur constant de valeur~$\Z$) pour essayer, dans ce cas, de décrire
une transformation naturelle explicite $\lambda_u$ comme ci-dessus.

On montre alors que si $u$ est asphérique en homologie et si $L_A : A
\to \Ch(\Ab)$ est un intégrateur sur $A$, alors le foncteur
$\prescript{A}{}{L}_B$ obtenu en appliquant degré par degré le
foncteur $({u}^{\op})_!^\ab$ à~$L_A$ est un intégrateur
sur~$B$, et qu'on obtient alors un diagramme commutatif à isomorphisme
près 
\[
\xymatrix{
\prefab{B} \ar[rr]^{u^*}="a" \ar[rd]_{{\prescript{A}{}{L}_B}_!}
\ar@{}"a";[rd]|{\simeq}
&&
\prefab{A} \ar[ld]^{{L_A}_!} \\
& \Ch(\Ab) & \pbox{.}
}
\]

Dans l'autre sens, la situation est moins favorable, puisque le
foncteur $u^*$ ne préserve pas les objets projectifs en général. 
On dégage toutefois des conditions pour que ce
soit le cas, qui s'appliquent bien à notre cadre. En particulier, on
montre que si $u$ est une fibration à fibres discrètes, alors le
foncteur~$L_A^B=u^*(L_B)$ est bien un intégrateur sur $A$. On montre alors que dans ce
cas, on obtient une transformation naturelle explicite s'insérant dans le
diagramme
\[
\xymatrix{
{\prefab{B}} 
\ar[rr]^{u^*} \ar[rd]_{{L_B}_!}^{}="c"
&& {\prefab{A}} 
\ar@{}"c"|(.4){}="d"|(.9){}="e"
\ar@2"d";"e"_{\lambda_u}
\ar[ld]^{{L^B_A}_!} \\
& {\Ch(\Ab)} & \pbox{,}
} 
\]
et que $u$ est un morphisme asphérique en homologie si et seulement si
$\lambda_u$ est un quasi-isomorphisme argument par argument. 

\medskip

En particulier, pour toute catégorie $B$ et pour tout préfaisceau~$F$
sur $B$, on peut définir de cette manière un intégrateur sur la
catégorie $\tranche{B}{F}$. Par exemple, si $F$ est un ensemble
simplicial, le procédé décrit ci-dessus montre qu'on peut calculer
l'homologie des préfaisceaux en groupes abéliens sur
$\tranche{\Delta}{F}$ comme l'homologie du complexe 
\[
\cdots \leftarrow \bigoplus_{\Delta_n \xrightarrow{u} F} X(\Delta_n,u)
\leftarrow \bigoplus_{\Delta_{n+1}\xrightarrow{u} F}
X(\Delta_{n+1},u) \leftarrow \cdots \pbox{.}
\]

De plus, on peut montrer grâce à cela que si $F$ est un ensemble
simplicial ayant le type d'homotopie du point, alors la catégorie
$\tranche{\Delta}{F}$ est également une catégorie pseudo-test
homologique. En d'autres termes, il existe des \emph{correspondances
de Dold-Kan locales} (toujours entendues comme des équivalences entre
les catégories dérivées). 

Ce constat nous amène alors à reproduire une partie des
schémas de la théorie des catégories test.

\section*{Catégories test homologiques}

Si $L_A : A \to \Ch(\Ab)$ est un foncteur, son extension de Kan à
gauche 
\[
{L_A}_! : \prefab{A}\to\Ch(\Ab) 
\]
possède un adjoint à droite que l'on note 
\[
L_A^* : \Ch(\Ab) \to \prefab{A} \mdvirg C \mapsto
(\Hom_{\Ch(\Ab})(L(-),C)) \pbox{.}
\]

Par exemple, en notant $ c : \Delta\to\Ch(\Ab)$ le foncteur envoyant
les objets de~$\Delta$ sur leur complexe normalisé, le théorème de
Dold-Kan affirme en fait que le couple de foncteurs $c_! \dashv c^*$ est une
équivalence de catégories. En particulier, le foncteur $c^*$ préserve
les équivalences faibles et induit un quasi-inverse au foncteur
$\Hf{\Delta}$.

On dit alors qu'un intégrateur $L_A$ sur $A$ est un \emph{intégrateur test
faible} si le foncteur $L_A^*$ envoie les quasi-isomorphismes sur des
équivalences test abéliennes de $A$, et si le foncteur induit entre
les catégories homotopiques 
\[
L_A^* : \Hotab \to \Hotab_A 
\]
est un quasi-inverse au foncteur $\Hf{A} : \Hotab_A \to \Hotab$.

On dit que $L_A$ est un \emph{intégrateur test local} si pour tout
objet $a$ de $A$, l'intégrateur
\[
  L_{\tranche{A}{a}}: \tranche{A}{a}\to A \xrightarrow{L_A} \Ch(\Ab)
  \]
évoqué dans la
section précédente, où la flèche de gauche désigne le foncteur
d'oubli, est un intégrateur test faible. Enfin, on dit que
$L_A$ est un \emph{intégrateur test} si $L_A$ est à la fois
un intégrateur test faible et test local.

\begin{defin*}
Une petite catégorie $A$ est une \emph{catégorie test
homologique faible} (resp. \emph{test homologique}) si il existe un
intégrateur test faible (resp. test) sur $A$.
\end{defin*}

Les exemples de catégories test de Whitehead énoncés précédemment sont
des catégories test homologiques. On montre aussi que la
sous-catégorie $\Delta'$ des monomorphismes de $\Delta$ est une
catégorie test homologique faible, ainsi que la catégorie des cubes
sans connexion $\cub$. En revanche, $\Delta'$ n'est pas une catégorie
test homologique, et on conjecture qu'il en est de même pour la
catégorie $\cub$.  On montre aussi que le produit d'une catégorie test
homologique faible avec une catégorie ayant le type d'homologie du
point est encore une catégorie test homologique faible, et que le
produit d'une catégorie test homologique locale avec une petite
catégorie quelconque est encore une catégorie test homologique locale.

Bien que nous ne disposions pas de caractérisations aussi simples que
pour les catégories test, il semble prometteur d'essayer de pousser
plus loin ces notions, et nous proposons à la lectrice et au lecteur
de les traiter comme des germes prometteurs. On montre toutefois le
résultat suivant.

\begin{theorem*}
Soit $A$ une petite catégorie. On suppose que le foncteur diagonal $A
\to A \times A$ est asphérique en homologie. Alors les conditions
suivantes sont équivalentes : \begin{enumerate}
\item $A$ est une catégorie test homologique faible;
\item $A$ est une catégorie test homologique.
\end{enumerate}
\end{theorem*}

On dit alors qu'une catégorie $A$ est une \emph{catégorie test
homologique stricte} si c'est une catégorie test homologique faible,
et si le foncteur $A \to A \times A$ est asphérique en homologie.

Dans le cas où $A$ est une catégorie test locale de
Whitehead, on verra que les notions de catégorie pseudo-test et de catégorie
test homologique faible sont équivalentes. On obtient finalement le
résultat suivant.
\begin{theorem*}
Toute catégorie test stricte de Whitehead est une catégorie test
homologique stricte.
\end{theorem*}
Ce théorème s'applique en particulier à la catégorie $\Delta$, la
catégorie cubique avec connexions $\cubc$, la catégorie $\Theta$ de
Joyal ainsi que les catégories $\Theta_n$ pour~$n\geq 1$, et les
catégories $\Xi_n$ pour tout entier~$n\geq1$. Nous ne connaissons pas
d'exemples de catégories test homologiques strictes qui ne soient pas des
catégorie de Whitehead.
\medskip

\section*{Les théorèmes de type Dold-Kan stricts}

Il existe plusieurs généralisations du théorème de Dold-Kan \og strict
\fg{} (par opposition à \og homotopique \fg{}), et nous n'essayons pas,
dans ce travail, d'en exhiber de nouvelles. On peut toutefois en citer
quelques unes. 

Le théorème de Dold-Puppe \cite{doldpuppe1961} affirme que pour toute
catégorie abélienne~$\A$, la catégorie des préfaisceaux simpliciaux
dans $\A$ est équivalente à celle des complexes de chaînes d'objets de
$\A$.

D'autres travaux établissent des équivalences de catégories de la
forme $\prefab{A}\simeq\Ch(\Ab)$ où $A$ est une petite catégorie,
notamment ceux de Brown et Higgins \cite{brown2003cubical} dans le cas
où $A$ est la catégorie cubique avec connexions. Bourn établit dans
\cite{bourn1990denormalization} une équivalence de catégories entre
la catégorie des $\omega$-groupoïdes stricts en groupes abéliens et
celle des complexes de chaînes, ce qui implique que la catégorie des
préfaisceaux en groupes abéliens sur la
catégorie~$\Grefl$ des globes réflexifs est aussi équivalente à la
catégorie des complexes de chaînes de groupes abéliens. On montrera
que ces deux exemples donnent bien lieu à des correspondances de
Dold-Kan homotopiques, en montrant que les foncteurs en jeu calculent
bien l'homologie au sens où on l'étudie dans cette thèse.

Une autre approche consiste à exhiber des équivalences de
catégories entre des catégories de préfaisceaux en groupes abéliens et
des généralisations des complexes de chaînes, comme le théorème de
Pirashvili~\cite{pirashvili2000dold} pour les préfaisceaux en groupes
abéliens sur la catégorie~$\Gamma$ de Segal. Lack et Street donnent dans
\cite{lack2014combinatorial} une grande classe d'équivalences de
catégories de la forme~$\Homi(A, P) \simeq \Homi_0(D,P)$ où $A$ est une
petite catégorie,~$P$ est une catégorie additive, $D$ est une
catégorie avec des morphismes nuls construite à partie de~$A$,
et~$\Homi_0(D,P)$ désigne la sous-catégorie des foncteurs~$D \to P$
préservant les morphismes nuls. On peut finalement citer le résultat de
\hbox{Gutiérrez}, Lukacs et Weiss \cite{gutierrez2011dendroidal} affirmant
que la catégorie des groupes abéliens dendroïdaux est équivalente à
celle des \og complexes
dendroïdaux \fg{}.

\section*{Plan de la thèse}

Le chapitre $1$ est principalement un chapitre préliminaire. Dans la
première section, on rappelle les points clés de la théorie des
catégories test en se basant sur le livre de Maltsiniotis
\cite{maltsiniotis2005}. La deuxième section est dédiée à l'homologie
des groupes abéliens simpliciaux. On tente de condenser tous les
résultats importants, et on énonce le théorème de Dold-Kan. Dans la
troisième section, on montre le théorème que Grothendieck appelle
\emph{théorème de Whitehead} pour la catégorie $\Delta$, et on
détaille le foncteur d'abélianisation des types d'homotopie décrit
dans cette introduction utilisant
l'intermédiaire de la catégorie $\Delta$. Dans la quatrième section,
on introduit le foncteur $\Hf{A} :
\prefab{A}\to\Hotab$ pour toute petite catégorie $A$, en montrant
comment il généralise l'homologie des groupes. On montre aussi comment
ce foncteur permet de définir un foncteur canonique associant à toute
petite catégorie son type d'homologie, sans utiliser l'intermédiaire
de la catégorie $\Delta$. Dans la cinquième section, on montre que le
foncteur~$\Hf{A}$ se comporte bien au sens de la théorie des
catégories test, c'est-à-dire qu'il associe bien à tout préfaisceau
d'ensembles l'homologie du nerf de sa catégorie des éléments. Dans la
sixième section, on introduit la catégorie $\Hotab_A$, et on définit
la notion de catégorie pseudo-test homologique.

Le chapitre $2$ est consacré à l'étude générale du foncteur $\Hf{A}$
défini au chapitre $1$. La première section est une section
préliminaire consistant en des rappels
d'algèbre homologique de base. Dans la deuxième section, on montre comment
calculer l'extension de Kan à gauche des foncteurs $L :A \to M$ où~$M$
est une catégorie additive cocomplète, afin d'obtenir un foncteur que
l'on notera $L_! : \prefab{A} \to M$ étendant le foncteur $L$. En
particulier, on montre qu'on peut calculer
cette extension de Kan en passant par l'intermédiaire de la catégorie
additive librement engendrée par $A$, appelée \emph{enveloppe
additive} de $A$. On montre aussi que l'extension de Kan $L_!$ d'un
foncteur~$L : A \to M$ admet toujours un adjoint à
droite, noté $L^*$. Dans la troisième section, on introduit le produit
tensoriel de foncteurs tel que Grothendieck le décrit dans
\emph{Pursuing Stacks}, qui nous permettra de dériver des foncteurs
$\prefab{A}\to\Ab$ dans la suite. Dans la quatrième section, on décrit
les objets projectifs de la catégorie~$\prefab{A}$ pour toute petite
catégorie $A$, et on montre
comment dériver les foncteurs~$\prefab{A}\to\Ab$ en utilisant le
produit tensoriel de foncteurs introduit dans la section précédente,
en utilisant des résolutions libres. Dans la cinquième section, on
introduit la notion d'intégrateur, permettant ainsi une description
du foncteur $\Hf{A} : \prefab{A}\to\Hotab$ utilisant le produit
tensoriel de foncteurs. Dans la sixième section, on
montre comment le formalisme des intégrateurs permet de calculer
l'homologie des catégories à coefficients dans $\Z$ d'une manière
combinatoire. On montre par exemple que la catégorie $\G$ des globes a
tous ses groupes d'homologie isomorphes à $\Z$, et en particulier
qu'elle n'est pas asphérique. Dans la septième section, on montre que la
donnée d'un intégrateur sur une petite catégorie $A$ permet aussi
calculer la limite inductive homotopique des foncteurs $A^{\op} \to
\Ch(\Ab)$. La huitième section est plus anecdotique : on décrit le
processus par lequel Grothendieck parvient, dans \emph{Pursuing
Stacks} à la notion d'intégrateurs, en cherchant à calculer la
cohomologie des préfaisceaux sur la catégorie opposée. On montre
toutefois qu'un intégrateur sur $A$ permet aussi de calculer la
cohomologie de la catégorie ${A}^{\op}$ à coefficients dans un
préfaisceau.

Le chapitre $3$ est consacré à la notion de morphisme asphérique en
homologie, c'est-à-dire des foncteurs commutant au foncteur $\Hf{A}$.
Dans la première section, on donne quelques rappels sur les
localisateurs fondamentaux, en se basant toujours sur le livre de
Maltsiniotis \cite{maltsiniotis2005}. Dans la deuxième section, on
introduit la classe $\W_\infty^\ab$ des foncteurs entre petites catégories
induisant un isomorphisme entre les groupes d'homologie, et on introduit
la notion de morphisme asphérique en homologie. Dans la troisième
section, on s'intéresse à la possibilité de transférer un intégrateur
le long d'un morphisme asphérique en homologie, permettant, étant
donné un foncteur asphérique en homologie~$u : A \to B$, d'obtenir une
transformation naturelle explicite entre l'homologie de $B$ à
coefficients dans $X$ et l'homologie de $A$ à coefficients dans le
préfaisceau $u^*X$. Dans la quatrième section, on introduit
les catégories de Whitehead, et on montre que $\Delta$ est une
catégorie de Whitehead. De plus, on montre que si $A$ est une catégorie
de Whitehead, et s'il existe un foncteur asphérique $A \to B$, alors
$B$ est également une catégorie de Whitehead. On utilise aussi la
notion de foncteur $A \to \pref{B}$ asphérique, et on montre que si un
tel foncteur existe, et si $A$ est une catégorie de Whitehead, alors
$B$ est également une catégorie de Whitehead. Dans la cinquième section,
on montre que pour toute catégorie test locale de Whitehead, la
catégorie $\prefab{A}$ peut être munie d'une structure de catégorie de
modèles dont les équivalences faibles sont les éléments de
$\W_A^{\ab}$ (ou donc de $\U^{-1}\W_A$), et que si $u : A \to B$ est un foncteur asphérique,
et si $A$ et $B$ sont des catégories test locales, alors le couple de
foncteurs adjoints~$u^* \dashv u_*^\ab$ est une équivalence de Quillen. On
montre aussi que pour toute catégorie test, il existe un zig-zag de
morphismes asphériques $A \rightarrow B \leftarrow \Delta$ où $B$ est une
catégorie test de Whitehead, permettant ainsi de démontrer que toute
catégorie test de Whitehead est une catégorie pseudo-test homologique.

Le chapitre $4$ introduit les notions de catégories test homologiques,
ainsi que ses variantes faibles, locales et strictes. Dans la première
section, on montre que toute catégorie peut être munie d'un
intégrateur canonique, et on trouve de cette manière une nouvelle
preuve de la formule de Bousfield-Kan pour les limites inductives
homotopiques. Dans la deuxième section, on introduit la notion de
catégorie test homologique faible, qui désigne les catégories
pseudo\nobreakdash-test homologiques dont un quasi-inverse au foncteur $\Hf{A}$
est induit par l'adjoint à droite de l'extension de Kan d'un intégrateur.
On montre que
toute catégorie test homologique faible a nécessairement le type
d'homologie du point, et on exhibe quelques caractérisation des
catégories test homologiques faibles, sans parvenir à une
caractérisation aussi simple que pour les catégories test faibles de
Grothendieck. Dans le cas où $A$ est une catégorie de
Whitehead, on montre que toute catégorie test est une catégorie test
homologique faible, et on montre également que le produit d'une
catégorie test homologique faible et d'une catégorie ayant l'homologie
du point est une catégorie test homologique faible. Dans la troisième
section, on montre que la donnée d'un intégrateur sur une catégorie
$A$ permet de définir un intégrateur sur les tranches $\tranche{A}{X}$
pour tout préfaisceau $X$ sur $A$. On introduit alors la
notion de catégorie test homologique locale, et de catégorie test
homologique. En revanche, on ne dispose pas pour l'instant de
caractérisation des catégories localement de Whitehead. On montre
aussi que la catégorie $\tranche{\Delta}{\nerf\Grefl}$, où $\Grefl$
désigne la catégorie des globes réflexifs, est une catégorie test,
test homologique, mais n'est pas une catégorie de Whitehead. Dans la
quatrième section, on introduit la notion de catégorie test
homologique stricte, qui fournit un grand nombre d'exemples de
catégories test homologiques.

Le chapitre $5$ a pour but de fournir des exemples de
théorèmes de Dold-Kan homotopiques obtenus grâce aux résultats de
cette thèse. Dans la première section, on donne une stratégie générale
pour construire des intégrateurs libres sur des catégories vérifiant
des axiomes permettant de définir une bonne notion de dimension. On
introduit la notion d'orientation et d'orientation asphérique, dont la
donnée permet de définir un intégrateur libre. Dans la deuxième
section, on montre que pour tout entier $n>0$, la
catégorie $\Delta^n$ est bien une catégorie test
homologique stricte, résultat du folklore que nous redémontrons avec
les outils élaborés dans cette thèse. Dans la troisième section, on montre que la
sous-catégorie $\Delta'$ des monomorphismes de $\Delta$ est une
catégorie test homologique faible, mais n'est pas une catégorie test
homologique locale. Les méthodes développées dans cette section
devraient pouvoir s'adapter à tous les cas dans lesquels on peut
définir un complexe normalisé comme quotient d'un complexe non
normalisé par les dégénérescences sans changer son homologie. Dans la
quatrième section, on évoque la situation des ensembles cubiques (en
tous cas, de trois des significations possibles de ce terme). On
montre que la catégorie cubique est une catégorie test homologique, et
que
la catégorie cubique avec connexions négatives (pour laquelle Brown et
Higgins ont prouvé dans \cite{brown2003cubical} un théorème de
Dold-Kan strict) et la catégorie cubique avec connexions positives et
négatives sont des catégories test homologiques strictes. On montre
aussi que la catégorie cubique avec connexions positives est une
catégorie de Whitehead. Dans la
cinquième section on traite l'exemple de la catégorie des globes et de
la catégorie des globes réflexifs. Cette dernière constitue un autre
exemple de catégorie munie d'un théorème de Dold-Kan strict,
établissant l'équivalence entre la catégorie des préfaisceaux abéliens
et celle des complexes de chaînes de groupes abéliens, qui est liée au
théorème de Bourn \cite{bourn1990denormalization} sur les
$\omega$-catégories strictes.
On montre que $\Grefl$, bien que n'étant pas une catégorie test, est
une catégorie test homologique faible. On montre aussi qu'elle n'est pas une
catégorie de Whitehead. On profite de cette section pour définir la
notion de $\omega$-catégorie stricte, qui nous servira pour les deux
derniers exemples. La sixième section est consacrée à la catégorie
$\Theta$ de Joyal. On montre que pour tout entier $n$, la
catégorie~$\Theta_n$ est une catégorie test homologique stricte de
Whitehead, ainsi que la catégorie $\Theta$. Dans la septième section,
on introduit la catégorie $\Xi_n$, obtenue comme l'image du foncteur
canonique~$\Delta^n \to \Theta_n$. On montre que $\Xi_n$ est une
catégorie test stricte, ainsi qu'une catégorie test homologique
stricte de Whitehead. On termine par la construction d'un intégrateur
explicite sur $\Xi_n$, pour tout entier $n\geq1$, utilisant la notion
de dimension introduite dans la section $\ref{secOrientation}$.

Quelques résultats dans cette thèse utilisent des éléments basiques
de théorie des dérivateurs, que l'on a regroupés dans l'annexe
\ref{annexeDerivateurs}. Les résultats en question sont précisément
indiqués au début de cette annexe. On y démontre également, en utilisant ce
formalisme, que les morphismes asphériques en homologie correspondent
aux foncteurs $u : A \to B$ tels que le foncteur $u^* :
\Homi({B}^{\op},\Ch(\Ab)) \to \Homi({A}^{\op},\Ch(\Ab))$ commute aux
limites inductives homotopiques, et que la classe $\W_\infty^\ab$ est
un localisateur fondamental. 

\section*{Notations}

Si $\C$ est une catégorie, on note $\Ob(\C)$ la classe des
objets de~$\C$ et~$\Fl(\C)$ la classe des morphismes de $\C$. Si $c$
et $c'$ sont deux objets de $\C$, on note~$\Hom_\C(c,c')$ l'ensemble
des morphismes de $c$ dans $c'$. On note ${\C}^{\op}$ la catégorie
opposée à $\C$. Si $C$ est une catégorie munie d'un objet final, on
le notera généralement $e_C$.

Étant données une petite catégorie $\C$ et une catégorie
$\mathcal{D}$, on note $\Homi(\C,\mathcal{D})$ la catégorie des
foncteurs de $\C$ vers $\mathcal{D}$, dont les morphismes sont les
transformations naturelles. On note 
\begin{align*}
\limind\nolimits_C^D : \Homi(\C,\mathcal{D}) \to \mathcal{D}
\hspace{1em}
\text{ et }  
\hspace{1em}
\limproj\nolimits_C^D : \Homi(\C,\mathcal{D}) \to \mathcal{D} 
\end{align*}
les foncteurs de limite inductive et projective, lorsqu'ils existent.

On note $\Cat$ la catégorie des petites catégories, dont on note $e$
l'objet terminal. On note $\Ens$ la catégorie des ensembles. 
Si $A$ est une petite catégorie, on note 
\[
\pref{A} = \Homi({A}^{\op}, \Ens) 
\]
la catégorie des préfaisceaux d'ensembles sur $A$, que l'on appellera
simplement préfaisceaux. Si $a$ est un objet de $A$, on notera
simplement $a$ le préfaisceau sur~$A$ représenté par $a$. 

Si $u : A \to B$ est un foncteur entre petites catégories, on note 
\[
u^* : \pref{B} \to \pref{A} 
\]
\notindex{$u^*$}%
le foncteur de précomposition par $u$. Ce foncteur admet un adjoint à
gauche noté $u_!$ ainsi qu'un adjoint à droite noté $u_*$. 
\notindex{$u_\bang$}%
\notindex{$u_*$}%

Si $\M$ est une catégorie et $\W$ est une classe de flèches de
$\M$, on note
\[
 \M \xrightarrow{\gamma} \W^{-1}\M 
\]
la localisation de Gabriel-Zisman de $(\M,\W)$, et on n'évoquera pas
les difficultés ensemblistes liées à sa définition, qu'on ne
rencontrera pas dans ce travail. Cette catégorie est caractérisée par la
propriété universelle suivante : pour toute catégorie $\C$, on a un
isomorphisme 
\[
\Homi(\W^{-1}\M,\C) \simeq \Homi^\W(\M,\C)
\]
où on a noté $\Homi^\W(\M,\C)$ la sous-catégorie pleine des foncteurs
qui envoient les éléments de $\W$ sur des isomorphismes de $\C$.
On dit que $\W^{-1}\M$ est la \emph{catégorie localisée} de $\M$ par
la classe de flèches $\W$. Le foncteur $\gamma : \M \to \W^{-1}\M$
est appelé \emph{foncteur de localisation} ou \emph{projection
canonique}.

\medskip

On note $\Ab$ la catégorie des groupes abéliens, et
\[
\prefab{A} = \Homi({A}^{\op}, \Ab) 
\]
la catégorie des préfaisceaux en groupes abéliens sur $A$.

Si $\mathcal{A}$ est une catégorie abélienne, on note
$\Chu(\mathcal{A})$ la catégorie des complexes de chaînes non bornés
\[
\cdots \leftarrow C_{n-1} \xleftarrow{d_{n}} C_{n}
\xleftarrow{d_{n+1}} C_{n+1} \leftarrow \cdots
\]
sur $\mathcal{A}$, que l'on écrira avec une différentielle de
degré~$-1$. On note $\Ch(\mathcal{A})$ la catégorie des complexes de
chaînes concentrés en degré positif, et $\Chm(\mathcal{A})$ la
catégorie des complexes concentrés en degré négatif, que l'on
représentera parfois comme des complexes de cochaînes 
\[
C_0 \xrightarrow{d_0} C_1 \xrightarrow{d_1} C_2 \to \cdots 
\]
concentrés en degré positif avec une différentielle de degré $1$.

On note $D(\A)$ la catégorie dérivée des complexes de chaînes d'objets
de $\A$, et on notera également $D_+(\A)$ (resp. $D^+(\A)$) la
catégorie dérivée concentrée en degré positif (resp. négatif) d'objets
de $\A$.

\chapter{Abélianisation des types d'homotopie}
\label{chapAbelianisation}
\section{Catégories test}
\label{secCatTest}

Cette section (mise à part les préliminaires topologiques et simpliciaux) est
entièrement basée sur le livre de Maltsiniotis
\cite{maltsiniotis2005}.

\paragr On note $\Top$ la catégorie des espaces topologiques et des
applications continues. On dit qu'un morphisme $f : X \to Y$ de $\Top$ est une
\ndef[équivalence faible!d'espaces topologiques]{équivalence faible d'homotopie} si le morphisme $\pi_0(f) :
\pi_0(X) \to \pi_0(Y)$ est une bijection, et si pour tout entier $n>0$
et tout élément $x\in X$, le morphisme
  \[
\pi_n(f,x) : \pi_n(X,x) \to \pi_n(Y, f(x))
  \]
est un isomorphisme.

On note $\W_\Top \subset \Fl(\Top)$ la classe des équivalences faibles
d'homotopie d'espaces topologiques. La catégorie des \ndef[type
d'homotopie]{types
d'homotopie} est la catégorie 
\[
\Hot_\Top = (\W_\Top^{-1}) \Top
\]
\notindex{$\Hot_\Top$}%
\notindex{$\W^{-1}\M$}%
obtenue en \ndef[localisation]{localisant} la
catégorie $\Top$ par les équivalences faibles d'homotopie.
Cela signifie qu'on dispose d'un foncteur noté 
\[
\gamma : \Top \to \Hot_\Top 
\]
envoyant les éléments de $\W_\Top$ sur des isomorphismes, et qui est
universel pour cette propriété, c'est-à-dire que pour toute catégorie
$\M$ et pour tout foncteur $F : \Top \to \M$ envoyant les équivalences
faibles d'homotopie sur des isomorphismes, il existe un unique
foncteur $G : \Hot_\Top \to \M$ faisant commuter le diagramme suivant
:
\[
\xymatrix{
\Top \ar[r]^{F} \ar[d]_{\gamma} & \M \\
\Hot_\Top \ar[ru]_{G} & \pbox{.}
} 
\]

\paragr\label{defDelta} On note $\Delta$ la sous-catégorie pleine de $\Cat$ dont les
\notindex{$\Delta$}%
\notindex{$\Delta_n$}%
objets sont les ensembles ordonnés
\[
\Delta_n =  \lbrace 0 , 1, \dots , n \rbrace \mdvirg n \in \N \pbox{.}
\]
La catégorie des \ndef[ensemble!simplicial]{ensembles simpliciaux}
désigne la catégorie $\pref{\Delta}$ des préfaisceaux d'ensembles sur
$\Delta$.

Pour $n\geq 1$ et $0 \leq i \leq n$, on note
\[
  \delta_i^n :\Delta_{n-1} \to \Delta_n
  \]
l'unique application injective n'atteignant pas la valeur $i$. On note
également, pour~$n \geq 0$ et~$0 \leq i \leq
{n}$,~
\[
  \sigma_i^n : \Delta_{n+1} \to \Delta_n
  \]
l'unique application surjective atteignant deux fois la valeur $i$. 

Si $X$ est un ensemble simplicial, on note simplement $X_n$ l'ensemble
$X(\Delta_n)$ des \emph{$n$-simplexes} de $X$ pour $n\geq 0$. On notera aussi
$d^n_i=X(\delta^n_i):X_n\to X_{n-1}$ et $s^n_i=X(\sigma^n_i):X_n\to
X_{n+1}$ les applications de \emph{faces} et de \emph{dégénérescences}.

\paragr On note $|-| : \EnsSimp \to \Top$ le foncteur de réalisation
géométrique. On appelle \emph{équivalences
faibles simpliciales} les morphismes d'ensembles simpliciaux dont
l'image par le foncteur de réalisation géométrique est une équivalence
faible d'espaces topologiques. On note $\W_\Delta$ la classe des équivalences faibles simpliciales de
$\EnsSimp$. On rappelle (voir \cite{milnor1957realization}) que le
foncteur de réalisation géométrique induit une équivalence de
catégories 
\[
\Hot_\Delta \xrightarrow{\simeq} \Hot_\Top 
\]
entre la catégorie $\Hot_\Delta$ obtenue en localisant la catégorie
$\EnsSimp$ par les équivalences faibles simpliciales, et la catégorie
$\Hot_\Top$. Un quasi-inverse est donné par le foncteur induit entre
les catégories localisée par le foncteur complexe singulier 
\begin{align*}
\mathsf{Sing} : \Top &\to \pref{\Delta} \\
T &\mapsto (\Delta_n \mapsto \Hom_\Top(|\Delta_n|,T)) \pbox{.}
\end{align*}
\notindex{$\mathsf{Sing}$}%
On renvoie à \cite{goerssjardine} pour plus de détails sur la théorie
de l'homotopie des ensembles simpliciaux.

\paragr\label{def:eqThomason} On peut montrer que l'inclusion $\Delta \hookrightarrow \Cat$
induit un foncteur pleinement fidèle 
\[
\nerf : \Cat \to \EnsSimp 
\]
appelé \ndef[nerf!simplicial]{nerf simplicial}, défini pour toute petite catégorie $C$
par 
\[
\nerf(C)_n = \Hom_\Cat(\Delta_n,C) \pbox{.} 
\]
\notindex{$\nerf(C)$}%
\notindex{$\W_\infty$}%
On note $\W_{\infty}$ la classe des morphismes de $\Cat$ envoyés par
le nerf sur des équivalences faibles simpliciales, et on appelle
\ndef[équivalence faible!de $\Cat$]{équivalences de Thomason}
\termindex{équivalence faible!de $\Cat$!de Thomason}ses éléments, ou encore
\emph{équivalences faibles} de $\Cat$.

\paragr On note $\Hot$ la catégorie obtenue en localisant la catégorie
$\Cat$ par les équivalences de Thomason 
\notindex{$\Hot$}%
\[
\Hot := \W_{\infty}^{-1} \Cat \pbox{.}
\]

Cette définition est justifiée par un résultat de Illusie
\cite[chapitre VI, théorème 3.3]{illusie1971cotangent}, qui affirme
que le foncteur nerf simplicial induit une équivalence de catégories 
\[
\Hot \xrightarrow{\nerf} \Hot_\Delta \mdvirg
\]
dont un quasi-inverse est donné par le foncteur associant à tout
ensemble simplicial sa \emph{catégorie des éléments}, que l'on va
définir au paragraphe \ref{defCatElements}. Dans le cadre de la
théorie de l'homotopie de Grothendieck, le couple $(\Cat,\W_\infty)$
est le \emph{modèle des types d'homotopie} privilégié.

\medskip
La motivation de Grothendieck, qui l'amène à introduire la notion de
catégorie test (voir
\cite{maltsiniotisPS} et~\cite{maltsiniotis2005}) est de caractériser
toutes les petites catégories $A$ telles que la catégorie~$\pref{A}$
modélise les types d'homotopie, généralisant ainsi la théorie de
l'homotopie simpliciale. Nous allons détailler ce que cela signifie
dans ce qui suit.

\paragr Soient $A$ et $B$ deux catégories. Si $u : A \to B$ est un
foncteur et $b$ est un objet de $B$, on note $\tranche{A}{b}$ la
\notindex{$\tranche{A}{b}$}%
catégorie dont les objets sont les couples
  \[
(a,f : u(a) \to b) \mdvirg a \in \Ob(A) \mdvirg f \in \Fl(B)
  \]
et dont les morphismes $(a,f) \to (a',f')$ sont les morphismes
$\varphi : a \to a'$ de $A$ faisant commuter le diagramme 
\[
\xymatrix{
u(a) \ar[rd]_{f} \ar[rr]^{u(\varphi)}  && u(a') \ar[ld]^{f'} 
\\ & b
}
\]
dans $B$.

Le résultat suivant constitue une propriété fondamentale de la
classe $\W_\infty$ des équivalences de Thomason, sur laquelle nous
reviendrons dans la section~\ref{secLocalisateursFondamentaux}. On
rappelle qu'on note $e$ l'objet final de $\Cat$.

\begin{theorem}[Quillen, théorème A]\label{theoremeA}
Soit $u : A \to B$ un morphisme de $\Cat$ tel que pour tout objet $b$
de $B$, le morphisme canonique $\tranche{A}{b} \to e$ soit dans
$\W_\infty$. Alors $u$ est un élément de $\W_\infty$.
\end{theorem}
\begin{proof}
Voir \cite[théorème A]{quillenktheory} pour la preuve originale.
\end{proof}

On introduit alors la terminologie suivante. 

\paragr\label{defAspherique} On dit qu'une
petite catégorie~$A$ est \ndef[catégorie!asphérique]{asphérique} si
l'unique morphisme~$A\to e$ est un élément de $\W_\infty$. On dit
qu'un morphisme $u : A \to B$ de~$\Cat$ est
\ndef[morphisme asphérique!de $\Cat$]{asphérique}\termindex{foncteur
asphérique!de $\Cat$ (voir \emph{morphisme asphérique)}} si pour tout objet
$b$ de $B$,
la catégorie $\tranche{A}{b}$ est une catégorie asphérique. 
Le théorème A de Quillen affirme donc que tout morphisme asphérique
de~$\Cat$ est une équivalence de Thomason.

\begin{prop}\label{propObjFinalOuInitialAspherique}
Toute catégorie admettant un objet final ou initial est asphérique.
\end{prop}
\begin{proof}
Si $A$ possède un objet final $e_A$, le foncteur $A\to e$ vers la
catégorie ponctuelle possède un adjoint à droite $e \to A$, associant à
l'unique objet de $e$ l'objet terminal $e_A$. Le foncteur nerf
simplicial envoie les transformations naturelles sur des homotopies
simpliciales, et les morphismes d'adjonction induisent donc une
équivalence d'homotopie entre les nerfs de $A$ et de $e$. Par
conséquent, le nerf de $A$ est contractile. On procède de même si $A$
possède un objet initial avec l'adjoint à gauche du foncteur $A\to e$.
\end{proof}

\begin{coro}
Tout foncteur admettant un adjoint à droite est asphérique. En
particulier, toute équivalence de catégories est asphérique.
\end{coro}
\begin{proof}
Voir \cite[proposition 1.1.9]{maltsiniotis2005}.
\end{proof}

\paragr\label{defCatElements} Si $A$ est une petite catégorie et $X$ est un préfaisceau sur $A$, on
note~$\tranche{A}{X}$ la \ndef{catégorie des éléments} de $X$, obtenue
en appliquant la construction introduite ci-dessus au plongement de Yoneda $A
\hookrightarrow \pref{A}$. \notindex{$\tranche{A}{X}$}
Concrètement, les objets de $\tranche{A}{X}$ sont les couples 
\[
(a,x : a \to X) \mdvirg a \in \Ob(A) \mdvirg x \in Xa  
\]
et les morphismes $(a,x) \to (a',x')$ sont les morphismes $\varphi : a
\to a'$ de $A$ faisant commuter le diagramme 
\[
\xymatrix{
{a} \ar[rr]^{\varphi} \ar[rd]_{x} && {a'} \ar[ld]^{x'} \\
& {X}
}  
\]
dans $\pref{A}$, c'est-à-dire tels que $x=X(\varphi)x'$.

\paragr\label{def:eqTest} Pour toute petite catégorie $A$, on note
\[
i_A : \pref{A} \to \Cat \mdvirg X \mapsto \tranche{A}{X}
\]
\notindex{$i_A$}%
le foncteur qui associe à un préfaisceau $X$ sa catégorie des éléments.
On note alors
\[
  \W_A \subset \Fl(\pref{A})
\]
\notindex{$\W_A$}%
la classe des morphismes $\varphi : X \to Y$ de $\pref{A}$ dont l'image 
\[
\tranche{A}{\varphi} : \tranche{A}{X} \to \tranche{A}{Y}
\]
par le foncteur
$i_A$ est un élément de~$\W_\infty$, et on appelle 
\ndef[équivalence faible!de préfaisceaux]{équivalences
faibles}, ou \emph{équivalences test} de $\pref{A}$ ses éléments. 
On dit qu'un préfaisceau $X$ sur $A$ est
\ndef[préfaisceau!asphérique]{asphérique} si son image par le
foncteur~$i_A$ est une catégorie asphérique, et on dit qu'un morphisme
de préfaisceaux sur $A$ est \ndef[morphisme asphérique!de
préfaisceaux]{asphérique} si son image par le foncteur $i_A$ est un
morphisme asphérique de $\Cat$.

\paragr\label{defPseudoTest} On note
\[
\Hot_A = \W_A^{-1}\pref{A} 
\]
\notindex{$\Hot_A$}%
la catégorie obtenue en localisant la catégorie $\pref{A}$ par les
équivalences test.
Le foncteur $i_A$ induit donc un foncteur à nouveau noté 
\[
i_A : \Hot_A \to \Hot \mdvirg
\]
et on dit que $A$ est une \ndef[catégorie!pseudo-test]{catégorie
pseudo-test} si ce foncteur est une équivalence de catégories. 

\begin{prop}[Grothendieck]\label{pseudoTestHomologiePoint}
Si $A$ est une catégorie pseudo-test, alors~$A$ est une
catégorie asphérique.
\end{prop}
\begin{proof}
On note $e_{\pref{A}}$ l'objet final de $\pref{A}$. On vérifie qu'on a
un isomorphisme 
  \[
i_A({e_{\pref{A}}}) \simeq A \pbox{.}
  \]
Ainsi, si $i_A$ est une équivalence de catégories, il préserve les
objets finaux, et~$A$ est donc nécessairement l'objet final de $\Hot$.
En vertu de 
\cite[lemme~1.3.6]{maltsiniotis2005}, l'objet final de $\Hot$ est
l'image de la catégorie ponctuelle $e$ par le fonteur de localisation
$\Cat \to \Hot$, et
cela signifie donc que le morphisme canonique~$A\to e$ est un élément de~$\W_{\infty}$.
\end{proof}

\paragr On peut montrer (voir \cite[proposition
1.2.2]{maltsiniotis2005}) que le foncteur $i_A$ admet un foncteur
adjoint à droite 
\[
i_A^* : \Cat \to \pref{A}
\]
\notindex{$i_A^*$}%
défini pour toute petite catégorie $C$ par 
\[
i_A^*(C) : a \mapsto \Hom_\Cat(\tranche{A}{a},C) \pbox{.}
\]

En général, ce foncteur ne préserve pas les équivalences faibles, et
n'induit pas un foncteur entre les catégories localisées. Grothendieck
introduit alors dans \emph{Pursuing Stacks} la terminologie suivante.

\begin{defin}
On dit que $A$ est une \ndef[catégorie!test!faible]{catégorie test faible} si 
\begin{enumerate}
\item le foncteur $i_A^*$ envoie les éléments de $\W_\infty$ sur des
équivalences faibles de~$A$;
\item les foncteurs $i_A$ et $i_A^*$ induisent des équivalences de
catégories 
\[
i_A : \Hot_A \to \Hot \mdvirg i_A^* : \Hot \to \Hot_A 
\]
quasi-inverses l'une de l'autre.
\end{enumerate}
\end{defin}

Les catégories test faibles sont alors caractérisées par la proposition
suivante : 

\begin{prop}[Grothendieck]\label{caracterisationTestFaible}
Soit $A$ une petite catégorie. Les conditions suivantes sont
équivalentes : \begin{enumerate}
\item $A$ est une catégorie test faible;
\item $A$ vérifie les conditions suivantes : \begin{enumerate}
\item $A$ est asphérique;
\item $i_A^*(\W_\infty) \subset \W_A$;
\end{enumerate}
\item pour toute petite catégorie $C$ possédant un objet final, le préfaisceau
$i_A^*(C)$ est asphérique.
\end{enumerate}
\end{prop}
\begin{proof}
Voir \cite[proposition 1.3.9]{maltsiniotis2005}.
\end{proof}

Grothendieck remarque que la propriété d'être une catégorie
test faible n'est pas une propriété locale. Par exemple, la
sous-catégorie $\Delta'$ des monomorphismes de $\Delta$ est une
catégorie test faible, mais pas ses
tranches~$\tranche{\Delta'}{\Delta_n}$ (voir~\cite[proposition~1.7.25
et remarque 1.7.26]{maltsiniotis2005}). Cela motive la définition
suivante.

\begin{defin}\label{defTestLocale}
On dit qu'une catégorie $A$ est une
\ndef[catégorie!test!locale]{catégorie test locale} si pour tout objet
$a$ de $A$, la catégorie $\tranche{A}{a}$ est une catégorie test
faible. On dit que $A$ est une \ndef[catégorie!test]{catégorie test}
si elle est à la fois une catégorie test locale et une catégorie test
faible.
\end{defin} 

\begin{prop}[Grothendieck]
Soit $A$ une petite catégorie. Les conditions suivantes sont
équivalentes : \begin{enumerate}
\item $A$ est une catégorie test;
\item $A$ est une catégorie test locale asphérique.
\end{enumerate}
\end{prop}
\begin{proof}
Voir \cite[remarque 1.5.4]{maltsiniotis2005}.
\end{proof}

\begin{prop}[Grothendieck]\label{prop:testLocalTranche}
Si $A$ est une catégorie test locale et si $F$ est un préfaisceau sur
$A$, alors la catégorie $\tranche{A}{F}$ est une catégorie test
locale. Si de plus $F$ est un préfaisceau asphérique, alors
$\tranche{A}{F}$ est une catégorie test.
\end{prop}
\begin{proof}
Voir \cite[remarque 1.5.4]{maltsiniotis2005}.
\end{proof}

\paragr\label{def:prefLocAspherique} On dit qu'un préfaisceau $F$ sur $A$ est
\ndef[préfaisceau!localement asphérique]{localement asphérique} si le
morphisme $F \to e_{\pref{A}}$ vers le préfaisceau terminal est
asphérique. Puisqu'on a un isomorphisme $i_A(e_{\pref{A}})\simeq A$,
cela revient à dire que le morphisme 
\[
\tranche{A}{F} \to A 
\]
est un morphisme asphérique de $\Cat$.

Grothendieck parvient alors à la caractérisation suivante,
particulièrement simple, des catégories test locales (et donc des
catégories test).
On note $\Delta_1$ la petite catégorie correspondant à l'ensemble
ordonné $\lbrace 0 \leq 1 \rbrace$.

\begin{theorem}[Grothendieck]
Soit $A$ une petite catégorie. Les conditions suivantes sont
équivalentes : \begin{enumerate}
\item $A$ est une catégorie test locale;
\item le préfaisceau $i_A^*(\Delta_1)$ est localement asphérique.
\end{enumerate}
\end{theorem}
\begin{proof}
Voir \cite[théorème 1.5.6]{maltsiniotis2005}.
\end{proof}

La dernière variation, la plus restrictive, autour du concept
de catégorie test est motivée par le fait que le type d'homotopie d'un
produit de préfaisceaux sur une petite catégorie n'est pas
nécessairement le produit des types d'homotopie, autrement dit que le
foncteur 
\[
i_A : \Hot_A \to \Hot 
\]
ne commute pas nécessairement aux produits. 
\begin{prop}[Grothendieck]\label{prop:totalementAspheriqueEquivalences}
Soit $A$ une petite catégorie. Les conditions suivantes sont
équivalentes : \begin{enumerate}
\item le foncteur $i_A : \pref{A}\to\Cat$ commute aux produits
binaires à équivalence faible près, c'est-à-dire que pour tous
préfaisceaux $X$ et $Y$ sur $A$, le morphisme canonique 
\[
\tranche{A}{(X\times Y)}\to \tranche{A}{X}\times \tranche{A}{Y} 
\]
est dans $\W_{\infty}$;
\item pour tous objets $a$ et $b$ de $A$, le produit $a\times b$ dans
$\pref{A}$ est un préfaisceau asphérique;
\item tout préfaisceau asphérique est localement asphérique;
\item le foncteur diagonal $A\to A\times A$ est asphérique.
\end{enumerate}
\end{prop}
\begin{proof}
Voir \cite[proposition 1.6.1]{maltsiniotis2005}.
\end{proof}

\paragr\label{def:totAspherique} On dit que $A$ est \ndef[catégorie!totalement
asphérique]{totalement asphérique} si elle asphérique, et si elle
vérifie une des conditions équivalentes de la proposition précédente.
On vérifie que les catégories totalement asphériques sont exactement
les catégories non vides vérifiant les conditions équivalentes de la
proposition précédente.

\begin{example}\label{exTotAspherique}
Toute catégorie non vide admettant des produits finis est donc
totalement asphérique. La catégorie $\Delta$ est également une
catégorie totalement asphérique (voir \cite[proposition
1.6.14]{maltsiniotis2005}). 
\end{example}

\begin{prop}[Grothendieck]\label{totAspheriqueMorphismeAspherique}
Si $A$ est une catégorie totalement asphérique, et si~$u : A\to B$ est
un foncteur asphérique, alors $B$ est une catégorie totalement
asphérique.
\end{prop}
\begin{proof}
Voir \cite[proposition 1.6.5]{maltsiniotis2005}.
\end{proof}

\paragr On dit que $A$ est une
\ndef[catégorie!test!stricte]{catégorie test stricte} si $A$ est une
catégorie test totalement asphérique.

\begin{prop}[Grothendieck]
Soit $A$ une catégorie totalement asphérique. Les conditions suivantes
sont équivalentes : \begin{enumerate}
\item $A$ est une catégorie test faible;
\item $A$ est une catégorie test locale;
\item $A$ est une catégorie test;
\item $A$ est une catégorie test stricte.
\end{enumerate}
\end{prop}
\begin{proof}
Voir \cite[proposition 1.6.6]{maltsiniotis2005}.
\end{proof}

\begin{prop}[Grothendieck]\label{propFonctAspheriqueSourceTotAsphTest}
Soient $A$ une catégorie totalement asphérique et $B$ une
catégorie test locale. Si $u : A \to B$ est un foncteur asphérique,
alors $A$ et~$B$ sont des catégories test strictes.
\end{prop}
\begin{proof}
Voir \cite[corollaire 1.6.11]{maltsiniotis2005}.
\end{proof}

\begin{example}\label{ex:test}
La catégorie $\Delta$ est une catégorie test stricte. Cela
implique en particulier que pour tout entier $n$, la catégorie
$\Delta^n$ est également une catégorie test stricte. La
sous-catégorie $\Delta'$ des monomorphismes de $\Delta$ est une
catégorie test faible, mais n'est pas une catégorie test locale (voir
\cite[proposition 1.7.25 et remarque~1.7.26]{maltsiniotis2005}).

Contrairement à l'intuition initiale de Grothendieck, tous les
exemples ne sont pas \og géométriques \fg{} : Maltsiniotis prouve par
exemple dans \cite[corollaire~1.6.10]{maltsiniotis2005} que toute petite
sous-catégorie pleine de $\Cat$ stable par produit, ne contenant pas
la catégorie vide, et contenant une catégorie ayant au moins deux
objets distincts, est une catégorie test stricte.

Toutefois, un grand nombre d'exemples géométriques existent. Cisinski
fournit une classe d'exemples de catégories test dans \cite[chapitre
8]{cisinskipref} comprenant entre autres la catégorie cubique (voir la
section \ref{secCubes}). Maltsiniotis prouve
dans~\cite{maltsiniotis2009cubique} que la catégorie cubique avec
connexions (voir la section \ref{secCubes}) est une catégorie test
stricte, ce
qui n'est pas le cas de la catégorie cubique.  
Ara, Cisinski et Moerdijk montrent dans
\cite{AraCisinskiMoerdijkDendroidalTest} que la catégorie dendroïdale $\Omega$
introduite par Moerdijk et Weiss dans \cite{moerdijk2007dendroidal}
est une catégorie test. Maltsiniotis et Cisinski démontrent dans
\cite{cisinski2011theta} que la catégorie $\Theta$ de Joyal (voir la
section \ref{secTheta}) est une catégorie test stricte. 
\end{example}

On introduit finalement une dernière notion liée aux catégories test.
Dans le cas des ensembles simpliciaux, un quasi-inverse au foncteur
$i_\Delta : \Hot_\Delta \to \Hot$ est induit par le foncteur nerf
simplicial (c'est ce que prouve Illusie dans \cite[chapitre VI,
théorème 3.3]{illusie1971cotangent}). 
La notion de \emph{foncteur test} généralise la notion de nerf
simplicial aux autres catégories, et permet de caractériser les
quasi-inverses au foncteur~$i_A : \Hot_A \to \Hot$.

\paragr\label{defFoncteurAspherique} Soit $A$ une petite catégorie et
$i: A \to \Cat$ un foncteur. On note 
\[
i^* : \Cat \to \pref{A} 
\]
le foncteur défini pour toute petite catégorie $C$ par 
\[
i^*(C) : a \mapsto \Homi(i(a),C) \pbox{.}
\]
On dit que~$i$ est un \ndef[foncteur
asphérique!de but $\Cat$]{foncteur asphérique} si les conditions
suivantes sont
satisfaites : \begin{enumerate}
\item pour tout objet $a$ de $A$, la catégorie $i(a)$ est asphérique;
\item pour qu'une petite catégorie $C$ soit asphérique, il faut et il
suffit que le préfaisceau $i^*(C)$ soit asphérique.
\end{enumerate}

\paragr\label{testFaiblei_AAspherique} Par exemple, on peut montrer qu'une petite catégorie $A$ est une
catégorie test faible si et seulement si le foncteur $i_A : A \to
\Cat$ est un foncteur asphérique (voir \cite[remarque
1.7.2]{maltsiniotis2005}).

\begin{prop}[Grothendieck]\label{lemmeFoncteursAspheriquesTriCommutatif}
Soient $A$ et $B$ deux petites catégories,~$u : A \to B$ un foncteur
et $j : B \to \Cat$ un foncteur tel que pour tout objet $b$ de $B$, la
catégorie $j(b)$ soit asphérique. On pose $i = ju : A \to \Cat$. Alors
: 
\begin{enumerate}
\item si $u$ est un morphisme asphérique de $\Cat$, alors $i$ est un
foncteur asphérique si et seulement si $j$ est un foncteur
asphérique;
\item si le foncteur $j$ est pleinement fidèle et si $i$ est un
foncteur asphérique, alors $u$ est un morphisme asphérique de $\Cat$
et $j$ est un foncteur asphérique.
\end{enumerate}
\end{prop}
\begin{proof}
Voir \cite[lemme 1.7.4]{maltsiniotisPS}.
\end{proof}

\begin{prop}[Grothendieck]\label{propEquivalencesFoncteurAspherique}
Soient $A$ une petite catégorie, et~$i : A \to \Cat$ un foncteur tel
que pour tout objet $a$ de $A$, la catégorie $i(a)$ soit asphérique. 
\begin{enumerate}
\item Les conditions suivantes sont équivalents : 
\begin{enumerate}
\item $i$ est un foncteur asphérique;
\item pour toute petite catégorie asphérique $C$, le préfaisceau $i^*(C)$ est
asphérique;
\item $A$ est asphérique, et un morphisme $w$ de $\Cat$ est une
équivalence faible si et seulement si le morphisme $i^*(w)$ est une
équivalence faible de préfaisceaux.
\end{enumerate}
\item Toutes ces conditions impliquent la condition suivante : 
\begin{enumerate}[resume]
\item $i^*(\W_\infty) \subset \W_A$, et le foncteur $i^* : \Hot \to
\Hot_A$ induit par $i^*$ est une équivalence de catégories
quasi-inverse à droite au foncteur $i_A : \Hot_A \to \Hot$.
\end{enumerate}
\item Si pour tout objet $a$ de $A$, la catégorie $i(a)$ admet un
objet final $e_a$, alors les conditions équivalentes $(a)-(c)$ sont encore
équivalentes à la condition : 
\begin{enumerate}[resume]
\item pour toute petite catégorie $C$ admettant un objet final, le
préfaisceau $i^*(C)$ est asphérique.
\end{enumerate}
De plus, si on note 
\[
\alpha : i_Ai^* \to \id_\Cat
\]
le morphisme défini pour toute catégorie $C$ par 
\[
\alpha_C : i_Ai^*(C) \to C \mdvirg (a, i(a) \xrightarrow{v} C) \mapsto
v(e_a) \mdvirg
\]
alors les conditions équivalentes $(a)-(c)$ et $(e)$ sont encore
équivalentes à la proposition 
\begin{enumerate}[resume]
\item pour toute petite catégorie, le foncteur $\alpha_C : i_Ai^*(C)
\to C$ est un morphisme asphérique de $\Cat$.
\end{enumerate}
\end{enumerate}
\end{prop}
\begin{proof}
Voir \cite[proposition 1.7.6]{maltsiniotis2005}.
\end{proof}

\paragr\label{def:foncteurTest} On dit alors qu'un foncteur $i : A \to \Cat$ est un
\ndef[foncteur test!faible]{foncteur test faible} si $A$ est une
catégorie test faible et si $i$ est un foncteur asphérique. On dit
que~$i$ est un \ndef[foncteur test!local]{foncteur test local} si pour
tout objet $a$ de $A$, le foncteur $\tranche{A}{a} \to \Cat$ induit
par $i$ est un
foncteur
test. Enfin, on dit que $i$ est un \ndef{foncteur test} s'il est à la
fois un foncteur test faible et un foncteur test local. Autrement dit,
$i$ est un foncteur test si $A$ est une catégorie test, et si les
foncteurs $i : A \to \Cat$ et $i|_{\tranche{A}{a}} :
\tranche{A}{a}\to\Cat$ sont asphériques, pour tout objet $a$ de $A$.

\begin{theorem}[Grothendieck]\label{propEquivalencesFoncteursTest}
Soient $A$ une petite catégorie et $i : A \to \Cat$ un foncteur tel
que pour tout objet $a$ de $A$, la catégorie $i(a)$ soit asphérique. 
\begin{enumerate}
\item Les conditions suivantes sont équivalentes : 
\begin{enumerate}
\item $i$ est un foncteur test local (et $A$ est donc une catégorie
test locale);
\item pour toute petite catégorie asphérique $C$, le préfaisceau
$i^*(C)$ est localement asphérique.
\end{enumerate}
\item De plus, si ces conditions sont satisfaites, alors les
conditions sont équivalentes :
\begin{enumerate}[resume]
\item $i$ est un foncteur test (et $A$ est donc une catégorie test);
\item $i$ est un foncteur asphérique;
\item $A$ est une catégorie asphérique.
\end{enumerate}
\item De plus, si pour tout objet $a$ de $A$, la catégorie $i(a)$
admet un objet final, alors les conditions équivalentes (a) et (b)
sont équivalentes à la condition suivante, en notant toujours
$\Delta_1$ la catégorie correspondant à l'ensemble ordonné $\lbrace 0
\leq 1 \rbrace$ : 
\begin{enumerate}[resume]
\item $i^*(\Delta_1)$ est un préfaisceau localement asphérique.
\end{enumerate}
\end{enumerate}
\end{theorem}
\begin{proof}
Voir \cite[théorème 1.7.13]{maltsiniotis2005}.
\end{proof}

\begin{remark}
Toutes les notions introduites dans cette section peuvent être
généralisées en remplaçant la classe $\W_\infty$ des équivalences de
Thomason par d'autres classes de morphismes de $\Cat$ vérifiant des
propriétés analogues à cette dernière, notamment une
généralisation du théorème A de Quillen (\ref{theoremeA}). On appellera
\emph{localisateurs fondamentaux} de telles classes de morphismes de
$\Cat$, et on les étudiera plus en profondeur dans la section
\ref{secLocalisateursFondamentaux}. 

Étant donné un localisateur fondamental $\W$ de $\Cat$, on peut ainsi
introduire la notion de $\W$-catégorie test (et toutes ses variantes).
On prouve alors (voir~\cite[section 3.3]{maltsiniotis2005}, ainsi
qu'en utilisant \cite{cisinski2004localisateur}) que toute catégorie
test (resp. faible, locale, stricte) est une~$\W$-catégorie test
(resp. faible, locale, stricte).
\end{remark}

\begin{remark}
Si la présence des ensembles simpliciaux semble encore trop prégnante
du fait de leur utilisation dans la définition des équivalences
de Thomason, on peut en fait s'en défaire de deux manières. D'une
part, les équivalences de Thomason coïncident avec les morphismes $u :
A \to B$ tels que le morphisme de topos induit $(u^*,u_*) : \pref{A} \to
\pref{B}$ est une équivalence
d'Artin-Mazur~(voir \cite{artin2006etale} et \cite{maltsiniotis2005}),
où $u^* : \pref{B} \to \pref{A}$ désigne le foncteur image inverse $X
\mapsto X\circ u$, et $u_*$ son adjoint à droite.
D'autre part, on peut caractériser de manière entièrement catégorique
la classe $\W_\infty$. Cette caractérisation est l'objet d'une
conjecture de Grothendieck démontrée par Cisinski dans
\cite{cisinski2004localisateur}, que l'on énoncera au paragraphe
\ref{locFondamentalMinimalThm}.
\end{remark}

\section{Homologie simpliciale et correspondance de Dold-Kan} \label{secDoldKan}

La catégorie $\Delta$, en plus d'être une catégorie test, est
particulièrement adaptée à la théorie de l'homologie. On rappelle dans
cette section comment calculer l'homologie des groupes abéliens
simpliciaux, c'est-à-dire des préfaisceaux en groupes abéliens sur
$\Delta$, avant d'énoncer le théorème de Dold-Kan.

\paragr  
On rappelle que si $X$ est un ensemble simplicial et $k$ est un entier
strictement positif, on dit qu'un $k$-simplexe de $X$ est dégénéré s'il est dans
l'image d'un des morphismes~$s_i^{k-1} : X_{k-1}\to X_k$ pour $0\leq i
\leq k-1$.

\paragr Si $\Delta_n$ est un objet de $\Delta$, on peut lui associer un
complexe de chaînes de groupes abéliens noté $c\Delta_n$ de la
manière suivante. On pose, pour $k\geq 0$,
\notindex{$c(\Delta_n)$}%
\[
(c\Delta_n)_k = \bigoplus_{\Delta_k \hookrightarrow \Delta_n}\Z
\]
le groupe abélien libre engendré par les $k$-simplexes non dégénérés
de~$\Delta_n$. On définit la différentielle en posant, pour tout
entier $k > 0$ et pour tout~$k$\nobreakdash-simplexe non dégénéré $\varphi$ de
$\Delta_n$, 
\[
d_k \langle\varphi\rangle  = \sum_{i=0}^k\left(-1\right)^i \langle
d_i^k\varphi \rangle
\]
en rappelant qu'on note $d_i^k : (\Delta_n)_k \to (\Delta_n)_{k-1}$
l'image par le préfaisceau $\Delta_n$ du morphisme $\delta_i^k :
\Delta_{k-1} \to \Delta_k$ de $\Delta$ pour $0 \leq i \leq k$.
On vérifie alors grâce aux relations simpliciales que le morphisme $d$
est bien une différentielle, c'est-à-dire que l'on a bien $d^2=0$.

\paragr Si $ u: \Delta_n \to \Delta_m$ est un morphisme de $\Delta$, on
définit un morphisme
\[
  c(u) : c\Delta_n \to c\Delta_m
\]
de la manière suivante. Si $\varphi : \Delta_k
\hookrightarrow~\Delta_n$ est un $k$-simplexe non dégénéré
de~$\Delta_n$, son image par le morphisme $c(u)_k : (c\Delta_n)_k \to
(c\Delta_m)_k$ est donnée par
\[
c(u)_k\langle\varphi\rangle = \begin{cases}
\langle u\circ \varphi\rangle & \text{ si } u\circ \varphi \text{ est
un monomorphisme;}
\\
0 & \text{ sinon.}
\end{cases} 
\]
Les relations simpliciales garantissent que $c(u)$ est un morphisme de
complexes, et qu'on a ainsi bien défini de cette manière un foncteur 
\[
c : \Delta \to \Ch(\Ab) \pbox{.}
\]

\paragr\label{adjdoldkan} On peut montrer (on le fera en détail dans la section
\ref{secExtKan}) que le foncteur $c$ défini ci-dessus induit un couple
de foncteurs adjoints 
\[
\adjpair{c_!}{\prefab{\Delta}}{\Ch}{c^*}
\]
où on a noté $\prefab{\Delta}$ la catégorie des préfaisceaux en
groupes abéliens sur $\Delta$. Si $C$ est un complexe de chaînes, le
\ndef[groupe abélien!simplicial]{groupe abélien simplicial} associé à
$C$ par le foncteur $c^*$ est 
  \[
c^*(C) : \Delta_n \mapsto \Hom_{\Ch}(c\Delta_n, C) 
  \]
\notindex{$c^*$}%
\notindex{${c_\bang^\ab}$}%
et $c_!$ est l'extension de Kan à gauche de $c : \Delta \to \Ch$ le long
du foncteur
\[
  \Delta \hookrightarrow \EnsSimp \xrightarrow{\Z}
  \prefab{\Delta} \mdvirg
\]
où on a noté $\Z$ le foncteur induit par le foncteur groupe abélien
libre.

\medskip
On va maintenant donner une description explicite de cet adjoint à gauche.

\paragr\label{defMooreComplex}
Si~$X$ est un groupe abélien simplicial, on commence par définir un
complexe de groupes abéliens $\mathsf{C}X$ appelé\footnote{Moore introduit en
fait le complexe normalisé (via la
formule d'intersection des noyaux) ainsi que le complexe non normalisé  dans
\cite{mooreAlgebraic}. Certain$\cdot$e$\cdot$s
auteur$\cdot$rice$\cdot$s appellent complexe de Moore le complexe
d'intersection des noyaux (on suit ici la
terminologie utilisée par Goerss et
Jardine dans \cite{goerssjardine}).} \ndef[complexe non normalisé]{complexe de
Moore} ou \emph{complexe non normalisé} associé à $X$ en posant, pour tout $k\geq 0$, 
\[
(\mathsf{C}X)_k = X_k
\] 
et, pour tout entier $k>0$, 
\[
d_k = \sum_{i=0}^k (-1)^i d_i^k : X_k \to X_{k-1}
\pbox{.}
\]
Les relations simpliciales garantissent que l'on a bien défini un
complexe de chaînes. 

On vérifie alors (toujours grâce aux relations
simpliciales) qu'en notant, pour tout entier positif $k$, $(DX)_k$ le sous-groupe de
$(\mathsf{C}X)_k$ engendré par les~$k$\nobreakdash-simplexes dégénérés de $X$, on obtient un
sous-complexe $DX$ de $\mathsf{C}X$, et on note $\dk X$ le complexe
quotient 
\[
\dk X = \mathsf{C}X / DX \pbox{.}
\]
\notindex{$\dk$}%
On a ainsi défini un foncteur $ \dk : \prefab{\Delta} \to \Ch(\Ab)$
appelé foncteur \ndef[complexe normalisé!d'un groupe abélien
simplicial]{complexe de chaînes normalisé associé}. 

\begin{proposition}[Kan]\label{normaliseExtKan}
Le foncteur $\dk : \prefab{\Delta}\to\Ch(\Ab)$ est adjoint à gauche au foncteur $
c^* : \Ch \to \prefab{\Delta}$ défini ci-dessus. Autrement dit, si $X$ est un
groupe abélien simplicial, on a un isomorphisme naturel 
\[
c_! X \simeq \dk X 
\]
dans $\Ch(\Ab)$.
\end{proposition}
\begin{proof}
Voir \cite[proposition 6.3]{kan1958functors} pour la preuve originale.
Il s'agit en fait d'un des premiers exemples de la construction de Kan
aujourd'hui appelée \emph{extension de Kan}.
On peut aussi simplement justifier que le foncteur~$\dk$ commute aux
limites inductives, et que sa restriction aux objets de $\Delta$
coïncide avec le foncteur $c$ (voir le paragraphe \ref{kanext}).
\end{proof}

\paragr On peut donner encore une autre description explicite du
foncteur complexe de chaînes normalisé de la manière suivante. Si $X$
est un groupe abélien simplicial, on note, pour tout entier
$k>0$, $(NX)_k$ le sous-groupe 
\[
(NX)_k =  \bigcap_{1 \leq i\leq k}\ker d^k_i
\]
de $(\mathsf{C} X)_k$. En utilisant à nouveau les relations
simpliciales, on vérifie qu'on obtient ainsi un sous-complexe 
\[
NX \hookrightarrow \mathsf{C} X 
\]
de $\mathsf{C} X$, en posant également $(NX)_0=X_0$. On remarque que
pour tout entier~$k > 0$, la différentielle de $\mathsf{C}X$ restreinte à $NX_k$
est simplement le morphisme de groupes abéliens $d^k_0$.

On obtient alors un nouveau foncteur noté
\[
  N : \prefab{\Delta}\to\Ch(\Ab) \mdvirg X \mapsto NX \pbox{.}
\]

\begin{prop}[Kan]\label{normaliseIntersectionNoyaux}
Si $X$ est un groupe abélien simplicial, il existe un isomorphisme naturel
\[
NX \simeq \dk X 
\]
dans $\Ch(\Ab)$.
\end{prop}
\begin{proof}
On renvoie à \cite[proposition 6.3]{kan1958functors} pour la preuve
originale de Kan, ou à \cite[théorème 2.1]{goerssjardine} pour un
exposé plus moderne.
\end{proof}

Un des intérêts de cette dernière expression est qu'elle permet de
faire le lien entre les groupes d'homologie et les groupes d'homotopie
d'un groupe abélien simplicial, d'une manière que l'on va maintenant
détailler. 

\paragr L'homologie d'un groupe abélien simplicial $X$ à coefficients
entiers est définie comme l'homologie du complexe $\dk X$. On note
alors, pour tout groupe abélien simplicial $X$ et pour tout entier
$n\geq 0$, 
\[
H_n(X,\Z) := H_n(\dk X) \pbox{.}
\]
En fait, le résultat suivant montre qu'on peut également choisir
d'utiliser le complexe non normalisé pour étudier le type
d'homologie d'un groupe abélien simplicial.

\begin{prop}[Eilenberg-Mac Lane]\label{DoldKanNormalisePasBesoin}
Pour tout groupe abélien simplicial $X$, la projection canonique
\[
  \mathsf{C}X \to \dk X
\]
est une équivalence d'homotopie de complexes de chaînes.
\end{prop}
\begin{proof}
On renvoie à \cite[théorème 4.1]{Eilenberg1953GroupsI}
pour la preuve originale, ou à 
\cite[chapitre III, théorème 2.4]{goerssjardine}.
\end{proof}

\begin{remark}
En utilisant le fait que le complexe normalisé est isomorphe à un
sous-complexe du complexe non normalisé, on peut montrer que pour tout
groupe abélien simplicial $X$, la suite exacte de complexes 
\[
0 \rightarrow DX \rightarrow \mathsf{C}X \rightarrow \dk X \rightarrow 0 
\]
est scindée. En particulier, le complexe $DX$ est donc un complexe
exact.
\end{remark}

La proposition \ref{DoldKanNormalisePasBesoin} constitue une propriété
importante de l'homologie des préfaisceaux abéliens sur la catégorie
$\Delta$. Par exemple, on verra que pour les préfaisceaux abéliens sur
la catégorie cubique, l'analogue du complexe non normalisé ne calcule
pas la bonne homologie (voir à la section
\ref{secCubes}).

\begin{proposition}[Moore]
L'ensemble simplicial sous-jacent à un préfaisceau en groupes sur
$\Delta$ est un complexe de Kan.
\end{proposition}
\begin{proof}
La preuve originale est dans \cite[chapitre II, théorème 2.2]{mooreAlgebraic}.
On renvoie également à \cite[lemme 3.4]{goerssjardine}.
\end{proof}

On rappelle qu'on peut décrire d'une manière combinatoire les groupes
d'homotopie d'un complexe de Kan (voir par exemple~\cite[chapitre~I,
section~7]{goerssjardine}). Il s'avère que
la description des groupes d'homotopie est encore plus simple dans le
cas des groupes abéliens simpliciaux (et des groupes simpliciaux), en
vertu du résultat suivant :

\begin{prop}[Moore]\label{homologieHomotopieNormalise}
On note $\U : \prefab{\Delta}\to\pref{\Delta}$ le foncteur d'oubli.
Si $X$ est un groupe abélien simplicial, alors
\begin{enumerate}
\item il existe une bijection naturelle $\pi_0(\U X) \simeq
H_0(N X)$;
\item pour tout entier $n\geq 1$, il existe un isomorphisme naturel de
groupes abéliens
\[
\pi_n(\U X,0) \simeq H_n(N X) \pbox{.}
\]
\end{enumerate}
\end{prop}
\begin{proof}
On trouve déjà une preuve de ce résultat dans un article non publié de
Moore~\cite[chapitre II, proposition 2.9]{mooreAlgebraic}, ainsi que
dans l'article de Kan \cite[lemme 7.7]{kan1958functors}. Voir aussi
\cite[chapitre III, corollaire 2.5]{goerssjardine}. 
\end{proof}

\begin{remark}
Moore prouve également dans \cite[chapitre II, proposition
2.9]{mooreAlgebraic} que pour tout groupe abélien simplicial $X$,
l'action du groupe $X_0$ sur l'ensemble simplicial sous-jacent à $X$
induit un isomorphisme, pour tout $0$-simplexe $x$ de $X$ et pour tout
entier $n\geq 1$,
\[
\pi_n(\U X, x) \simeq \pi_n(\U X, 0) \pbox{.}
\]
\end{remark}

\paragr Puisque la catégorie $\Delta$ est une catégorie test, il peut
sembler naturel de se demander si tout type d'homologie peut être
réalisé comme l'homologie d'un groupe abélien simplicial. En d'autres
termes, si $C$ est un complexe de chaînes, existe-t'il un groupe
abélien simplicial $X$ et un quasi-isomorphisme $c_! X \to C$ ? 

Non seulement la
réponse à la question précédente est affirmative, mais il s'avère que
l'on peut en fait réaliser tout \emph{complexe} comme le complexe
normalisé d'un groupe abélien simplicial. C'est ce qu'affirme le
théorème de Dold-Kan.

\begin{theorem}[Dold, Kan]\label{doldkan}
Les foncteurs intervenant dans l'adjonction
\[
\adjpair{c_!}{\prefab{\Delta}}{\Ch(\Ab)}{c^*}
\]
sont des équivalences de catégories quasi-inverses l'une de l'autre.
\end{theorem}
\begin{proof}
La preuve originale est due à Kan \cite[théorèmes 8.1
et~8.2]{kan1958functors}, et on trouvera la preuve de Dold dans
\cite[théorème 1.9]{dold1958homology}. Les deux preuves montrent en
fait que le foncteur $N : \prefab{\Delta} \to \Ch(\Ab)$ est un
quasi-inverse au foncteur $c^*$. La forme énoncée ici utilise aussi
les proposition \ref{normaliseIntersectionNoyaux} et
\ref{normaliseExtKan}.
\end{proof}

\paragr Une des conséquences de ce résultat, quand on l'associe à la proposition
\ref{homologieHomotopieNormalise}, est que, étant donné un groupe
abélien $M$ et un entier $n \geq 1$, la réalisation géométrique de
l'ensemble simplicial
\[
X = \U c^*(M[n]) \mdvirg
\]
où on a noté $M[n]$ le complexe de chaînes concentré en degré $n$ de
valeur $M$, est un espace d'Eilenberg-Mac Lane de type $K(M,n)$. Cela
\notindex{$K(\pi,n)$}%
signifie que $X$ est connexe, et que pour tout entier $k >0$ et pour
tout $0$-simplexe $x$ de $X$, on a
\[
\pi_k(X,x)=
\begin{cases}
M & k= n \pbox{;}\\
0 & k\neq n \pbox{.}
\end{cases}
\]

\paragr Le foncteur $c^* : \Ch(\Ab) \to \prefab{\Delta}$, quant à lui,
peut être décrit de la manière suivante. Si $C$ est un complexe de
chaînes de groupes abéliens, on pose, pour tout entier $n\geq0$, 
\notindex{$\Gamma$}%
\[
\Gamma(C)_n = \bigoplus_{\Delta_n \twoheadrightarrow \Delta_i} C_i
\mdvirg
\]
où la somme est indexée par les épimorphismes $\Delta_n \to \Delta_i$
de $\Delta$.
Pour tout morphisme $u : \Delta_n \to \Delta_m$, on définit un
morphisme de groupes abéliens~$\Gamma(C)_m \to \Gamma(C)_n$ de la
manière suivante. Si
$\varphi : \Delta_m \twoheadrightarrow \Delta_i$ est un épimorphisme,
on écrit la factorisation du morphisme $\varphi \circ u$ sous la forme
d'un épimorphisme suivi d'un monomorphisme 
\[
\xymatrix{
\Delta_n \ar[r]^{u} \ar@{>>}[rd]_{\varphi_u} & \Delta_m \ar@{>>}[r]^{\varphi} & \Delta_i
\\
& \Delta_j \ar@{^(->}[ru]_{\psi_u} & \pbox{.}
}
\]
Pour tout élément $x$ de $C_i$, le morphisme $\Gamma(C)(u) :
\Gamma(C)_m \to \Gamma(C)_n$ envoie
alors le couple $(\varphi,x)$ sur le couple $(\varphi_u, x_u)$ avec
\[
x_u= \begin{cases}
x & \text{ si } j=i \pbox{;}\\
dx & \text{ si } j=i-1 \text { et } \psi_u=\delta_0^i\pbox{;}\\
0 & \text{ sinon.} 
\end{cases}
\]
On peut alors vérifier que $\Gamma(C)$ est bien un groupe abélien
simplicial.

\begin{prop}
Si $C$ est un complexe de chaînes, il existe un isomorphisme naturel 
\[
c^*(C) \simeq \Gamma(C)
\]
dans $\prefab{\Delta}$.
\end{prop}
\begin{proof}
Dans leur exposé du théorème de Dold-Kan, Goerss et Jardine
\cite[chapitre III, corollaire 2.3]{goerssjardine} prouvent que
$\Gamma$ est un quasi-inverse au foncteur $N$, ce qui implique en
particulier que les foncteurs $\Gamma$ et $c^*$ sont naturellement
isomorphes.
\end{proof}

\section{Le foncteur d'abélianisation}
La correspondance entre les groupes d'homologie et les groupes
d'homotopie d'un groupe abélien simplicial est une propriété
remarquable, que Grothendieck décide dans \emph{Pursuing Stacks}
\cite{pursuingstacks} de mettre en avant en vue d'un
traitement axiomatique de la \emph{linéarisation des types
d'homotopie}. Nous allons détailler ce que cela signifie dans cette
section.

\paragr
Si $ A $ est une petite catégorie, on note 
\[ \Whf{A} : \pref{A} \to \prefab{A} \] 
le foncteur qui associe à un préfaisceau d'ensembles le préfaisceau en
groupes abéliens libres associé. Si $X$ est un préfaisceau d'ensembles
sur $A$, on a donc 
\begin{align*} 
\Wh{A}{X} : A^{\op} &\to \Ab
\\ a & \mapsto \Z^{(Xa)} \pbox{.}
\end{align*}
\notindex{$\Wh{A}{X}$}%
Grothendieck appelle ce foncteur le \ndef{foncteur de
Whitehead}
\footnote{Une recherche de la référence précise à laquelle
pense Grothendieck en attribuant ce résultat à Whitehead n'a pas
abouti à un résultat particulièrement concluant. Il est possible qu'il
fasse référence au résultat de J.H.C Whitehead selon lequel toute
équivalence faible entre complexes de Kan est une équivalence
d'homotopie. Nous avons tout de même choisi de garder ici sa
terminologie pour faciliter la comparaison avec \emph{Pursuing Stacks}.},
en raison du résultat suivant dans le cas où $A$ est la catégorie
$\Delta$ :

\begin{prop}\label{theoremeWhitehead}
Le foncteur $\Whf{\Delta} : \pref{\Delta} \to \prefab{\Delta}$
envoie les équivalences faibles simpliciales sur des 
morphismes dont le morphisme d'ensembles simpliciaux sous-jacent est
une équivalence faible simpliciale.
\end{prop}
\begin{proof}
Le théorème de Whitehead affirme que toute équivalence faible d'homotopie
entre CW-complexes est une équivalence d'homotopie. Or, la
réalisation géométrique d'un ensemble simplicial est un
$CW$\nobreakdash-complexe, donc
si~$u : X \to Y$ est une équivalence faible simpliciale, le morphisme
induit entre les réalisations géométriques~$|u| : |X| \to |Y|$
est donc une équivalence d'homotopie de $\Top$, et induit donc un
quasi-isomorphisme entre les complexes d'homologie singulière
\[
\mathsf{H}^{\mathsf{sing}}(|X|,\Z) \xrightarrow{\simeq}
\mathsf{H}^{\mathsf{sing}}(|Y|,\Z) \pbox{.} \]
Un résultat de Milnor \cite[lemme 5]{milnor1957realization} énonce
qu'on a un isomorphisme naturel entre l'homologie d'un ensemble
simplicial et l'homologie singulière de sa réalisation topologique. En
utilisant la proposition \ref{DoldKanNormalisePasBesoin} ainsi que la
proposition~\ref{homologieHomotopieNormalise}, on obtient finalement
un isomorphisme naturel en $X$
\[
\pi_n(\U\Wh{\Delta}{X},x) \simeq \mathsf{H}^\mathsf{sing}_n(|X|,\Z)
\]
pour tout $0$-simplexe $x$ de $X$ et pour tout entier $n\geq 1$, ainsi
qu'une bijection naturelle
\[
  \pi_0(\U\Wh{\Delta}{X}) \simeq H^\mathsf{sing}_0(|X|) \mdvirg
\]
qui
nous permet de conclure que si $u$ est une équivalence faible
simpliciale, il en est de même pour le morphisme $\U\Wh{\Delta}{u}$.
\end{proof}

\paragr Si $A$ est une petite catégorie et $X$ est un préfaisceau sur
$A$, on notera, quand aucune ambigüité n'est possible, simplement 
\[
\Z^{(X)} = \Wh{A}{X}
\]
le préfaisceau en groupes abéliens sur $A$ librement engendré par $X$.

\paragr\label{paragr:commutativiteu^*WhetU} Si $u : A \to B$ est un
foncteur entre petites catégories, on rappelle qu'on note 
\[
u^* : \pref{B} \to \pref{A} 
\]
le foncteur défini pour tout préfaisceau $X$ sur $B$ par 
\[
u^*X : a \mapsto Xu(a) \pbox{.}
\]
Puisque le foncteur $u^*$ est un adjoint à droite, il préserve les
produits et envoie donc les groupes abéliens internes à $\pref{B}$ sur
des groupes abéliens internes à $\pref{A}$. En d'autres termes,
puisque les groupes abéliens internes aux catégories de préfaisceaux
coïncident avec les préfaisceaux en groupes abéliens, $u^*$
induit un foncteur
\[
(u^*)^\ab : \prefab{B} \to \prefab{A} 
\]
\notindex{$(u^*)^\ab$}%
que l'on notera simplement $u^*$ lorsqu'aucune confusion n'en découle,
et on obtient le diagramme commutatif
\[
\xymatrix{
{ \prefab{B} } \ar[r]^{ u^* } \ar[d]_-{ \U } & { \prefab{A} } \ar[d]^-{ \U} \\
{ \pref{B} } \ar[r]_{ u^* } & { \pref{A} } \pbox{.}
} 
\]
De plus, on remarque que le diagramme
\[
\xymatrix{
\prefab{B} \ar[r]^{u^*} & \prefab{A} \\
\pref{B} \ar[u]^{\Whf{B}} \ar[r]_{u^*} & \pref{A} \ar[u]_{\Whf{A}}
} 
\]
est également commutatif.

\paragr\label{WhAaa} On rappelle que si $a$ est un objet de $A$, on
note également $a$ son image dans $\pref{A}$ par le plongement de
Yoneda. On notera donc $\Wh{A}{a}$ son image dans $\prefab{A}$,
c'est-à-dire le préfaisceau en groupes abéliens 
\begin{align*}
\Wh{A}{a} :  A^{\op} &\to \Ab \\
x &\mapsto \Z^{\left(\Hom_A(x,a)\right)} \pbox{.}
\end{align*}
\notindex{$\Wh{A}{a}$}%
On dira alors qu'un tel préfaisceau est \emph{représentable} : c'est le
préfaisceau en groupes abéliens \emph{représenté par} $a$. On rappelle
la version additive du lemme de Yoneda.

\begin{prop}
Soit $\A$ une petite catégorie additive. Si $X$ est un
préfaisceau additif en groupes abéliens sur $\A$, alors pour tout
objet $a$ de $\A$, on a un isomorphisme naturel 
\[
\Hom_{\prefab{{\A}}}(\Wh{\A}{a}, X) \simeq Xa \pbox{.}
\]
\end{prop}
\begin{proof}
On renvoie par exemple à \cite[proposition 1.3.7]{borceuxHandbook2}. 
\end{proof}

\begin{coro}
Si $A$ est une petite catégorie et $X$ est un préfaisceau en
groupes abéliens sur $A$, alors, pour tout objet $a$ de $A$, on a un
isomorphisme naturel
\[
\Hom_{\prefab{A}}(\Wh{A}{a}, X) \simeq Xa \pbox{.}
\]
\end{coro}
\begin{proof}
Il suffit d'appliquer le lemme de Yoneda additif, en remarquant que la
catégorie des préfaisceaux abéliens sur $A$ coïncide avec la catégorie
des préfaisceaux additifs sur la catégorie préadditive libre sur $A$,
c'est-à-dire la catégorie munie des mêmes objets que $A$, et dont les
groupes abéliens de morphisme sont les groupes abéliens libres sur les
ensembles de morphismes de $A$.
\end{proof}

\begin{remark}
Le foncteur $\Whf{A} : A\to\prefab{A}$ est fidèle, mais n'est pas
plein. 
\end{remark}

Nous allons maintenant détailler ce que Grothendieck entend, dans
\cite{pursuingstacks}, par \emph{linéarisation} (ou
\emph{abélianisation}) des types d'homotopie, et énoncer la
correspondance de Dold-Kan homotopique pour $\Delta$, que l'on va 
chercher à généraliser aux autres catégories test dans cette thèse.

\label{sectypehomologieassocie}

\paragr On note 
\[
  \Hotab = \W_{qis}^{-1}\Ch(\Ab)
\]
\notindex{$\Hotab$}%
la catégorie dérivée $D_+(\Ab)$ des groupes abéliens, obtenue en
localisant la catégorie~$\Ch(\Ab)$ par la classe des
quasi-isomorphismes. On appellera \ndef[type
d'homologie]{types d'homologie} les objets de $\Hotab$. 

\paragr\label{HotabUDelta} On note aussi
\[
  \Hotab_{\Delta}^{\U} = \Big( \U^{-1}\W_{\Delta} \Big)^{-1} \prefab{\Delta}
\]
la catégorie obtenue en localisant $\prefab{\Delta}$ par les
morphismes qui sont des équivalences faibles simpliciales après oubli
de la structure de groupe. 

La \emph{correspondance de Dold-Kan homotopique} pour la catégorie
$\Delta$, qui est un corollaire de la correspondance de Dold-Kan \og
stricte \fg{} et de la proposition~\ref{homologieHomotopieNormalise},
s'énonce alors de la manière suivante.

\begin{prop}
Le foncteur complexe normalisé $ c_! : \sAb \to \Ch $ et son adjoint
à droite $c^* : \Ch(\Ab) \to \sAb$ induisent des 
équivalences de catégories quasi-inverses l'une de l'autre
\[
c_! : \Hotab_{\Delta}^{\U} \xrightarrow{\simeq} \Hotab \mdvirg
c^* : \Hotab \xrightarrow{\simeq} \Hotab^{\U}_{\Delta}
\]
\end{prop}
\begin{proof}
On sait grâce à la proposition \ref{homologieHomotopieNormalise} et à
la proposition~\ref{normaliseIntersectionNoyaux} que le
foncteur $c_!$ préserve et reflète les équivalences faibles.
On peut alors conclure grâce au théorème de Dold-Kan que le foncteur
$c^*$ préserve également les équivalences faibles. Les foncteurs
$c_!$ et $c^*$ induisent donc bien des foncteurs entre les
catégories localisées, et ces foncteurs induits sont encore des
équivalences de catégories.
\end{proof}

\paragr On obtient alors, grâce à la proposition
\ref{theoremeWhitehead}, le
diagramme commutatif à isomorphisme près 
\begin{equation}\label{diagrammeabelianisationdelta}
\xymatrix{
\Cat \ar[r]^{\nerf} \ar[d]^{} &\pref{\Delta} \ar[r]^{\Whf{\Delta}} \ar[d]^{} &
\prefab{\Delta} \ar[r]^{\dk}_{\simeq} \ar[d]^{} & \Ch \ar[d]^{}  \\
\Hot \ar[r]^{\simeq} & \Hot_{\Delta} \ar[r]_{\Whf{\Delta}} &
\Hotab_{\Delta}^{\U}
\ar[r]^{\simeq}_{\dk} & \Hotab \pbox{,} 
} 
\end{equation}
où les flèches verticales correspondent aux localisations. 
Le foncteur de Whitehead peut donc être vu comme un foncteur
\emph{d'abélianisation des types d'homotopie}, qui associe à un type
d'homotopie son type d'homologie, en utilisant le vocabulaire de la
catégorie test $\Delta$. 

La question naturelle que pose Grothendieck dans \cite{pursuingstacks}
est alors la suivante : \emph{existe-t'il un foncteur canonique} 
\[
\Whf{} : \Hot \to \Hotab 
\]
ne faisant pas intervenir les particularités de la catégorie $\Delta$
? Par ailleurs, si~$A$ est une petite catégorie, on peut
introduire la catégorie 
\[
\Hotab_A^{\U} = \Big( \U^{-1}\W_A \Big)^{-1} \prefab{A}
\]
où $ \W_A $ désigne la classe des équivalences test de $ A $, et se
poser la question de l'existence d'un foncteur canonique 
\[
\Hf{A} : \Hotab_{A}^{\U} \to \Hotab \mdvirg
\]
de sorte que le carré suivant ait un sens, et soit commutatif à
isomorphisme naturel près : 
\[
\xymatrix{
\Hot_A \ar[r]^-{\Whf{A}} \ar[d]_{i_A} & \Hotab_A^{\U} \ar[d]^{H_A} \\
\Hot \ar[r]_-{\Whf{}} & \Hotab \pbox{.}
} 
\]
L'existence de la flèche horizontale du haut, dans le cas de $\Delta$,
provient de la coïncidence des groupes d'homotopie et des groupes
d'homologie des groupes abéliens simpliciaux. La commutativité de ce
carré revient à demander à ce que le foncteur $H_A$ envoie
les préfaisceaux représentables de la forme $\Wh{A}{a}$ sur le type
d'homologie de la catégorie $\tranche{A}{a}$.

 La stratégie que nous allons explorer est la suivante : 
\begin{enumerate}
  \item Définir un foncteur canonique $\Hf{A} : \prefab{A} \to \Hotab$
(\ref{defH_A}).
  \item Définir une classe d'\emph{équivalences faibles abéliennes} de
  $\prefab{A}$ en posant
\[
f \in \Wab_A \iff H_A(f)\in \Fl\big(\Hotab\big) \text{ est un
isomorphisme} \pbox{.}
\]
On posera alors (\ref{secPseudoTest})
\[
\Hotab_A = \big(\Wab_A\big)^{-1}\prefab{A} \pbox{.}
\]
  \item Explorer sous quelles conditions les catégories 
\[
\Hotab_A^{\U} \text{ , }\Hotab_A
\]
coïncident. Par exemple, dans le cas de la catégorie
$\Delta$, on va montrer que les classes de flèches 
\[
\U^{-1}\W_{\Delta} \mdvirg \W_\Delta^\ab \subset \Fl\big(\sAb\big)
\]
coïncident (voir proposition \ref{DeltaWhitehead}), et les
deux catégories localisées associées coïncident. 
On dira que $\Delta$ est une \emph{catégorie de Whitehead}. En
revanche, il existe des petites catégories qui ne sont pas des
catégories de Whitehead, comme la catégorie des globes réflexifs
$\Grefl$
(\ref{globespaswhitehead}).
\end{enumerate}

\section{Homologie et colimite homotopique}

Dans cette section, on va introduire le foncteur $\Hf{A}$ annoncé
ci-dessus pour toute petite catégorie $A$. On va voir qu'il permet à
son tour de définir un foncteur associant à toute petite catégorie son
homologie à coefficients entiers, sans passer par l'intermédiaire de la
catégorie $\Delta$.

\paragr Si $G$ est un groupe, on note $BG$ le groupoïde associé, c'est
à dire la catégorie munie d'un unique objet $*$ avec $\Hom_{BG}(*,*) =
G$. On rappelle que se donner un foncteur $A : BG \to \Ab$ correspond
à se donner un $G$\nobreakdash-module~$A$, c'est-à-dire un $\Z[G]$-module, où
$\Z[G]$ désigne l'anneau du groupe $G$. On remarque que la limite
inductive d'un tel foncteur $A : BG \to \Ab$
n'est autre que le groupe abélien
\[
  A_G = A / I
\]
des coinvariants, où $I$ est le sous-groupe de $A$ engendré par les
éléments de type $g \cdot x - x$ pour $x\in A$ et $g \in G$. Le calcul de l'homologie $\H{G}{A}$ de~$G$
à coefficients dans $A$ est classiquement exprimé comme le foncteur
dérivé à gauche du foncteur des coinvariants. Autrement dit, on peut
énoncer la proposition suivante.

\begin{prop}
Si $A : BG \to \Ab$ est un $G$-module, alors on a un
isomorphisme canonique naturel
\[
\H{G}{A} \simeq L\limind A 
\]
dans $\Hotab$.
\begin{remark}
On note ici $\H{G}{A}$ l'objet de la catégorie $\Hotab$, c'est-à-dire
le complexe vu à quasi-isomorphisme près, qui calcule l'homologie de
$G$ à coefficients dans $A$, et pas
le groupe abélien gradué des groupes d'homologie. Par exemple, si
$P_{\bullet} \to A$ est une résolution projective, alors le complexe
$\limind{P_{\bullet}}$ est un représentant de $\H{G}{A}$ dans
$\Hotab$. On notera au besoin $\operatorname{H}_i(G,A)$
le~$i$\nobreakdash-ième groupe d'homologie de ce complexe. 
\end{remark}
\end{prop}

Cette observation justifie en partie le choix que l'on va faire, en
définissant le foncteur $\Hf{A}$ comme le foncteur de limite inductive
homotopique. Avant cela, on introduit quelques préliminaires relatifs
aux catégories de modèles.

\begin{prop}[Quillen]
La catégorie $\Ch(\Ab)$ des complexes de chaînes
concentrés en degré positif admet une structure de catégorie de
modèles combinatoire dite \emph{projective}, dont \begin{itemize}
\item les équivalences faibles sont les quasi-isomorphismes;
\item les fibrations sont les morphismes $f_\bullet$ de complexes tels
que pour tout entier $n > 0$, $f_n$ est un morphisme surjectif de
groupes abéliens; %
\item les cofibrations sont les monomorphismes dont le conoyau est un
complexe de groupes abéliens projectifs (c'est-à-dire, libres).
\end{itemize} 
\end{prop}
\begin{proof}
Voir \cite[chapitre II, pages 4.11 et 4.12]{quillenhomotalg}. Il
s'agit en fait d'un corollaire du théorème de Dold-Kan
(\ref{doldkan}), mais on peut le prouver directement, comme le font
Goerss et Jardine dans \cite[chapitre~III, corollaire~2.10]{goerssjardine}.
\end{proof}

\paragr\label{colimiteHomotopiqueDiagrammeCh} 
Puisque la structure de catégorie de modèles définie ci-dessus est
combinatoire, si~$A$ est une petite catégorie, on peut munir la catégorie 
\[
\Homi\big({A}^{\op}, \Ch(\Ab)\big) 
\]
d'une structure de catégorie de modèles dite \emph{projective}
(voir par exemple \cite[théorème 11.6.1]{hirschhorn2003model})
dont les équivalences faibles (resp. les fibrations) sont les
morphismes de foncteurs~$\alpha : F \to G$ tels que
pour tout objet $a$ de $A$, la flèche~$\alpha_a : Fa \to Ga$ est un
quasi-isomorphisme (resp. une fibration). 

\paragr Si $A$ est une petite catégorie, on note
\[
\DerHotab(A) = \W^{-1}\Homi\Big(A^{\op}, \Ch(Ab)\Big)
\]
\notindex{$\DerHotab(A)$}%
la catégorie homotopique associée, obtenue en localisant la catégorie
des foncteurs ${A}^{\op}\to \Ch(\Ab)$ par les morphismes de foncteurs
qui sont des quasi-isomorphismes argument par argument. On remarque
qu'on a un isomorphisme canonique $\DerHotab(e) \simeq \Hotab$, en
notant $e$ la catégorie ponctuelle.

Le foncteur diagramme constant 
\begin{align*}
\Ch(\Ab) &\to \Homi({A}^{\op},\Ch(\Ab)) \\
C &\mapsto \big(a \mapsto C\big)
\end{align*}
envoie les quasi-isomorphismes sur des équivalences faibles et
préserve les fibrations, et est donc un foncteur de Quillen à droite
pour cette structure de catégorie de modèles. Son adjoint à gauche 
\[
\limind\nolimits_{{A}^{\op}}^{\Ch(\Ab)} :
\Homi\big({A}^{\op},\Ch(\Ab)\big) \to \Ch(\Ab)
\]
est donc un foncteur de Quillen à gauche. Par conséquent, il admet un
foncteur dérivé à gauche noté 
\[
\hocolim\nolimits_{A^{\op}}^{\Hotab} : \DerHotab(A) \to \Hotab
\]
\notindex{$\hocolim$}%
appelé foncteur \ndef{limite inductive homotopique}, ou encore
\emph{colimite homotopique}. Ce foncteur est défini par la propriété
universelle suivante : on a un~$2$-carré 
\[
\xymatrix@=4em{
\Homi({A}^{\op}, \Ch(\Ab)) \ar[r]^-{\limind_{{A}^{\op}}^{\Ch(\Ab)}}
\ar[d]_{\gamma} & \Ch(\Ab)
\ar[d]^{\gamma} \\
\DerHotab(A) 
\ar@{}[ru]|(0.4){}="a"|(0.6){}="b"
\ar@{=>}"a";"b"^{\alpha_A} \ar[r]_-{\hocolim_{{A}^{\op}}^{\Hotab}} & \Hotab
} 
\]
où on a noté $\gamma$ les foncteurs de localisation, qui est
universel au sens suivant : pour tout foncteur $F:\DerHotab(A)\to\Hotab$
et pour tout morphisme de foncteurs $\beta : F \circ \gamma \to \gamma \circ
\limind_{{A}^{\op}}$, il existe une unique transformation naturelle
\[
  \eta : F \to \hocolim_{{A}^{\op}}^{\Hotab}
\]
telle qu'on ait $\beta = \alpha(\eta\star \gamma)$.

\begin{remark}
On peut aussi justifier l'existence du foncteur limite inductive
homotopique en utilisant les outils de l'algèbre homologique, puisque
la catégorie des foncteurs à valeurs dans une catégorie abélienne est
une catégorie abélienne, et que le foncteur limite inductive est exact
à droite. 
\end{remark}

\paragr\label{defH_A} Si $A$ est une petite catégorie et $X$ est un
préfaisceau en groupes abéliens sur $A$, on peut
voir $X$ comme un préfaisceau en complexes de chaînes de groupes
abéliens en le composant
avec l'inclusion $\Ab \hookrightarrow \Ch(\Ab)$ associant à un groupe
abélien $A$ le complexe concentré en degré $0$ associé. On dispose
ainsi d'un foncteur 
\[
\prefab{A} \xrightarrow{i} \Homi(A^{\op}, \Ch(\Ab)) \pbox{.}
\]
On définit alors le foncteur $\Hf{A}$ en posant
\begin{align*}
\Hf{A} : \prefab{A} &\to \Hotab \\
X &\mapsto \H{A}{X} := \hocolim_{A^{\op}}^{\Hotab}(iX) \pbox{.}
\end{align*}
\notindex{$\Hf{A}$}%
\notindex{$\H{A}{X}$}%
On dira que $\H{A}{X}$ est \ndef[homologie d'une catégorie à
coefficients!dans un préfaisceau]{l'homologie de $A$ à coefficients
dans le préfaisceau~$X$}, ou plus simplement
\emph{l'homologie du préfaisceau} $X$. En d'autres termes, le
foncteur $\Hf{A}$ est défini comme la composée des foncteurs
\[
\prefab{A} \hookrightarrow \Homi\big({A}^{\op}, \Ch(\Ab)\big)
\xrightarrow{\gamma} \DerHotab(A)
\xrightarrow{\hocolim_{{A}^{\op}}^{\Hotab}} \Hotab \pbox{.}
\]
\begin{example}
Par construction, dans le cas où $A$ est le groupoïde~$BG$ associé à
un groupe $G$, on a bien $\H{BG}{M}=\H{G}{M}$ pour
tout $G$\nobreakdash-module~$M$ vu comme un
préfaisceau abélien sur $BG$.
\end{example}

On peut définir, à partir de ce foncteur, un foncteur associant à
toute petite catégorie son homologie à coefficients dans $\Z$, de la
manière suivante.

\paragr Si $u : A \to B$ est un foncteur entre petites catégories, il
induit un foncteur 
\[
u^* : \Homi({B}^{\op},\Ch(\Ab)) \to \Homi({A}^{\op},\Ch(\Ab))  
\]
qui préserve les quasi-isomorphismes argument par argument, et induit
donc un foncteur
\[
\DerHotab(B) \xrightarrow{u^*} \DerHotab(A)
\xrightarrow{\hocolim_{{A}^{\op}}^{\Hotab}} \Hotab \pbox{.}
\]
On obtient de plus un morphisme canonique
\[
\hocolim_{{A}^{\op}}^{\Hotab} \circ\,\, u^* \circ \gamma
\xrightarrow{\alpha_A \star u^*} \gamma \circ
\limind\nolimits_{{A}^{\op}}^{\Ch(Ab)} \circ\,\, u^*
\to \gamma\circ \limind\nolimits_{{B}^{\op}}^{\Ch(\Ab)} \mdvirg
\]
en notant $\gamma$ les foncteurs de localisation, et où la deuxième
flèche est obtenue par la propriété universelle de la limite
inductive. Par propriété universelle de la limite inductive
homotopique, on obtient donc un morphisme canonique 
\[
\hocolim_{{A}^{\op}}^{\Hotab} \circ \,\,u^* \xrightarrow{\eta}
\hocolim_{{B}^{\op}}^{\Hotab} \pbox{.}
\]

\paragr\label{foncteurHomologieZ} Si $A$ est une petite catégorie, on note $\Z_A$ le préfaisceau
constant de valeur $\Z$ sur $A$. Pour simplifier les notations, on
notera souvent, pour toute petite catégorie $A$,
\[
\H{A}{\Z} := \H{A}{\Z_A} 
\]
\notindex{$\H{A}{\Z}$}%
lorsqu'aucune confusion n'en découle. Si $u : A \to B$ est un foncteur entre petites catégories, on
remarque qu'on a $u^*\Z_B=\Z_A$. 

Comme expliqué au paragraphe précédent, on obtient alors un morphisme
canonique 
\[
\H{A}{\Z} = \hocolim_{{A}^{\op}}^{\Hotab} u^* \Z_B \to
\hocolim_{{B}^{\op}}^{\Hotab} \Z_B = \H{B}{\Z} \pbox{.}
\]
On a ainsi défini un foncteur
\begin{align*}
\Cat &\to \Hotab \\
A &\mapsto \H{A}{\Z} 
\end{align*}
associant à toute petite catégorie son homologie, au sens introduit
précédemment, à coefficients dans le préfaisceau constant de valeur
$\Z$. On va montrer dans la proposition
\ref{homologieZ=HomologieSinguliere} que ce foncteur coïncide avec le
foncteur envoyant une petite catégorie sur le type d'homologie de son
nerf simplicial, et par conséquent qu'il envoie les équivalences de Thomason sur
des isomorphismes. 

\section{Homologie et catégorie des éléments}
On va maintenant montrer que la notion d'homologie introduite dans la
section précédente est la \og bonne notion \fg{} d'homologie au sens
des catégories test, c'est-à-dire qu'elle associe aux préfaisceaux
d'ensembles l'homologie du nerf de leur catégorie des éléments.

On commence par souligner que le foncteur $\Hf{A}$ est en
fait l'exact analogue abélien du foncteur $i_A$, du point de vue des
foncteurs dérivés, c'est-à-dire de rappeler que ce dernier calcule
également la limite inductive homotopique des diagrammes d'ensembles
simpliciaux discrets.

\begin{prop}[Quillen]\label{structureCatModKanQuillen}
La catégorie $\pref{\Delta}$ peut être munie d'une structure de
catégorie de modèles combinatoire dite de Kan-Quillen, dont 
les équivalences faibles sont les équivalences faibles
simpliciales et
les cofibrations les monomorphismes. 
\end{prop}
\begin{proof}
Voir \cite[chapitre II, section 3, théorème 3]{quillenhomotalg}.
\end{proof}

\paragr En utilisant le même raisonnement qu'au paragraphe
\ref{colimiteHomotopiqueDiagrammeCh}, si $A$ est une petite catégorie,
on peut donc justifier l'existence d'un foncteur 
\[
\hocolim_{{A}^{\op}}^{\Hot_{\Delta}} : \DerHot(A) \to \Hot_\Delta \mdvirg
\]
dérivé à gauche du foncteur 
\[
\limind\nolimits_{{A}^{\op}}^{\pref{\Delta}} : \Homi({A}^{\op},\pref{\Delta})
\to \pref{\Delta}
\]
\notindex{$\DerHot(A)$}%
où on a noté $\DerHot(A)$ la catégorie obtenue en localisant la
catégorie des foncteurs ${A}^{\op}\to\pref{\Delta}$ par les morphismes
de foncteurs qui sont des équivalences faibles simpliciales argument
par argument. En particulier, on a un isomorphisme $\DerHot(e) \simeq
\Hot_\Delta$.

Si $X$ est un préfaisceau d'ensembles sur $A$, on peut le voir comme
un préfaisceau en ensembles simpliciaux discrets 
\[
\tilde{X} : A^{\op} \rightarrow \Ens \hookrightarrow \Cat
\xrightarrow{\nerf} \EnsSimp \mdvirg
\]
où on note $\nerf$ le nerf simplicial. 

\begin{prop}[Thomason]\label{i_Ahocolim}
En gardant les notations du paragraphe précédent, si $X$ est un
préfaisceau sur une petite catégorie $A$, on a un isomorphisme naturel 
\[
\nerf i_A(X) \simeq \hocolim_{A^{\op}}^{\Hot_{\Delta}}\tilde{X} 
\]
dans $\Hot_\Delta$.
\end{prop}
\begin{proof}
Voir \cite[théorème 1.2]{thomason1979homotopy}.
\end{proof}

Pour montrer que le foncteur
\[
  \Hf{A} : \prefab{A} \to \Hotab
\]
introduit au paragraphe \ref{defH_A} envoie les préfaisceaux
d'ensembles sur le type d'homologie de leur catégorie des éléments, 
on va utiliser le fait que le foncteur associant à un
ensemble simplicial son complexe normalisé commute aux limites
inductives homotopiques, ce qui découle du lemme suivant.

\begin{lemme}\label{complexeNormaliseLibresQuillen}
En munissant les catégories $\pref{\Delta}$ et $\Ch(\Ab)$ des
structures de catégories de modèles des propositions
\ref{colimiteHomotopiqueDiagrammeCh} et
\ref{structureCatModKanQuillen}, le foncteur
\[
\pref{\Delta}\xrightarrow{\Whf{\Delta}}\prefab{\Delta}\xrightarrow{\dk}\Ch(\Ab)
\]
est un foncteur de Quillen à gauche.
\end{lemme}
\begin{proof}
Voir par exemple \cite[lemme 4.1.4]{guetta2021homologie}.
\end{proof}

\begin{prop}\label{homologiePrefaisceauxLibres}
Pour toute petite catégorie $A$, le diagramme suivant est commutatif à
isomorphisme naturel près : 
\begin{equation}\label{diag-colimitehomotopiquelibres}
\xymatrix{
\pref{A} \ar[r]^{i_A} \ar[rrd]_{\Whf{A}} & \Cat \ar[r]^{\nerf} & \Hot_{\Delta}
\ar[r]^{\Whf{\Delta}} & \Hotab^{\U}_{\Delta} \ar[r]^{\dk} & \Hotab \\
&& \prefab{A} \ar[rru]_{\Hf{A}} && \pbox{.}
} 
\end{equation}
\end{prop}
\begin{proof}
Si $X$ est un préfaisceau abélien sur $A$, on a la
chaîne d'isomorphismes naturels suivante dans $\Hotab$, en gardant les
notations de la proposition \ref{i_Ahocolim} : 
\begin{align*}
\dk\Whf{\Delta}\nerf i_A(X) 
&\simeq
\dk\Whf{\Delta}\hocolim_{{A}^{\op}}^{\Hot_{\Delta}} \tilde{X}
& (\ref{i_Ahocolim})
\\
&\simeq \hocolim_{{A}^{\op}}^{\Hotab} \dk \circ \Whf{\Delta} \circ \tilde{X}
& (\ref{complexeNormaliseLibresQuillen})
\end{align*}
Pour conclure, il suffit alors de remarquer que le diagramme de
complexes 
\[
{A}^{\op} \xrightarrow{\tilde{X}}  
\pref{\Delta} \xrightarrow{\Whf{\Delta}} 
\prefab{\Delta} \xrightarrow{\dk} 
\Ch(\Ab)
\]
correspond au diagramme associant à tout objet $a$ de $A$ le complexe
concentré en degré $0$ de valeur $\Z^{(Xa)}$, c'est-à-dire l'image de
$X$ par le foncteur $\Whf{A}$ suivi de l'inclusion en complexe de
chaînes concentré en degré $0$. On obtient finalement :
\begin{align*}
\dk\Whf{\Delta}\nerf i_A(X) 
&\simeq \hocolim_{{A}^{\op}}^{\Hotab} \dk \circ \Whf{\Delta} \circ \tilde{X}
\\
&\simeq \hocolim_{{A}^{\op}}^{\Hotab}(i\Wh{A}{X}) \\
&= \H{A}{\Wh{A}{X}} \pbox{.} \qedhere
\end{align*}
\end{proof}

\begin{coro}
Si $X$ est un ensemble simplicial, il existe un isomorphisme
canonique naturel
\[
\dk \Wh{\Delta}{X} \simeq \H{\Delta}{\Wh{\Delta}{X}}
\]
dans $\Hotab$.
\end{coro}
\begin{proof}
On applique la proposition précédente à la petite catégorie~$\Delta$.
Illusie prouve (en attribuant le résultat à Quillen) dans \cite[chapitre VI,
théorème~3.3]{illusie1971cotangent} que les foncteurs $\nerf$ et $i_{\Delta}$
induisent des équivalences de catégories quasi-inverses
\[
\nerf : \Hot \to \Hot_{\Delta} \mdvirg
i_{\Delta} : \Hot_{\Delta} \to \Hot \mdvirg
\]
et la commutativité du diagramme \ref{diag-colimitehomotopiquelibres}
implique donc bien le résultat énoncé.
\end{proof}

\begin{remark}
On prouvera qu'on a bien un isomorphisme naturel 
\[
\dk X \simeq \H{\Delta}{X} 
\]
dans $\Hotab$ pour tout groupe abélien simplicial $X$ au paragraphe
\ref{ex:IntegrateurLibreDelta}. Ce résultat est en fait déjà connu
de Bousfield et Kan \cite[chapitre XII,
paragraphe~5.5]{bousfield1972homotopy}, mais n'est pas un prérequis
dans ce travail.
\end{remark}

\paragr On rappelle que pour toute petite catégorie $A$, on note
$e_{\pref{A}}$ l'objet final de la catégorie $\pref{A}$. On remarque
alors que le préfaisceau abélien $\Whf{A}(e_{\pref{A}})$ coïncide avec
le préfaisceau constant de valeur $\Z$ sur $A$, que l'on note $\Z_A$
ou simplement $\Z$.

\begin{coro}\label{homologieZ=HomologieSinguliere}
Pour toute petite catégorie $A$, on a un isomorphisme naturel
\[
\H{A}{\Z} \simeq \dk\Wh{\Delta}{\nerf{A}}
\]
dans $\Hotab$. En particulier, le foncteur 
\begin{align*}
\Cat &\to \Hotab \\
A &\mapsto \H{A}{\Z_A}
\end{align*}
envoie les équivalences de Thomason sur des isomorphismes de $\Hotab$,
et induit donc un foncteur noté 
\[
\Whf{} : \Hot \to \Hotab 
\]
appelé \ndef[foncteur de Whitehead!absolu]{foncteur de Whitehead absolu}.
\end{coro}
\begin{proof}
Il suffit de remarquer que $i_A(e_{\pref{A}}) \simeq A$ avant
d'utiliser la commutativité du diagramme
\ref{diag-colimitehomotopiquelibres}. 
\end{proof}

\begin{coro}
Soit $A$ une petite catégorie. Alors \begin{enumerate}
\item pour tout préfaisceau $X$ sur $A$, on a un isomorphisme naturel
\[
\H{A}{\Wh{A}{X}} \simeq \dk\Wh{\Delta}{\nerf(\tranche{A}{X})} 
\]
dans $\Hotab$;
\item pour tout objet $a$ de $A$, on a un isomorphisme naturel 
\[
\H{A}{\Wh{A}{a}} \simeq \Z \leftarrow 0 \leftarrow 0 \leftarrow \cdots
\]
dans $\Hotab$.
\end{enumerate}
\end{coro}
\begin{proof}
La première proposition est simplement une reformulation de la
commutativité du diagramme \ref{diag-colimitehomotopiquelibres}. Pour
la deuxième proposition, il suffit de remarquer que puisque la
catégorie $\tranche{A}{a}$ admet $(a,\id_a)$ comme objet final, c'est
une catégorie asphérique en vertu de la proposition
\ref{propObjFinalOuInitialAspherique}. La première
proposition permet alors de conclure.
\end{proof}

\section{Catégories pseudo-test homologiques}\label{secPseudoTest}
Maintenant que l'on dispose d'un foncteur canonique, pour toute petite
catégorie $A$, associant un type d'homologie aux préfaisceaux en groupes
abéliens sur $A$, on peut partir à la \og quête des modèles des types
d'homologie \fg{}. Pour l'instant, seul l'équivalent de la notion de
catégorie pseudo-test (\ref{defPseudoTest}) a un sens. On introduira au
chapitre \ref{chapCatTestHomologiques} des équivalents aux notions de
catégories test, test faibles, locales et strictes.

\begin{defin}
On dit qu'un morphisme $ \varphi : X \to Y $ de préfaisceaux en
groupes abéliens sur $A$ est une \ndef[équivalence faible!de préfaisceaux!abéliens]{équivalence
test abélienne}, ou \emph{équivalence faible}, si son image 
\[
\H{A}{\varphi} : \H{A}{X} \to \H{A}{Y} 
\]
par le foncteur $\Hf{A}$ (\ref{defH_A}) est un isomorphisme dans $\Hotab$. 
\end{defin}

\paragr\label{defWab_A} On note $\Wab_A$ la classe des équivalences
test abéliennes de
\notindex{$\Wab_A$}%
$\prefab{A}$, et on note 
\[
\Hotab_A = \big(\Wab_A\big)^{-1} \prefab{A}
\]
\notindex{$\Hotab_A$}%
la catégorie obtenue en localisant la catégorie $\prefab{A}$ par les
équivalences test abéliennes. Le foncteur $\Hf{A}$ induit donc un
foncteur encore noté
\[
\Hf{A} : \Hotab_A \to \Hotab 
\]
et en vertu de la commutativité du diagramme \ref{diag-colimitehomotopiquelibres}, pour toute petite
catégorie $A$, on a un diagramme commutatif à isomorphisme près
\[
\xymatrix@=4em@C=4em{
\Hot_A \ar[r]^{\Whf{A}} \ar[d]_{i_A} & \Hotab_A \ar[d]^{\Hf{A}} \\
\Hot \ar[r]_{\Whf{}} & \Hotab \pbox{.}
} 
\]

On peut alors se poser la question de caractériser les petites
catégories $A$ telles que la flèche de droite est une équivalence de
catégories. La conjecture que Grothendieck énonce dans
\cite{pursuingstacks} est que si la flèche de gauche est une
équivalence de catégorie ($A$ est pseudo-test), il en est de même pour
la flèche de droite.
\begin{defin}\label{def:pseudoTestHomologique}
On dit qu'une petite catégorie $A$ est une
\ndef[catégorie!pseudo-test homologique]{catégorie pseudo-test
homologique} si le morphisme 
\[
\Hf{A} : \Hotab_A \to \Hotab  
\]
est une équivalence de catégories.
\end{defin}
Comme pour les catégories pseudo-test (\ref{defPseudoTest}), le
préfixe \og pseudo\fg{} souligne le fait qu'on n'impose pas de
contrainte au quasi-inverse du foncteur $\Hf{A}$. On introduira des
notions plus fortes (plus restrictives) dans le chapitre
\ref{chapCatTestHomologiques}.

\begin{example}\label{doldkancolimitehomotopique}
La catégorie $\Delta$ est une catégorie pseudo-test homologique. Cette
assertion revient à dire que le foncteur complexe normalisé $\dk$ calcule la
colimite homotopique, résultat déjà connu de Bousfield et Kan
(voir~\cite[chapitre XII, paragraphe 5.5]{bousfield1972homotopy}. On
donnera également une preuve alternative de ce résultat
en~\ref{ex:IntegrateurLibreDelta}.

On va aussi montrer de cette manière que la catégorie des globes
réflexifs~$\Grefl$
(bien qu'elle ne soit pas une catégorie pseudo-test) et la catégorie
des cubes, qui disposent toutes les deux d'un théorème de Dold-Kan
strict, sont des catégories pseudo-test homologiques (voir aux sections
\ref{secGlobes} et \ref{secCubes}).
\end{example}

\chapter{Homologie des préfaisceaux}\label{chapHomologiePrefaisceaux}

\section{Rappels d'algèbre homologique}

\paragr Si $\A$ est une catégorie abélienne, la catégorie
$\Chu(\A)$ des complexes de chaînes non bornés sur $\A$ est une
catégorie abélienne. La
catégorie des \ndef{complexe double} sur $\A$ est la catégorie
$\Chu(\Chu(\A))$. Autrement dit, un complexe double
$C_{\bullet,\bullet}$ dans $\A$ est un diagramme de la forme
\[
\xymatrix@=1.5em{
& \vdots\ar[d]&\vdots\ar[d]&\vdots\ar[d]
\\
\cdots
&C_{i-1,j+1} \ar[l]\ar[d]&C_{i,j+1}\ar[d]\ar[l]&C_{i+1,j+1}\ar[d]\ar[l]
&\cdots\ar[l]
\\
\cdots
&C_{i-1,j} \ar[l]\ar[d]&C_{i,j}\ar[d]\ar[l]&C_{i+1,j}\ar[d]\ar[l]
&\cdots\ar[l]
\\
\cdots
&C_{i-1,j-1} \ar[l]\ar[d] &C_{i,j-1}\ar[d]\ar[l] &C_{i+1,j-1}\ar[d]\ar[l]
&\cdots\ar[l]
\\
& \vdots&\vdots&\vdots
} 
\]
dans $\A$, où chaque ligne et chaque colonne est un complexe de
chaînes dans~$\A$, et où tous les carrés sont commutatifs. 
Si $C_{\bullet,\bullet}$ est un complexe double sur $\mathcal{A}$ de
différentielles
\[
  d^H : C_{\bullet,\bullet} \to C_{\bullet-1,\bullet} \mdvirg
  d^V : C_{\bullet,\bullet} \to C_{\bullet,\bullet-1} \mdvirg
\]
et si $\mathcal{A}$ possède des sommes directes ou des produits
infinis, on pose, pour tout entier $n\in \Z$, 
\[
\Tot^{\oplus}(C)_n = \bigoplus_{i+j=n}C_{i,j} \mdvirg 
\Tot^\pi(C)_n = \prod_{i+j=n} C_{i,j} \mdvirg
\]
\notindex{$\Tot^\oplus(C)$}%
\notindex{$\Tot^\pi(C)$}%
\notindex{$\Tot(C)$}%
la différentielle étant donnée pour tout entier $n$ et tout couple
$(i,j)$ d'entiers tels que $i+j=n$ par  
\[
d_{n}|_{C_{i+j}} = d^H + (-1)^id^V \pbox{.}
\]

On dit qu'un complexe double $C$ est \emph{de diagonale finie} si pour
tout entier~$n\in \Z$, l'ensemble des couples d'entiers $(i,j)$ tels que $i+j = n$ est
fini. Si $C$ est un tel complexe, les complexes $\Tot^\oplus(C)$ et
$\Tot^\pi(C)$ coïncident, et on notera parfois simplement~$\Tot(C)$ \emph{le}
complexe total de $C$.
\begin{prop}\label{prop:lemmeComplexeDouble}
Soit $u : C_{\bullet,\bullet}\to D_{\bullet,\bullet}$ un morphisme de
complexes doubles d'objets de $\A$. On suppose que $C$ et $D$ sont
deux complexes de diagonales finies. Si pour tout entier $n$, le morphisme
\[
  u_{n,\bullet} : C_{n,\bullet}\to D_{n,\bullet}
\]
est un quasi-isomorphisme, alors le morphisme $\Tot(u) : \Tot(C) \to
\Tot(D)$ est un quasi-isomorphisme.
\end{prop}
\begin{proof}
Voir \cite[théorème 12.5.4]{kashiwara2005categories}.
\end{proof}

\begin{coro}
Soit $C$ un complexe double de $\A$ de diagonale finie. Si chaque ligne de
$C$ est exacte, alors $\Tot(C)$ est un complexe exact.
\end{coro}

\paragr On rappelle qu'on dit qu'un foncteur $F : \A \to \A'$ entre
deux catégories abéliennes est \emph{exact à droite} si il commute aux
limites inductives \emph{finies}, et \emph{exact à gauche} si il
commute aux limites projectives finies.

On dit qu'un objet $P$ de $\A$ est \emph{projectif} si le foncteur 
\[
\Hom_\A(P,-) : \A \to \Ab
\]
est exact à droite. On peut montrer qu'il suffit pour cela que ce dernier
préserve les épimorphismes. On dit que la catégorie $\A$ possède \emph{assez
d'objets projectifs} si pour tout objet $M$ de $\A$, il existe un objet
projectif $P$ et un épimorphisme~$P \twoheadrightarrow M$.

\paragr Si $M$ est un objet de $\A$, une \emph{résolution projective} de
$M$ dans $\A$ est la donnée d'un complexe $P_\bullet$ d'objets projectifs de
$\A$ et d'un quasi-isomorphisme 
\[
P_\bullet \to M[0] \in \Fl(\Ch(\A))
\]
où on note $M[0]$ le complexe concentré en degré $0$ de valeur $M$. Un exercice classique consiste à montrer que si $A$ admet
assez d'objets projectifs, alors tout objet de $\A$ admet une
résolution projective. De plus, on peut montrer que les résolutions
projectives sont uniques à équivalence d'homotopie près (voir par
exemple \cite[théorème 2.2.6]{weibel1994Introduction}).

\paragr Étant donné un foncteur additif $F : \A \to
\A'$, on note également 
\[
F : \Ch(\A) \to \Ch(\A') 
\]
le foncteur obtenu en appliquant degré par degré le foncteur $F$ aux
complexes de chaînes en degré positif d'objets de $\A$. 

Si $\A$ possède assez d'objets projectifs et si $F$ est un foncteur
additif exact à droite, on peut montrer l'existence de son foncteur
dérivé à gauche, noté
\[
LF : D_+(\A) \to D_+(\A') \mdvirg
\]
\notindex{$LF$}%
et caractérisé par la propriété universelle suivante. On dispose
d'un diagramme 
\[
\xymatrix{
\Ch(\A) \ar[d]_{\gamma}  \ar[r]^{F} & \Ch(\A') \ar[d]^{\gamma} \\
D_+(\A) \ar[r]_{LF} 
\ar@{}[ur]|(.35){}="aa"|(.65){}="zz"
\ar@2"aa";"zz"^{\alpha} & D_+(\A')
}
\]
où les flèches verticales sont les foncteurs de localisation, qui est
universel au sens où
pour tout foncteur $G : D_+(\A) \to D_+(\A')$ muni d'une
transformation naturelle $\beta : G \circ \gamma \to \gamma \circ F$,
il existe une unique transformation naturelle~$\eta : G \to LF$ telle
que $\beta = \alpha(\eta\star\gamma)$. 

Si $M$ est un objet de $\A$, on peut calculer la valeur de $LF$ en $M$
de la manière suivante : pour toute résolution projective $P_\bullet$
de $M$, on a un isomorphisme canonique
\[
LF(M) \simeq F(P_\bullet)
\]
dans $D_+(\A')$. 

On renvoie à la littérature classique (par exemple
\cite{weibel1994Introduction} et \cite{kashiwara2005categories}) pour
plus de détails sur les foncteurs dérivés, ainsi que pour les notions
duales de résolutions injectives, et de foncteurs dérivés à droite
pour les foncteurs exacts à gauche.

\section{Extensions de Kan et préfaisceaux abéliens}\label{secExtKan}
Dans toute cette section, on fixe une petite catégorie $A$. On va
montrer que pour toute catégorie additive cocomplète $\M$, l'extension de Kan à gauche des foncteurs $A\to \M$ le long du
foncteur
\[
  A \hookrightarrow \pref{A} \xrightarrow{\Whf{A}}\prefab{A}
\]
existe, et peut être calculée en passant par l'intermédiaire de la
catégorie additive libre sur $A$.

\paragr \label{Add(F)} 
Si $A$ est une petite catégorie, on définit son \ndef[enveloppe
additive]{enveloppe additive}, notée~$\Add(A)$, par la
propriété universelle suivante : pour toute
catégorie additive~$M$, on a un isomorphisme naturel
\[ 
\Homadd(\Add(A), M) \simeq \Homi(A, M ) 
\] 
\notindex{$\Homadd$}%
\notindex{$\Add(A)$}%
où $\Homadd(A,M)$ désigne la sous-catégorie pleine des foncteurs qui commutent aux
sommes directes finies. 
Cela signifie que $\Add(A)$ est une catégorie
munie d'un foncteur $A \to \Add(A)$, et que pour toute catégorie
additive $M$ et pour tout foncteur $F : A\to M$, il existe un unique
foncteur commutant aux sommes directes finies noté $\Add(F) :
\Add(A) \to M$ tel que le diagramme 
\[
\xymatrix{
A \ar[d]^{} \ar[r]^{F} & M \\
\Add(A) \ar[ru]_{\Add(F)} 
} 
\]
soit commutatif. 

\paragr On peut réaliser l'enveloppe additive de $A$ comme
la sous-catégorie pleine de $\prefab{A}$ formée des sommes directes
\emph{finies} de préfaisceaux en groupes abéliens représentables,
c'est-à-dire des préfaisceaux de la
forme 
\[
\bigoplus_{i\in I}\Wh{A}{a_i} : x \mapsto \bigoplus_{i\in
I}\Z^{(\Hom_A(x,a_i))} 
\]
où $I$ est un ensemble \emph{fini} (on renvoie au paragraphe
\ref{WhAaa} pour la notation~$\Whf{A}$). En effet, si $M$ est une catégorie
additive et si~$F : A \to M$ est un foncteur, on peut définir le
foncteur 
\[
\Add(F) : \Add(A) \to \M 
\]
\notindex{$\Add(F)$}%
en posant, pour tout préfaisceau abélien de la forme
$\bigoplus_{i\in I} \Wh{A}{a_i}$,
\[
\Add(F) \big(\bigoplus_{i\in I} \Wh{A}{a_i}\big) = \bigoplus_{i\in I}
F(a_i) \pbox{.}
\]
Un morphisme dans $\prefab{A}$ de la forme 
\[
  \bigoplus_{i\in I} \Wh{A}{a_i} \to \bigoplus_{j\in J}
  \Wh{A}{b_j}
\]
correspond par le lemme de Yoneda à
la donnée, pour tout $i\in I$, d'un élément du groupe abélien $\bigoplus_{j\in
J}\Z^{(\Hom_A(a_i,b_j))}$, et donc à la donnée d'une matrice 
\[
  u = 
\begin{bmatrix}
u_{0,0} & u_{0,1} & \cdots & u_{0,n} \\
u_{1,0} & u_{1,1} & \cdots & u_{1,n} \\
\cdots & \cdots & \cdots & \cdots \\
u_{m,0} & u_{n,1} & \cdots & u_{m,n} \\
\end{bmatrix} 
\mdvirg
\begin{aligned}
 m&= |J| \\
 n&= |I|
\end{aligned}
\]
où chaque coefficient $u_{j,i}$ est une combinaison linéaire de
morphismes $a_i \to b_j$ dans $A$. En étendant le foncteur $F$ par
linéarité et en l'appliquant à chaque coefficient, on obtient alors
une nouvelle matrice dont le coefficient $(j,i)$ est dans
$\Z^{(\Hom_M(F(a_i),F(b_j))}$, c'est-à-dire un morphisme
\[
\Add(F)(u) : \bigoplus_{i\in I}F(a_i) \to
\bigoplus_{j\in J}F(b_j) \pbox{.}
\]
dans $\M$. 

La composition de deux morphismes $M$ et $N$ dans $\Add(A)$ est donnée
par la multiplication matricielle $N\times M$, et on vérifie sans
difficultés que $\Add(F)$ est bien un foncteur.

\paragr
L'enveloppe additive de $A$ va nous permettre de décrire l'extension
de Kan à gauche le long du foncteur
$\Whf{A} : A \to \prefab{A}$, et donc d'étendre un foncteur
$F:A\to M$ en un foncteur $F_!:
\prefab{A}\to M$. On a déjà vu qu'on peut étendre le morphisme $F$
à l'enveloppe additive de $A$, et il reste à étendre ce nouveau
morphisme aux préfaisceaux en groupes abéliens de sorte que le
diagramme suivant soit commutatif :
\[
\xymatrix@C=4em{
A \ar[d]_{\Whf{A}} \ar[rd]^{F}\\
\Add(A) \ar@{^{(}->}[d]^{} \ar[r]^{\Add(F)} & M \\  
\prefab{A} \ar[ru]_{F_!} & \pbox{.}
} 
\]

\begin{prop}\label{densiteAdd}
La catégorie $\Add(A)$ est dense dans $\prefab{A}$.
\end{prop}
\begin{proof}
Il s'agit de montrer que pour tout préfaisceau en groupes abéliens $X
\in \prefab{A}$, on a un isomorphisme naturel
\[
X \simeq \!\!\!\!\!\underset{L \in
\Add(A)/X}{\limind\nolimits^{\prefab{A}}}\!\!\!\!L \pbox{.}
\]
Soit alors un cône inductif de sommet $Y$ pour le foncteur
$\tranche{\Add(A)}{X}\to\prefab{A}$, c'est-à-dire la donnée d'un objet $Y$ de
$\prefab{A}$ muni, pour tout objet $(L, L\xrightarrow{s}X)$ de
$\tranche{\Add(A)}{X}$, d'un morphisme $u_{s} : L \to Y$ de sorte que
pour tout morphisme $f:(L,s)\to(L',s')$, le triangle 
\[
\xymatrix@=2em{
{L} \ar[rr]^{f} \ar[rd]_{u_{s}} && {L'} \ar[ld]^{u_{s'}} \\
& {Y}
} 
\]
soit commutatif. On doit montrer qu'il existe un unique morphisme
$\tilde{u}:X\to Y$ de $\prefab{A}$ de sorte que pour tout $(L,s)\in
\tranche{\Add(A)}{X}$, le triangle 
\[
  \label{diag-DensiteEnveloppeAdditive}
\xymatrix@=2em{
  {L} \ar[rr]^{s} \ar[rd]_{u_{s}} && {X} \ar[ld]^{\tilde{u}} \\
  & {Y}
}
\]
soit commutatif. Pour l'unicité, on remarque que la commutativité des
triangles de la forme
\[
\xymatrix@=2em{
{\Wh{A}{a}} \ar[rr]^{s} \ar[rd]_{u_{s}} && {X} \ar[ld]^{\tilde{u}} \\
& {Y}
}  
\]
impose, par le lemme de Yoneda additif, de définir le morphisme $\tilde{u}$ de la manière suivante :
pour tout objet $a$ de $A$, et pour tout objet $s$ de $Xa \simeq
\Hom(\Wh{A}{a},X)$, on pose $\tilde{u}(s)=u_{s}\in
\Hom(\Wh{A}{a},Y)\simeq Ya$. 
On vérifie alors
que si $(L,s)=(\bigoplus_{i\in I}\Wh{A}{a_i},(s_i)_{i\in I})$, le triangle
\[
\xymatrix@=2em{
{\bigoplus_i\Wh{A}{a_i}} \ar[rr]^-{(s_i)_i} \ar[rd]_{u_{s}} && {X}
\ar[ld]^{\tilde{u}} \\
& {Y}
}  
\]
est bien commutatif : il suffit de
vérifier que la restriction à $\Wh{A}{a_i}$ de $u_{s}$ est bien le
morphisme $u_{s_i}$, ce qui est
impliqué par la commutativité du triangle
\[
\xymatrix@=2em{
{\Wh{A}{a_i}} \ar[rr]^-{in_i} \ar[rd]_{s_i} && {L} \ar[ld]^{s} \\
& {X} & \pbox{.} } 
\]
Il reste à vérifier que, pour tout objet $a$ de $A$, $\tilde{u}_a : Xa \to Ya$
est bien un morphisme de groupes. Si $x_1$ et $x_2$ sont deux éléments de
$Xa$, on a un triangle commutatif 
\[
\xymatrix@=2em{
{\Wh{A}{a}} \ar[rr]^-{\Delta} \ar[rd]_{x_1+x_2} &&
{\Wh{A}{a}\oplus\Wh{A}{a}} \ar[ld]^{(x_1,x_2)} \\
& {X}
} 
\]
en notant $\Delta$ la diagonale, qui induit donc un triangle
commutatif 
\[
\xymatrix@=2em{
{\Wh{A}{a}} \ar[rr]^-{\Delta} \ar[rd]_{\tilde{u}(x_1+x_2)} &&
{\Wh{A}{a}\oplus\Wh{A}{a}}
\ar[ld]^[l]{u_{(x_1,x_2)}} \\
& {Y}
} 
\]
et puisque $u_{(x_1,x_2)} =
(u_{x_1},u_{x_2})$, on peut finalement conclure grâce à la
commutativité du triangle 
\[
\xymatrix@=2em{
{\Wh{A}{a}} \ar[rr]^-{\Delta} \ar[rd]_{\tilde{u}(x_1)+\tilde{u}(x_2)} &&
{\Wh{A}{a}\oplus\Wh{A}{a}} \ar[ld]^{(u_{x_1},u_{x_2})} \\
& {Y}
} 
\]
que $\tilde{u}$ est bien un morphisme de groupes, ce qui termine la
démonstration.
\end{proof}

\paragr\label{kanext} 
Soit $\M$ une catégorie additive cocomplète. 
Étant donné un foncteur~$F : A \to \M$, la proposition
précédente permet de l'étendre en un foncteur~$F_{!} : \prefab{A} \to
\M$ en posant
\[
  F_{!} X = \!\!\!\sideset{}{^\M}\limind_{L \in
  \Add(A)/X}\!\!\!\Add(F)(L) \pbox{.}
\]
On va montrer dans ce qui suit que ce foncteur commute aux limites inductives.

\begin{remark}
Pour tout objet $L$ de $\Add(A)$, on a un isomorphisme naturel 
\[ 
  F_! (L) \simeq \Add(F)(L) \mdvirg 
\] 
dans $\M$, puisque la catégorie~$\tranche{\Add(A)}{L}$ admet
$(L,\id_L)$ comme objet final.
Le diagramme 
\[
\xymatrix@C=4em{
A \ar[d]_{\Whf{A}} \ar[rd]^{F}\\
\Add(A) \ar@{^{(}->}[d]^{} \ar[r]^{\Add(F)} & \M \\  
\prefab{A} \ar[ru]_{F_!}
} 
\]
est donc bien commutatif.
\end{remark}

\begin{prop}\label{kanextadj}
Soit $\M$ une catégorie additive cocomplète. Pour tout foncteur $F :
A \to \M$, on a un couple de foncteurs adjoints 
\[
  F_! : \prefab{A} \to \M \mdvirg F^* : \M \to \prefab{A} 
\]
où $F^*$ est défini pour tout objet $x$ de $\M$ par 
\begin{align*} 
F^*(x) : A^{\op} &\to \Ab \\ a & \mapsto \Hom_\M(F(a), x) \pbox{.}
\end{align*}
\end{prop}
\begin{proof}
Soient $x$ un objet de $\M$ et $X$ un préfaisceau abélien sur $A$. 
On a alors une suite d'isomorphismes naturels :
\begin{align*}
\Hom_\M\left(F_!\,X, x \right)
& \simeq \Hom_\M\Big(\!\!\!\!\!
\sideset{}{^\M}\limind_{L \in \tranche{\Add(A)}{X}} \!\!\!\Add(F)(L),\, x
\Big)
\\ & \simeq \!\!\sideset{}{^\Ab}\limproj_{\substack{L \in \tranche{\Add(A)}{X} \\
\mathclap{L = \oplus_i\Wh{A}{a_i}}}} 
\Hom_\M\Big(\bigoplus_i F(a_i),\, x\Big)
\\ & \simeq \!\!\sideset{}{^\Ab}\limproj_{\substack{L \in \Add(A)/X\\ 
\mathclap{L = \oplus_i\Wh{A}{a_i}}}} \prod_i F^*(x)(a_i)
\\ & \simeq \!\!\sideset{}{^\Ab}\limproj_{\substack{L \in \tranche{\Add(A)}{X} \\ 
\mathclap{L = \oplus_{i}\Wh{A}{a_i}}}} \prod_i
\Hom_{\prefab{A}}\big(\Wh{A}{a_i}, F^*(x)\big)
\\ & \simeq \Hom_{\prefab{A}}\Big(\sideset{}{^{\prefab{A}}}\limind_{L \in
\Add(A)/X}\!\!\!L,\,\,
F^*(x)\Big)
\\ & \simeq \Hom_{\prefab{A}} \left(X, F^*(x)\right) \pbox{.}\qedhere
\end{align*}
\end{proof}

\begin{prop}\label{kanextEqCat}
Pour toute catégorie additive cocomplète $\M$, le foncteur de restriction 
\[ 
\Homi_!(\prefab{A},\M) \to \Homi(A,\M) \mdvirg F \mapsto F \circ \Whf{A}
\mdvirg
\]
où $\Homi_!(A,\M)$ désigne la
\notindex{$\Homi_\bang$}%
sous-catégorie pleine des foncteurs qui commutent aux petites limites
inductives,
est une équivalence de catégories dont un quasi-inverse est donné par 
\[
\Homi(A,\M) \to \Homi_!(\prefab{A},\M) \mdvirg F \mapsto F_! \pbox{.}
\]
\end{prop}
\begin{proof}
Par construction, si on a un foncteur $F :A \to \M$ et un objet $a$ de
$A$, alors
$F_! \circ \Wh{A}{a} \simeq F(a)$. Inversement, si $F : \prefab{A} \to
\M$ est un foncteur commutant aux petites limites inductives, il suffit
de vérifier que les foncteurs $F$ et $\Add(F\circ\Whf{A})$ coïncident
sur les représentables, ce qui est immédiat. 
\end{proof}

\section{Copréfaisceaux en groupes abéliens} 
Grothendieck introduit dans \emph{Pursuing Stacks}
\cite{pursuingstacks} un formalisme permettant d'étudier les
propriétés homologiques des catégories de préfaisceaux en groupes
abéliens. La construction qu'il obtient est équivalente à ce qui est
connu sous le nom de produit tensoriel de foncteurs
(voir par exemple~\cite[chapitre IX, section 6]{maclane}), à la différence que
celle-ci est entièrement exprimée dans le langage des préfaisceaux (et
donc des modèles des types d'homotopie): un foncteur d'une petite
catégorie $A$ dans
$\Ab$ sera traité comme un préfaisceau sur la catégorie opposée à $A$,
qu'il note simplement $B$. L'idée est de considérer qu'au lieu d'un
simple jeu de dualité formel, on dispose de deux catégories $A$ et $B$
vérifiant la \emph{propriété} $B={A}^{\op}$, ou plus particulièrement
que les préfaisceaux abéliens sur $B$ coïncident avec les foncteurs commutant
aux limites inductives~$\prefab{A} \to \Ab$.

\bigskip
\emph{
Dans toute cette section, on fixe une petite catégorie $A$, et on note
$B = A^{\op}$ la catégorie opposée. 
}

\paragr \label{couplageOdot}
Le foncteur de restriction
associe alors à un foncteur~$F:\prefab{A}\to\Ab$ le préfaisceau
$F\circ\Whf{A} : a \mapsto F(\Wh{A}{a})$ sur $B$. Inversement, étant
donné un préfaisceau $Y$ sur $B$, on peut l'étendre grâce à la
proposition \ref{kanext} en un foncteur
commutant aux limites inductives $Y_! : \prefab{A}\to A$.

\begin{prop}\label{prefabBeqcat}
Le foncteur de restriction 
\[\Homi_!(\prefab{A}, \Ab) \to \prefab{B}\]
est une équivalence de catégories, dont un quasi-inverse est donné par 
\begin{align*} 
\prefab{B} &\to \Homi_!(\prefab{A}, \Ab) \\
Y &\mapsto Y_! \pbox{.}
\end{align*}
\end{prop}
\begin{proof}
Puisque $\prefab{B}=\Homi(A,\Ab)$, il s'agit simplement d'un cas
particulier de la proposition \ref{kanextEqCat}
\end{proof}

Cette remarque nous permet d'introduire un couplage entre la catégorie
des préfaisceaux abéliens sur $A$ et celle des préfaisceaux abéliens
sur $B$, noté 
\[
- \odot_A - : \prefab{A} \times \prefab{B} \to \Ab 
\] 
\notindex{$-\odot -$}%
\termindex{produit tensoriel de foncteurs}%
défini par le diagramme ci-dessous : 
\[ 
\xymatrix{
{\prefab{A}\times\prefab{B} } \ar[r]^{- \odot_A -}
\ar[d]^{\simeq}_{\id_{\prefab{A}}\times(-)_! } & \Ab \\
{\prefab{A}\times\Hom_!(\prefab{A},\Ab)} \ar[ru]_{\ev} & \pbox{.} }
\]
Explicitement, si $X$ est un préfaisceau abélien sur $A$ et $Y$ est un
préfaisceau abélien sur $B = {A}^{\op}$, on a :
\[ 
X \odot_A Y := Y_! (X) 
\simeq \sideset{}{^\Ab}\limind_{\substack{
L\in\tranche{\Add(A)}{X}}}\!\!\!\Add(Y)(L)\pbox{.}
\]

\begin{prop}\label{prop:adjCoprefAb}
Pour tout préfaisceau abélien $Y$ sur $B={A}^{\op}$, on a un couple de foncteurs
adjoints
\[ - \odot_A Y \ : \prefab{A} \to \Ab \mdvirg Y^* : \Ab \to \prefab{A} \]
où l'adjoint à droite $Y^* : \Ab \to \prefab{A}$ est défini pour tout
groupe abélien $M$ par 
\begin{align*} 
Y^* (M) : A^{\op} &\to \Ab \\
a &\mapsto \Hom_{\Ab}(Y(a),M) \pbox{.}
\end{align*}
\end{prop}
\begin{proof}
Il s'agit d'un cas particulier de la proposition \ref{kanextadj}.
\end{proof}

\begin{remark}\label{remarqueHomologieFoncteurs}
On peut donner une expression générale au bifoncteur
\[
  - \odot_A - : \prefab{A} \times \prefab{B} \to \Ab
\]
en utilisant le formalisme des cofins, et on retrouve la notion de
\emph{produit tensoriel de foncteurs}, défini pour $F :
{A}^{\op}\to\Ab$ et $G:A\to\Ab$ par 
\[
F\otimes_A G = \int^{b \in B} Fb \otimes Gb \pbox{.}
\]
Ce formalisme est notamment utilisé pour les calculs d'homologie des
foncteurs (voir par exemple \cite{pirashvili2003introduction} et
\cite[appendice A]{djament2010homology}). Les
deux points de vue sont équivalents, et nous choisissons ici d'exposer le
formalisme introduit par Grothendieck dans \cite{pursuingstacks}.
\end{remark}

\begin{remark}
La catégorie $A$ (et donc $B={A}^{\op}$) étant fixée, on peut
intervertir le rôle de $A$ et $B$ et définir un autre couplage, noté
\[
- \odot_B - : \prefab{B} \times \prefab{A} \to \Ab 
\]
de la même manière que précédemment, c'est-à-dire en utilisant le
diagramme suivant :
\[ 
\xymatrix{
{\prefab{B}\times\prefab{A} } \ar[r]^{- \odot_B -}
\ar[d]^{\simeq}_{\id_{\prefab{B}}\times(-)_!} & \Ab \\
{\prefab{B}\times\Hom_!(\prefab{B},\Ab)} \ar[ru]_{ev} & \pbox{.} }
\]
\end{remark}

\begin{prop}\label{kanextsym}
Si $X$ est un préfaisceau abélien sur $A$ et $Y$ est un préfaisceau
abélien sur $B={A}^{\op}$, il existe un isomorphisme canonique naturel
\[ 
X \odot_A Y \simeq Y \odot_B X 
\]
dans $\Ab$.
En d'autres termes, si $X$ est un préfaisceau abélien sur $A$ et si~$Y :
A \to \Ab$ est un foncteur, les deux groupes abéliens suivants sont
isomorphes :
\begin{itemize}
\item l'extension de Kan à gauche de $Y$  le long de
$\Whf{A} : A \to \prefab{A}$
évaluée en $X$; 
\item l'extension de Kan à gauche de $X$ (que l'on traite comme un
foncteur défini sur les représentables de $\prefab{B}$) le long de
$\Whf{B} : B \to \prefab{B}$ évaluée en $Y$. 
\end{itemize}

\end{prop}
\begin{proof}
Le foncteur \[ - \odot_A - : \prefab{A} \times \prefab{B} \to \Ab \]
commute aux petites limites inductives en chaque variable. En effet,
on a vu que pour $Y \in \prefab{B}$, le foncteur $ - \odot_A Y :
\prefab{A} \to \Ab$ est un adjoint à gauche. De plus, pour $X \in
\prefab{A} $, le foncteur $ X \odot_A -: \prefab{B} \to \Ab$ est la
composée de l'équivalence de catégories \ref{prefabBeqcat} et du
foncteur d'évaluation, qui commutent tous les deux aux limites
inductives. La situation est identique pour le foncteur~$- \odot_B -$,
qui commute aussi aux limites inductives en chaque variable. 
Il suffit donc de vérifier que ces deux bifoncteurs coïncident sur les
représentables. Or, pour tout objet $a$ de $A$ et tout objet $b$ de
$B$, on a 
\[ \Wh{A}{a} \odot_A \Wh{B}{b} \simeq \Z^{(\Hom_B(a,b))} \mdvirg B={A}^{\op}\]
et
\[  \Wh{B}{b} \odot_B \Wh{A}{a} \simeq \Z^{(\Hom_A(b,a))} \mdvirg
A={B}^{\op} \pbox{.} \qedhere\]
\end{proof}

On notera donc simplement $ - \odot - $ par la suite quand le contexte le permet.

\paragr On peut donner une description explicite de la counité
intervenant dans l'adjonction de la proposition
\ref{prop:adjCoprefAb}. Si $Y$ est un préfaisceau abélien sur~$B$, on
peut écrire $Y$ comme la limite inductive 
\[
Y \simeq \limind_{\substack{L \to Y \\ \mathclap{L \in \Add(B)}}} L \pbox{.}
\]
Pour tout groupe abélien $M$, on a alors des isomorphismes canoniques
naturels 
\[
Y^*(M) \odot Y 
\simeq \limind_{\substack{L \to Y \\\mathclap{ L \in \Add(B)}}}
\Add(Y^*(M))(L)
\simeq
\limind_{\substack{L\to Y \\\mathclap{L =
\bigoplus\limits_{i\in I}\Wh{B}{b_i}}}}\Hom(Y(b_i),M) \pbox{.}
\]
On peut alors définir un morphisme $Y^*(M)\odot Y \to M$ de la manière
suivante. Pour tout morphisme $u : L = \bigoplus_i\Wh{B}{b_i} \to Y$, et
pour tout morphisme~$f_j : Y(b_j)\to M$ avec $j\in I$, on dispose d'un
morphisme
\[
\Wh{B}{b_j} \to \bigoplus_{i\in I}\Wh{B}{b_i}\xrightarrow{u} Y 
\]
correspondant donc à un objet de $Y(b_i)$. On peut alors associer à un
tel objet son image par $f_j$ dans $M$. On vérifie que le morphisme
ainsi défini est bien naturel, et qu'il s'agit bien du morphisme
d'adjonction annoncé.

\section{Enveloppe additive infinie}
\label{defAddinf}
On sait maintenant que si $F : \prefab{A}\to \Ab$ est un foncteur
commutant aux limites inductives, on a un isomorphisme 
\[
  F \simeq
  -\odot Y : \prefab{A} \to \Ab
  \] 
  en notant $Y$ la restriction de $F$ aux objets de $A$. Pour
  dériver à gauche un tel foncteur, on peut donc chercher à dériver le
  produit tensoriel de foncteurs, c'est-à-dire introduire un
  bifoncteur 
  \[
  - \Lodot - : D_+(\prefab{A}) \times D_+(\prefab{B}) \to \Hotab
  \mdvirg B={A}^{\op} \pbox{.}
  \]
L'idée est alors d'utiliser la symétrie prouvée en \ref{kanextsym} afin de
montrer que pour dériver ce bifoncteur, il suffit de le dériver en une
variable. On pourra donc utiliser une résolution de $Y$ fixée pour
\og calculer \fg{} la valeur de $LF(X)$ pour tout objet $X$ de
$\prefab{A}$. Pour cela, on a besoin d'utiliser des résolutions
libres, c'est-à-dire des résolutions dans l'enveloppe additive. Mais
l'enveloppe additive finie ne fournit que des résolutions libres de
type fini, qui n'existent pas en toute généralité : on va donc
introduire l'enveloppe additive infinie, qui fournit des résolutions
libres, qui ne sont pas nécessairement de type fini.

\paragr Si $A$ est une petite catégorie, on définit son
\ndef[enveloppe additive!infinie]{enveloppe additive infinie} notée
$\Addinf(A)$, de manière similaire à
\notindex{$\Addinf(A)$}%
l'enveloppe additive finie, par la propriété universelle suivante :
pour toute catégorie additive $\M$ possédant des coproduits, on a un isomorphisme naturel
\[
\Homaddinf(\Addinf(A), \M) \simeq \Homi(A,\M) 
\]
\notindex{$\Homaddinf$}%
où $\Homaddinf$ désigne la catégorie des foncteurs commutant aux
coproduits. 

On peut la réaliser comme
la sous-catégorie pleine de $\prefab{A}$ formée des coproduits de
représentables, c'est-à-dire des préfaisceaux de la
forme 
\[
\bigoplus_{i\in I} \Wh{A}{a_i} : x \mapsto \bigoplus_{i\in
I}\Z^{(\Hom_A(x,a_i))}
\]
où cette fois, l'ensemble $I$ peut être infini. 

Un morphisme entre
deux préfaisceaux de cette forme
\[
\bigoplus_{i\in I} \Wh{A}{a_i} \to \bigoplus_{j\in J}\Wh{B}{b_j}
\]
est alors une matrice infinie dont les coefficients sont des
combinaisons linéaires de morphismes $a_i\to b_j$ de $A$, et où chaque colonne
possède seulement un nombre fini de coefficients non nuls. 

\bigskip
\emph{Dans toute cette section, on fixe une petite catégorie $A$, et
on note $B={A}^{\op}$ la catégorie opposée.}
\bigskip

Ainsi, si $F : A \to \M$ est un foncteur et si $\M$ est une catégorie
additive possédant des coproduits,
on peut encore une fois l'appliquer par linéarité à chaque
coefficient pour obtenir une matrice où chaque colonne possède
toujours seulement un nombre fini de coefficients non nuls, et le
foncteur~$F$ induit donc bien un foncteur 
\label{Addinf(F)}
\begin{align*}
\Addinf(F) : \Addinf(A) &\to \M\\
\bigoplus_{i\in I} \Wh{A}{a_i} &\mapsto \bigoplus_{i\in I}
F(a_i) \pbox{.}
\end{align*}
\begin{prop}\label{kanextpourEnveloppeAdditive}
Si $\M$ est une catégorie additive possédant des coproduits et
$F:A\to \M$ est un foncteur, alors l'extension de Kan à gauche le long
de $F$ (voir \ref{kanext}) coïncide avec le foncteur $\Addinf(F)$ :
pour tout objet $L$ de $\Addinf(A)$, on a un isomorphisme canonique
naturel
  \[
F_!(L)\simeq\Addinf(F)(L) \pbox{.}
  \]
\notindex{$\Addinf(F)$}%
\end{prop}
\begin{proof}
C'est un corollaire immédiat de la propriété universelle de
l'enveloppe additive infinie.
\end{proof}

\begin{prop}
Les objets de $\Addinf(A)$ sont projectifs dans $\prefab{A}$.
\end{prop}
\begin{proof}
Un morphisme $X \to Y$ de $\prefab{A}$ est un épimorphisme si et
seulement si son image par le foncteur $\U : \prefab{A}\to \pref{A}$
est un épimorphisme. Ainsi, si $a$ est un objet de $A$ et $f : X\to Y$
est un épimorphisme de $\prefab{A}$, 
l'existence du relèvement de $f$ dans le diagramme de gauche
ci-dessous
\[
\xymatrix{
& X \ar@{>>}[d]^{f} \\
\Wh{A}{a} \ar[r]^{} \ar@{-->}[ru]^{}  & Y & \mdvirg
} 
\xymatrix{
& \U X \ar@{>>}[d]^{\U f} \\
a \ar[r]^{} \ar@{-->}[ru]^{} & \U Y
} 
\]
est équivalente à l'existence d'un relèvement de $\U(f)$
dans le diagramme de droite. On vérifie immédiatement que
le morphisme $\U(f)$ admet bien un relèvement, et donc que $\Wh{A}{a}$
est un objet projectif de $\prefab{A}$.

Puisque les
coproduits d'objets projectifs sont projectifs, on obtient bien que
tous les objets de $\Addinf(A)$ sont projectifs.
\end{proof}
\begin{prop}
Tout objet de $\prefab{A}$ admet une résolution projective par des
objets de $\Addinf(A)$. 
\end{prop}
\begin{proof}
Si $X$ est un préfaisceau en groupes abéliens sur $A$, on pose
\[
  L_X=\bigoplus_{a\in\Ob(A)}\bigoplus_{x\in Xa}\Wh{A}{a}
\]
et on définit un morphisme $\varphi : L_X \to X$ sur chaque composante
$(a,x)$ par~$\varphi_{(a,x)}=x:\Wh{A}{a}\to X$, qui est
évidemment un épimorphisme. En itérant le processus à partir du noyau
de $\varphi$, on peut ainsi construire une résolution de $X$ par des
objets de $\Addinf(A)$.
\end{proof}
On a en fait une caractérisation générale des objets projectifs de la
catégorie~$\prefab{A}$ : il s'agit de l'enveloppe karoubienne de la
catégorie $\Addinf(A)$. 
\begin{prop}\label{projectifPrefab}
Les objets projectifs de la catégorie $\prefab{A}$ sont les rétractes
de coproduits de préfaisceaux représentables.
\end{prop}
\begin{proof}
Un rétracte d'un objet projectif est toujours projectif.
Réciproquement, si $X \in \prefab{A}$ est un objet projectif, alors on
obtient une section de l'épimorphisme $L_X\to X$ utilisé dans la
preuve de la proposition précédente.
\end{proof}

\begin{coro}
\label{coro:fonctPreserveProj}
Soient $C$ et $D$ deux petites catégories. On suppose que~$F :
\prefab{C}\to\prefab{D}$ est un foncteur tel que pour tout
objet $L$ de $\Addinf(C)$, le préfaisceau $F(L)$ est un objet
projectif de $\prefab{D}$. Alors $F$ préserve les objets projectifs.
\end{coro}
\begin{proof}
Si $P$ est un objet projectif de $\prefab{C}$, alors $P$ est un
rétracte d'un préfaisceau $L$ dans $\Addinf(C)$. Par fonctorialité,
$F(P)$ est un rétracte de $F(L)$, et est donc projectif si $F(L)$ est
projectif.
\end{proof}

\paragr
Considérons maintenant un foncteur $F : \prefab{A} \to \Ab$ commutant
aux limites inductives, que l'on souhaite dériver à gauche. On note
\[
  Y = F \circ \Whf{A} 
\]
le préfaisceau en groupes abéliens sur $B={A}^{\op}$ obtenu par
restriction de $F$ à~$A$. On rappelle (voir la proposition
\ref{prefabBeqcat}) que le foncteur
\[
  - \odot Y : \prefab{A} \to \Ab
\]
coïncide alors avec le foncteur $F$. Dériver à gauche le foncteur $F$
revient donc à dériver à gauche le foncteur $ - \odot Y$, et on peut
alors utiliser la symétrie démontrée en \ref{kanextsym} de la manière
suivante.

\begin{prop}\label{produittensorielderive}
Soient $X$ un préfaisceau abélien sur $A$ et $Y$ un préfaisceau
abélien sur $B={A}^{\op}$. Les foncteurs 
\begin{align*}
X \odot - : \prefab{B} &\to \Ab  \\
- \odot Y : \prefab{A} &\to \Ab
\end{align*}
admettent chacun un foncteur dérivé à gauche 
\begin{align*}
L(X \odot -) : \D_+(\prefab{B}) &\to \Hotab \mdvirg \\
L(- \odot Y) : \D_+(\prefab{A}) &\to \Hotab
\end{align*}
et on a un isomorphisme naturel
\[
L(X \odot -)(Y) \simeq L(- \odot Y)(X) 
\]
dans $\Hotab$.
\end{prop}
\begin{proof}
Il s'agit d'utiliser un argument classique d'algèbre homologique (voir par
exemple \cite[2.7]{weibel1994Introduction}). Considérons deux
résolutions projectives~$M_{\bullet}\xrightarrow{\epsilon} X$ et
$L_{\bullet} \xrightarrow{\eta} Y$ par
des objets de $\Addinf(A)$ et de $\Addinf(B)$ respectivement. On peut
alors former les complexes de groupes abéliens $M_{\bullet}\odot Y$,
$X\odot L_{\bullet}$, ainsi que le complexe double $M_{\bullet}\odot
L_{\bullet}$ obtenu en appliquant terme à terme le couplage $\odot$.
On va montrer que ces trois complexes sont quasi-isomorphes.

En considérant les complexes $M_{\bullet}\odot Y$ et $X\odot L_{\bullet}$
comme des complexes doubles concentrés en la première colonne ou en la
première ligne, les résolutions $\epsilon$ et $\eta$ fournissent des morphismes de
complexes doubles qui à leur tour induisent des morphismes de
complexes 
\[
\Tot (M_\bullet \odot Y) \xleftarrow{M \odot \eta} \Tot(M_\bullet
\odot L_\bullet)
\xrightarrow{\epsilon\odot L} \Tot(X \odot L_\bullet)
\]
dont on peut montrer qu'ils sont des quasi-isomorphismes
sous condition que les foncteurs 
\[
- \odot L_j \mdvirg M_i \odot - \mdvirg i,j\in \N
\]
soient exacts. 
Par symétrie, il suffit de montrer que $N \odot -$ est exact pour tout
objet $N=\bigoplus_{i\in I}\Wh{A}{a_i}$ dans $\Addinf(A)$. Mais
puisque, pour tout préfaisceau abélien $Y$ sur $B$, on a
\[
\left(\bigoplus_{i\in I}\Wh{A}{a_i}\right)\odot Y = \bigoplus_{i\in I}Y(a_i)
\]
et que les foncteurs d'évaluation et de somme directe sont exacts, on
peut bien conclure.
\end{proof}

\paragr On note alors (en notant toujours~$B$ la catégorie opposée à
$A$)
\[
- \Lodot - : D_+(\prefab{A}) \times D_+(\prefab{B}) \to \Hotab 
\]
le bifoncteur défini grâce à la proposition
\ref{produittensorielderive}. 

\begin{coro}\label{calculfoncteursderives}
Soit $ F : \prefab{A} \to \Ab $ un foncteur commutant aux limites
inductives. On note $Y$ le préfaisceau $F\circ \Whf{A}$ sur
$B={A}^{\op}$. Alors il existe un isomorphisme canonique naturel 
\[
LF(X) \simeq X \Lodot Y 
\]
dans $\Hotab$. En particulier, si $L$ est une résolution projective de
$Y$ dans $\prefab{B}$, alors on a un isomorphisme naturel
\[
LF(X) \simeq X \odot L_\bullet 
\]
dans $\Hotab$.
\end{coro}
\begin{proof}
C'est immédiat en utilisant la proposition \ref{prefabBeqcat} qui affirme
entre autres que le foncteur $- \odot Y$ coïncide avec le foncteur
$F$, et la proposition \ref{produittensorielderive} qui affirme qu'il
suffit de dériver le produit $\odot$ en une seule variable.
\end{proof}

En choisissant des résolutions à valeurs dans l'enveloppe additive
infinie, on
obtient une expression particulièrement simple des foncteurs dérivés
dans ce contexte.

\begin{coro}
On garde les notations du corollaire précédent. Alors
si $L$ est une résolution de $Y$ dans l'enveloppe additive
infinie de $B$, on a un isomorphisme naturel, pour tout préfaisceau
abélien $X$ sur $A$,
\[
LF(X) \simeq \Addinf(X)(L_{\bullet}) \mdvirg
\]
ce qui signifie qu'en notant, pour tout entier $n\geq 0$,
\[
L_n = \bigoplus_{i \in I_n} \Wh{B}{b_i} \mdvirg
\]
alors on a un isomorphisme naturel dans $\Hotab$
\[
LF(X) \simeq \bigoplus_{i \in I_{0}}X(b_i)
\leftarrow \bigoplus_{i \in I_{1}}X(b_i)
\leftarrow \bigoplus_{i \in I_{2}}X(b_i)
\leftarrow \cdots
\]
où la différentielle est obtenue en appliquant additivement le
préfaisceau $X$ sur~$A$  (vu comme 
foncteur~$B \to \Ab$) à la différentielle du complexe $L$ d'objets
de $\Addinf(B)$.
\end{coro}
\begin{proof}
En vertu de la proposition \ref{kanextpourEnveloppeAdditive}, il
s'agit simplement d'utiliser la description du foncteur~$\Addinf(X) :
\Addinf(B) \to \Ab$ donnée en \ref{Addinf(F)}.
\end{proof}

\section{Intégrateurs}\label{secintegrateurs}

\bigskip
\emph{Dans toute cette section, on fixe une petite catégorie $A$, et
on note $B={A}^{\op}$ la catégorie opposée.}
\bigskip

On s'intéresse maintenant au cas où le foncteur
que l'on souhaite dériver à gauche est le foncteur 
\[
\limind\nolimits_{{A}^{\op}}^\Ab : \prefab{A} \to \Ab \pbox{.}
\]
Par l'équivalence de catégories \ref{prefabBeqcat}, ce foncteur (qui
commute évidemment aux limites inductives) correspond au préfaisceau
abélien sur $B = {A}^{\op}$ obtenu par restriction aux objets de $A$
\[
  Y = \limind\nolimits_{{A}^{\op}}^\Ab \circ \Whf{A} \pbox{.}
\]
On va voir que $Y$ coïncide avec le préfaisceau abélien $\Z_B$
constant de valeur $\Z$ sur $B$, c'est-à-dire, en notant $e_{\pref{B}}$ l'objet final de la
catégorie $\pref{B}$, au préfaisceau abélien~$\Wh{B}{e_{\pref{B}}}$ sur $B$.

\begin{prop}\label{prop:limIndProdTensZ}
Pour tout préfaisceau en groupes abéliens $X$ sur ${A}$, on a un
isomorphisme canonique naturel de groupes abéliens
\[
\limind\nolimits_{{A}^{\op}}^\Ab X \simeq X \odot \Z_B \pbox{.}
\]
\end{prop}
\begin{proof}
Si $X$ est un préfaisceau en groupes abéliens sur $A$, on a une chaine
d'isomorphismes naturels dans $\Ab$ : 
\begin{align*}
X \odot \Z_B 
& \simeq X \odot \Whf{B} \left( e_{\pref{B}} \right)\\
& \simeq X \odot \Whf{B} \left( \limind\nolimits_{b\in B}^{\pref{B}} b \right)
\\ & \simeq X \odot \limind\nolimits_{b \in B}^{\prefab{B}} \Wh{B}{b} 
\\ & \simeq \limind\nolimits_{b \in B}^{\Ab} \left( X \odot \Wh{B}{b} \right)
\\ & \simeq \limind\nolimits_{b \in B}^{\Ab} X(b) \pbox{.} \qedhere
\end{align*}
\end{proof}

\paragr Ainsi, pour exprimer le foncteur
\[
  \Hf{A} : \prefab{A}\to\Hotab  
\]
qui a été défini au paragraphe
\ref{defH_A} comme foncteur dérivé à gauche du
foncteur limite inductive, il suffit donc en vertu de
\ref{produittensorielderive} de trouver une résolution projective du
préfaisceau constant de valeur $\Z$ sur $B$. 

Une expression particulièrement simple sera possible si on trouve une
telle résolution dans $\Addinf(B)$, d'après le corollaire
\ref{calculfoncteursderives}, et d'autant plus dans le cas où une
telle résolution est dans $\Add(B)$, l'enveloppe additive
\emph{finie}. 

\begin{defin} \label{defIntegrateurs}
En notant toujours $B$ la catégorie
opposée à $A$, on appelle \ndef[intégrateur]{intégrateur}
\footnote{Dans \emph{Pursuing Stacks}, Grothendieck commence par
appeler intégrateurs sur $A$ les résolutions projectives du préfaisceau
constant de valeur $\Z$ dans la catégorie $\prefab{A}$, qu'il décide
de nommer plus loin \emph{cointégrateurs} sur $A$ (c'est-à-dire,
intégrateurs sur ${A}^{\op}$). La transition précise se situe à~\cite[§103]{pursuingstacks}.}
sur $A$ une
résolution projective~$L_{\bullet}$ du préfaisceau $\Z_B$ constant de
valeur $\Z$ dans la catégorie abélienne $\prefab{B}$. On introduit de
plus la terminologie suivante \footnote{Grothendieck utilise dans
\emph{op. cit.} les termes d'intégrateur \emph{spécial} et
\emph{quasi-spécial}, que nous avons choisi de remplacer par les plus
évocateurs \emph{libre} et libre de
type fini.}: 
\begin{itemize}
\item si, pour tout $n$, le préfaisceau $L_n$ est dans $\Addinf(B)$,
on dit que $L$ est un intégrateur \ndef[intégrateur!libre]{libre};
\item si, pour tout $n$, $L_n$ est dans $\Add(B)$, on dit que $L$ est
un intégrateur \ndef[intégrateur!libre de type fini]{libre de type fini}.
\end{itemize}
\end{defin}
En d'autres termes, un intégrateur libre (resp. libre de type fini)
est une résolution du foncteur $A\to\Ab$ constant de valeur $\Z$ par
des sommes (resp. sommes finies) de foncteurs représentables.

\begin{example}
Dans le cas où la catégorie $B$ est le groupoïde $BG$ associé à un groupe
$G$, la donnée d'un intégrateur correspond à la donnée d'une
résolution projective du $G$-module $\Z$ par des modules projectifs,
et un intégrateur libre (resp. libre de type fini) correspond à une
résolution par des modules libres (resp. libres de type fini).
\end{example}

\begin{remark}
En complément de la remarque \ref{remarqueHomologieFoncteurs}, on
précise que l'homologie, au sens où on l'étudie ici, correspond donc à
l'homologie des foncteurs à coefficients entiers, telle qu'étudiée par
exemple dans \cite{pirashvili2003introduction} et
\cite{djament2010homology}.
\end{remark}

\begin{prop}\label{propExpressionHomologieDerivateurLibre}
Pour tout intégrateur $L$ sur $A$, on a un isomorphisme naturel dans
$\Hotab$
\[
\H{A}{X} \simeq L_!(X) \pbox{.}
\]
\notindex{$L_\bang$}%
Si $L$ est un intégrateur libre sur $A$, alors, en
notant pour tout entier $n\geq 0$
\[
L_n = \bigoplus_{i \in I_{n}} \Wh{B}{a_i} \mdvirg
\]
on a un isomorphisme naturel dans $\Hotab$:
\begin{align*}
\H{A}{X} 
&\simeq \bigoplus_{i \in I_{0}} X(a_i) \leftarrow \bigoplus_{i \in I_{1}}
X(a_i) \leftarrow \bigoplus_{i \in I_{2}} X(a_i) \leftarrow \cdots
\end{align*}
\end{prop}
\begin{proof}
En vertu de la proposition \ref{prop:limIndProdTensZ}, le foncteur
dérivé à gauche du foncteur limite inductive coïncide avec le produit
tensoriel dérivé évalué en~$\Z_B$. On applique alors simplement le
corollaire \ref{calculfoncteursderives}.
\end{proof}

\begin{example}\label{ex:IntegrateurLibreDelta}
On considère le complexe de chaînes 
\[
L= L_{\Delta} = 
   \Wh{{\Delta}^{\op}}{\Delta_0} \leftarrow
   \Wh{{\Delta}^{\op}}{\Delta_1} \leftarrow
   \Wh{{\Delta}^{\op}}{\Delta_2} \leftarrow \cdots  
\]
\notindex{$L_\Delta$}%
dans $\prefab{{\Delta}^{\op}}$, où la différentielle est donnée par la somme
alternée des faces. Concrètement, pour tout objet $\Delta_n$ de
$\Delta$, on a 
\[
L(\Delta_n) = \Z^{(\Delta_0 \to \Delta_n)} \leftarrow
\Z^{(\Delta_1 \to \Delta_n)}\leftarrow \Z^{(\Delta_2 \to \Delta_n)}
\leftarrow \cdots
\]
où on note $\Z^{(\Delta_k \to \Delta_n)}$ le groupe libre engendré par
\notindex{$\Z^{(a\to b)}$}%
$\Hom_{\Delta}(\Delta_k, \Delta_n)$, et dont la différentielle est le
morphisme défini pour $k\geq 1$ par 
\begin{align*}
d_k = \sum_{i=0}^{k}(-1)^i d_i^{k} : L(\Delta_n)_{k} \to
L(\Delta_n)_{k-1} \pbox{.}
\end{align*}
On va montrer que $L$ est un intégrateur libre de type fini, que l'on
appellera \emph{intégrateur libre standard} sur $\Delta$. 

En effet, puisque toutes les composantes de $L$ sont dans l'enveloppe
additive de $\Delta^{\op}$, il suffit de vérifier que $L$ est une
résolution du préfaisceau constant de valeur $\Z$ sur $\Delta^{\op}$.
Un argument classique pour montrer cela (sans
utiliser la correspondance de Dold-Kan) consiste à construire un
rétracte par déformation du complexe~$L\Delta_n$ sur le complexe
$\Z[0]$. On
note
\[
  r : L\Delta_n \to \Z[0]
\]
l'unique morphisme envoyant les générateurs de dimension $0$ sur
l'unité $1\in \Z$. On construit une section 
 \[
s_0 : \Z[0] \to L\Delta_n
 \]
de $r$ définie par le $0$-simplexe minimal de $\Delta_n$.
Explicitement, $s_0$ est le morphisme envoyant, en degré $0$,
l'élément $1 \in \Z$ sur le $0$-simplexe $\langle 0 \rangle : \Delta_0
\to \Delta_n$. On a bien $r\circ s_0 = \id_{\Z[0]}$, et on va
construire une famille de morphismes \[
 h_k : \Z^{(\Delta_k \to \Delta_n)} \to \Z^{(\Delta_{k+1} \to
 \Delta_n)} \mdvirg k\geq 0
 \]
de sorte que $h = (h_k)$ soit une homotopie de $s_0\circ r$
vers l'identité, c'est-à-dire que l'on ait 
\begin{align*}
\begin{cases}
h_{k-1}d_{k} + d_{k+1}h_k = \id & \mdvirg k > 1 \\
d_1h_0 = \id - s_0\circ r \pbox{.}
\end{cases}
\end{align*}

Pour tout entier $k\geq 0$, on définit $h_k$ de la manière suivante :
si
\[
  \varphi = \langle i_0 \to i_1 \to \cdots \to i_k\rangle
\]
est un $k$-simplexe de
$\Delta_n$ avec $0\leq i_0 \leq i_1 \leq \cdots \leq i_k \leq k$, on définit
$h_k(\varphi)$ comme le $k+1$-simplexe 
\[
  h_k(\varphi) = \langle 0 \to i_0 \to i_1 \to \cdots \to i_k\rangle
  \pbox{.}
\]
Par construction, pour tout $0$-simplexe $\langle i \rangle$ de $\Delta_n$, on a 
\[
d_1h_0(i) = d \langle 0 \to i \rangle = (\id - s_0r)(i) \mdvirg
\]
et pour tout entier $k \geq 1$, on a
\begin{align*}
d_i^{k+1}h_k = 
\begin{cases}
  h_{k-1}d_{i-1}^k & \mdvirg 1 \leq i \leq k \\
  \id & \mdvirg i=0\pbox{.}
\end{cases} 
\end{align*}
On obtient alors, pour tout entier $k > 1$,
\begin{align*}
h_{k-1}d_{k} + d_{k+1}h_k 
&= 
\sum_{i=0}^{k} (-1)^i h_{k-1}d_i^{k}
+
\sum_{j=0}^{k+1} (-1)^j d_j^{k+1}h_k
\\
&= 
\sum_{i=0}^{k} (-1)^i h_{k-1}d_i^{k}
+
\id + \sum_{j=1}^{k+1} (-1)^{j} h_{k-1}d_{j-1}^k 
\\
&=
\sum_{i=0}^{k} (-1)^i h_{k-1}d_i^{k}
+
\id + \sum_{j=0}^{k} (-1)^{j+1} h_{k-1}d_{j}^k 
\\
&= \id \pbox{.}
\end{align*}
On a bien montré $h$ est une homotopie de $s_0 \circ r$ vers
l'identité. Puisque $s_0$ est une section de $r$, ces deux morphismes
sont donc des équivalences d'homotopie. En particulier, le
complexe~$L_\Delta$ est donc bien un intégrateur libre sur $\Delta$.
\end{example}

On obtient donc une nouvelle preuve du résultat bien connu suivant.
\begin{coro}\label{propHomologieDeltaMoore}
Pour tout préfaisceau abélien $X$ sur $\Delta$, on a un isomorphisme
naturel 
\[
\H{\Delta}{X} \simeq X_0 \leftarrow X_1 \leftarrow X_2 \leftarrow \cdots
\]
entre l'homologie de $\Delta$ à coefficients dans $X$ et le complexe
non normalisé associé (\ref{defMooreComplex}) à $X$.
\end{coro}
\begin{proof}
Il suffit de remarquer que le foncteur ${L_\Delta}_!$ coïncide avec le
foncteur associant à un groupe abélien simplicial son complexe de
chaînes non normalisé, ce qui est une conséquence immédiate de la
proposition \ref{propExpressionHomologieDerivateurLibre}.
\end{proof}

On obtient grâce à l'exemple précédent une nouvelle preuve du fait que
la catégorie $\Delta$ est une catégorie pseudo-test homologique, puisqu'on sait
que le foncteur ${L_{\Delta}}_! : \Hotab_{\Delta} \to \Hotab$ coïncide
avec le foncteur complexe normalisé, qui est une équivalence de
catégories. Le foncteur $\Hf{\Delta} : \Hotab_{\Delta} \to \Hotab$
est donc bien une équivalence de catégories.

\bigskip
Dans ce qui suit, on va montrer comment construire des intégrateurs
sur le produit de deux petites catégories.

\paragr Si $X$ et $Y$ sont deux préfaiseaux abéliens sur une petite
catégorie~$C$, on note 
\[
X \otimes Y : c \mapsto Xc \otimes Yc 
\]
le préfaisceau abélien sur $C$ obtenu en appliquant le produit
tensoriel de groupes abéliens argument par argument.

\paragr On fixe deux petites catégories $C$ et $D$. On définit un
foncteur 
\begin{align*}
\boxtimes : \prefab{C} \times \prefab{D} &\to \prefab{C\times D} \\
(X,Y) &\mapsto \pr_1^*(X)\otimes\pr_2^*(Y)
\end{align*}
\notindex{$\boxtimes$}%
où 
\[
\xymatrix{
& C \times D \ar[ld]_{\pr_1}  \ar[rd]^{\pr_2}\\
C && D
} 
\]
désignent les projections. On remarque que si $c$ est un objet de $C$
et $d$ est un objet de $D$, on a alors 
\[
\Wh{C}{c}\boxtimes\Wh{D}{d} = \Wh{C\times D}{c,d} \pbox{.} 
\]

En particulier, si $L_C$ est un objet de $\Addinf(C)$ et $L_D$ est un
objet de~$\Addinf(D)$, alors $L_C\boxtimes L_D$ est un objet de
$\Addinf(C\times D)$. En passant aux catégories opposées, on obtient
alors le résultat suivant : 

\begin{prop}\label{integrateurProduit}
Si $L_C$ et $L_D$ sont deux intégrateurs sur $C$ et $D$, alors le
foncteur 
\begin{align*}
L_{C\times D} : C\times D &\to \Ch(\Ab)\\
(c,d) &\mapsto L_C(c)\otimes L_D(d)
\end{align*}
est un intégrateur sur $C\times D$.
\end{prop}
\begin{proof}
On vérifie que $L_{C\times D}$ est obtenu en appliquant degré par
degré le foncteur 
\[
  \boxtimes : \prefab{{C}^{\op}}\times
  \prefab{{D}^{\op}}\to \prefab{{C^{\op}\times D^{\op}}}
\]
au couple $(L_C,L_D)$. On peut alors utiliser une légère variation
du corollaire~\ref{coro:fonctPreserveProj} : si $L_C$ et $L_D$ sont
des intégrateurs libres, alors on sait que $L_{C\times D}$ est
un complexe dans $\Addinf({C}^{\op}\times {D}^{\op})$. Si $L_C$ et
$L_D$ sont deux intégrateurs quelconques, alors il existe deux
intégrateurs libres $F_C$ et $F_D$ sur $C$ et $D$ tels que $L_C$ est
un rétracte de $F_C$ et $L_D$ est un rétracte de $F_D$, et $L_{C\times
D}$ est un rétracte de $F_C\boxtimes F_D$.

Dans tous les cas, $L_{C\times D}$ est bien un complexe de
préfaisceaux projectifs. De plus, on a des
quasi-isomorphismes~$L_C \to \Z_{{C}^{\op}}$ et $L_D\to
\Z_{{D}^{\op}}$, qui induisent donc, par projectivité, un
quasi-isomorphisme $L_{C\times D}\to \Z_{{(C\times D)}^{\op}}$.
\end{proof}

 \begin{example}
Si $X$ est un groupe abélien bisimplicial, on
peut calculer l'homologie de $\Delta\times\Delta$ à coefficients dans
$X$ par le complexe total du complexe double obtenu en
appliquant colonne par colonne le foncteur complexe normalisé ou non
normalisé. 
\end{example}

\paragr\label{integrateurNormalise} Bien qu'on puisse toujours trouver un intégrateur libre sur
n'importe quelle petite catégorie $A$, on utilisera fréquemment des
intégrateurs qui ne sont pas libres : l'exemple
phare est donné par la correspondance de \hbox{Dold-Kan}. On rappelle que
pour tout groupe abélien simplicial $X$, la suite exacte de complexes
\[
0 \rightarrow \mathsf{D}X \rightarrow \mathsf{C}X \rightarrow \dk X
\rightarrow 0 \mdvirg
\]
où on a noté $\mathsf{D}X$ le sous-complexe des simplexes dégénérés du complexe
non normalisé~$\mathsf{C}X={L_\Delta}_!X$ associé à $X$ (voir
la section \ref{secDoldKan}), est scindée. En particulier, pour tout
objet $\Delta_i$ de $\Delta$, le préfaisceau abélien sur $\Delta^{\op}$
\[
  c_i : \Delta_n \mapsto \Z^{(\Delta_i \hookrightarrow \Delta_n)} =
  \dk(\Delta_n)_i
\]
est un rétracte de $\Wh{\Delta^{\op}}{\Delta_i}$, et $c : \Delta \to
\Ch(\Ab)$ est donc bien un complexe de préfaisceaux projectifs sur
${\Delta^{\op}}$. 

De plus, le rétracte par déformation de l'exemple
\ref{ex:IntegrateurLibreDelta} passe au quotient et fournit un
rétracte par déformation de $c\Delta_n$ sur~$c\Delta_0$, ce qui prouve
effectivement que $c$ est un intégrateur sur $\Delta$.

Pour montrer que le foncteur $\dk : \prefab{\Delta}\to \Ch(\Ab)$
calcule bien la colimite homotopique, sans utiliser le fait qu'il est
quasi-isomorphe au complexe non normalisé (ce qui par exemple est faux
en remplaçant la catégorie $\Delta$ par la catégorie cubique
$\square$), on peut alors simplement remarquer que puisque $\dk$ est
un adjoint à gauche, il commute aux
limites inductives, et il est donc isomorphe à l'extension de Kan à gauche de
sa restriction $c : \Delta \to \Ch(\Ab)$. Plus généralement, on peut
utiliser la proposition suivante : 
\begin{prop}\label{propIntegrateurCocontinuExtKan}
Soit $F : \prefab{A}\to \Ch(\Ab)$ un foncteur commutant aux limites
inductives dont la restriction $F\circ \Whf{A} : A \to \Ch(\Ab)$
soit un intégrateur sur $A$. Alors pour tout $X\in\prefab{A}$, on a un
isomorphisme naturel 
\[
FX \simeq \H{A}{X} 
\]
dans $\Hotab$.
\end{prop}
\begin{proof}
Il suffit de remarquer que si $F$ commute aux limites inductives, il
est isomorphe au foncteur $(F\circ\Whf{A})_!$.
\end{proof}
On utilisera cette proposition pour montrer que les autres théorèmes
de Dold-Kan stricts (pour les groupes abéliens cubiques dans la
section \ref{secCubes} et globulaires en~\ref{secGlobes}) donnent bien
lieu à des exemples de catégories pseudo-test homologiques
(\ref{def:pseudoTestHomologique}).

\begin{remark}
On a utilisé à plusieurs reprises le fait que pour toute petite catégorie $B$, la donnée d'un complexe de chaînes
$L_{\bullet}\in\Ch(\prefab{B})$ correspond à la donnée d'un foncteur
$L: {B}^{\op}\to \Ch(\Ab)$, et on
utilisera sans les différencier les deux points de vue. La notation
$L_{\bullet}$ permet d'utiliser le vocabulaire de l'algèbre
homologique, tandis que la notation fonctorielle permet d'utiliser les
outils des catégories de modèles. Aussi, on préfèrera dans la suite
introduire des foncteurs $A\to\Ch(\Ab)$ plutôt que des complexes
dans~$\Ch(\prefab{{A}^{\op}})$.
\end{remark}

\begin{prop}\label{propAdjIntegrateurs}
Soit $L : A \to \Ch(\Ab)$ un intégrateur sur $A$. Le foncteur 
\[
 -\odot L_{\bullet} : \prefab{A}\to\Ch(\Ab)
\]
admet comme adjoint à droite le foncteur
\notindex{$L^*$}%
\begin{align*}
L^* : \Ch(\Ab)&\to \prefab{A} \\
C &\mapsto \big( a \mapsto \Hom_{\Ch(\Ab)}(L(a),C)\big) \pbox{.}
\end{align*}
\end{prop}
\begin{proof}
C'est un cas particulier de la proposition \ref{kanextadj}. Il faut
seulement vérifier qu'on a bien $L_! = -\odot L_{\bullet}$,
c'est-à-dire que l'extension de Kan à gauche (\ref{kanext}) d'un
complexe de préfaisceaux se calcule degré par degré, ce qui est
impliqué par le fait que les limites inductives se calculent degré par
degré.
\end{proof}
On retrouve par exemple l'adjoint à droite au foncteur complexe
normalisé sur $\prefab{\Delta}$ détaillé dans la section
\ref{secDoldKan}. De plus, le foncteur associant à un groupe abélien
simplicial son complexe de chaînes non normalisé 
\[
{L_\Delta}_! : \prefab{\Delta}\to \Ch(\Ab) 
\]
admet pour adjoint à droite le foncteur 
\[
L_\Delta^* : \Ch(\Ab)\to \prefab{\Delta} 
\]
défini comme dans la proposition ci-dessus.

\paragr Il n'est pas clair, pour un intégrateur $L$ quelconque sur une
catégorie $A$, que le foncteur $L^*$ défini ci-dessus envoie les
quasi-isomorphismes sur des équivalences abéliennes, et donc qu'il
induise un foncteur 
\[
L^* : \Hotab \to \Hotab_A  \pbox{.}
\]
On verra dans la proposition \ref{coroQuillenIntegrateur} que si $A$ est
une catégorie test locale de Whitehead (\ref{defWhitehead}), alors
tout intégrateur libre vérifie cette propriété. 

\section{Homologie des petites catégories} \label{homologieCatIntegrateurs}
On rappelle (voir \ref{homologieZ=HomologieSinguliere}) que l'homologie des petites
catégories peut être calculée en utilisant la catégorie test $\Delta$,
c'est-à-dire qu'on a un isomorphisme naturel, pour toute petite catégorie $A$,
\[
\H{A}{\Z} \simeq \H{\Delta}{\Z^{(\nerf A)}} 
\]
dans $\Hotab$. On va détailler ici une manière de décrire le
foncteur d'abélianisation (ou, dans la terminologie de Grothendieck,
le \og foncteur de Whitehead absolu \fg{})
\[
\Whf{} : \Hot \to \Hotab 
\]
à l'aide du formalisme des intégrateurs. On a défini ce foncteur comme
associant à un catégorie $A$ son homologie à coefficients dans $\Z$
\begin{align*}
\Whf{}:\Hot &\to \Hotab\\
A &\mapsto \H{A}{\Z} 
\end{align*}
et on peut en donner une description concrète grâce aux intégrateurs.
Si $A$ est une petite catégorie et $L_A$ est un intégrateur libre sur
$A$ (\ref{defIntegrateurs}), alors on sait qu'on
a un isomorphisme naturel, pour tout préfaisceau $X$ sur $A$,
\[
\H{A}{X} \simeq \bigoplus_{i\in I_0} Xa_i \leftarrow \bigoplus_{i\in
I_1} Xa_i \leftarrow \cdots 
\]
dans $\Hotab$, où on note, pour tout entier $n\geq 0$, $L_n = \bigoplus_{i\in I_n}
\Wh{{A}^{\op}}{a_i}$. Pour tout entier $n>0$, la différentielle $L_n \to
L_{n-1}$ est donnée par une matrice de la forme 
\[
  m = 
\begin{bmatrix}
u_{0,0} & u_{0,1} & \cdots & u_{0,s} \\
u_{1,0} & u_{1,1} & \cdots & u_{1,s} \\
\cdots & \cdots & \cdots & \cdots \\
u_{r,0} & u_{r,1} & \cdots & u_{r,s} \\
\end{bmatrix} 
\mdvirg
\begin{aligned}
 s&= |I_n| \\
 r&= |I_{n-1}|
\end{aligned}
\]
où chaque coefficient $u_{i,j}$ est dans le groupe libre
$\Z^{\left(\Hom_A(a_i, a_j)\right)}$ avec $i\in I_{n-1}$ et~$j\in I_{n}$. 

Appliquer additivement le
préfaisceau constant de valeur $\Z$ donne donc un morphisme 
\[
\Z^{(I_{n-1})} \xleftarrow{d} \Z^{(I_n)}
\]
où $d$ est donné par la matrice 
\[
\begin{bmatrix}
\epsilon_{0,0} & \epsilon_{0,1} & \cdots & \epsilon_{0,s} \\
\epsilon_{1,0} & \epsilon_{1,1} & \cdots & \epsilon_{1,s} \\
\cdots & \cdots  & \cdots & \cdots \\
\epsilon_{r,0} & \epsilon_{r,1} & \cdots & \epsilon_{r,s} \\
\end{bmatrix} 
\begin{aligned}
\in \M_{r,s}(\Z)
\end{aligned}
\]
et où pour tout $(i,j)\in I_{n-1}\times I_n$, l'entier $\epsilon_{i,j}$ est
l'image de $u_{i,j}$ par l'augmentation, c'est-à-dire l'unique morphisme
$\Z^{(\Hom_A(a_i,a_j))} \to \Z$ de valeur $1$ sur les générateurs.
\begin{example}
Dans le cas de l'intégrateur libre standard sur $\Delta$
(voir l'exemple ~\ref{ex:IntegrateurLibreDelta})
\[
L = \Wh{{\Delta}^{\op}}{\Delta_0} \xleftarrow{d_0 - d_1} 
\Wh{{\Delta}^{\op}}{\Delta_1} \xleftarrow{d_0 - d_1 + d_2}
\Wh{{\Delta}^{\op}}{\Delta_2} \leftarrow \cdots \mdvirg
\]
appliquer additivement le préfaisceau constant de valeur $\Z$ donne le
complexe
\[
\Z \xleftarrow{1-1=0} \Z \xleftarrow{1-1+1=1} \Z \xleftarrow{1-1+1-1=0} \cdots 
\]
qui a bien l'homologie du point, comme prévu puisque $\Delta$ admet un
objet final.
\end{example}

\begin{example}\label{exGlobesHomologie}
La catégorie des globes non réflexifs est la petite catégorie notée
$\G$ engendrée par le graphe 
\[
\xymatrix{
\D_0 \ar@/^/[r]^\sigma \ar@/_/[r]_\tau & 
\D_1 \ar@/^/[r]^\sigma \ar@/_/[r]_\tau & 
\D_2 \ar@/^/[r]^\sigma \ar@/_/[r]_\tau & 
\cdots
     \ar@/^/[r]^\sigma \ar@/_/[r]_\tau & 
\D_i \ar@/^/[r]^\sigma \ar@/_/[r]_\tau & 
\cdots
}
\]
soumise aux relations coglobulaires
\[
\sigma\circ \sigma = \tau \circ \sigma \mdvirg \sigma\circ \tau = \tau
\circ \tau \pbox{.}
\]

On reviendra plus en détail sur la catégorie $\G$ dans la section
\ref{secGlobes}, et sur son utilisation pour la définition des
$\omega$-catégories strictes. 

On peut se représenter géométriquement l'objet $\D_i$ comme le graphe
engendrant le disque $i$\nobreakdash-catégorique, dont voici les
quatre premiers : 
\shorthandoff{;}
\begin{align*}
\xymatrix@=4em{
\bullet &
\bullet \ar[r] & \bullet & 
\bullet \ar@/^1.5pc/[r]^{}|{\vphantom{X}}="a"
\ar@/_1.5pc/[r]^{}|{\vphantom{X}}="b"
\ar@2"a";"b" & \bullet &
\bullet \ar@/^1.5pc/[r]_{}="a" \ar@/_1.5pc/[r]^{}="b" & \bullet
\ar@/_1pc/@2"a";"b"_{}="c"
\ar@/^1pc/@2"a";"b"^{}="d"
\ar@3"c";"d"
} \pbox{.} 
\end{align*}
\shorthandon{;}%

On va alors montrer que le complexe de préfaisceaux sur $\G^{\op}$
\[
L_{\G} = \Wh{{\G}^{\op}}{\mathbb{D}_0} \xleftarrow{t-s}
\Wh{{\G}^{\op}}{\mathbb{D}_1} \xleftarrow{t-s} 
\Wh{{\G}^{\op}}{\mathbb{D}_2} \xleftarrow{t-s} 
\cdots
\]
est un intégrateur libre sur $\G$, où on note $t$ (resp. $s$) la
précomposition par $\sigma$ (resp. par $\tau$). En effet, on vérifie grâce aux
relations coglobulaires que pour tous objets
$\D_i$, $\D_n$ de $\G$, on a 
\[
\Hom_{\G}(\D_i, \D_n) = \begin{cases}
\lbrace \sigma^{n-i},\tau^{n-i} \rbrace & \text{si } i<n \pbox{;}\\
\lbrace \id \rbrace &\text{si } i=n \pbox{;}\\
\emptyset &\text{sinon} \pbox{.}
\end{cases} 
\]
On a donc un isomorphisme de complexes
\[
L_{\G}(\D_n) \simeq
\Z^2 
\xleftarrow{ \begin{bsmallmatrix*}[r]
-1 & -1 \\ 1 & 1
\end{bsmallmatrix*}} 
\Z^2 
\xleftarrow{\begin{bsmallmatrix*}[r]
-1 & -1 \\ 1 & 1
\end{bsmallmatrix*}} 
\cdots 
\leftarrow
\Z^2
\xleftarrow{ \begin{bsmallmatrix*}[r]
-1 \\ 1 
\end{bsmallmatrix*}}
\Z \leftarrow 0 \leftarrow \cdots
\]
qui montre que $L_{\G}(\D_n)$ a bien le type d'homologie du point, et
que le morphisme d'augmentation $L_{\G}(\D_n) \to \Z[0]$ est un
quasi-isomorphisme.

Le complexe $L_\G$ est donc un intégrateur libre sur $\G$, 
et l'homologie de $\G$ est obtenue en appliquant additivement le
préfaisceau constant $\Z_\G$ de valeur~$\Z$ sur $\G$ à $L_{\G}$. On
obtient alors le complexe 
\[
{L_\G}_!(\Z_\G) = 
\Z \xleftarrow{0}  
\Z \xleftarrow{0}  
\Z \xleftarrow{0}  
\cdots
\]
ce qui signifie donc qu'on a, pour tout entier $i\geq0$
\[
\mathsf{H}_i(\G;\Z) = \Z \pbox{.}
\]
\end{example}

\section{Hyperhomologie}

\bigskip
\emph{Dans toute cette section, on fixe une petite catégorie $A$, et
on note $B={A}^{\op}$ la catégorie opposée.}
\bigskip

On va construire un analogue enrichi de l'adjonction
\[
\xymatrix{
{\prefab{A}} \ar@<.5em>[r]^-{L_!}|-{}="a" 
& {\Ch(\Ab)}
\ar@<.5em>[l]^-{L^*}_-{}="b"
\ar@{}"a";"b"|{\perp}
}
\]
obtenue au paragraphe \ref{propAdjIntegrateurs} pour tout intégrateur
$L$ sur $A$, et on va montrer que celle-ci est
toujours homotopique. On va voir que celle-ci permet de calculer la
limite homotopique inductive des diagrammes de complexes de chaînes.

\paragr On rappelle que la catégorie $\Ch(\Ab)$ admet une structure de
catégorie monoïdale fermée, qui peut être décrite de la manière
suivante. Notons~$T_0$ le complexe $\Z[0]$ et, pour tout $k\geq1$, 
\[
T_k = 0 \leftarrow \cdots \leftarrow 0 \leftarrow \Z \xleftarrow{\id}
\Z \leftarrow 0 \leftarrow \cdots 
\]
où les deux copies de $\Z$ sont en degrés $k$ et $k-1$.
Pour $X$ et $Y$ deux complexes de chaînes en degré
positif, on définit alors le complexe 
\[
\Homi_{\Ch(\Ab)}(C,D)_n = \Hom_{\Ch(\Ab)}(C \otimes T_n, D)
\]
\notindex{$\Homi_{\Ch(\A)}$}%
où la différentielle est induite par le morphisme canonique $T_k \to
T_{k+1}$. On peut alors vérifier qu'on a des isomorphismes 
\begin{align*}
\Homi(C,D)_n &\simeq \prod_i\Hom(C_i, D_{i+n}) \mdvirg n>0 \\
\Homi(C,D)_0 &\simeq \Hom_{\Ch(\Ab)}(C,D) \pbox{.}
\end{align*}

On peut encore voir apparaitre cette construction à partir du bifoncteur 
\begin{align*}
\Chm(\Ab)^{\op} &\times \Ch(\Ab) \to
\Ch(\Ab) \\
(C^{\bullet}&, D_{\bullet}) \mapsto
\Tot^{\pi}\left(\Hom(C^{\bullet},D_{\bullet})\right)
\end{align*}
appliqué à deux complexes de chaînes $C$ et $D$ concentrés en degré
positif, en posant $C^{n} = C_{-n}$ et $d^n = d_{-n} : C^{n-1}\to C^n$,
qui donne alors \[
\Tot^\pi \left(\Hom_{\Ab}(C^{\bullet},D_{\bullet})\right) =
\Homi_{\Ch(\Ab)}(C,D) \pbox{.}
\]

\paragr \label{def:LSouligne}
On rappelle que le couplage défini au paragraphe \ref{couplageOdot} 
\[
  - \odot - : \prefab{A}\times\prefab{B} \to \Ab 
\]
admet un foncteur dérivé à gauche
\[
- \Lodot - : D_+(\prefab{A}) \times D_+(\prefab{B}) \to D_+(\Ab) 
\]
que nous n'avons jusqu'ici considéré qu'en l'appliquant à des éléments de
$\prefab{A}$, c'est-à-dire à des complexes concentrés en degré $0$. 

Si  $X_{\bullet}$
est un complexe dans $\prefab{A}$ et $Y_{\bullet}$ est un complexe
dans $\prefab{B}$, on peut former le complexe double 
\[
X_{\bullet} \odot Y_{\bullet} \in \Ch(\Ch(\Ab)) 
\]
dont peut considérer le complexe total
$\Tot^{\oplus}(X_{\bullet}\odot Y_{\bullet})$. 
Si $L : A \to \Ch(\Ab)$ est un foncteur, on note alors 
\begin{align*}
\underline{L}_! : \Ch(\prefab{A})&\to \Ch(\Ab) \\
X_{\bullet} &\mapsto \Tot^{\oplus}(X_{\bullet} \odot L_{\bullet})
\end{align*}
ce foncteur. On a déjà montré dans la preuve de la proposition
\ref{produittensorielderive} que si $L$ est un intégrateur sur $A$,
alors le foncteur $\underline{L}_!$ préserve les quasi-isomorphismes. 

\begin{prop}
Si $L : A \to \Ch(\Ab)$ est un foncteur, alors le foncteur
\[
  \underline{L}_! : \Homi({A}^{\op},\Ch(\Ab)) \to \Ch(\Ab) 
    \]
\notindex{$\underline{L}_\bang$}%
\notindex{$\underline{L}^*$}%
admet pour adjoint à droite le foncteur
\[
\underline{L}^* : \Ch(\Ab) \to \Homi({A}^{\op},\Ch(\Ab)) 
\]
défini pour tout complexe de chaînes $C$ par 
\[
\underline{L}^*(C) : a \mapsto \Homi_{\Ch(\Ab)}(La, C) \pbox{.}
\]
De plus, 
en notant $t : \Ch(\Ab) \to \Ab$ le foncteur de troncation \og bête
\fg{} des~$0$\nobreakdash-chaînes, adjoint à droite de l'inclusion $\Ab
\hookrightarrow \Ch(\Ab)$ en degré $0$, les diagrammes suivants sont
commutatifs :
\[
\xymatrix{
\prefab{A} \ar[r]^{L_!} \ar[d]_{i} & \Ch(\Ab)\\
\Homi({A}^{\op},\Ch(\Ab)) \ar[ru]_{\underline{L}_!} & \mdvirg
} 
\xymatrix{
\prefab{A} & \ar[l]_{L^*} \Ch(\Ab) \ar[ld]^{\underline{L}^*} \\
\Homi({A}^{\op},\Ch(\Ab)) \ar[u]^{t} &\pbox{.}
} 
\]
\end{prop}
\begin{proof}

Pour tout foncteur $L : A \to \Ch(\Ab)$ vu comme un complexe de
préfaisceaux $L_{\bullet}$ sur $B = {A}^{\op}$,
le foncteur \[
\Tot^{\oplus}(- \odot L_{\bullet}) : \Ch(\prefab{A}) \to \Ch(\Ab)
\]
admet comme adjoint à droite le foncteur 
\[
\Tot^{\pi}(L_{\bullet}^*(-_{\bullet})) : \Ch(\Ab) \to
\Ch(\prefab{A})
\]
où on traite $L$ comme un objet de $\Chu(\prefab{B})$ en posant
$L^n=L_{-n}$ et $d^n =d_{-n} : L^{n-1} \to L^n$. On obtient donc que
pour tout complexe de chaînes $C$, on a \[
\Tot^\pi(L_{\bullet}^*(C_{\bullet})) = \Tot^\pi (\Hom(L^{\bullet}(-),
C_{\bullet})) = \Homi_{\Ch(\Ab)}(L(-), C) \pbox{.}
\]
La commutativité des triangles résulte alors simplement d'une
composition de foncteurs adjoints.
\end{proof}

\paragr  
On rappelle (voir la section \ref{colimiteHomotopiqueDiagrammeCh})
qu'on note~$\DerHotab(A)$ la catégorie localisée \[
\W^{-1}\Homi({A}^{\op},\Ch(\Ab)) 
\]
où on note $\W$ la classe des transformations naturelles qui sont des
quasi-isomorphismes argument par argument. Sous l'isomorphisme entre
la catégorie des
foncteurs~${A}^{\op}\to\Ch(\Ab)$ et la catégorie des complexes de préfaisceaux
abéliens sur~$A$, la catégorie $\DerHotab(A)$ coïncide avec la catégorie dérivée de la
catégorie abélienne des préfaisceaux en groupes abéliens sur $A$, c'est-à-dire qu'on a un
isomorphisme canonique 
\[
\DerHotab(A) \simeq D_+(\prefab{A}) \pbox{.}
\]

On va maintenant voir que si $L$ est un intégrateur
sur $A$, le foncteur $\underline{L}_!$ défini ci-dessus calcule la
limite inductive homotopique
\[
\DerHotab(A) \xrightarrow{\hocolim_{{A}^{\op}}^{\Hotab}} \Hotab
\pbox{.}
\]

On remarque d'abord que si $L$ est un intégrateur sur $A$, on dispose
d'un
quasi\nobreakdash-isomorphisme $\epsilon : L \to \Z_B$, en notant $\Z_B$
le préfaisceau constant de valeur $\Z$ sur $B={A}^{\op}$ vu comme un complexe
concentré en degré $0$. Mais~$\epsilon$ est
nécessairement une équivalence d'homotopie, puisque
$L$ et $\Z_B$ sont des complexes d'objets projectifs de $\prefab{B}$.
\begin{prop}\label{prop:LsoulignePreserveEqf}
Pour tout intégrateur~$L$ sur $A$,
le diagramme 
\[
\xymatrix{
{ \Ch(\Ab) } \ar[r]^-{ \underline{L}^* } \ar[d]_{  } & { \Homi({A}^{\op},\Ch(\Ab)) } \ar[d]^{  } \\
{ \Hotab } \ar[r]_{ \const } & { \DerHotab(A) }\pbox{,} }
\]
où les flèches verticales désignent les foncteurs de localisation, est
commutatif à isomorphisme naturel près. En particulier,
le foncteur
\[
  \underline{L}^* : \Ch(\Ab) \to \Homi({A}^{\op}, \Ch(\Ab))
  \]
préserve les équivalences faibles.
\end{prop}
\begin{proof}
Pour tout objet $a$ de $A$, le morphisme $\epsilon_a : La \to
\Z[0]$ est une équivalence d'homotopie. Cela implique que pour tout
complexe $C$, le morphisme 
\[
\underline{\epsilon}^*(C) : \Homi_{\Ch(\Ab)}(\Z_B,C) \to
\underline{L}^*(C)
\]
est une équivalence d'homotopie argument par argument. En d'autres
termes, on dispose d'une transformation naturelle du foncteur
diagramme constant 
\[
\const : \Ch(\Ab) \to \Homi({A}^{\op},\Ch(\Ab))
\]
vers le foncteur $\underline{L}^*$ qui est une équivalence faible
argument par argument (c'est-à-dire que pour tout complexe $C$, le
morphisme $\const_C \to \underline{L}^*(C)$ est un quasi-isomorphisme
argument par argument).
Cela implique en particulier que le foncteur $\underline{L}^*$ envoie
les quasi-isomorphismes sur des quasi-isomorphismes arguments par
arguments.
\end{proof}

\begin{prop}\label{colimiteHomotopiqueDegreparDegre}
Pour tout intégrateur $L$ sur $A$, le diagramme  
\[
\xymatrix{
\Homi({A}^{\op},\Ch(\Ab)) \ar[d]^{}  
\ar[r]^-{{\underline{L}}_!} 
& \Ch(\Ab) \ar[d]^{} \\
\DerHotab(A) \ar[r]_{\hocolim_{{A}^{\op}}^{\Hotab}}  
&\Hotab
} 
\]
est commutatif à isomorphisme naturel près, où les flèches verticales
correspondent aux foncteurs de localisation.
\end{prop}
\begin{proof}
On sait que la catégorie des foncteurs ${A}^{\op}\to \Ch(\Ab)$
peut être munie d'une structure de catégorie de modèles dite
\emph{projective} (voir \ref{colimiteHomotopiqueDiagrammeCh}), pour
laquelle le foncteur diagramme constant est un foncteur de Quillen à
droite. Soit alors $L$ un intégrateur sur $A$. On va montrer que
$\underline{L}^*$ est également un foncteur de Quillen à droite pour
cette structure. 

Pour cela, il suffit de montrer que pour tout objet $a$ de $A$, le
foncteur 
  \begin{align*}
\underline{L}^*(-)(a) : \Ch(\Ab) &\to \Ch(\Ab) 
\\
C &\mapsto \Homi_{\Ch(\Ab)}(La, C)
  \end{align*}
est un foncteur de Quillen à droite. Or, on vérifie rapidement que pour tout objet $a$
de $A$, $La$ est un complexe de groupes abéliens projectifs, et donc
que~$La$ est
un objet cofibrant de $\Ch(\Ab)$. Puisque la structure de catégorie de
modèles sur $\Ch(\Ab)$ est monoïdale (voir par exemple
\cite[lemme 1.5]{jardine2003presheaves}), le foncteur $La \otimes - : \Ch(\Ab)
\to \Ch(\Ab)$ est un foncteur de Quillen à gauche, et son adjoint à
droite $\underline{L}^*(-)(a)$ est donc bien un foncteur de Quillen à
droite.

On peut alors conclure de la manière suivante : les foncteurs
$\underline{L}_!$ et $\underline{L}^*$ préservent les
équivalences faibles. Par adjonction des
foncteurs dérivés, l'adjonction $\underline{L}_!
\dashv \underline{L}^*$ induit donc une adjonction
après localisation. Or, il résulte de la preuve de la proposition
\ref{prop:LsoulignePreserveEqf} que le foncteur $\underline{L}^*$
coïncide, après localisation, avec le foncteur dérivé à droite du
foncteur diagramme constant. Le foncteur~$\underline{L}_!$ coïncide
donc, après localisation, au foncteur dérivé à gauche de l'adjoint à
gauche du foncteur diagramme constant, qui est le foncteur de limite
inductive homotopique.
\end{proof}

\paragr En d'autres termes, pour tout complexe $X$ de préfaisceaux abéliens
sur $A$ on a un isomorphisme naturel
\[
 \hocolim\nolimits_{{A}^{\op}}^{\Hotab} X \simeq \underline{L}_!(X_\bullet) 
\]
dans $\Hotab$, ce qui peut se traduire en disant que
\emph{la colimite homotopique des foncteurs~${A}^{\op}\to \Ch(\Ab)$ se
calcule degré par degré}.

On notera alors, si $X_\bullet$ est un complexe de préfaisceaux en
groupes abéliens sur $A$, 
\[
\H{A}{X_\bullet} := \hocolim_{{A}^{\op}}^{\Hotab} X
\]
\notindex{$\H{A}{X_\bullet}$}%
l'\ndef[homologie d'une catégorie à coefficients!dans un complexe de préfaisceaux]{hyperhomologie} de $A$ à coefficients dans le complexe de
préfaisceaux $X_\bullet$, qui coïncide donc avec le foncteur \og
hyperdérivé\fg{} du
foncteur
\[
- \odot \Z_B : \prefab{A} \to \Ab \pbox{.}
\]

\begin{prop}
Soient $X_{\bullet}$ et $Y_{\bullet}$ deux complexes de préfaisceaux
abéliens sur $A$, et $\alpha : X_\bullet \to Y_\bullet$ un morphisme.
Si, pour tout entier $n\geq 0$, $\alpha_n$ est dans $\Wab_A$, alors le
morphisme 
\[
\H{A}{\alpha} : \H{A}{X_\bullet} \to \H{A}{Y_\bullet} 
\]
est un isomorphisme dans $\Hotab$.
\end{prop}
\begin{proof}
\emph{A priori}, pour que deux complexes de préfaisceaux $X_\bullet$
et~$Y_\bullet$ sur $A$ aient la même limite inductive homotopique, il
faut qu'on ait un quasi-isomorphisme $X_\bullet \to Y_\bullet$,
c'est-à-dire un morphisme $X_\bullet \to Y_\bullet$ dans~$\Ch(\prefab{A})$ qui soit un
quasi-isomorphisme argument par argument. Considérons alors un
intégrateur libre $L$ sur $A$.
Le morphisme $\alpha$ induit alors par naturalité un morphisme de
complexes doubles 
\[
X_{\bullet} \odot L_{\bullet} \xrightarrow{\alpha} Y_{\bullet} \odot
L_{\bullet}
\]
et par hypothèse, le morphisme $X_n \odot L_{\bullet} \to Y_n
\odot L_{\bullet}$
est un quasi-isomorphisme pour tout $n\geq0$. La proposition
\ref{prop:lemmeComplexeDouble} permet d'affirmer que le morphisme
induit entre les complexes totaux est un quasi-isomorphisme.
\end{proof}

\paragr Dans ce contexte, on peut aussi comparer facilement les
adjoints à droites induits par différents intégrateurs. En effet, si
$L$ et $L'$ sont deux intégrateurs sur $A$, alors il existe une
équivalence d'homotopie entre $L$ et $L'$ (c'est un corollaire du
théorème de comparaison, puisqu'un intégrateur est une résolution
projective du préfaisceau $\Z_B$. Voir par exemple~\cite[2.2.6]{weibel1994Introduction}). Puisque le foncteur
$\Homi(-,C)$ envoie les homotopies sur des homotopies, on peut énoncer
le résultat suivant.
\begin{prop}
Soient $L$ et $L'$ deux intégrateurs sur une petite catégorie $A$. Pour tout complexe $C$,
il existe une équivalence d'homotopie naturelle
\[
\underline{L}^*(C) \to \underline{L'}^*(C)
\]
dans $\Ch(\prefab{A})$. En particulier, il existe un
isomorphisme dans $\Hotab$
\[
\underline{L}_!\underline{L}^*(C) \simeq
\underline{L}_!\underline{L'}^*(C) \pbox{.}
\]
\end{prop}

\section{Dualité}

La manière dont Grothendieck parvient, dans \cite{pursuingstacks}, à
la définition des intégrateurs est en fait duale de celle que nous avons
exposé dans la section précédente. Plus habitué à la cohomologie des
faisceaux qu'à l'homologie, son idée est de chercher un moyen d'étudier
les propriétés cohomologiques de la catégorie opposée $B={A}^{\op}$,
c'est-à-dire de trouver un complexe de chaines $L_{\bullet}$ fixé dans
$\Addinf(B)$ tel que pour tout préfaisceau abélien $Y$ sur $B$, on ait
un isomorphisme 
\[
R\Gamma_B(Y) \simeq \Hom_{\prefab{B}}^{\bullet}(L_{\bullet},Y)
\]
dans $D^+(\Ab)$, où le foncteur $\Gamma : \prefab{B}\to \Ab$ est le
foncteur section globale. Vient alors le slogan \og les
cointégrateurs sur $B$ coïncident avec les intégrateurs sur
$A={B}^{\op}$\fg, que nous allons expliquer dans cette section.

\paragr
Si $B$ est une petite catégorie (oublions pour l'instant la
catégorie~$A={B}^{\op}$), on dispose d'un foncteur \[
\Gamma_B : \prefab{B} \to \Ab \mdvirg Y \mapsto \Hom_{\prefab{B}}(\Z_B,Y)
\]
appelé le foncteur \emph{section globale}. Ce foncteur est exact à
gauche (puisque le foncteur $\Hom(\Z_B,-)$ commute aux limites projectives), et
admet donc un foncteur dérivé à droite 
\[
R\Gamma_B :  D^+(\prefab{B})\to D^+(\Ab) \pbox{.}
\]

Si $Y$ est un préfaisceau abélien sur $B$, la \emph{cohomologie de
$B$ à coefficients dans $Y$} est définie comme la cohomologie du
complexe $R\Gamma_B(Y)$ : on note alors 
\[
\coH{B}{Y} := R\Gamma_B(Y) 
\]
le complexe, vu dans la catégorie $D^+(\Ab)$, calculant la cohomologie
de $B$ à coefficients dans $Y$.

\begin{prop}
Pour toute petite catégorie $B$ et tout préfaisceau abélien $Y$ sur $B$, on a un isomorphisme
canonique naturel
\[
  \Gamma_B(Y) \simeq \limproj_B Y
\]
dans $\Ab$.\end{prop}
\begin{proof}
Si $Y$ est un préfaisceau abélien sur $B$, en utilisant le fait que
$\Z_B=\Wh{B}{e_{\pref{B}}}$, et que les foncteurs $\Whf{B}$ et
$\Hom(-,Y)$ commutent aux limites inductives,
on obtient la chaîne d'isomorphismes suivante :
\begin{align*}
\Hom(\Z_B, Y) &\simeq \Hom\big(\limind_B \Wh{B}{b}, \,Y\big) \\
&\simeq \limproj_B \Hom(\Wh{B}{b},\, Y)\\
&\simeq \limproj_B Y \pbox{.} \qedhere
\end{align*}
\end{proof}

\begin{coro}
Pour toute petite catégorie $B$ et tout préfaisceau abélien~$Y$ sur $B$, on a un isomorphisme
naturel dans $D_+(\Ab)$
\[
\coH{B}{Y} \simeq R\limproj_B Y \pbox{.}
\]
\end{coro}

\begin{example}
En reprenant l'exemple du groupoïde associé à un groupe~$G$, si $M:BG
\to \Ab$ est un $G$-module, alors la limite projective de $M$
est isomorphe au groupe des invariants sous l'action de $G$ : 
\[
\limproj\nolimits_{BG} M \simeq M^G \pbox{,}
\]
et la cohomologie du groupe $G$ à coefficients dans $M$ correspond
bien à la limite projective homotopique de $M$ vu comme un foncteur
$BG \to \Ab$.
\end{example}

\begin{prop}
Soit $B$ une petite catégorie et $L$ un objet de $\Add(B)$. Alors pour tout préfaisceau
abélien~$Y$ sur ${B}$, on a un isomorphisme canonique naturel 
\[
\Hom_{\prefab{B}}(L,Y) \simeq \Add(Y)(L). 
\]
En notant $L=\bigoplus_{i} \Wh{B}{b_i}$, on a donc 
\[
\Hom_{\prefab{B}}(L,Y) \simeq \bigoplus_{i} Y(b_i) \pbox{.}
\]
\end{prop}
\begin{proof}
On peut le vérifier directement en utilisant le fait que le foncteur
$\Hom$ commute aux produits en la première variable, et en utilisant
le lemme de Yoneda.
\end{proof}

\begin{prop}
Soit $B$ une petite catégorie et $L$ un intégrateur sur la catégorie $A={B}^{\op}$. Alors pour
tout préfaisceau abélien $Y$ sur $B$, on a un isomorphisme naturel
\[
R\Gamma_B(Y) \simeq \Hom_{\prefab{B}}^{\bullet}(L_{\bullet}, Y) \pbox{.}
\]
Ainsi, si $L$ est un intégrateur libre de type fini, alors en notant,
pour tout entier $n\geq 0$,
$L_n=\bigoplus_{i\in I_n}\Wh{B}{b_i}$, on a un isomorphisme naturel
\[
R\Gamma_B(Y) \simeq \cdots \to \bigoplus_{i\in I_n} Y(b_i) \rightarrow
\bigoplus_{i\in I_{n+1}} Y(b_i) \to \cdots
\]
dans $D^+(\Ab)$.
\end{prop}
En d'autres termes, la donnée d'un intégrateur sur la catégorie
$A={B}^{\op}$, exprimant les propriétés homologiques de $\prefab{A}$,
coïncide avec la donnée de ce que Grothendieck appelle
\emph{cointégrateur}
\footnote{On renvoie à la note de bas de page de la définition
\ref{defIntegrateurs} pour la confusion possible entre intégrateurs et
cointégrateurs dans \emph{Pursuing Stacks} \cite{pursuingstacks}.}
sur
$B$, exprimant les propriétés cohomologiques de $\prefab{B}$.

\begin{example}
En considérant l'intégrateur libre de type fini
\[
L_{\Delta} = \Wh{{\Delta}^{\op}}{\Delta_0} \leftarrow
\Wh{{\Delta}^{\op}}{\Delta_1} \leftarrow \Wh{{\Delta}^{\op}}{\Delta_2}
\cdots
\]
sur $\Delta$, c'est-à-dire le foncteur $\Delta \to \Ch(\Ab)$ associant
à tout objet $\Delta_n$ le complexe de chaînes 
\[
\Z^{(\Delta_0 \to \Delta_n)} \leftarrow
\Z^{(\Delta_1 \to \Delta_n)} \leftarrow
\Z^{(\Delta_2 \to \Delta_n)} \leftarrow \cdots \mdvirg
\]
le complexe opposé associé est le complexe de
préfaisceaux sur $\Delta$ qui associe à un objet $\Delta_n$ de
$\Delta$ le complexe de cochaînes 
\[
\Z^{(\Delta_n \to \Delta_0)} \rightarrow
\Z^{(\Delta_n \to \Delta_1)} \rightarrow
\Z^{(\Delta_n \to \Delta_2)} \rightarrow \cdots
\]
où la différentielle est donnée par la somme alternée des faces. Alors
si $X$ est un préfaisceau sur ${\Delta}$ et $Y$ est un préfaisceau sur
${\Delta^{\op}}$, on a 
\begin{align*}
\H{\Delta}{X} &\simeq X_0 \leftarrow X_1 \leftarrow X_2 \leftarrow
\cdots \\
\coH{\Delta^{\op}}{Y} &\simeq Y_0 \rightarrow Y_1 \rightarrow Y_2 \rightarrow \cdots 
\end{align*}
où les différentielles sont toutes les deux données par la somme
alternée des faces.
\end{example}
\begin{prop}\label{dualiteAddFinie}
Soit $B$ une petite catégorie et $A={B}^{\op}$ la catégorie opposée.
On a alors une équivalence de catégories 
\[
\gamma_A : \Add(A)^{\op} \xrightarrow{\simeq} \Add(B) \mdvirg
\]
qui induit donc une équivalence de catégories
\[
\Ch(\Add(B)) \simeq \Chm(\Add(A)) \pbox{.}
\]

\end{prop}
\begin{proof}
Par la propriété universelle de l'enveloppe additive, il suffit,
pour définir $\gamma_A$, de définir un foncteur $A \to \Add(B)^{\op}$,
puisque cette dernière est une catégorie additive. On pose alors, pour
tout objet $a$ de $A$,
\[
\gamma_A(a) = \Wh{B}{{a}^{\op}} = \Z^{(\Hom_A(a,-))} \pbox{.}
\]
On obtient alors un foncteur 
\[
\gamma_A : \Add(A)^{\op} \to \Add(B) \mdvirg \Wh{A}{a} \mapsto
\Wh{B}{{a}^{\op}}
\]
dont l'inverse est donné par le foncteur $\gamma_B$ défini
identiquement.
\end{proof}

\begin{remark}
On ne peut pas reproduire cet argument dans le cas de l'enveloppe
additive infinie ! En effet, on a utilisé la propriété universelle de
l'enveloppe additive qui permet d'étendre de manière unique les
foncteurs~$A\to M$ pour toute catégorie additive $M$. Pour l'enveloppe
additive infinie, on a besoin que $M$ soit une catégorie additive
admettant des sommes indexées par des ensembles quelconques. Or la
catégorie $\Addinf(B)^{\op}$ ne possède pas de telles sommes (car
l'enveloppe additive infinie de possède pas de produits infinis). 
\end{remark}

Ainsi, si $L_{\bullet}$ est un intégrateur libre de type fini sur une
petite catégorie $A$,
on peut lui associer le complexe de cochaînes de préfaisceaux sur $A$
\[
J^{\bullet} = \gamma(L)^{\bullet} \in \Chm(\Add(A)) 
\]
et on remarque alors qu'on a, pour tout préfaisceau abélien $Y$ sur
$B={A}^{\op}$,
\[
\Hom(L_{\bullet},Y) \simeq \gamma(L)^{\bullet}\odot Y \pbox{.}
\]

\begin{prop}
Si $Y$ est un préfaisceau abélien sur $B$ et $L$ est un intégrateur
libre de type fini sur $A={B}^{\op}$, on a un isomorphisme naturel  
\[
R\Gamma_B(Y) \simeq \gamma(L)^{\bullet} \odot_A Y
\]
dans $D_+(\Ab)$.
\end{prop}
\begin{proof}
Si $L$ est un intégrateur libre \emph{de type fini}, il suffit
d'utiliser les isomorphismes naturels explicités ci-dessus
\[
R\Gamma_B(Y)\simeq\Hom(L_\bullet,Y) \simeq \gamma(L)^{\bullet}\odot Y
\pbox{.}
\qedhere
\]
\end{proof}
Il peut alors sembler tentant de calculer la cohomologie de $B$ à
coefficients dans le complexe $L \in \Ch(\prefab{B})$, ou ce qui
revient au même, l'homologie de~$A={B}^{\op}$ à coefficients dans
$J=\gamma(L)$. On obtient alors le complexe double de groupes abéliens
\[
  \H{A}{J^{\bullet}} = J^{\bullet} \odot L_{\bullet} =
R\Gamma_B(L_{\bullet}) = \coH{B}{L_{\bullet}}
\]
dont chaque ligne est alors quasi-isomorphe à $\H{A}{\Z}$ et chaque colonne
est quasi-isomorphe à $\coH{B}{\Z}$. On a explicitement le complexe
double suivant : 
\begin{align*}
\xymatrix{
J^0 \odot L_0 \ar[r]^{} & J^1 \odot L_0 \ar[r] & J^2\odot L^0
\ar[r]^{} & \cdots \\
J^0 \odot L_1 \ar[u]^{} \ar[r]^{} & J^1 \odot L_1 \ar[u]^{} \ar[r] &
J^2\odot L_1 \ar[u]^{} \ar[r] & \cdots \\
J^0 \odot L_2 \ar[u]^{} \ar[r]^{} & J^1 \odot L_2 \ar[u]^{} \ar[r] &
J^2\odot L_2 \ar[u]^{} \ar[r] & \cdots \\
\vdots \ar[u]^{} &
\vdots \ar[u]^{} &
\vdots \ar[u]^{} & \pbox{.}
}
\end{align*}
Par exemple, en appliquant cette construction
à l'intégrateur $L_{\Delta}$,
on obtient le complexe double 
\[
\xymatrix{
\Z^{(\Delta_0 \to \Delta_0)} \ar[r]^{}  & 
\Z^{(\Delta_0 \to \Delta_1)} \ar[r]^{}  & 
\Z^{(\Delta_0 \to \Delta_2)} \ar[r]^{}  & \cdots  \\
\Z^{(\Delta_1 \to \Delta_0)} \ar[r]^{} \ar[u]^{} & 
\Z^{(\Delta_1 \to \Delta_1)} \ar[r]^{} \ar[u]^{} & 
\Z^{(\Delta_1 \to \Delta_2)} \ar[r]^{} \ar[u]^{} & \cdots \\
\Z^{(\Delta_2 \to \Delta_0)} \ar[r]^{} \ar[u]^{} & 
\Z^{(\Delta_2 \to \Delta_1)} \ar[r]^{} \ar[u]^{} & 
\Z^{(\Delta_2 \to \Delta_2)} \ar[r]^{} \ar[u]^{} & \cdots \\
\vdots \ar[u] & 
\vdots \ar[u] & 
\vdots \ar[u] & 
} 
\]
dont les différentielles sont toutes données par les sommes alternées
des faces.

\paragr 
On rappelle (\ref{propAdjIntegrateurs}) que pour tout intégrateur $L$
sur une petite catégorie $A$, le foncteur $ - \odot L_{\bullet} : \prefab{A} \to
\Ch(\Ab)$ admet comme adjoint à droite le foncteur $L^* :
\Ch(\Ab) \to \prefab{A}$, défini pour tout complexe $C$ par 
\begin{align*}
L^*(C) : {A}^{\op} &\to \Ab \\
a &\mapsto \Hom_{\Ch(\Ab)}(La,C ) \pbox{.}
\end{align*}

On dira plus loin (\ref{defTestHomologiqueFaible}) qu'un intégrateur $L$
sur $A$ est test faible si cet adjoint à droite envoie les
quasi-isomorphismes sur des éléments de $\Wab_A$, et si le foncteur
induit est une équivalence de catégories quasi-inverse au foncteur~
\[
  L_! = - \odot L : \Hotab_A \to \Hotab \mdvirg
\]
c'est-à-dire si la counité de l'adjonction $L_! \dashv L^*$ induit un
isomorphisme naturel, pour tout complexe $C$, 
\[
L^*(C)\odot L \simeq C \pbox{.} 
\]

On peut alors remarquer que le terme de gauche ci-dessus correspond
au complexe obtenu en appliquant terme à terme le foncteur
$\Hom_{\Ch(\Ab)}(-,C)$ aux colonnes de ce complexe double.
Intuitivement, cela correspond à
considérer le \og $C$-dual de l'homologie de $A$ à valeurs dans le
complexe $J^\bullet$ \fg. 

Une piste qui semble prometteuse pour caractériser les intégrateurs
libres de type fini test faibles, mais qui en est toujours une au
moment de l'écriture de cette thèse, est d'utiliser le résultat
suivant.

\begin{prop}
Si $L$ est un intégrateur de type fini sur une petite catégorie $A$,
alors, pour tout complexe $C\in \Ch(\Ab)$, on a un isomorphisme  
\[
L^*(C) \odot L \simeq \Hom_{\Ch(\Ab)}(J^{\bullet}\odot L, C) 
\]
dans $\Hotab$, où on note $J^\bullet=\gamma(L_\bullet)$ l'image de~$L$
par l'équivalence de catégories~$\gamma : \Ch(\Add({A}^{\op}) \simeq
\Chm(\Add(A))$ de la proposition \ref{dualiteAddFinie}.
\end{prop}
\begin{proof}
Pour tout entier $n$, en notant $L_n = \bigoplus_{i\in
I_n}\Wh{B}{a_i}$ on a alors $J^n = \bigoplus_{i\in I_n}\Wh{A}{a_i}$ et
peut conclure de la manière suivante : 
\begin{align*}
(L^*(C)\odot L)_n &\simeq \Add(L^*(C))(L_n) \\
&\simeq \bigoplus_{i\in I_n}(\Hom(La_i, C) \\
&\simeq \Hom(\bigoplus_{i\in I_n}La_i, C) \\
& \simeq \Hom(\bigoplus_{i\in I_n}\Wh{A}{a_i}\odot L, C)\\
& \simeq \Hom(J^n \odot L, C) \qedhere
\end{align*}
\end{proof}

\chapter{Asphéricité en homologie}\label{chapAsphericite}
\section{Localisateurs fondamentaux}
\label{secLocalisateursFondamentaux}

On introduit dans cette section la notion de localisateur fondamental,
dont l'un des principaux intérêts est de fournir une caractérisation
interne à $\Cat$ de la classe $\W_\infty$ des équivalences de
Thomason. L'idée de Grothendieck en introduisant cette notion est de
construire une axiomatique basée sur le théorème A de Quillen
\cite[théorème A]{quillenktheory}, et de montrer que $\W_\infty$ est
la plus petite classe de~$\Fl(\Cat)$ satisfaisant ces axiomes. Ce
résultat sera en fait démontré par Cisinski dans
\cite{cisinski2004localisateur}. On va, dans la section
\ref{secMorphismesAsphEnHomologie}, étudier l'analogue abélien de
$\W_\infty$, et on montrera que celle-ci forme également un
localisateur fondamental. Cela nous permettra entre autres de
caractériser les foncteurs commutant aux limites inductives
homotopiques, en utilisant quelques énoncés de l'annexe
\ref{annexeDerivateurs}.

Nous renvoyons au livre de Maltsiniotis \cite{maltsiniotis2005}, sur
lequel nous nous sommes largement appuyés, pour un exposé détaillé sur
les localisateurs
fondamentaux.

\paragr \label{defFaibleSaturation} Soit $\W$ une classe de flèches de
$\Cat$. On dit que $\W$ est \ndef[partie faiblement
saturée]{faiblement saturée} si elle vérifie les conditions suivantes
: 
\begin{enumerate}
\item $\W$ contient les isomorphismes;
\item si, dans un triangle commutatif de $\Cat$, deux morphismes sur trois sont
dans $\W$, alors le troisième morphisme est dans $\W$;
\item si $i:A\to B$ et $r:B\to A$ sont deux morphismes tels que
$ri=\id_A$ et~$ir \in \W$, alors $r\in \W$ (ce qui implique aussi que
$i\in\W$).
\end{enumerate}

On peut remarquer que si $F : \Cat \to \M$ est un foncteur, alors la
classe des foncteurs dont l'image par $F$ est un isomorphisme de $\M$
vérifie toujours ces propriétés. Il s'avère que cette axiomatique, assez
minimale, suffit pour caractériser~$\W_\infty$ comme la plus petite
classe \og raisonnable \fg{} de $\Fl(\Cat)$ vérifiant le théorème A
de Quillen, comme on va l'énoncer au théorème
\ref{locFondamentalMinimalThm}.

\paragr \label{defAspheriqueLocalisateurFondamental}
Soit $\W$ une partie faiblement saturée de $\Cat$.
On dit qu'une petite catégorie $A$ est
\ndef[catégorie!asphérique!relativement à un localisateur
fondamental]{$\W$-asphérique} si le morphisme~$p_A:A\to e$ est dans
$\W$, où on note encore $e$ la catégorie ponctuelle. On dit qu'un
foncteur $u:A\to B$ est~\ndef[morphisme asphérique!de
$\Cat$!relativement à un localisateur fondamental]{$\W$-asphérique} si
pour tout objet $b$ de~$B$, le morphisme \[
  \tranche{u}{b} : \tranche{A}{b} \to \tranche{B}{b} \mdvirg
(a, ua \xrightarrow{\varphi} b) \mapsto (ua, \varphi)
\]
est dans $\W$.

\paragr\label{morphismeAspheriqueAuDessus}Étant donné un morphisme
$u:A\to B$ de $\Cat$ au-dessus d'une petite catégorie $C$,
c'est-à-dire un triangle commutatif de la forme 
\[
\xymatrix{
{A} \ar[rr]^{u} \ar[rd]_{v} && {B} \ar[ld]^{w} \\
& {C}
} 
\]
dans $\Cat$, on peut définir, pour tout objet $c$ de $C$, un morphisme
\begin{align*}
\tranche{A}{c} &\xrightarrow{\tranche{u}{c}} \tranche{B}{c} \\ 
(a,\varphi) 
&\mapsto (ua, \varphi) \pbox{.}
\end{align*}
Si $\W$ est une partie faiblement saturée de $\Cat$, on dit que $u$
est \ndef[morphisme asphérique!de $\Cat$!au-dessus d'une petite
catégorie]{$\W$-asphérique au-dessus de $C$} si pour tout objet $c$ de
$C$, le morphisme $\tranche{u}{c}$ est dans $\W$. 

\begin{defin}\label{defLocFondamental}
Soit $\W$ une classe de flèches de $\Cat$. On dit que $\W$ est un
\ndef{localisateur fondamental} si elle vérifie les propriétés
suivantes : 
\begin{enumerate}
\item $\W$ est faiblement saturée;
\item toute petite catégorie admettant un objet final est
$\W$-asphérique;
\item tout morphisme $\W$-asphérique au-dessus d'une petite catégorie
est dans~$\W$.
\end{enumerate}
\end{defin}

\begin{example}
Le théorème A de Quillen \cite[théorème A]{quillenktheory} implique
que la classe~$\W_{\infty}$ des équivalences de Thomason (\ref{def:eqThomason}) est un
localisateur fondamental. D'autres exemples de localisateurs
fondamentaux incluent les classes~$\W_n$ des morphismes induisant des
isomorphismes sur tous les groupes d'homotopie jusqu'au degré $n$ pour
$n\in \N$ et une bijection sur le $\pi_0$ pour tout point de base. 
\end{example}

On renvoie à \cite{maltsiniotis2005} pour un exposé approfondi sur la
théorie de l'homotopie de~$\Cat$ utilisant le formalisme des
localisateurs fondamentaux. Les propositions suivantes sont des
conséquences directes des axiomes \ref{defLocFondamental}, et on
renvoie à la première section de \emph{op. cit.} pour les preuves.

\medskip
\emph{
On fixe, jusqu'à la fin de cette section, un localisateur fondamental $\W$.
}
\medskip

\begin{prop}[Grothendieck]\label{morphismesWAspheriquesEquivalences}
Soit $u:A\to B$ un foncteur. Les propositions suivantes
sont équivalentes : \begin{enumerate}
\item $u$ est un morphisme $\W$-asphérique de $\Cat$;
\item pour tout objet $b$ de $B$, le morphisme 
\[
\tranche{u}{b}: \tranche{A}{b}\to \tranche{B}{b}
\]
est dans $\W$;
\item pour tout objet $b$ de $B$, la catégorie $\tranche{A}{b}$ est
$\W$-asphérique.
\end{enumerate}
\end{prop}

\begin{prop}[Grothendieck]
Tout foncteur entre petites catégories admettant un adjoint à droite
est $\W$\nobreakdash-asphérique. En particulier, toute équivalence de
petites catégories est $\W$\nobreakdash-asphérique (et donc, est un
élément de~$\W$).
\end{prop}

\begin{prop}[Grothendieck]
Toute petite catégorie admettant un objet initial
est~$\W$\nobreakdash-asphérique. \end{prop}

\begin{prop}[Grothendieck]
Soit $A \xrightarrow{u} B \xrightarrow{v} C$ un couple de morphismes
composables de $\Cat$. Si $u$ est $\W$-asphérique, alors $vu$ est
$\W$\nobreakdash-asphérique si et seulement si $v$ est $\W$-asphérique.
\end{prop}

On rappelle (voir la section \ref{secCatTest}) que si $A$ est une
petite catégorie, on note 
\[
i_A : \pref{A}\to\Cat \mdvirg X \mapsto \tranche{A}{X}
\]
le foncteur associant à tout préfaisceau sur $A$ sa catégorie des éléments.
On a déjà vu dans la proposition \ref{i_Ahocolim} que pour toute
petite catégorie $A$, le foncteur~$i_A$ calcule la limite inductive homotopique des
préfaisceaux d'ensembles, vus comme des préfaisceaux en petites
catégories discrètes. On dispose donc d'un modèle \emph{strict}
particulièrement pratique pour calculer la limite inductive homotopique. 

Dans ce cadre, le formalisme des morphismes asphériques que l'on
vient d'introduire permet de caractériser les morphismes commutant aux
limites inductives homotopiques, comme on va maintenant l'expliquer.

\paragr\label{def:PrefaisceauAspherique} Si $A$ est une petite catégorie, on dit qu'un préfaisceau
$X$ sur~$A$ est~\ndef[préfaisceau!asphérique!relativement à un
localisateur fondamental]{$\W$\nobreakdash-asphérique} si son image par le
foncteur $i_A$
est une catégorie $\W$\nobreakdash-asphérique, c'est-à-dire si le morphisme
$p_{\tranche{A}{X}}: \tranche{A}{X} \to e$, où $e$ désigne la
catégorie ponctuelle, est dans $\W$. On dit qu'un morphisme de
préfaisceaux~$f : X\to Y$ sur $A$ est une
\ndef[équivalence faible!de préfaisceaux!relativement à
un localisateur fondamental]{$\W$\nobreakdash-équivalence} si
son image par le foncteur~$i_A$ est un élément de $\W$. 
On dit également qu'un préfaisceau $X$ sur $A$ est
\ndef[préfaisceau!localement asphérique!relativement à un localisateur
fondamental]{localement $\W$-asphérique} si l'image du morphisme~$X
\to e_{\pref{A}}$ par le foncteur $i_A$ est un morphisme
$\W$-asphérique de~$\Cat$.

\paragr 
Si $A$ est une petite catégorie et $u :A\to B$ est un foncteur, on
dispose, pour tout préfaisceau $X$ sur ${B}$, d'un foncteur
\begin{align*}
\tranche{A}{u^*(X)}&\xrightarrow{\tranche{u}{X}} \tranche{B}{X} \\
(a, a\xrightarrow{x}u^*(X)) &\mapsto (u(a), u(a) \xrightarrow{x} X)
\end{align*}
et on vérifie qu'on a ainsi défini un morphisme de foncteurs
$\lambda_u : i_Au^* \to i_B$, s'insérant donc dans le diagramme 
\begin{equation}\label{diag:morphismeWAspheriqueHomotCofinal}
\xymatrix{
{\pref{B}} 
\ar[rr]^{u^*} \ar[rd]_{i_B}^{}="c"
&& {\pref{A}} 
\ar@{}"c"|(.4){}="d"|(.9){}="e"
\ar@2"d";"e"_{\lambda_u}
\ar[ld]^{i_A} \\
& {\Cat} &\pbox{.}
} 
\end{equation}

On peut alors montrer la proposition suivante, en utilisant
exclusivement les propriétés du foncteur $i_A$ et des localisateurs
fondamentaux.

\begin{prop}[Grothendieck]\label{propMorphismesAspheriquesLocFond}
Soit $u:A\to B$ un morphisme de $\Cat$. Les conditions suivantes sont
équivalentes : \begin{enumerate}
\item u est un morphisme $\W$-asphérique de $\Cat$;
\item pour qu'un préfaisceau $X$ sur $B$ soit $\W$-asphérique, il faut
et il suffit que~$u^*(X)$ soit un préfaisceau $\W$-asphérique sur $A$;
\item pour tout objet $b$ de $B$, le préfaisceau $u^*(b)$ est
$\W$-asphérique;
\item pour tout préfaisceau $X$ sur $B$, le morphisme $\lambda_{u,X}$
dans le diagramme~\ref{diag:morphismeWAspheriqueHomotCofinal} est dans
$\W$.
\end{enumerate}
De plus, ces conditions impliquent la propriété suivante : 
\begin{enumerate}[resume]
\item pour
qu'un morphisme $f : X\to Y$ soit une $\W$-équivalence de $\pref{B}$, il faut et il
suffit que $u^*(f)$ soit une $\W$-équivalence de $\pref{A}$.
\end{enumerate}
\end{prop}
\begin{proof}
Voir \cite[proposition 1.2.9]{maltsiniotis2005}.
\end{proof}

Pour terminer cette courte introduction aux localisateurs fondamentaux,
on énonce le résultat suivant,
conjecturé par Grothendieck et démontré par
Cisinski dans \cite{cisinski2004localisateur}, qui fournit une
caractérisation de la classe $\W_{\infty}$ des équivalences de
Thomason d'une manière entièrement catégorique, et en particulier
n'utilisant pas le nerf simplicial. 

\begin{theorem}[Grothendieck, Cisinski]
\label{locFondamentalMinimalThm}
La classe $\W_{\infty}$ des équivalences de Thomason est le
localisateur fondamental minimal : si $\W$ est un
localisateur fondamental, alors $\W_{\infty}\subset \W$.
\end{theorem}
\begin{proof}
Voir \cite{cisinski2004localisateur}.
\end{proof}

\section{Morphismes asphériques en homologie}
\label{secMorphismesAsphEnHomologie}
On rappelle qu'on dispose d'un foncteur 
\[
\Whf{} : \Cat \to \Hotab \mdvirg A \mapsto \H{A}{\Z}
\]
défini au paragraphe \ref{homologieZ=HomologieSinguliere}, et que pour
toute petite catégorie $A$, on a un isomorphisme naturel
\[
\H{A}{\Z} \simeq \H{\Delta}{\Z^{(\nerf A)}} 
\]
dans $\Hotab$, où on a noté $\Z^{(\nerf A)}$ le groupe abélien
simplicial libre sur le nerf de $A$. En particulier, ce
foncteur envoie les éléments de $\W_\infty$ sur des isomorphismes.

\paragr On dit
qu'un foncteur $u : A \to B$ est une \ndef[équivalence!faible!de
$Cat$!abélienne]{équivalence faible abélienne} de $\Cat$ si son
image par le foncteur $\Whf{}$ est un isomorphisme dans $\Hotab$. On
note $\W_\infty^\ab$ la casse des équivalences faibles abéliennes de
$\Cat$.
\notindex{$\W_{\infty}^\ab$}%

\begin{prop}
\label{localisateurFondamentalAbelien}
$ \W_\infty^\ab $ est un localisateur fondamental.
\end{prop}
\begin{proof}
On renvoie à l'annexe \ref{annexeWabLocFond}.
\end{proof}

On introduit alors la terminologie suivante, traduisant simplement le
fait qu'on travaille maintenant avec le localisateur fondamental
$\W_\infty^\ab$.

\paragr \label{defAspheriqueEnHomologie} On dit qu'une
petite catégorie $A$ est \ndef[catégorie!asphérique!en homologie]{asphérique en
homologie} si le morphisme canonique 
\[
p_A : A \to e 
\]
est un élément de $\W_\infty^\ab$, où on note toujours $e$ la
catégorie ponctuelle. On dit qu'un morphisme
$u:A\to B$ de~$\Cat$ est \ndef[morphisme asphérique!de $\Cat$!en
homologie]{asphérique en homologie} si pour tout objet $b$ de $B$, la
catégorie $\tranche{A}{b}$ est asphérique en homologie. 

\begin{remark}
On réservera toujours l'adjectif \og asphérique \fg{} pour parler de
$\W_\infty$-asphéricité.
\end{remark}

\begin{prop}\label{propAspheriqueImpliqueAsphHomologie}
On a $\W_\infty \subset \W_\infty^\ab$. En particulier, tout morphisme asphérique est
asphérique en homologie, et toute catégorie asphérique est asphérique
en homologie.
\end{prop}
\begin{proof}
Cela découle immédiatement du corollaire \ref{homologieZ=HomologieSinguliere}.
\end{proof}

La proposition \ref{propMorphismesAspheriquesLocFond} nous permet de
décrire les morphismes $\W_\infty^\ab$\nobreakdash-asphériques comme les
morphismes homotopiquement cofinaux au sens, \emph{a priori}
insuffisant, où un foncteur~$u : A \to B$ est asphérique en homologie
si et seulement la transformation naturelle $\lambda_u$ dans le
diagramme
\[
\xymatrix{
{\pref{B}} 
\ar[rr]^{u^*} \ar[rd]_{i_B}^{}="c"
&& {\pref{A}} 
\ar@{}"c"|(.4){}="d"|(.9){}="e"
\ar@2"d";"e"_{\lambda_u}
\ar[ld]^{i_A} \\
& {\Cat} 
} 
\]
est dans $\W_\infty^\ab$ argument par argument. Autrement dit, en
utilisant la proposition \ref{homologiePrefaisceauxLibres}, si pour
tout préfaisceau~$X$ sur $B$, on a un isomorphisme naturel
\[
\H{A}{\Z^{(u^*X)}} \simeq \H{B}{\Z^{(X)}}
\]
dans $\Hotab$. La proposition suivante affirme que cela suffit pour
obtenir la cofinalité homotopique pour tous les préfaisceaux en
groupes abéliens.

\begin{prop}
\label{morphismesAsphEnHomologieEquivalences}
\label{morphismesAsphEnHomologieRefleteEqf}
Soit $u:A\to B$ un morphisme de $\Cat$. Les conditions suivantes sont
équivalentes : \begin{enumerate}
\item $u$ est un morphisme asphérique en homologie;
\item pour tout objet $b$ de $B$, la catégorie $\tranche{A}{b}$ a le
type d'homologie du point;
\item pour tout diagramme $X: {A}^{\op}\to \Ch(\Ab)$, le morphisme
canonique 
\[
\hocolim^{\Hotab}_{{A}^{\op}}u^*(X) \to
\hocolim^{\Hotab}_{{B}^{\op}} X 
\]
est un isomorphisme dans $\Hotab$;
\item le diagramme 
\[
\xymatrix{
{\prefab{B}} \ar[rr]^{u^*} \ar[rd]_{\Hf{B}} && {\prefab{A}}
\ar[ld]^{\Hf{A}} \\
& {\Hotab}
} 
\]
est commutatif à isomorphisme naturel près.
\end{enumerate}
De plus, toutes les conditions précédentes impliquent la condition
\begin{enumerate}[resume]
\item pour tout morphisme $f$ de $\prefab{B}$, on a l'équivalence 
\[
f \in \W_B^\ab \iff u^*(f) \in \W_A^\ab \quad (\ref{defWab_A})\pbox{.}
\]

\end{enumerate}

\end{prop}
\begin{proof}
L'équivalence de $(1)$ et $(2)$ est un cas particulier de la
proposition \ref{morphismesWAspheriquesEquivalences}. L'équivalence
de $(2)$ et $(3)$ provient de la caractérisation de~$\W_\infty^\ab$
comme les équivalences pour le dérivateur $\DerHotab$ décrite dans
l'annexe~\ref{appendix:DerHotab}. 
L'implication $(3)\implies(4)$ est
immédiate, ainsi que l'implication $(4)\implies(5)$.
Démontrons que $(4)\implies(2)$. Pour tout objet $b$ de $B$, on a des
isomorphismes naturels \begin{align*}
\H{\tranche{B}{b}}{\Z} &\simeq \H{B}{\Wh{B}{b}} & (\ref{homologieZ=HomologieSinguliere})\\
&\simeq \H{A}{u^*\Wh{B}{b}} & (\text{par hypothèse}) \\
&\simeq \H{A}{\Wh{A}{u^*(b)}} & (\ref{paragr:commutativiteu^*WhetU})\\
&\simeq \H{\tranche{A}{u^*(b)}}{\Z} &
(\ref{homologieZ=HomologieSinguliere})
\end{align*}
et on peut conclure, en remarquant qu'on a
$\tranche{A}{b}=\tranche{A}{u^*(b)}$, et que puisque la catégorie $\tranche{B}{b}$ a
un objet final, elle a le type d'homologie du point, ce qui entraine
que $\tranche{A}{b}$ a également le type d'homologie du point.
\end{proof}

\begin{coro}
Soient $A$ une catégorie, $L_A$ un intégrateur~(\ref{defIntegrateurs})
sur~$A$, et~$u:A\to B$ un morphisme asphérique en homologie. Alors pour tout
préfaisceau $X\in \prefab{B}$, on a un isomorphisme canonique naturel
\[
\H{B}{X} \simeq {L_A}_! u^*(X)
\]
dans $\Hotab$.
\end{coro}
\begin{remark}
On montrera dans la proposition \ref{coro:integrateurInduit} que le
foncteur~${L_A}_!u^*$ de la proposition précédente est également
l'extension de Kan à gauche d'un intégrateur.
\end{remark}

\begin{example}\label{exTotalementAspheriqueEilenbergZilber}
On s'intéressera particulièrement aux catégories $A$ telles que le
foncteur diagonal $\delta :A\to A\times A$ est asphérique en
homologie (on dira que $A$ est totalement asphérique en homologie, voir section
\ref{secCatTestHomStrictes}). Cela implique qu'on peut calculer
l'homologie de $A\times A$ à coefficients dans un préfaisceau abélien
$X$ en utilisant la diagonale : on a un isomorphisme naturel 
\[
\H{A\times A}{X} \simeq \H{A}{\delta^* X} 
\]
dans $\Hotab$. En particulier, un morphisme entre deux préfaisceaux
abéliens sur~$A\times A$ est une équivalence faible abélienne si et
seulement si son image par le foncteur diagonal est une équivalence
faible abélienne de $\prefab{A}$. La catégorie~$\Delta$ est un exemple
de catégorie totalement asphérique (\ref{exTotAspherique}), et est
donc totalement asphérique en homologie. 

En utilisant l'intégrateur produit (\ref{integrateurProduit}) pour
calculer l'homologie des groupes abéliens bisimpliciaux, on
retrouve une preuve de la généralisation par Dold et Puppe du théorème
d'Eilenberg-Zilber~\cite{eilenbergZilber1953products,doldpuppe1961}:
si $X$ est un groupe abélien bisimplicial, on a un quasi-isomorphisme 
\[
\Tot(\dk X) \simeq \dk \delta(X)
\]
entre le complexe total du complexe double obtenu en appliquant terme
à terme le foncteur complexe de chaînes normalisé à $X$, et le complexe
normalisé associé à la diagonale de $X$. On peut ainsi voir les catégories
totalement asphériques en homologie comme des catégories vérifiant une
généralisation du théorème d'Eilenberg-Zilber.
\end{example}

\section{Induction et coinduction}
Le problème principal auquel nous sommes confrontés pour abélianiser la
théorie des catégories test est qu'on ne dispose pas
d'un \og meilleur choix \fg{} pour remplacer le foncteur $i_A$, bien
qu'on dispose d'un intégrateur canonique que l'on définira dans la
section \ref{defIntegrateurBousKan}. Aussi, étant donnés un foncteur
$u : A \to B$ asphérique en homologie et deux intégrateurs $L_A$ et
$L_B$ respectivement sur $A$ et $B$, l'espoir d'obtenir une
transformation naturelle explicite s'insérant dans le diagramme 
\[
\xymatrix{
\prefab{B} \ar[rd]_{{L_B}_!}|{}="b" 
\ar[rr]^{u^*}="c" 
&& 
\prefab{A}
\ar@{}"b"|(.5){}="x"|(.8){}="y"
\ar@2"x";"y"_{\lambda_u}
\ar[ld]^{{L_A}_!}_(0.5){}="a" \\
& \Ch(\Ab) 
} 
\]
et telle que $\lambda_u$ soit un quasi-isomorphisme argument par
argument est déraisonnable. La stratégie que nous allons exposer dans
cette section consiste à 
transférer un intégrateur le long d'un foncteur asphérique. On pourra
alors, sous certaines conditions, expliciter une telle transformation
naturelle. 

\paragr Soient $C$ et $D$ deux petites catégories, et $v : C \to D$ un
foncteur. On rappelle que 
le foncteur $v$ induit alors un foncteur 
\[
\prefab{D}\xrightarrow{v^*}\prefab{C} 
\]
admettant un adjoint à gauche $v^\ab_!$ et un adjoint à droite
$v_*^\ab : \prefab{C}\to\prefab{D}$. L'adjoint à droite $v_*^\ab$ est
\notindex{$u_\bang^\ab$}%
\notindex{$u_*^\ab$}%
induit par le foncteur $v_*: \pref{C}\to \pref{D}$, puisqu'il préserve
les produits et envoie donc les groupes abéliens internes à $\pref{C}$
sur des groupes abéliens internes à $\pref{D}$. On obtient donc un
diagramme commutatif 
\[
\xymatrix{
\prefab{C} \ar[d]_{v_*^\ab} \ar[r]^{\U} & \pref{C} \ar[d]^{v_*} \\
\prefab{D} \ar[r]_{\U} & \pref{D} \pbox{.}
} 
\]
qui donne une description complète du comportement des extensions de
Kan à droite lors du processus d'abélianisation. Pour décrire
l'adjoint à gauche~$v_!^\ab : \prefab{C}\to\prefab{D}$, 
on peut utiliser la proposition~\ref{kanextadj} : si $L$ est un objet
de $\Addinf(C)$, l'enveloppe additive infinie de $C$ définie en
\ref{defAddinf}, alors on a 
\[
v_!^\ab(L) = \Addinf(v)(L) \pbox{,}
\]
ce qui signifie qu'en notant $L=\bigoplus_{i}\Wh{C}{c_i}$ où les $c_i$
sont des objets de $C$, on a 
\[
v_!^\ab(L) = \bigoplus_i \Wh{D}{v(c_i)} \pbox{.}
\]
En particulier, on obtient le diagramme commutatif suivant : 
\[
\xymatrix{
\pref{C} \ar[r]^{\Whf{C}} \ar[d]_{v_!} & \prefab{C} \ar[d]^{v_!^\ab} \\
\pref{D} \ar[r]_{\Whf{D}}  & \prefab{D} \pbox{.}
} 
\]

Si $X$ est un préfaisceau abélien sur $C$, on peut alors (voir
\ref{densiteAdd}) le décomposer comme une limite inductive de
préfaisceaux dans l'enveloppe additive finie  
\[
X \simeq \limind_{\substack{L\to X \\\mathclap{ L \in \Add(C)}}} L
\]
et on obtient un isomorphisme canonique naturel 
\[
v_!^\ab(X) \simeq \limind_{\substack{L\to X \\\mathclap{ L \in \Add(C)}}}
\Add(v)(L) \pbox{.}
\]

\begin{lemme}
Si $v:C\to D$ est un foncteur, le foncteur $v_!^\ab :
\prefab{C}\to\prefab{D}$ envoie les objets projectifs sur des objets
projectifs.
\end{lemme}
\begin{proof}
C'est un cas particulier du corollaire \ref{coro:fonctPreserveProj},
puisque le foncteur $v_!^\ab$ préserve les préfaisceaux
représentables. On peut aussi utiliser le fait que son adjoint à
droite $v^*$ est un foncteur exact.
\end{proof}

\medskip
\emph{On fixe, jusqu'au corollaire \ref{coro:integrateurInduit}, deux petites catégories $A$ et $B$, un
intégrateur~$L_A$ sur~$A$, et un foncteur $u : A \to B$.}
\bigskip

On va produire, à partir de $L_A$, un
intégrateur sur la petite catégorie $B$. Le morphisme $u$ induit un foncteur
\[
{u}^{\op} : {A}^{\op}\to {B}^{\op}
\]
et on va montrer que si $u$ est asphérique en homologie, alors le
complexe de préfaisceaux sur ${B}^{\op}$ obtenu à partir de $L_A$ en
appliquant degré par degré le foncteur
\[
  ({u}^{\op})_!^\ab : \prefab{{A}^{\op}}\to \prefab{{B}^{\op}}
  \]
est un intégrateur sur $B$.

\begin{remark}
On sait déjà que si $u$ est asphérique en homologie, on a un triangle
commutatif à isomorphisme près
\[
\xymatrix{
{\Hotab_B} \ar[rr]^{u^*} \ar[rd]_{\Hf{B}} && {\Hotab_A}
\ar[ld]^{\Hf{A}} \\
& {\Hotab}
}  
\]
et des isomorphismes canoniques naturels dans $\Hotab$ 
\begin{align*}
\Hf{A} \simeq {L_A}_! \mdvirg \Hf{B} \simeq {L_A}_!u^* \pbox{.}
\end{align*}
Le but de ce qui suit est alors de montrer que le morphisme
${L_A}_!u^*$ est bien l'extension de Kan à gauche d'un intégrateur sur
$B$. Pour cela, on va montrer qu'il existe un isomorphisme 
\[
{L_A}_!u^* \simeq ({u}^{\op})_!^\ab ({L_A}_\bullet) 
\]
dans la catégorie $\Ch(\Ab)$.
\end{remark}

\paragr On sait que les préfaisceaux représentables de $B$ ont le type
d'homologie du point. Ainsi, si $u$ est un morphisme asphérique en
homologie, le foncteur~$\prescript{A}{}{L}_B : B \to \Ch(\Ab)$ défini
par le diagramme commutatif
\[
\xymatrix{
B \ar[rd]_{\prescript{A}{}{L}_B} \ar[r]^{\Whf{B}} & \prefab{B} \ar[r]^{u^*} &
\prefab{A} \ar[ld]^{(L_A)_!} 
\\
& \Ch(\Ab)
} 
\]
envoie les objets de $B$ sur des complexes de chaînes ayant
l'homologie du point. 

On va maintenant prouver que~$\prescript{A}{}{L}_B$ est bien un
complexe d'objets projectifs de~$\prefab{{B}^{\op}}$.

\begin{lemme}\label{lemme:extKanNaturalité}
Le diagramme 
\[
\xymatrix{
\Homi_!(\prefab{A},\Ab) \ar[r]^{-\circ u^*} &
\Homi_!(\prefab{B},\Ab) \ar[d]^{- \circ \Whf{B}} \\
\prefab{{A}^{\op}} \ar[u]^{(-)_!} 
\ar[r]_{({u}^{\op})_!^\ab} &
\prefab{{B}^{\op}}
}
\]
est commutatif à isomorphisme canonique naturel près, où les flèches
verticales désignent les foncteurs intervenant dans l'équivalence de
catégories de la proposition \ref{prefabBeqcat}.

\end{lemme}
\begin{proof}
Une limite inductive de foncteurs cocontinus est
un foncteur cocontinu, et les limites inductives dans
$\Homi_!(\prefab{A},\Ab)$ se caculent donc point par point.
Cela implique que le foncteur $-\circ u^*$ commute aux limites
inductives, et donc que tous les foncteurs de ce diagramme commutent aux
limites inductives. Il suffit donc de vérifier la commutativité du
diagramme pour des
préfaisceaux abéliens représentables de ${A}^{\op}$. 

Or, si $a$ est un
objet de~${A}^{\op}$, alors on a  
\[
({u}^{\op})_!^\ab(\Wh{{A}^{\op}}{a})=\Wh{{B}^{\op}}{{u}^{\op}(a)} \mdvirg
\]
qui est une manière compliquée de désigner le foncteur 
\begin{align*}
B &\to \Ab \\
b &\mapsto \Z^{(\Hom_{{B}^{\op}}(b,{u}^{\op}(a)))} = 
\Z^{(\Hom_B(ua,b))} \pbox{.}
\end{align*}

Par ailleurs, l'image de $\Wh{{A}^{\op}}{a}$ par la flèche verticale de gauche est
le foncteur~$\prefab{A}\to\Ab$ défini par 
\[
(\Wh{{A}^{\op}}{a})_! = \ev_a : X \mapsto X{a} \pbox{,}
\]
et son image par le foncteur $- \circ u^*$ est donc le foncteur
\[
\ev_a \circ u^* = \ev_{u(a)} : Y \mapsto Yu(a) \pbox{.}
\]
Finalement, son image dans $\prefab{{B}^{\op}}=\Hom(B,\Ab)$ est donc
bien le foncteur 
\[
b \mapsto \Wh{B}{b}(ua) = \Z^{(\Hom_B(ua,b))} \pbox{.} \qedhere
\]
\end{proof}

\paragr Puisque $-\circ\Whf{B}$ est une équivalence de catégories, le
lemme précédent implique donc que le diagramme 
\[
\xymatrix{
\Homi_!(\prefab{A},\Ab) \ar[r]^{-\circ u^*} &
\Homi_!(\prefab{B},\Ab)  \\
\prefab{{A}^{\op}} \ar[u]^{(-)_!} 
\ar[r]_{({u}^{\op})_!^\ab} &
\prefab{{B}^{\op}} \ar[u]_{(-)_!} 
}
\]
est commutatif à isomorphisme près. En d'autres
termes, si $L$ est un préfaisceau abélien sur ${A}^{\op}$ et si $Y$
est un préfaisceau abélien sur $B$, alors on a
un isomorphisme canonique naturel de groupes abéliens
\[
Y \odot_B ({u}^{\op})_!^\ab(L) \simeq u^*(Y) \odot_A L 
\]
(on renvoie au paragraphe \ref{couplageOdot} pour la définition du
couplage $\odot$).
\begin{coro}\label{coro:integrateurInduit}
Soit $L_A$ un intégrateur sur $A$ et $u:A\to B$ un foncteur
asphérique en homologie. Alors le foncteur $\prescript{A}{}{L}_B : B
\to \Ch(\Ab)$
est un intégrateur sur $B$.
\end{coro}
\begin{proof}
Pour tout entier $n\geq 0$, le préfaisceau abélien $(\prescript{A}{}{L}_B)_n$
sur~${B}^{\op}$ est bien projectif, puisque le lemme
\ref{lemme:extKanNaturalité} implique qu'il est l'image par le
foncteur $({u}^{\op})_!^\ab$ de l'objet projectif $(L_A)_n$ de
$\prefab{{A}^{\op}}$. De plus, $u$ est asphérique en homologie, donc
pour tout objet $b$ de $B$, le complexe ${L_A}_!u^*\Wh{B}{b}$ a l'homologie du point. On
en conclut que $\prescript{A}{}{L}_B$ est bien une résolution
projective du préfaisceau constant de valeur $\Z$.
\end{proof}

On remarque de plus que, par construction, le diagramme 
\[
\xymatrix{
{\prefab{B}} \ar[rd]_{{\prescript{A}{}{L}_B}_!} \ar[rr]^{u^*} && {\prefab{A}}
\ar[ld]^{{L_A}_!} \\
& {\Ch(\Ab)}
}  
\]
est commutatif à isomorphisme canonique naturel près.

\bigskip 
\paragr \label{notationsInduction}\emph{On fixe maintenant, et jusqu'à la fin de
cette section, deux petites catégories~$A$ et $C$, un
intégrateur~$L_A$ sur~$A$ et un foncteur $u : C\to A$}.
\bigskip

On va s'intéresser à des conditions pour que le foncteur
\[
  {L^A_C} : C \to \Ch(\Ab)
\]
défini par le diagramme commutatif 
\[
\xymatrix{
{C} \ar[rr]^{u} \ar[rd]_{L_C^A} && {A} \ar[ld]^{L_A} \\
& {\Ch(\Ab)}
}  
\]
soit un intégrateur sur $C$. 

\paragr Il est immédiat que si $c$ est un objet de $C$,
alors le complexe~$L_C^A(c)$ a bien l'homologie du point. Il n'est pas
clair en revanche que~$L_C^A$ soit un complexe de préfaisceaux
projectifs sur ${C}^{\op}$. On peut cependant exhiber des conditions
suffisantes pour que ce soit le cas.

\begin{prop}\label{conditionsTransfertIntegrateur}
On suppose que le foncteur
\[
  (u^{\op})^* : \prefab{A^{\op}}\to \prefab{C^{\op}}
\]
envoie les
préfaisceaux représentables de $\prefab{A^{\op}}$ sur des objets
projectifs. Alors le foncteur $L_C^A$ défini ci-dessus
est un intégrateur sur $C$. En particulier, pour tout préfaisceau
abélien $X$ sur $C$, on a un isomorphisme canonique naturel
\[
\H{C}{X} \simeq \H{A}{u_!^\ab X} 
\]
dans $\Hotab$, et le foncteur $u_!^\ab$ préserve donc les équivalences
faibles abéliennes.
\end{prop}
\begin{proof}
Le foncteur $L_C^A$, vu comme un complexe de préfaisceaux sur
${A}^{\op}$, est obtenu à partir de $L_A$ en appliquant degré par
degré le foncteur 
\[
{(u^{\op})}^* : \prefab{{A}^{\op}}\to \prefab{{C}^{\op}} \pbox{.}
\]
On sait de plus que $({u}^{\op})^*$ est un foncteur additif. Puisque
par hypothèse, ce foncteur envoie les préfaisceaux représentables sur
des objets projectifs, il envoie les objets de $\Addinf({A}^{\op})$
sur des objets projectifs. On peut donc appliquer le corollaire
\ref{coro:fonctPreserveProj}, qui nous permet de conclure que le
foncteur $({u}^{\op})^*$ préserve les objets projectifs. Si $L_A$ est
un intégrateur sur $A$, alors $L_C^A$ est donc bien une
résolution projective du préfaisceau constant de valeur $\Z$ sur ${C}^{\op}$.
\end{proof}

\begin{example}
En revenant au cas particulier de l'homologie des groupes, si $u :
H \to G$ est un morphisme de groupes, alors le morphisme
\[
  u^* : \Hom(BG, \Ab) \to \Hom(BH,\Ab)
\]
est traditionnellement appelé le foncteur de \emph{restriction}. Son
adjoint à gauche 
\[
  u_!^\ab : \Hom(BH,\Ab) \to \Hom(BG,\Ab)
\]
est alors le foncteur \emph{d'induction}, et l'adjoint à droite 
\[
u_*^\ab : \Hom(BH,\Ab) \to \Hom(BG,\Ab)
\]
est
le foncteur de \emph{coinduction}. Un résultat classique attribué à
Shapiro (voir par exemple~\cite[6.3.2]{weibel1994Introduction})
affirme alors que si $u$
est une inclusion de sous-groupe, alors il existe un 
quasi-isomorphisme $\H{BG}{u_!^\ab(M)} \simeq \H{BH}{M}$ pour
tout~$H$-module $M$. On peut voir ce résultat comme un cas particulier
de la commutativité du diagramme \ref{diagShapiroGen} et du fait que
si $u$ est l'inclusion d'un sous-groupe, alors $u^*(\Z G)$ est une
somme directe indexée par les éléments de~$\tranche{G}{H}$ du $H$-module $\Z
H$. Étant donné une résolution projective $L_\bullet$ de $\Z$ par
des $G$-modules (c'est-à-dire un intégrateur sur $BG$), le complexe
$u_!^\ab(L_\bullet)$ est donc un intégrateur sur $BH$. 
\end{example}

On peut exhiber une condition suffisante pour appliquer la
proposition~\ref{conditionsTransfertIntegrateur}, que l'on rencontrera
souvent dans la suite.

\paragr On rappelle qu'un morphisme $p : C \to A$ est une
\ndef{fibration à fibres discrètes} si pour tout objet~$c$ de $C$ et
tout morphisme $f : a \to p(c)$ de $A$, il existe un unique
morphisme~$u\in\Fl(C)$ de but $c$ tel que $p(u)=f$.

\begin{prop}\label{inductionFibrationFibDiscrete}
Si le foncteur $u : C \to A$ est une fibration à fibres discrètes,
alors le foncteur $L_C^A$ défini ci-dessus est bien un intégrateur.
\end{prop}
\begin{proof}
Si $a$ est un objet de $A$, on veut montrer que le
préfaisceau en groupes abéliens sur $C$
\[
  ({u}^{\op})^*\Wh{{A}^{\op}}{a}
\]
  est un objet projectif de
$\prefab{{C}^{\op}}$. On peut écrire plus simplement ce préfaisceau
sur ${C}^{\op}$ comme le foncteur 
\[
C \to \Ab \mdvirg c \mapsto \Z^{\left(\Hom_A(a,u(c))\right)} \pbox{.}
\]
Puisque $u$ est une fibration à fibres discrètes, alors pour tout
morphisme $f : a \to u(c)$, il existe un unique morphisme $c_f \to c$
de $C$ dont l'image par $u$ est le morphisme $f$.
On a alors un isomorphisme 
\[
\Z^{(\Hom_A(a,u(c)))} \simeq \bigoplus_{c_f \in C_a}
\Z^{(\Hom_C(c_f,c))}
\]
où on a noté $C_a$ la fibre de $u$ au-dessus de $a$.
En repassant aux catégories opposées, on a exhibé un isomorphisme
\[
({u}^{\op})^*\Wh{{A}^{\op}}{a} \simeq \bigoplus_{c\in
C_a}\Wh{{C}^{\op}}{c} \pbox{.}
\]
En particulier, $({u}^{\op})^*$ envoie les représentables sur des
objets de $\Addinf({C}^{\op})$, et envoie donc les préfaisceaux
projectifs sur ${A}^{\op}$ sur des préfaisceaux projectifs sur
${C}^{\op}$.
\end{proof}

\begin{coro}
\label{coro:inductionTranche} Si $A$ est une petite catégorie et $F$ est un
préfaisceau d'ensembles sur $A$, alors le
foncteur~$L_{\tranche{A}{F}}$ défini par le diagramme commutatif 
\[
\xymatrix{
\tranche{A}{F} \ar[rd]_{L_{\tranche{A}{F}}} \ar[r]^{\alpha_F} & A
\ar[d]^{L_A} \\
& \Ch(\Ab)
} 
\]
est un intégrateur sur $\tranche{A}{F}$.
On a donc un isomorphisme
naturel, pour tout préfaisceau en groupes abéliens $X$ sur $\tranche{A}{F}$, 
\[
\H{\tranche{A}{F}}{X} \simeq \H{A}{{\alpha_F}_!^\ab(X)} \pbox{,}
\]
et le foncteur 
\[
{\alpha_F}_!^\ab : \prefab{\tranche{A}{F}} \to \prefab{A} 
\]
préserve donc les équivalences faibles.
\end{coro}
\begin{proof}
Il s'agit simplement d'une reformulation de la proposition précédente,
puisque les fibrations à fibres discrètes correspondent exactement aux
morphismes de projections $\tranche{A}{F}\to A$.
\end{proof}

\begin{remark}
Il découle de la preuve de la proposition précédente que si $L_A$ est un
intégrateur libre, alors $L_{\tranche{A}{F}}$ est également un
intégrateur libre.
\end{remark}

\bigskip
\emph{On garde les notations précédentes, et on suppose désormais que le foncteur~$L_C^A$ définit bien un intégrateur sur $C$}.
\bigskip

\paragr Grâce à la commutativité à isomorphisme près du diagramme 
\begin{equation}\label{diagShapiroGen}
\xymatrix{
{\prefab{C}} \ar[rr]^{u_!^\ab} \ar[rd]_{{L_C^A}_!} &&
{\prefab{A}} \ar[ld]^{{L_A}_!} \\
& {\Ch(\Ab)} & \pbox{,}
} 
\end{equation}
on définit une transformation naturelle 
\[
\xymatrix{
\prefab{A} \ar[rd]_{{L_A}_!}^{}="b" \ar[rr]^{u^*}="c" 
&& \prefab{C}
\ar@{}"b"|(.5){}="x"|(.8){}="y"
\ar@2"x";"y"_{\lambda_u}
\ar[ld]^{{L^A_C}_!}_(0.5){}="a" \\
& \Ch(\Ab) 
} 
\]
de la manière suivante : on note $\epsilon : u_!^\ab u^* \to
\id_{\prefab{C}}$ le morphisme d'adjonction, et on définit le
morphisme $\lambda_u$ en posant 
\[
\lambda_u : {L_C^A}_!u^* \xrightarrow{\simeq} {L_A}_!u_!^\ab u^*
\xrightarrow{{L_A}_!(\epsilon)} {L_A}_! \pbox{.}
\]

\begin{prop}\label{integrateurInduitMorphismeAspheriqueTF}
Le morphisme $u$ est asphérique en homologie \emph{si et seulement si} le morphisme
$\lambda_u$ défini ci-dessus est un quasi-isomorphisme argument par
argument.
\end{prop}
\begin{proof}
On renvoie à l'annexe
\ref{appendix:morphismeHotabAspheriqueCounite} pour la preuve du fait
suivant : un foncteur
$u : A \to B$ est asphérique en homologie si et seulement
si pour tout préfaisceau abélien $X$ sur $B$ et pour toute résolution
projective $p : P_\bullet \to u^*(X)$ dans $\Ch(\prefab{A})$, la
flèche oblique dans le diagramme commutatif
\[
\xymatrix@C=4em{
\underline{L_A}_!\underline{u}_!^\ab P_\bullet
\ar[rd]^{} 
\ar[r]^{{\underline{L_A}}_!\underline{u}_!^\ab(p)} 
& 
{L_A}_!u_!^\ab u^*(X) \ar[d]^{{L_A}_!(\epsilon)} 
\\
& {L_A}_!(X)
}
\]
est un isomorphisme dans $\Hotab$, où on a noté $\underline{u}_!^\ab$
et $\underline{L_A}_!$ les foncteurs obtenus en appliquant degré par
degré les foncteur ${L_A}_!$ et $u_!^\ab$ aux complexes de chaînes de
préfaisceaux. 
Mais puisque $L_A\circ u = L^A_{C}$ est un intégrateur
sur~$C$, le foncteur~${\underline{L_A}}_!\underline{u}_!^\ab$ calcule
la limite inductive homotopique des complexes de préfaisceaux sur $C$
(voir \ref{colimiteHomotopiqueDegreparDegre}), et il préserve donc
les quasi-isomorphismes. La flèche horizontale du diagramme ci-dessus
est donc un quasi-isomorphisme, et, par deux-sur-trois, la flèche
oblique est donc un quasi-isomorphisme \emph{si et seulement si} il en
est de même pour la flèche verticale, qui n'est autre que le
morphisme $(\lambda_u)_X$.
\end{proof}

Ce résultat sera en particulier utile pour exhiber des \emph{intégrateurs
test locaux} (voir \ref{defIntegrateurTestLocal}).

\begin{prop}\label{integrateurCoinduitCommuteCounite}
En gardant les notations de~\ref{notationsInduction}, 
on note
\[
  \epsilon^A : {L_A}_!{L_A}^*\to\id_{\prefab{A}} \mdvirg
  \epsilon^C : {L^A_C}_!{L^A_C}^*\to\id_{\prefab{C}}
\]
les morphismes d'adjonction (voir \ref{propAdjIntegrateurs}). Alors pour tout complexe de groupes
abéliens $X$, le diagramme 
\[
\xymatrix@C=4em{
{{L^A_C}_!{L^A_C}^*(X)} \ar[rd]_{\epsilon^C} \ar@{=}[r] &
{{L^A_C}_!u^*L_A^*(X)} \ar[r]^-{\lambda_{u,L_A^*(X)}} &
{{L_A}_!L_A^*(X)} \ar[ld]^{\epsilon^A} \\
& {X}
} 
\]
dans $\Ch(\Ab)$ est commutatif.
\end{prop}
\begin{proof}
Il s'agit simplement du théorème de composition des adjonctions, que
l'on trouvera dans \cite[chapitre IV, 8.3]{maclane}. En particulier,
on n'a pas besoin de supposer que $L_C^A$ est un intégrateur sur $C$,
ni que $L_A$ est un intégrateur sur $A$.
\end{proof}

\section{Catégories de Whitehead}\label{secWhitehead}
Étant donnée une petite catégorie $A$, on rappelle que deux classes de morphismes de
$\prefab{A}$ sont candidates pour le rôle d'équivalences faibles : la
classe~$\Wab_A$ des morphismes dont l'images par le foncteur 
\[
\Hf{A} : \Hotab_A \to \Hotab 
\]
est un isomorphisme, introduite dans la section \ref{secPseudoTest}, et la classe $\U^{-1}\W_A$ des
morphismes dont le morphisme sous-jacent de préfaisceaux d'ensembles
sur $A$ est une équivalence test (\ref{secCatTest}).

\begin{defin}\label{defWhitehead}
On dit qu'une petite catégorie $A$ est une \ndef[catégorie!de
Whitehead]{catégorie de
Whitehead} si les classes de flèches $\Wab_A$ et $\U^{-1}\W_A$
coïncident.
\end{defin} 

L'exemple fondamental est donné par le théorème de
Dold-Kan.
\begin{prop}\label{DeltaWhitehead}
La catégorie $\Delta$ est une catégorie de Whitehead.
\end{prop}
\begin{proof}
Si $f : X\to Y$ est un morphisme de $\sAb$, on sait grâce à la
proposition \ref{homologieHomotopieNormalise} qu'on a, pour entier
$n>0$, un carré commutatif 
\[
\xymatrix{
{\pi_n(\U X, 0)} \ar[r]^{\pi_n(\U f)} \ar[d]_{\simeq}  & {\pi_n(\U Y, 0)}
\ar[d]^{\simeq}  \\
\mathsf{H}_n(X,\Z) \ar[r]_{\mathsf{H}_n(f)} & \mathsf{H}_n(Y,\Z)
}
\]
où $\mathsf{H}_n(-,\Z)$ désigne l'homologie du complexe normalisé, et
de même en degré~$0$. Puisque le complexe normalisé calcule bien
l'homologie au sens du foncteur $\Hf{\Delta}$ (voir
proposition~\ref{propHomologieDeltaMoore}), on peut alors conclure par
deux-sur-trois que le morphisme $f$ est dans $\Wab_{\Delta}$ si et seulement si $\U
f$ est dans $\W_{\Delta}$ (en utilisant également le théorème d'Illusie-Quillen
\cite[chapitre VI, théorème~3.3]{illusie1971cotangent} permettant
d'affirmer que les équivalences test sur $\Delta$ coïncident avec les
équivalences faibles simpliciales).
\end{proof}

Notre stratégie pour exhiber des catégories de Whitehead sera de les
comparer à $\Delta$, grâce à la proposition suivante :

\begin{prop}\label{aspheriqueSourceWhitehead}
Soit $A$ une catégorie de Whitehead. Si $u : A\to B$ est un foncteur
asphérique, alors $B$ est une catégorie de Whitehead.
\end{prop}
\begin{proof}
Pour tout morphisme $u:A\to B$, on a un
diagramme commutatif 
\begin{equation}\label{diag:commUu*}
\xymatrix{
\pref{B} \ar[r]^{u^*} \ar[d]_{\Whf{B}} & \pref{A} \ar[d]^{\Whf{A}} \\
\prefab{B} \ar[r]_{u^*} \ar[d]_{\U} & \prefab{A} \ar[d]^{\U} \\
\pref{B} \ar[r]_{u^*} & \pref{A}  
}
\end{equation}
et on rappelle que tout morphisme asphérique est asphérique en
homologie. Ainsi, si $u$ est asphérique, on a la chaîne d'équivalences
suivante : 
\begin{align*}
f \in \Wab_B &\iff u^*(f) \in \Wab_A
& (\text{$u$ asphérique en homologie}) \\
&\iff \U u^*(f) \in \W_A & (\text{$A$ de Whitehead})\\
&\iff u^* \U f \in \W_A &\text{(commutativité)}\\
&\iff \U f \in \W_B &\text{($u$ asphérique)}
\end{align*}
ce qui montre bien que $B$ est une catégorie de Whitehead.
\end{proof}

\begin{remark}
On a bien besoin que le foncteur $u$ soit asphérique, c'est-à-dire
$\W_\infty$-asphérique : l'asphéricité en
homologie ne suffit pas ici.
\end{remark}

\begin{remark}
Si $A$ est une petite catégorie et $u : \Delta \to A$ est un morphisme
asphérique, alors on a en fait une propriété \emph{a priori} plus
forte, établissant un isomorphisme entre les groupes d'homotopie et
les groupes d'homologie des préfaisceaux abéliens sur $A$. En effet,
on sait (voir la proposition \ref{propMorphismesAspheriquesLocFond})
que sous cette condition, pour tout préfaisceau d'ensembles $X$ sur
$A$, le morphisme 
\begin{align*}
\tranche{u}{X} : \tranche{\Delta}{u^*(X)} &\to \tranche{A}{X} \\
(\Delta_n, u(\Delta_n) \xrightarrow{x} X) &\mapsto (u(\Delta_n), x)
\end{align*}
est une équivalence de Thomason (\ref{def:eqThomason}).
Ainsi, si $X$ est un préfaisceau en groupes abéliens sur $A$, on
dispose de la chaîne d'isomorphismes naturels de groupes abéliens
suivante, pour tout entier $n\geq1$ :
\begin{align*}
\mathsf{H}_n(A,X) 
&\simeq \mathsf{H}_n(\Delta, u^* X) & \text{(asphéricité en homologie de $u$)} \\
&\simeq \pi_n(\U u^*(X), 0) & \text{(proposition
\ref{homologieHomotopieNormalise})}\\
&\simeq \pi_n(u^*\U(X),0) &\text{(commutativité de \ref{diag:commUu*})} \\
&\simeq \pi_n\left(\nerf i_{\Delta}{u^*\U(X)},(\Delta_0,0)\right) \\
&\simeq \pi_n\left(\nerf i_{A}{\U X}, (u(\Delta_0),0)\right)
&\text{(asphéricité de $u$)} 
\end{align*}
où l'avant dernier isomorphisme provient du théorème d'Illusie-Quillen
\cite[chapitre VI, théorème 3.3]{illusie1971cotangent}. On vérifie
qu'on obtient de la même manière une bijection sur le $\pi_0$. 
\end{remark}

La proposition précédente implique par exemple que la catégorie $\Theta_n$ de
Joyal est une catégorie de Whitehead, pour tout entier $n\geq1$
(\ref{secTheta}). En revanche,
elle ne suffit pas pour prouver que $\Theta$ est de Whitehead. Pour
cela, on utilise une variation autour de la notion d'asphéricité.
\footnote{Dans \emph{Pursuing Stacks}, Grothendieck introduit la notion de
\emph{structures d'asphéricité} dont cette définition (comme celle
d'asphéricité liée à un localisateur fondamental) est un cas
particulier.}

\paragr\label{defFoncteurAspheriqueÂ} Soit $\W$ un localisateur
fondamental. On dit qu'un \ndef[foncteur asphérique!de but une
catégorie de préfaisceaux]{foncteur}
\[
  j: B \to \pref{A}
\]
est $\W$-\emph{asphérique} si il vérifie les conditions suivantes : \begin{enumerate}
\item pour tout objet $b$ de $B$, le préfaisceau $j(b)$ sur ${A}$ est
$\W$-asphérique (\ref{def:PrefaisceauAspherique});
\item pour tout objet $a$ de $A$, le préfaisceau ${j^*(a)}$ sur ${B}$
est $\W$-asphérique;
\item pour tout objet $b$ de $B$, le préfaisceau $j^*j(b)$ sur ${B}$ est
$\W$-asphérique.
\end{enumerate}

On note en particulier qu'un morphisme $u : B \to A$ de $\Cat$ 
est $\W$-asphérique si et seulement si le foncteur
\[
  B \xrightarrow{u} A \hookrightarrow \pref{A}
\]
est $\W$-asphérique.

\paragr\label{paragrFoncteurAspheriquediag} Si $j:B\to \pref{A}$ est un foncteur, on considère la
sous-catégorie pleine~$C$ de $\pref{A}$ contenant les représentables
et les images par $j$ des objets de $B$. On dispose alors d'un
diagramme commutatif 
\begin{equation}\label{diagMorphismesAspheriquesStructureÂ}
\xymatrix{
& \pref{A} \\
B \ar[ru]^{j} \ar[r]_{j'} & C \ar@{^{(}->}[u]^{i} & A \ar@{_{(}->}[lu]_{h} \ar[l]^{h'} 
} 
\end{equation}
où on a noté $h$ le plongement de Yoneda.
On vérifie alors qu'on obtient un nouveau diagramme commutatif 
\begin{equation}\label{diagMorphismesAspheriquesStructureÂ2}
\xymatrix{
& \pref{A} \ar[ld]_{j^*} \ar[d]^{i^*} \ar[rd]^{\id} \\
\pref{B} & \pref{C} \ar[l]^{j'^*} \ar[r]_{h'^*} & \pref{A}
} 
\end{equation}
où le sens de \og l'étoile en haut \fg{} de $j $ et $i$ est
légèrement différent de celui de~$j'$ et $h'$ : les foncteurs $j^*$ et
$i^*$ désignent les foncteurs nerfs associés à $i$ et $j$. Par
exemple, $i^*$ est le foncteur
\begin{align*}
i^* : \pref{A} &\to \pref{C}\\
X &\mapsto \big( c \mapsto \Hom_{\pref{A}}(i(c),X) \big)
\end{align*}
qui est un adjoint à droite, son adjoint à gauche étant le foncteur de
réalisation associé (c'est-à-dire l'extension de Kan à gauche de $i$ le long
du plongement de Yoneda).

\begin{prop}\label{propFoncteurAspheriqueZigZagMorphAsph}
Dans le diagramme \ref{diagMorphismesAspheriquesStructureÂ}, si $j$
est un foncteur $\W$-asphérique, alors $j'$ et $h'$ sont des morphismes
$\W$-asphériques de $\Cat$.
\end{prop}
\begin{proof}
On rappelle (voir \ref{propMorphismesAspheriquesLocFond}) que pour montrer
qu'un morphisme~$u:A\to B$ de $\Cat$ est asphérique, il suffit de
vérifier que pour tout objet $b$ de $B$, le préfaisceau $u^*(b)$ est
asphérique, c'est-à-dire que $i_Au^*(b)$ est une catégorie asphérique.
On fixe alors un objet $c$ de $C$. Puisque le foncteur~$i$ est pleinement
fidèle, on a un isomorphisme canonique naturel $i^*i(c) \simeq c$
dans~$\pref{C}$, et on obtient donc un isomorphisme naturel
\[
  h'^*(c) \simeq h'^*i^*i(c) \simeq i(c) \pbox{.}
\]
On distingue alors deux cas : si $c$ est un objet $a$ de $A$,
alors~$i(c)$ est un préfaisceau représentable de $\pref{A}$. Si $c$
est de la forme $j'(b)$ pour un objet $b$ de $B$, alors~$i(c)=j(b)$.
Dans les deux
cas, $h'^*(c)$ est bien un préfaisceau asphérique, ce qui prouve que
$h'$ est un morphisme asphérique de $\Cat$. 

On procède de la même manière pour montrer que le morphisme $j$ est
asphérique : on a cette fois un isomorphisme
\[
j'^*(c) \simeq j'^*i^*i(c) \simeq j^*i(c)
\]
et on distingue à nouveau deux cas : soit $c$ est un objet de $A$,
soit $c$ est de la forme $j'(b)$ pour un objet $b \in B$, ce qui
implique que $j^*ij'(b) = j^*j(b)$. Dans les deux cas, l'asphéricité
de $j$ implique que~$j'$ est un morphisme asphérique de $\Cat$.
\end{proof}

\paragr Ainsi, si $j:B\to \pref{A}$ est un foncteur asphérique, on obtient un
zig-zag de morphismes asphériques 
\[
B \to C \leftarrow A 
\]
de $\Cat$. Cela sera utile par la suite, puisqu'on verra
(\ref{foncteurAspheriqueEqQuillenAb}) que les morphismes asphériques
entre catégories étant à la fois test locales et de Whitehead induisent des équivalences de Quillen
entre les catégories de préfaisceaux abéliens.

\begin{prop}\label{propFoncteurAspheriqueSourceWhitehead}
Soient $B$ une catégorie de Whitehead et $A$ une petite catégorie. Si $j
: B \to \pref{A}$ est un foncteur asphérique
au sens du paragraphe \ref{defFoncteurAspheriqueÂ}, alors $A$ est une catégorie de
Whitehead.
\end{prop}
\begin{proof}
On gardes les notations du paragraphe
\ref{paragrFoncteurAspheriquediag}. En particulier, on a un diagramme
commutatif 
\[
\xymatrix{
& \pref{A} \ar[ld]_{j^*} \ar[d]^{i^*} \ar[rd]^{\id} \\
\pref{B} & \pref{C} \ar[l]^{j'^*} \ar[r]_{h'^*} & \pref{A} 
} 
\]
où $C$ est la sous-catégorie pleine de $\pref{A}$ contenant les
représentables et les images par $j$ des objets de $B$. Puisque les foncteurs $i^*$
et ${h'}^*$ sont des adjoints à droite, ils
préservent les groupes abéliens internes à $\pref{A}$ et $\pref{C}$ et
induisent donc des morphismes 
\[
i^*_{\ab} : \prefab{A}\to\prefab{C} \mdvirg h'^*_{\ab} : \prefab{C} \to
\prefab{A}
\]
de sorte que, dans le diagramme 
\[
\xymatrix{
\pref{A} \ar[d]_{\Whf{A}} \ar[r]^{i^*}  &
\pref{C} \ar[d]^{\Whf{C}} \ar[r]^{h'^*} &
\pref{A} \ar[d]^{\Whf{A}} \\
\prefab{A} \ar[d]_{\U_{A}} \ar[r]_{i^*_{\ab}}
\ar@{}[rd]|{(\star)}
&
\prefab{C} \ar[d]^{\U_{C}} \ar[r]_{h'^*_{\ab}} &
\prefab{A} \ar[d]^{\U_{A}} \\
\pref{A}  \ar[r]_{i^*} &
\pref{C}  \ar[r]_{h'^*} &
\pref{A} \pbox{,}
} 
\]
les carrés de droite ainsi que le carré marqué $\star$ en bas à gauche
sont commutatifs. Puisqu'on a 
\[
h'^*i^* = \id \mdvirg h'^*_{\ab}\,i^*_{\ab} = \id \mdvirg
\]
le carré extérieur est aussi commutatif. 

Puisque $j'$ est un
morphisme asphérique, la proposition \ref{aspheriqueSourceWhitehead}
implique que $C$ est une catégorie de
Whitehead (\ref{aspheriqueSourceWhitehead}). De plus,~$h'$ est aussi
un morphisme asphérique de $\Cat$. Ainsi, si $f$ est un
morphisme de~$\pref{A}$, on peut conclure par la chaîne d'équivalences
suivante :
\begin{align*}
f \in \Wab_A &\iff h'^*_{\ab}\,i^*_{\ab}(f) \in \Wab_A 
&(h'^*_{\ab}i^*_{\ab} = \id)
\\&\iff i^*_{\ab}(f) \in \Wab_C  &\text{ (asphéricité de $h'$)}
\\&\iff \U_C i^*_{\ab}(f) \in \W_C &\text{ ($C$ est de Whitehead)}
\\&\iff i^* \U_A(f) \in \W_C &\text{ (commutativité de $\star$)}
\\&\iff h'^*i^*\U_A(f) \in \W_A     &\text{ (asphéricité de $h'$)}
\\&\iff \U_A(f)\in\W_A &(h'^*i^*=\id) \pbox{.}
\end{align*}
\end{proof}

Ce résultat nous permettra de montrer que la catégorie $\Theta$ de Joyal est une
catégorie de Whitehead pseudo-test homologique (\ref{secTheta}). 

\begin{prop}
Soit $A$ une petite sous-catégorie pleine de $\Cat$. On suppose que
$A$ est dense dans $\Cat$, et que le foncteur d'inclusion $i : A \hookrightarrow
\Cat$ est un foncteur test (en particulier, $A$ est une catégorie test
par \ref{def:foncteurTest}). Alors $A$
est une catégorie de Whitehead.
\end{prop}
\begin{proof}
On note $k : \Delta \hookrightarrow \Cat$ le foncteur d'inclusion
($k^*$ est donc le foncteur nerf simplicial). On
va montrer que le foncteur
\[
j : \Delta \xhookrightarrow{k} \Cat \xrightarrow{i^*} \pref{A} 
\]
est un foncteur asphérique. Puisque $i$ est un foncteur test, alors pour
tout objet~$\Delta_n$ de $\Delta$, 
$j(\Delta_n)$ est un préfaisceau asphérique sur $A$. Par ailleurs, par
pleine fidélité de $i$, le foncteur $i^*i : A
\to \pref{A}$ est isomorphe à l'inclusion canonique $A \hookrightarrow
\pref{A}$, ce qui permet d'obtenir, pour tout objet $a$ de $A$ et pour tout
objet $\Delta_n$ de $\Delta$, une chaîne de bijections naturelles
\begin{align*}
j^*(a)_n &= \Hom_{\pref{A}}(i^*k(\Delta_n),a)  \\
&\simeq\Hom_{\pref{A}}\left(i^*k(\Delta_n),i^*i(a)\right)
&(\text{pleine fidélité de } i) \\
& \simeq\Hom_{\Cat}\left(k(\Delta_n),i(a)\right) & (\text{densité de } i)\\
&\simeq k^*i(a)_n
\end{align*}
qui montre qu'on a un isomorphisme $j^*(a) \simeq k^*i(a)$. Puisque
$i(a)$ est une catégorie asphérique et $k$ est un foncteur test, $j^*(a)$ est donc un préfaisceau
asphérique sur $\Delta$. Il reste à montrer que pour tout objet
$\Delta_n$ de $\Delta$, le préfaisceau~$j^*j(\Delta_n)$ est
asphérique. On montre pour cela que ce dernier est isomorphe au
nerf de $\Delta_n$, grâce à la chaîne de bijections naturelles
suivante, pour tout objet $\Delta_m$ de $\Delta$ :
\begin{align*}
j^*j(\Delta_n)_m &= \Hom_{\pref{A}}(j(\Delta_m),j(\Delta_n)) \\
&=\Hom_{\pref{A}}(i^*k(\Delta_m), i^*k(\Delta_n)) \\
&\simeq \Hom_\Cat(k(\Delta_m),k(\Delta_n)) & (\text{densité de } i)\\ 
&\simeq k^*k(\Delta_n)_m \pbox{.}
\end{align*}
Puisque $k$ est un foncteur test, $k^*k(\Delta_n)$ est bien un
préfaisceau asphérique sur $\Delta$. Ainsi, le foncteur $j : \Delta \to \pref{A}$ est asphérique, et la
proposition \ref{propFoncteurAspheriqueSourceWhitehead} permet de
conclure que $A$ est une catégorie de Whitehead.
\end{proof}

\begin{remark}
En utilisant la notion de segment séparant (voir
\cite[section 1.4]{maltsiniotis2005}), on peut montrer que toute
sous-catégorie pleine dense de $\Cat$, totalement asphérique et formée
de catégories asphériques est une catégorie test de Whitehead. En effet,
si $A$ est une telle catégorie, il suffit, pour reproduire la preuve
de la proposition précédente, de vérifier que le foncteur~$j : \Delta \to \pref{A}$ envoie les objets de $\Delta$ sur des
préfaisceaux asphériques. Par densité de $A$, il existe un
objet $I$ de $A$ et un épimorphisme $I \to \Delta_1$ dans $\Cat$. En
choisissant un objet $e_0$ (resp. $e_1$) de $I$ au-dessus de
$0$ (resp. $1$), on construit un segment
séparant représentable~$(I,e_0,e_1)$ dans $\pref{A}$ muni d'un
morphisme de segments~$(I,e_0,e_1) \to j(\Delta_1,0,1)$.
Puisque $j$ commute aux produits,
alors pour tout objet
$\Delta_n$ de $\Delta$, $j(\Delta_n)$ est contractile relativement au
segment $j(\Delta_1,0,1)$ et donc relativement au segment $I$. Puisque
$A$ est totalement asphérique et que $I$ est représentable, $I$ est
localement asphérique. Par~\cite[lemme 1.4.6]{maltsiniotis2005},
le préfaisceau $j(\Delta_n)$ est donc bien asphérique, et on peut
conclure comme dans la preuve de la proposition précédente.
En particulier, la sous-catégorie pleine de $\Cat$ dont les objets
sont les cubes~$\cub_n = (\Delta_1)^n$ pour~$n\geq0$ est une catégorie de Whitehead. On étudiera les catégories
cubiques plus en profondeur dans la section \ref{secCubes}.
\end{remark}

\section{Structure de catégorie de modèles abélienne}
Dans \cite{cisinskipref}, Cisinski montre que toute catégorie
test locale peut être munie d'une structure de catégorie de modèles à
engendrement cofibrant, dont les équivalences faibles sont les
équivalences test, et les cofibrations sont les monomorphismes. Dans
le cas des catégories de Whitehead, la classe des équivalences faibles
abéliennes coïncide avec la classe des morphismes qui
deviennent des équivalences test après oubli de la structure de
groupe. On montre dans cette section que dans ce cas, on peut définir
une structure de catégorie de modèles sur $ \prefab{A} $ par un
transfert de la structure de Cisinski sur $\hat{A}$.

\paragr 
Si $\M$ est une catégorie localement présentable et $I$ est une classe
de morphismes de $\M$, on note~$l(I)$ (resp. $r(I)$) la classe des
morphismes vérifiant la propriété de relèvement à gauche (resp. à
droite) par rapport à $I$. On appelle \emph{saturation} de $I$
l'ensemble~$l(r(I))$. On dit qu'une catégorie de modèles $(\M, \W_\M,
\textit{Cof}_\M,\textit{Fib}_\M)$ est à \emph{engendrement cofibrant}
s'il existe deux ensembles~$I,~J\subset\Fl \M$ tels que
$\textit{Cof}_\M$ (resp. $\textit{Cof}_\M\cap \W$) soit la saturation de $I$
(resp. de $J$). On dit alors que le couple~$(I,J)$ \emph{engendre} la
structure de catégorie de modèles.  

\begin{theorem}[Crans]\label{ThCransTransfert}
Soit $(\M, \W_\M, \textit{Cof}_\M,\textit{Fib}_\M)$ une catégorie de
modèles fermée engendrée par un couple $(I,J)$. Soient $\M'$ une catégorie localement
présentable et $G : \M \to \M'$ un foncteur
admettant un adjoint à droite~$D$. On suppose qu'on a 
\[
D(l(r(GJ)))\subset \W_\M \pbox{.} 
\]
Alors la catégorie $\M'$ admet une structure de catégorie de modèles
fermée engendrée par le couple $(GI, GJ)$ dont les équivalences faibles sont les
éléments de $D^{-1}(\W_\M)$ et les fibrations sont les éléments de
$D^{-1}(\textit{Fib}_\M)$.
\end{theorem}
\begin{proof}
On renvoie à \cite[théorème 3.3]{Crans1995Quillen}.
\end{proof}

L'outil principal pour vérifier la deuxième condition est le
corollaire suivant de l'argument du petit objet.
\begin{prop}
Soit $\M$ une catégorie localement présentable et $I$ un
ensemble de flèches de $\M$. Alors la saturation $l(r(I))$ coïncide
avec la classe des rétractes de composées transfinies d'images
directes d'éléments de $I$.
\end{prop}
\begin{proof}
Voir par exemple \cite[corollaire 10.5.22]{hirschhorn2003model}.
\end{proof}

La structure de catégorie de modèles à laquelle nous allons appliquer
ces résultats est la suivante :

\begin{theorem}[Cisinski]\label{structureCisinskiTestLocal}
Soit $A$ une catégorie test locale. On peut munir $\pref{A}$ d'une
structure de catégorie de modèles fermée à engendrement cofibrant dont les
équivalences faibles sont les éléments de $\W_A$, et les cofibrations
sont les monomorphismes.
\end{theorem}
\begin{proof}
Voir \cite[corollaire 4.2.18]{cisinskipref}.
\end{proof}

\paragr Nous allons montrer dans ce qui suit que si $A$ est une
catégorie test locale de Whitehead, on peut appliquer le théorème
\ref{ThCransTransfert} au couple de foncteurs adjoints 
\[
\xymatrix{
{\prefab{A}} \ar@<.5em>[r]^{\U}|{}="a" 
& {\pref{A}}
\ar@<.5em>[l]^{\Whf{A}}_{}="b"
\ar@{}"a";"b"|{\perp}
} \pbox{.}
\]

\begin{prop}
Soient $A$ une petite catégorie et $ i : X \to Y $ un monomorphisme dans
$\hat{A}$. Si $L : A \to \Ch(\Ab)$ est un intégrateur libre sur $A$, alors
le morphisme 
\[
L_!(\Wh{A}{i}) : \Wh{A}{X} \to \Wh{A}{Y} 
\]
est un monomorphismes à conoyau projectif.
\end{prop}
\begin{proof}
On doit montrer que pour tout entier positif $n$, le morphisme de
groupes abéliens
\[
(L_n)_!\Wh{A}{i}:(L_n)_!\Wh{A}{X} \to (L_n)_!\Wh{A}{Y}
\]
est un monomorphisme dont le conoyau est projectif. Puisque $L$ est un
intégrateur libre, on peut écrire
\[
L_n = \bigoplus_{i\in I_n}\Wh{{A}^{\op}}{a_i} \pbox{(\text{voir
\ref{defIntegrateurs}})}
\]
et, par additivité, il suffit de montrer que si $a$ est un objet de
$A$, alors le morphisme 
\[
\Wh{{A}^{\op}}{a}_!\Wh{A}{X} \to \Wh{{A}^{\op}}{a}_!\Wh{A}{Y}
\]
est un monomorphisme à conoyau projectif. Mais ce morphisme est
simplement le morphisme de groupes abéliens libres
\[
\Z^{(i_a)} : \Z^{(Xa)} \to \Z^{(Ya)} 
\]
qui est un monomorphisme dont le conoyau est le groupe abélien libre
engendré par l'ensemble $Y(a) \setminus \im i$.
\end{proof}

\begin{coro}\label{homologielibrecofibtriviale}
Soient $A$ une catégorie test locale, $j: X \to Y$ une cofibration
triviale de $\pref{A}$ pour la structure de Cisinski, et~$L : A
\to \Ch(\Ab)$ un intégrateur libre. Alors le morphisme $L_!(\Wh{A}{j})$ est une
cofibration triviale dans $\Ch(\Ab)$.
\end{coro}
\begin{proof}
On a déjà montré en \ref{homologiePrefaisceauxLibres} que le foncteur
$L_! \circ \Whf{A}$ envoie les équivalences test sur des
quasi-isomorphismes. La proposition précédente permet donc de
conclure. 
\end{proof}

\begin{thm}\label{structurecatmodeleabelienne}
Soit $A$ une catégorie de Whitehead test locale
(c'est-à-dire une catégorie test locale telle que~$\Wab_A =
U^{-1}\W_A$). Alors la catégorie~$\prefab{A} $ peut être munie d'une
structure de catégorie de modèles à engendrement cofibrant dont les
équivalences faibles sont les éléments de $ \Wab_A $, et les
fibrations sont les morphismes dont l'image par le foncteur d'oubli $U
: \prefab{A} \to \hat{A}$ sont des fibrations pour la structure de
Cisinski (\ref{structureCisinskiTestLocal}). 
\end{thm}
\begin{proof}
On va appliquer le théorème \ref{ThCransTransfert} au foncteur 
\[
\Whf{A} : \pref{A}\to\prefab{A} 
\]
en munissant $\pref{A}$ de la structure de catégorie de modèles de Cisinski.
Puisque la catégorie $\prefab{A}$ est localement présentable, il
suffit de montrer que si $J$ est un ensemble générateur des
cofibrations triviales de $\pref{A}$, alors la classe $l(r(\Z J))$ est
incluse dans $\U^{-1}\W_A$.

Soit $L : A \to \Ch$ un intégrateur libre sur $A$. Puisque $A$ est une catégorie de
Whitehead, on a
\[
\U^{-1}\W_A = L_!^{-1}\W_{qis}
\]
et il revient donc au même de montrer que la classe $L_!(l(r(\Z J)))$
est incluse dans celle des quasi-isomorphismes de complexes de
chaînes. Puisque le foncteur~$L_! : \prefab{A} \to \Ch(\Ab)$ est un
adjoint à gauche (voir la proposition \ref{propAdjIntegrateurs}), il
commute aux images directes et aux compositions transfinies (ainsi
qu'aux rétractes), et on a donc 
\[
L_!(l(r(\Z J))) = l(r(L_!\Z(J))) \pbox{.}
\]
Or on a montré précédemment que les éléments de $L_!\Z(J)$ sont des
cofibrations triviales de $\Ch(\Ab)$, et cette classe est stable par
images directes, composées transfinies et rétractes. Finalement,
$l(r(\Z J))$ est donc inclus dans~$\Wab_A$, ce qui conclut
la preuve.
\end{proof}

\begin{remark}\label{cofibGenStructureMod}
Si $I$ désigne un ensemble générateur des cofibrations de $\pref{A}$
et $J$ un ensemble générateur des cofibrations triviales, alors les
ensembles~$\Z I$ (resp.~$\Z J$) sont générateurs des cofibrations
(resp. cofibrations triviales) dans~$\prefab{A}$. De plus, puisque les
cofibrations de $\pref{A}$ sont les monomorphismes et que le foncteur
$\Whf{A}$ est un foncteur de Quillen à gauche, tous les objets libres
(ainsi que les coproduits d'objets libres) sont cofibrants.
\end{remark}

\begin{example}
On sait que $\Delta$ est une catégorie de Whitehead, et
on retrouve la structure de catégorie de modèles usuelle sur
$\prefab{\Delta}$. On obtient de plus une structure de catégorie de
modèles sur $\prefab{A}$ dès que~$A$ est une catégorie test locale
munie d'un foncteur asphérique $\Delta \to A$
(\ref{aspheriqueSourceWhitehead}) ou d'un foncteur asphérique $\Delta
\to \pref{A}$ au sens du paragraphe
\ref{propFoncteurAspheriqueSourceWhitehead}. En particulier, on verra
qu'obtient une structure de catégorie de modèles sur la catégorie des
préfaisceaux abéliens sur la catégorie $\Theta$ de Joyal
(voir \ref{secTheta}).
\end{example}

\begin{prop}\label{propIntegrateurQuillenGauche}
Soient $A$ une catégorie de Whitehead test locale et~$L$ un
intégrateur tel que $L_! \circ \Whf{A} : \pref{A}\to\Ch(\Ab)$ envoie
les monomorphismes sur des monomorphismes à conoyau projectif (c'est
en particulier le cas si $L$ est un intégrateur libre).
Alors le foncteur 
\[
L_! : \prefab{A} \to \Ch(\Ab) 
\]
est un foncteur de Quillen à gauche. En particulier, $A$ est une
catégorie pseudo-test
homologique (\ref{def:pseudoTestHomologique}) si et seulement si $L_!$
est une équivalence de Quillen.
\end{prop}
\begin{proof}
La condition sur $L$ garantit que le foncteur $L_! \circ \Whf{A}$ est
un foncteur de Quillen à gauche. Comme expliqué dans la remarque
\ref{cofibGenStructureMod}, on peut choisir des ensembles générateurs
des cofibrations et des cofibrations triviales dont tous les éléments
sont dans l'image du foncteur $\Whf{A}$, ce qui permet de conclure.
\end{proof}

On introduit alors la terminologie suivante.

\paragr\label{defIntegrateurQuillen} Si $A$ est une catégorie de
Whitehead test locale, on dit qu'un intégrateur $L$ sur $A$ est un
\ndef[intégrateur!de Quillen]{intégrateur de Quillen} si le foncteur $L_! :
\prefab{A}\to\Ch(\Ab)$ est un foncteur de Quillen à gauche. En vertu
de la proposition \ref{propIntegrateurQuillenGauche},
il suffit pour cela que le foncteur $L_! \circ \Whf{A} :
\pref{A}\to\Ch(\Ab)$ envoie les monomorphismes sur des monomorphismes
à conoyaux projectifs. En particulier, tout intégrateur libre est un
intégrateur de Quillen, ce qui implique que toute catégorie test
locale de Whitehead peut être munie d'un intégrateur de Quillen.

\medskip

On a vu dans la section \ref{secintegrateurs} que si $L : A \to
\Ch(Ab)$ est un intégrateur sur~$A$, alors le foncteur
\[
L_!: \prefab{A} \to \Ch(\Ab)
\]
calcule effectivement l'homologie de $A$ à coefficients dans les
préfaisceaux abéliens. En particulier, si $L$ et $L'$ sont deux
intégrateurs, alors ils induisent le même foncteur $\Hf{A} :
\Hotab_A\to\Hotab$. Pour leurs adjoints à droite \[
L^*, L'^* : \Ch(\Ab)\to\prefab{A}\mdvirg 
\]
sans hypothèses sur $A$, il est difficile d'affirmer quoi que ce soit, y
compris le fait qu'ils préservent les équivalences faibles. 

Dans le cas où $A$ est une catégorie de Whitehead test locale,
l'existence de la structure de catégorie de modèles assure que les
foncteurs dérivés d'un couple de foncteurs adjoints sont encore des
foncteurs adjoints.

\begin{coro}\label{coroQuillenIntegrateur}
Si $A$ est une catégorie de Whitehead test locale, alors pour tout
intégrateur de Quillen $L :A\to\Ch(\Ab)$ sur $A$, le foncteur
\[ 
L^* : \Ch(Ab) \to \prefab{A} \mdvirg C \mapsto \Hom_{\Ch(\Ab)}(L(-), C)
\]
préserve les équivalences faibles. Il induit donc un foncteur 
\[
L^* : \Hotab \to \Hotab_A \pbox{.}
\]
\end{coro}
\begin{proof}
Si $L_!$ est un foncteur de Quillen à gauche, son adjoint à droite
$L^*$ est un foncteur de Quillen à droite, et on peut conclure en
remarquant que tous les complexes de chaînes sont fibrants dans
$\Ch(\Ab)$.
\end{proof}

Autrement dit, si $f : C\to D$ est un quasi-isomorphisme de complexes
de chaînes de groupes abéliens, alors le morphisme
\[
  \H{A}{L^*(f)} : \H{A}{L^*(C)}\to \H{A}{L^*(D)}
\]
est un isomorphisme dans $\Hotab$.

\begin{prop}\label{comparaisonAdjointsDroite}
Soient $A$ une catégorie de Whitehead test locale, et~$L$, $L'$ deux
intégrateurs de Quillen sur $A$. Alors, pour tout complexe $C$, il
existe un isomorphisme naturel 
\[
L^*(C) \simeq L'^*(C) 
\]
dans la catégorie $\Hotab_A$.
\end{prop}
\begin{proof}
La condition sur $L$ et $L'$ garantit que les foncteurs~$L_!$ et
$L'_!$ sont des foncteurs de Quillen à gauche. Or, il existe un
isomorphisme naturel  
\[
L_! X \simeq L'_!X \simeq \H{A}{X} 
\]
dans $\Hotab$ pour tout préfaisceau abélien $X$ sur ${A}$. On peut
alors appliquer le théorème d'adjonction des foncteurs dérivés (voir
\cite[théorème 3, p 4.5]{quillenhomotalg}), qui implique alors que les foncteurs dérivés
à droite de $L^*$ et $L'^*$ sont naturellement isomorphes. Mais
puisque ces deux foncteurs préservent les équivalences faibles, on
obtient bien le résultat annoncé.
\end{proof}

\begin{coro}
Notons $c : \Delta_n \mapsto \dk(\Delta_n)$ l'intégrateur défini par
le complexe normalisé intervenant dans la correspondance de Dold-Kan
et $L$ l'intégrateur non normalisé. Alors pour tout complexe de
chaînes $C$, les complexes 
\begin{align*}
L_!c^*(C) = (\Hom(c\Delta_\bullet, C)) &\mdvirg L_!L^*(C) =
(\Hom(L\Delta_\bullet, C)) \\
\dk c^*(C) = (\Hom(c\Delta_\bullet, C))^{nd} &\mdvirg \dk L^*(C) =
(\Hom(L\Delta_\bullet, C))^{nd}
\end{align*}
sont quasi-isomorphes.
\end{coro}
\begin{proof}
C'est un cas particulier de la proposition
\ref{comparaisonAdjointsDroite}, puisque $c$ est un intégrateur de
Quillen (voir \ref{complexeNormaliseLibresQuillen}) et $L$ est un
intégrateur libre.
\end{proof}

Dans ce qui suit, on va montrer que si $u: A \to B$ est un morphisme
asphérique entre deux catégories test locales de Whitehead, alors le
foncteur~
\begin{align*}
  (u^*)^\ab : \prefab{B}&\to\prefab{A}\\
  X &\mapsto X \circ u 
\end{align*}
est l'adjoint à gauche d'une équivalence de Quillen.

\begin{prop}[Cisinski]\label{StructureCisinskiMorphAspheriquesEqQuillen}
Soient $A$ et $B$ deux catégories test locales, et $u : A \to B$ un
foncteur asphérique. Alors les foncteurs 
\[
u^* : \pref{B}\to \pref{A}\mdvirg u_* : \pref{A}\to \pref{B} 
\]
forment une équivalence de Quillen pour les structures de catégorie de
modèles définies en \ref{structureCisinskiTestLocal}. 
\end{prop}
\begin{proof}
Voir \cite[proposition 6.4.29]{cisinskipref}
\end{proof}

\paragr On rappelle que pour tout foncteur $u : A \to B$ entre petites
catégories, puisque $u^*$ et $u_*$ sont des adjoints à droite, ils
préservent les groupes abéliens et induisent donc des morphismes que
l'on notera ici
\[
(u^*)^\ab : \prefab{B} \to \prefab{A} \mdvirg u_*^\ab :
\prefab{A} \to
\prefab{B}
\]
de sorte que les trois carrés
\begin{align*}
\xymatrix{
\prefab{B} \ar[d]_{\U} \ar[r]^{(u^*)^\ab} \ar@{}[rd]|{(\star)} &
\prefab{A} \ar[d]^{\U} \\
\pref{B} \ar[r]_{u^*} & \pref{A} \pbox{,}
} \qquad
\xymatrix{
\prefab{B} \ar[d]_{\U} \ar@{}[rd]|{(\circledast)}  
& \prefab{A} \ar[l]_{u_*^\ab}  \ar[d]^{\U} \\
\pref{B}  & \pref{A} \ar[l]^{u_*} \pbox{,} 
} \qquad
\xymatrix{
\prefab{B} \ar[r]^{(u^*)^\ab} & \prefab{A} \\
\pref{B} \ar[u]^{\Whf{B}} \ar[r]_{u^*} & \pref{A}
\ar[u]_{\Whf{A}} 
}
\end{align*}
soient commutatifs. On les notera simplement (et on l'a déjà fait plus
tôt) $u^*$ et $u_*$
lorsqu'aucune confusion ne sera possible.

\begin{prop}\label{foncteurAspheriqueEqQuillenAb}
Soient $A$ une catégorie de Whitehead test locale, $B$ une catégorie
test locale et $u:A\to B$ un foncteur asphérique. Alors :
\begin{enumerate}
\item $B$ est une catégorie de Whitehead et, en particulier, on peut munir
$\prefab{B}$ de la structure de catégorie de modèles définie en
\ref{structurecatmodeleabelienne};
\item les foncteurs
\[
  (u^*)^\ab : \prefab{B} \to \prefab{A}\mdvirg u_*^\ab : \prefab{A} \to
  \prefab{B}
\]
forment une équivalence de Quillen.
\end{enumerate}
\end{prop}
\begin{proof}
On a déjà montré le premier point dans la
proposition~\ref{aspheriqueSourceWhitehead}. 
Pour le deuxième point, on commence par rappeler que les cofibrations
et les cofibrations triviales de $\prefab{B}$ sont engendrées par
l'image par $\Whf{B}$ des cofibrations et cofibrations triviales de
la structure de Cisinski~(\ref{structureCisinskiTestLocal}) sur~$\pref{B}$. Pour vérifier que le foncteur
$(u^*)^\ab$ est un foncteur de Quillen à gauche, il suffit donc de
vérifier qu'il envoie les éléments de l'image par $\Whf{B}$ des
cofibrations de $\pref{B}$ sur des cofibrations de $\prefab{A}$.
Or les foncteurs~$\Whf{A}$ et $u^*$ sont des foncteurs de Quillen à gauche, et on peut donc
conclure de la commutativité du carré de droite ci-dessus que le morphisme~$(u^*)^\ab
: \prefab{B}\to \prefab{A}$ est bien un foncteur de Quillen à gauche.

Soient alors $X$ un préfaisceau abélien sur ${B}$ et $Y$ un objet
fibrant de $\prefab{A}$.  On remarque alors que $\U(Y)$ est un objet
fibrant de $\pref{A}$, et on rappelle que tout objet de $\pref{B}$ est
cofibrant. En utilisant la proposition
\ref{StructureCisinskiMorphAspheriquesEqQuillen}, on peut alors
conclure grâce à la chaîne d'équivalences suivante 
\begin{align*}
& f : (u^*)^\ab X \to Y \in \Wab_A \\
&\iff \U f : \U(u^*)^\ab X \to \U Y \in \W_A & (A \text{ de
Whitehead})
\\
&\iff \U f : u^*\U X \to \U Y \in \W_A &(\text{commutativité de }\star)
\\
&\iff (\U f)^\sharp : \U X \to u_*\U Y \in \W_B
&(\ref{StructureCisinskiMorphAspheriquesEqQuillen})
\\
&\iff \U(f^\sharp) : \U X \to \U u_*^\ab Y \in \W_B
&(\text{commutativité de }\circledast)
\\
&\iff f^\sharp : X \to u_*^{\ab} Y \in \Wab_B &(B\text{ de Whitehead})
\end{align*}
que l'adjonction $(u^*)^\ab \dashv u_*^\ab$ est bien une équivalence de Quillen.
\end{proof}

\begin{remark}
Il existe des hypothèses plus faibles qui garantissent que les
foncteurs 
\[
u^* : \pref{B}\to \pref{A} \mdvirg u_* : \pref{A}\to \pref{B} 
\]
induits par $u : A \to B$ forment une adjonction de Quillen pour les
structures de Cisinski : il
suffit que le foncteur $u$ soit
$\W_{\infty}$-localement constant \cite[paragraphe 6.4.1]{cisinskipref}. Cela
signifie que pour tout morphisme $b \to b'$ de $B$, le morphisme induit 
\[
\tranche{A}{b} \to \tranche{A}{b'} 
\]
est dans $\W_\infty$. Sous
ces hypothèses, et si $A$ et $B$ sont des catégories de Whitehead, la
même preuve que ci-dessus montre que les foncteurs 
\[
(u^*)^\ab : \prefab{B}\to \prefab{A} \mdvirg u_*^\ab : \prefab{A}\to \prefab{B} 
\]
forment aussi une adjonction de Quillen. Cisinski montre dans
\emph{loc. cit.} que si~$u$ est
un morphisme $\W_{\infty}$-localement constant, alors l'adjonction
$u^*\dashv u_*$ est une équivalence de Quillen si et seulement si $u$
est asphérique.
\end{remark}

\begin{coro}\label{coro:fonctAspheriqueTestLocalWhiteheadPsTestHomEq}
Soient $A$ une catégorie de Whitehead test locale, $B$ une catégorie
test locale et $u: A \to B$ un foncteur asphérique. Alors $A$ est une
catégorie pseudo-test homologique si et seulement si il en est de même
pour~$B$.
\end{coro}
\begin{proof}
Puisque $u$ est asphérique et que $A$ est une catégorie de Whitehead,
la catégorie $B$ est une catégorie de Whitehead par la proposition
\ref{aspheriqueSourceWhitehead}. Ainsi, $\prefab{B}$ peut aussi être munie
de la structure de catégorie de modèles du théorème
\ref{structurecatmodeleabelienne},
et, grâce à la proposition \ref{foncteurAspheriqueEqQuillenAb}, les
foncteurs~$u^* : \prefab{B} \to \prefab{A}$ et $u_* : \prefab{A} \to
\prefab{B}$ forment une équivalence de Quillen. Puisque $u$ est
asphérique, il est asphérique en homologie et on peut conclure par
deux-sur-trois dans le diagramme commutatif à isomorphisme naturel
près
\[
\xymatrix{
\Hotab_B \ar[rr]^{u^*}_{\simeq} \ar[rd]_{\Hf{B}} && \Hotab_{A}
\ar[ld]^{\Hf{A}} \\
& \Hotab 
} 
\]
que $\Hf{A}$ est une équivalence de catégories si et seulement si
$\Hf{B}$ est une équivalence de catégories.
\end{proof}

On va maintenant montrer le résultat principal de cette section, à
savoir que toute catégorie test de Whitehead est une catégorie
pseudo-test homologique. On va en fait essentiellement montrer que
pour toute catégorie test $A$, il existe un zig-zag de morphismes
asphériques $A \to B \leftarrow \Delta$, où $B$ est une catégorie test
stricte.

\paragr On rappelle (\ref{defFoncteurAspherique}) qu'on
dit qu'un foncteur~$i : B \to \Cat$ est un foncteur asphérique (resp.
localement asphérique) si, pour toute petite catégorie asphérique $C$,
le préfaisceau $i^*(C)$ est asphérique (resp. localement asphérique).
De plus, l'existence d'un foncteur localement asphérique~$i : B \to
\Cat$ garantit que $B$ est une catégorie test locale.

\begin{theorem}\label{thm:TestWhiteheadPseudoTestHomologique}
Toute catégorie test de Whitehead est une catégorie pseudo-test
homologique.
\end{theorem}
\begin{proof}
On fixe une petite catégorie $A_0$ test stricte, de Whitehead, et
pseudo-test homologique (par
exemple, on peut prendre pour~$A_0$ la catégorie~$\Delta$). Notons $B$
la sous-catégorie pleine de $\Cat$ dont les objets sont les images par
le foncteur $i_A$ des objets de $A$ et les images par $i_{A_0}$ des
objets de $A_0$. On obtient alors le diagramme commutatif 
\[
\xymatrix{
& \Cat 
\\
A \ar[r]_{u} \ar[ru]^{i_A} &
B \ar@{^{(}->}[u]^{i} &
A_0 \ar[l]^{v} \ar[lu]_{i_{A_0}} 
} 
\]
dans $\Cat$. On va montrer que $u$ et $v$ sont des morphismes
asphériques de~$\Cat$ et que $B$ est une catégorie test locale de
Whitehead, afin d'appliquer le
corollaire~\ref{coro:fonctAspheriqueTestLocalWhiteheadPsTestHomEq} aux
foncteurs $u$ et $v$. 

La catégorie $A_0$ est une catégorie test faible, et le foncteur
$i_{A_0}$ est donc un
foncteur asphérique (\ref{testFaiblei_AAspherique}). Puisque $i$ est
un foncteur pleinement fidèle, la proposition
\ref{lemmeFoncteursAspheriquesTriCommutatif} implique donc que $i$ est
un foncteur asphérique, et que $v$ est un morphisme asphérique de
$\Cat$. Le même raisonnement montre que $i_A$ est un foncteur
asphérique et que $u$ est donc un morphisme asphérique de
$\Cat$. 

Puisque $v$ est un morphisme asphérique et que $A_0$ est une
catégorie de Whitehead, la proposition \ref{aspheriqueSourceWhitehead}
permet d'affirmer que $B$ est une catégorie de Whitehead. 

Montrons que $B$ est un catégorie test locale. L'asphéricité de $v$,
et le fait que $A_0$ est une catégorie totalement asphérique, implique
que $B$ est également une catégorie totalement asphérique
(\ref{totAspheriqueMorphismeAspherique}). 
Par la proposition \ref{prop:totalementAspheriqueEquivalences}, on
sait donc qu'un préfaisceau sur
$B$ est asphérique si et seulement si il est localement asphérique.
Soit alors $C$ une petite catégorie asphérique. Puisque~$i$ est un
foncteur asphérique, $i^*(C)$ est un préfaisceau asphérique, et est
donc localement asphérique. Ainsi, en vertu de la proposition
\ref{propEquivalencesFoncteursTest}, $B$ est une catégorie test locale
(elle est donc même test stricte) et $i$ est un foncteur test.

On peut alors conclure de la manière suivante : les catégories $A$,
$B$ et $A_0$ sont des catégories test locales de Whitehead, et on
dispose d'un zig-zag de morphismes asphériques 
\[
A \xrightarrow{u} B \xleftarrow{v} A_0 \pbox{.}
\]
En vertu du corollaire
\ref{coro:fonctAspheriqueTestLocalWhiteheadPsTestHomEq}, $A$ est donc une
catégorie pseudo-test homologique si et seulement si il en est de même
pour $A_0$, ce qui termine la preuve.
\end{proof}

Le corollaire suivant résume la situation.

\begin{coro}\label{coro:ResumeCatPseudoTest}
Soient $A$ une catégorie test et $A_0$ une catégorie de Whitehead. On
suppose qu'il existe un foncteur asphérique $A_0 \to \pref{A}$ au sens du
paragraphe \ref{defFoncteurAspheriqueÂ}. Alors $A$ est une catégorie
pseudo-test homologique de Whitehead. 

En particulier, si il existe un morphisme asphérique $A_0
\to A$ de $\Cat$, alors~$A$ est une catégorie pseudo-test homologique de
Whitehead.

\end{coro}
\begin{proof}
On a montré dans la proposition
\ref{propFoncteurAspheriqueSourceWhitehead} que sous cette condition,
la catégorie $A$ est une catégorie de Whitehead, ce qui permet de
conclure grâce au théorème précédent. De plus, si $A_0 \to A$ est un
morphisme asphérique de $\Cat$, alors $A_0 \to A \hookrightarrow \pref{A}$ est un foncteur
asphérique au sens du paragraphe~\ref{defFoncteurAspheriqueÂ}.
\end{proof}

On appliquera souvent ce corollaire dans le cas ou $A_0$
est la catégorie $\Delta$. 
En particulier, on pourra prouver de cette manière que les catégories
$\Theta_n$ de Joyal ainsi que la catégorie $\Theta$ sont pseudo-test
homologiques (\ref{secTheta}).

\chapter{Catégories test homologiques} \label{chapCatTestHomologiques}

\section{L'intégrateur de Bousfield-Kan}

\paragr\label{defIntegrateurBousKan} Si $A$ est une petite catégorie,
on dispose d'un intégrateur (\ref{defIntegrateurs}) canonique, qu'on appellera
\ndef[intégrateur de Bousfield-Kan]{intégrateur de Bousfield-Kan} sur $A$, défini par la composition 
\[
\Lbk{A}:A \hookrightarrow \pref{A} \xrightarrow{i_A} \Cat \xrightarrow{\nerf} \EnsSimp
\xrightarrow{\mathsf{C}} \Ch(\Ab)
\]
où on a noté $\mathsf{C}$ le foncteur associant à un ensemble
simplicial $X$ son complexe non normalisé (\ref{defMooreComplex}).

En effet, pour tout objet $a$ de $A$, la catégorie $\tranche{A}{a}$
est asphérique, ce qui implique que le complexe $\Lbk{A}(a)$ a le type
d'homologie du point. De plus, pour tout entier $n\geq 0$, on a  
\[
\Lbk{A}(a)_n = \!\!\!\!\!\!\!\!\!
\bigoplus_{a_0 \to a_1 \to \cdots \to a_n} \!\!\!\!\!\!\!\!\!
\Z^{(a_n \to a)} \mdvirg
\]
ce qui implique en particulier que $\Lbk{A}$ est un intégrateur libre
sur $A$.
On peut alors exprimer son extension de Kan à gauche de la
manière suivante : pour tout préfaisceau $X\in \prefab{A}$, on a 
\[
\Lbk{A}_!(X)_n = \bigoplus_{\Delta_n\xrightarrow{u}A}X(u(n)) \pbox{.}
\]

\begin{remark}
On peut aussi utiliser le foncteur $\dk : \prefab{\Delta} \to
\Ch(\Ab)$ pour définir une version normalisée de l'intégrateur de
Bousfield-Kan. Nous privilégions ici le choix de l'intégrateur non
normalisé, puisque ce dernier est un intégrateur libre.
\end{remark}

On retrouve le résultat de Bousfield-Kan sur le calcul des limites
inductives homotopiques introduit dans
\cite[chapitre XII, section 5.5]{bousfield1972homotopy}, dont la proposition
\ref{propExpressionHomologieDerivateurLibre} donne une nouvelle
preuve.

\begin{prop}[Bousfield-Kan]
Si $X$ est un préfaisceau abélien sur $A$, alors on a un
isomorphisme canonique naturel 
\[
\H{A}{X} \simeq \Lbk{A}_!(X) 
\]
dans $\Hotab$.
\end{prop}
\begin{proof}
C'est un cas particulier de la proposition
\ref{propExpressionHomologieDerivateurLibre}.
\end{proof}

\paragr On va détailler une autre manière de voir apparaître la
formule de Bousfield-Kan. Si $A$ est une petite
catégorie, on dispose d'un foncteur 
\[
\pref{A} \xrightarrow{i_A} \Cat \xrightarrow{\nerf} \pref{\Delta} 
\]
qu'on \cite[§92, p.327]{pursuingstacks} pourrait être tenté d'utiliser afin
d'obtenir un foncteur canonique $\prefab{A}\to\prefab{\Delta} \simeq
\Ch(\Ab)$. Hélas ce foncteur ne commute pas aux produits, et ne
préserve donc pas les groupes abéliens. En revanche, $N\circ i_A$
commute aux produits fibrés (voir \cite[lemme
1.3.5]{cisinski2004localisateur}), ce qui implique que le foncteur 
\[
\tranche{\pref{A}}{e_{\pref{A}}} \to \tranche{\Cat}{A} \to \tranche{\pref{\Delta}}{\nerf A}
\]
commute aux produits, et donc qu'il envoie les préfaisceaux
abéliens sur $A$ sur des groupes abéliens internes à
$\tranche{\EnsSimp}{\nerf A}$. Pour décrire ces derniers, on peut
utiliser l'isomorphisme bien connu
\[
\tranche{\pref{\Delta}}{\nerf A} \xrightarrow{\simeq}
\pref{\tranche{\Delta}{\nerf A}} \mdvirg
\]
défini de la manière suivante : si $(F,\varphi)$ est un ensemble
simplicial au dessus de $\nerf A$, on lui associe le préfaisceau
$\tilde{F}$ sur $\tranche{\Delta}{\nerf A}$ défini sur les
objets~$(\Delta_n, \Delta_n \xrightarrow{u} A)$ de
$\tranche{\Delta}{\nerf A}$ par le carré cartésien
\[
\xymatrix{
\tilde{F}(\Delta_n,u) \ar[d]^{} \ar[r]^{} &F_n \ar[d]^{\varphi_n} \\
\cdot \ar[r]_-{u} & \nerf A_n \pbox{.}
}
\]

Par cette isomorphisme, un groupe abélien interne à
$\tranche{\pref{\Delta}}{\nerf A}$ correspond donc à un préfaisceau en
groupes abéliens sur $\tranche{\Delta}{\nerf A}$.

On obtient alors un foncteur $u_0 : \pref{A} \to
\pref{\tranche{\Delta}{\nerf A}}$ qui commute aux produits, et qui
définit donc un foncteur 
\[
u_0 : \prefab{A} \to \prefab{\tranche{\Delta}{\nerf A}} \pbox{.}
\]

\paragr\label{oubliDelta/NAaspherique} Pour toute petite catégorie $A$, on dispose d'un foncteur 
\begin{align*}
\sigma : \tranche{\Delta}{\nerf A} &\to A \\
(\Delta_n, u) &\mapsto u(n) \pbox{,}
\end{align*}
et on sait que $\sigma$ est un morphisme
asphérique, puisque l'inclusion $\Delta \hookrightarrow \Cat$ est
un foncteur test et que $\Delta_n$ possède un objet final pour tout
entier $n\geq 0$~\cite[proposition 1.7.4]{maltsiniotis2005}. De plus, on
vérifie que le morphisme $u_0$ introduit ci-dessus coïncide avec
l'image inverse du morphisme $\sigma$, c'est-à-dire qu'on a un
isomorphisme canonique 
\[
u_0(X) \simeq \sigma^*(X)
\]
dans $\prefab{\tranche{\Delta}{\nerf A}}$pour tout préfaisceau abélien
$X$ sur $A$.

\paragr
On dispose aussi d'un foncteur d'oubli
\[
  \alpha : \tranche{\Delta}{\nerf A} \to \Delta \mdvirg (\Delta_n, u)
  \mapsto \Delta_n
  \]
et on peut définir un intégrateur sur la catégorie
$\tranche{\Delta}{\nerf A}$ en appliquant le procédé expliqué au
paragraphe \ref{notationsInduction} à l'intégrateur
libre standard~$L_\Delta$ sur~$\Delta$
(voir~\ref{ex:IntegrateurLibreDelta}) : en vertu de la proposition
\ref{inductionFibrationFibDiscrete}, puisque $\alpha$ est une
fibration à fibres discrètes, le foncteur $L_{\tranche{\Delta}{\nerf
A}}$ défini par le diagramme commutatif 
\[
\xymatrix{
\tranche{\Delta}{\nerf A} \ar[rd]_{L_{\tranche{\Delta}{\nerf A}}} \ar[r]^{\alpha} & \Delta
\ar[d]^{L_{\Delta}} \\
& \Ch(\Ab)
} 
\]
est bien un intégrateur sur $\tranche{\Delta}{\nerf A}$. Par
conséquent, si $X$ est un préfaisceau en groupes abéliens sur
$\tranche{\Delta}{\nerf A}$, on a un isomorphisme naturel
\[
\H{\tranche{\Delta}{\nerf A}}{X} \simeq \H{\Delta}{\alpha_!^\ab X} \pbox{.}
\]

On remarque de plus que l'intégrateur $L_{\tranche{\Delta}{\nerf A}}$
est un intégrateur libre, et que pour tout préfaisceau abélien $X$ sur
$\tranche{\Delta}{\nerf A}$, on a un isomorphisme naturel
\[
{L_{\tranche{\Delta}{\nerf A}}}_! (X) \simeq 
\bigoplus_{\Delta_0 \xrightarrow{u} A} X(\Delta_0, u) 
\leftarrow
\bigoplus_{\Delta_1 \xrightarrow{u} A} X(\Delta_1, u) 
\leftarrow
\cdots 
\]
dans $\Ch(\Ab)$.

\paragr On obtient finalement un foncteur 
\begin{align*}
\prefab{A}\xrightarrow{\sigma^*} \prefab{\tranche{\Delta}{\nerf A}}
\xrightarrow{\alpha_!^{\ab}} \prefab{\Delta} 
\end{align*}
et on remarque alors que le diagramme suivant est commutatif à
isomorphisme près :
\[
\xymatrix@C=3em{
\prefab{A} \ar[r]^-{\sigma^*} \ar[rd]_{\Lbk{A}_!}
& \prefab{\tranche{\Delta}{\nerf A}} \ar[r]^-{\alpha_!^\ab} 
\ar[d]|(.4){{L_{\tranche{\Delta}{\nerf A}}}_!} 
& \prefab{\Delta} \ar[ld]^{{L_\Delta}_!} \\
& \Ch(\Ab) & \pbox{.}
} 
\]
On retrouve donc de cette manière la formule de Bousfield-Kan. En
outre, on obtient le diagramme commutatif à isomorphisme près 
:
\begin{equation}\label{diag:cubetranchedelta}
\xymatrix@=2em{
  & \Hot_A
  \ar[rr]^{} \ar[dd]|(.49){\vphantom{X}}^(.3){\sigma^*}
  && \Hotab_A \ar[dd]^(.3){\sigma^*}
\\
\pref{A} \ar[rr]^(.3){} 
         \ar[dd]_(.3){\sigma^*}
         \ar[ur]
&& \prefab{A} 
         \ar[dd]^(.3){\sigma^*} 
         \ar[ur]
\\
  &\Hot_{\tranche{\Delta}{\nerf A}} 
            \ar[rr]|(.48){\phantom{X}}
            \ar[dd]|(.51){\vphantom{X}}^(.3){L\alpha_!}
  && \Hotab_{\tranche{\Delta}{\nerf A}} \ar[dd]^(.3){\alpha_!^\ab}
\\
\pref{\tranche{\Delta}{\nerf A}} 
         \ar[rr]_(.3){}
         \ar[dd]_(0.3){\alpha_!}
         \ar[ur]
&& \prefab{\tranche{\Delta}{\nerf A}} 
         \ar[dd]^(0.3){\alpha_!^\ab}
         \ar[ur]
\\ 
  &\Hot_\Delta \ar[rr]|(.48){\phantom{X}}
  && \Hotab_\Delta \ar[r]^{\simeq}_{\Hf{\Delta}} & \Hotab
\\
\pref{\Delta} \ar[rr]_(.3){} \ar[ur]
&& \prefab{\Delta} \ar[ur] 
} 
\end{equation}
où les flèches horizontales correspondent aux foncteurs
d'abélianisation et les flèches obliques correspondent aux
localisations. En effet, le cube du haut est commutatif à isomorphisme
naturel près puisque le
morphisme $\sigma$ est asphérique, et le cube du bas est commutatif en
vertu de la proposition \ref{inductionFibrationFibDiscrete}.

\section{Catégories test homologiques faibles}

Cette section, ainsi que toutes celles de ce chapitre, sont à lire principalement comme une
motivation à développer plus en profondeur la théorie des
\og catégories test homologiques\fg{}. L'espoir est évidemment d'arriver à une
caractérisation de ces dernières, dans l'esprit de la théorie des
catégories test, ou de parvenir à démontrer (ou infirmer) que toute
catégorie test est une catégorie test homologique.
Une telle caractérisation est pour l'instant hors de portée, mais il nous
a semblé intéressant de pointer certains aspects de la théorie des
catégories test qui se transposent bien au cadre abélien. 
En particulier, la notion de catégorie test homologique locale a un
sens, et toute catégorie totalement asphérique en homologie qui est
également test homologique faible est aussi test homologique locale. 

Cette \og théorie \fg{} est cependant loin d'être mature, et on essaiera
de souligner autant que possible les briques identifiées comme
manquantes.

\paragr\label{defIntegrateurTranche} On rappelle que si $L$ est un intégrateur
(\ref{defIntegrateurs}) sur une petite
catégorie~$A$, on dispose de morphismes d'adjonctions 
\[
\epsilon : {L}_!L^* \to \id_{\Ch(\Ab)} \mdvirg 
\eta : \id_{\prefab{A}} \to L^*{L}_!
\]
où ${L}_!$ est l'extension de Kan à gauche de $L$ étudiée en
\ref{secintegrateurs} et $L^*$ est défini pour tout complexe $C$ par
\[
L^*(C) : a \mapsto \Hom_{\Ch(\Ab)}(L(a), C) \pbox{.}
\]

\begin{defin}
\label{defTestHomologiqueFaible}
On dit qu'un intégrateur $L$ sur $A$ est un
\emph{intégrateur test homologique faible} (ou plus simplement,
\ndef[intégrateur!test!faible]{intégrateur test faible}) si il vérifie les conditions suivantes : 
\begin{enumerate}
\item les foncteurs ${L}_!$ et $L^*$ préservent les équivalences
faibles;
\item pour tout complexe $C$, le morphisme $\epsilon_C$ est un
quasi-isomorphisme, et pour tout préfaisceau en groupes
abéliens $X$ sur $A$, le morphisme $\eta_X$ est dans $\Wab_A$.
\end{enumerate}
 On dit qu'une petite catégorie $A$ est une \ndef[catégorie!test homologique!faible]{catégorie test
homologique faible} s'il existe un intégrateur test faible sur $A$.
\end{defin}

\paragr Autrement dit, $L$ est un intégrateur test homologique faible si le
foncteur 
\[
L^* : \Ch(\Ab) \to \prefab{A}
\]
envoie les quasi-isomorphismes sur des équivalences faibles de
$\prefab{A}$, et si le foncteur induit, encore noté
\[
L^* : \Hotab \to \Hotab_A \mdvirg
\]
est une équivalence de catégories quasi-inverse au foncteur 
\[
L_! : \Hotab_A \to \Hotab \pbox{,}
\]
les isomorphismes $L_!L^* \simeq \id_{\Hotab}$ et $L^*L_! \simeq
\id_{\Hotab_A}$ étant induits par les morphismes d'adjonction. En
particulier, toute catégorie test homologique faible est une catégorie
pseudo-test homologique (\ref{def:pseudoTestHomologique}).

Grâce à la proposition suivante, on va pouvoir montrer que les axiomes
dans cette définition sont largement redondants.

\begin{remark}
Le choix que nous avons fait dans la définition ci-dessus peut être
questionné. En fait, on pourrait \emph{a priori} développer deux
théories différentes : celle que nous exposons ici, ne privilégiant
aucun choix d'intégrateur, et une autre qui consisterait à
caractériser les petites catégories $A$ telles que l'intégrateur de
Bousfield-Kan $\Lbk{A}$ est un intégrateur test faible. 
\end{remark}

\begin{prop}[Grothendieck]
Considérons deux foncteurs adjoints~$F : \M\to\M'$ et $G:\M'\to\M$, et
\[
\epsilon : FG\to \id_{\M'} \mdvirg \eta : \id_{\M} \to GF 
\]
les morphismes d'adjonction. Soient $\W$ (resp. $\W'$) une classe de
flèches faiblement saturée (\ref{defFaibleSaturation}) de $\M$ (resp. $\M'$). Alors les conditions
suivantes sont équivalentes :
\begin{enumerate}
\item $\W\!=\!F^{-1}(\W')$ et pour tout objet $a'$ de $\M'$, le morphisme
$\epsilon_{a'}$ est dans $\W'$;
\item $\W'\!=\!G^{-1}(\W)$ et pour tout objet $a$ de $\M$, le morphisme
$\eta_a$ est dans $\W$. 
\end{enumerate}
De plus, chacune de ces conditions implique la condition 
\begin{enumerate}[resume]
\item $F(\W) \subset \W'$, $G(\W') \subset \W$ et les foncteurs 
\[
\bar{F} : \W^{-1}\M \to \W'^{-1}\M' \mdvirg \bar{G} : \W'^{-1}\M' \to
\W^{-1}\M 
\]
sont des équivalences de catégories quasi-inverses l'une de l'autre.
\end{enumerate}
\end{prop}
\begin{proof}
Voir \cite[lemme 1.3.8]{maltsiniotis2005}.
\end{proof}

On obtient donc une caractérisation un peu plus simple des catégories
test homologiques faibles, et donc une condition suffisante pour
qu'une petite catégorie $A$ soit pseudo-test homologique
(\ref{def:pseudoTestHomologique}).

\begin{coro}\label{coro:caracIntegrateurTestFaible}
Soient $A$ une petite catégorie, et~$L$ un intégrateur sur~$A$. Les
conditions suivantes sont équivalentes : \begin{enumerate}
\item $L$ est un intégrateur test faible, et $A$ est une catégorie test
homologique faible;
\item pour tout objet $C$ de $\Ch(\Ab)$, le morphisme d'adjonction
\[
\epsilon_C : {L}_!L^*(C) \to C 
\]
est un quasi-isomorphisme;
\item pour tout objet $C$ de $\Ch(\Ab)$, le morphisme $\epsilon$
induit un isomorphisme 
\[
\H{A}{L^*(C)} \simeq C 
\]
dans $\Hotab$.

\end{enumerate}
\begin{proof}
On applique la proposition précédente au couple de foncteurs adjoints
défini par l'intégrateur $L$. Puisque $L$ est un intégrateur, on a
bien $\Wab_A = {L}_!^{-1}(\W_{qis})$.
\end{proof}

\end{coro}
\begin{remark}\label{etsionavaitlabonnecondition}
Dans le cas des catégories test, on montre
(voir \cite[proposition 1.3.9]{maltsiniotis2005})
que $A$ est une catégorie test faible si et seulement si pour toute
catégorie asphérique $C$, la catégorie $i_Ai_A^*(C)$ est asphérique.
Malgré nos efforts dans cette direction, nous n'avons pas d'équivalent
de cette condition : la situation serait grandement simplifiée si l'on
parvenait à montrer qu'un intégrateur~$L$ est test homologique faible
si et seulement si pour tout complexe exact (resp. ayant l'homologie
du point), le complexe ${L}_!L^*(C)$ est exact (resp. a
l'homologie du point). La propriété dont l'équivalent nous manque est
l'observation que si~$c$ est un objet d'une petite catégorie $C$,
alors $\tranche{i_Ai_A^*(C)}{c} \simeq i_Ai_A^*(\tranche{C}{c})$.
\end{remark}

\begin{example}\label{ex:c_deltaTestFaible}
On sait que le foncteur $\dk : \prefab{\Delta}\to \Ch(\Ab)$ est une
équivalence de catégories, et que sa restriction à $\Delta$ 
\[
c_\Delta : \Delta \to \Ch(\Ab) 
\]
est un
intégrateur. Ce dernier est donc bien un intégrateur test faible.
\end{example}

\begin{prop}\label{prop:testHomFaibleImpliqueAspherique}
Si $A$ est une catégorie test homologique faible, alors~$A$ est une
catégorie asphérique en homologie (\ref{defAspheriqueEnHomologie}).
\end{prop}
\begin{proof}
Notons $\Z[0]$ le complexe concentré en degré $0$ de valeur~$\Z$.
Si $L_A$ est un intégrateur test faible sur $A$, alors on a un
quasi-isomorphisme 
\[
{L_A}_!L_A^*(\Z[0]) \to \Z[0] 
\]
et on va montrer que le terme de gauche coïncide avec l'homologie de $A$ à
coefficients dans $\Z$. On remarque que si $C$ est un complexe de
chaînes, alors on a un isomorphisme 
\[
\Hom_{\Ch(\Ab)}(C,\Z[0]) \simeq \Hom_\Ab(H_0(C),\Z) \pbox{.}
\]
Par ailleurs, on sait qu'on dispose d'un quasi-isomorphisme de $L_A$ vers le foncteur
constant de valeur $\Z$ sur $A$, et par conséquent, le foncteur 
\[
H_0(L_A) : A \to \Ab 
\]
est isomorphe au foncteur constant de valeur $\Z$ sur $A$. Le préfaisceau
\[
L_A^*(\Z[0]) : a \mapsto \Hom_{\Ch(\Ab)}(L_A(a),\Z[0])
\]
est donc isomorphe au préfaisceau constant de valeur
$\Hom(\Z,\Z)$, et donc au préfaisceau constant de valeur $\Z$ sur
$A$.
Son homologie coïncide donc avec l'homologie de $A$
à coefficients dans le préfaisceau constant de valeur $\Z$,
c'est-à-dire avec l'homologie de $A$ à coefficients dans
$\Z$, en vertu de la proposition~\ref{homologieZ=HomologieSinguliere}.
Ainsi, $A$ a nécessairement le type d'homologie du point.
\end{proof}

\begin{prop}\label{prop:intCoinduitAsphEqTestHomFaible}
Soient $u : C \to A$ un foncteur, et $L_A$ un intégrateur sur~$A$. On
suppose que le foncteur $L_C : C \to \Ch(\Ab)$ défini par le diagramme
commutatif 
\[
\xymatrix{
{C} \ar[rr]^{u} \ar[rd]_{L_C} && {A} \ar[ld]^{L_A} \\
& {\Ch(\Ab)}
}  
\]
est un intégrateur sur $C$. Si $u$ est un morphisme asphérique en
homologie~(\ref{defAspheriqueEnHomologie}), alors $L_A$ est un
intégrateur test faible si et seulement si $L_C$ est un
intégrateur test faible.
\end{prop}
\begin{proof}
On a déjà montré dans la proposition
\ref{integrateurCoinduitCommuteCounite} que sous ces conditions, on
dispose d'un diagramme commutatif, pour tout complexe de chaînes $X$,
\[
\xymatrix{
{{L_C}_!L_C^*(X)} \ar[rd]_{\epsilon^C_X} \ar@{=}[r] &
{{L_C}_!u^*L_A^*(X)} \ar[r]^{\lambda_X} &
{{L_A}_!L_A^*(X)} \ar[ld]^{\epsilon^A_X} \\
& {D}
} 
\]
où le morphisme $\lambda_X$ est un quasi-isomorphisme
(\ref{integrateurInduitMorphismeAspheriqueTF}). Ainsi, par
deux-sur-trois, le corollaire \ref{coro:caracIntegrateurTestFaible}
permet d'affirmer que $L_A$ est un intégrateur test faible si et seulement
si $L_C$ est un intégrateur test faible.
\end{proof}

\begin{remark}
On renvoie à la proposition \ref{conditionsTransfertIntegrateur} pour des conditions
suffisantes pour que le foncteur $L_C$ ci-dessus soit bien un
intégrateur. Une de ces conditions est que le morphisme $u$ 
soit une fibration à fibres discrètes.
\end{remark}

\begin{coro}\label{trancheDeltaNerfTestHomologiqueFaible}
Pour toute petite catégorie asphérique $A$, la
catégorie~$\tranche{\Delta}{NA}$ est une catégorie test homologique
faible.
\end{coro}
\begin{proof}
Puisque $\Delta$ est une catégorie totalement asphérique, en vertu de
la proposition \ref{prop:totalementAspheriqueEquivalences}, tout
préfaisceau asphérique est localement asphérique. Ainsi, si $A$ est une
catégorie asphérique, alors $NA$ est un préfaisceau localement
asphérique sur $\Delta$, ce qui signifie que que le foncteur~$\alpha :
\tranche{\Delta}{NA}\to \Delta$
intervenant dans le diagramme~\ref{diag:cubetranchedelta} est asphérique. 

De plus, le morphisme $\alpha$ est une fibration à fibres discrètes,
et on peut donc, en vertu de la remarque précédente, appliquer la proposition
\ref{prop:intCoinduitAsphEqTestHomFaible} à l'intégrateur
$L_{\tranche{\Delta}{\nerf A}}$ défini par le diagramme commutatif, 
\[
\xymatrix{
\tranche{\Delta}{NA} \ar[rr]^{\alpha}
\ar[rd]_{L_{\tranche{\Delta}{\nerf A}}} 
&& \Delta \ar[ld]^{c_\Delta} 
\\ & \Ch(\Ab) 
} 
\]
où on a noté $c_\Delta$ l'intégrateur correspondant au complexe normalisé.
Puisque~$c_\Delta$ est un intégrateur test faible, il en est donc de même
pour l'intégrateur induit (qui correspond à la version normalisée de l'intégrateur de
Bousfield-Kan défini au paragraphe \ref{defIntegrateurBousKan}), ce qui implique que
$\tranche{\Delta}{NA}$ est une catégorie test homologique faible.
\end{proof}

La situation se trouve grandement simplifiée dans le cas où $A$ est
une catégorie de Whitehead, grâce à la structure de catégorie de
modèles introduite dans le théorème \ref{structurecatmodeleabelienne}.

\begin{lemme}\label{CounitePasBesoin}
Soient $A$ une catégorie test locale de Whitehead (\ref{defWhitehead})
et~$L$ un intégrateur de Quillen (\ref{defIntegrateurQuillen})
sur~$A$. Les conditions suivantes sont équivalentes : 
\begin{enumerate}
\item $L$ est un intégrateur test faible;
\item il existe un isomorphisme naturel $L_!L^* \simeq \id_{\Hotab}$.
\end{enumerate}
\end{lemme}
\begin{proof}
Puisque $A$ est une catégorie test locale de Whitehead et que $L$ est
un intégrateur de Quillen, le foncteur $L_!$ est un
foncteur de Quillen à gauche pour la structure de catégorie de modèles
sur $\prefab{A}$ définie au
théorème~\ref{structurecatmodeleabelienne}. Puisque tous les complexes
de groupes abéliens sont fibrants, le foncteur $L^*$ préserve les équivalences
faibles, et on sait donc que les foncteurs
\[
  L_! : \Hotab_A \to \Hotab \mdvirg L^* : \Hotab \to \Hotab_A
\]
forment un couple de foncteurs adjoints. On peut alors appliquer la
proposition~\cite[1.1.1]{johnstone2022elephant} 
qui affirme que la counité d'une adjonction $\epsilon : FG \to \id$
est un isomorphisme si et seulement si il existe un isomorphisme
naturel $FG \simeq\id$.
\end{proof}

\begin{prop}\label{WhiteheadIntegrateursTousPareil}
Si $A$ est une catégorie test locale de Whitehead, alors les propositions suivantes
sont équivalentes : 
\begin{enumerate}
\item il existe un intégrateur de Quillen test faible sur $A$;
\item tout intégrateur de Quillen sur $A$ est un intégrateur test faible.
\end{enumerate}
\end{prop}
\begin{proof}
Si $L$ et $L' : A \to \Ch(\Ab)$ sont deux intégrateurs, on sait que
les foncteurs $\Hotab_A \to \Hotab$ induits entre les catégories
localisées par les foncteurs $L_!$ et $L'_!$ coïncident. Puisque $L$
et $L'$ sont des intégrateurs de Quillen, les foncteurs $L^*$ et
${L'}^*:\Ch(\Ab)\to\prefab{A}$ 
sont des foncteurs de Quillen à droite, et préservent donc les
équivalences faibles.
Par adjonction, les foncteurs $L^*$ et~${L'}^* : \Hotab \to
\Hotab_A$ coïncident. On a donc un isomorphisme~$L_!L^* \simeq
L'_!L'^*$. Ainsi, si la counité associée à $L'$ est un
quasi-isomorphisme, on obtient un isomorphisme naturel~$L_!L^* \simeq
\id_{\Hotab}$, et le lemme~\ref{CounitePasBesoin} permet de
conclure.
\end{proof}

\begin{example}
La catégorie $\Delta$ est une catégorie test locale de Whitehead (voir
\ref{ex:test} et \ref{DeltaWhitehead}) et l'intégrateur $c_\Delta$
obtenu par restriction du foncteur complexe normalisé est un
intégrateur de Quillen test faible
(voir~\ref{complexeNormaliseLibresQuillen} et
\ref{ex:c_deltaTestFaible}). Par conséquent, tous les intégrateurs de
Quillen sur~$\Delta$ sont des
intégrateurs test faibles. En particulier, l'intégrateur non
normalisé~$L_\Delta$, qui est un intégrateur libre, est donc un
intégrateur test faible. Le foncteur 
\[
L_\Delta^* : \Ch(\Ab)\to\prefab{\Delta} \mdvirg C \mapsto 
( \Delta_n \mapsto \Hom({L_\Delta}(\Delta_n),C) )
\]
envoie donc les quasi-isomorphismes sur des équivalences faibles de
$\prefab{\Delta}$, et induit par passage aux catégories localisées
un quasi-inverse au foncteur~$\Hf{\Delta}$.
\end{example}

\begin{prop}\label{prop:WhiteheadPseudoTestSsiTestFaible}
Soit $A$ une petite catégorie de Whitehead test locale. Alors les conditions
suivantes sont équivalentes : \begin{enumerate}
\item $A$ est une catégorie pseudo-test homologique; 
\item $A$ est une catégorie test homologique faible.
\end{enumerate}
\end{prop}
\begin{proof}
On a déjà vu que toute catégorie test homologique faible est une
catégorie pseudo-test homologique. Soient alors $A$ une catégorie
pseudo-test homologique, test locale, et de Whitehead, et $L$ un
intégrateur libre sur $A$. On sait que le foncteur $L_! :
\prefab{A}\to\Ch(\Ab)$ est une équivalence de Quillen pour la
structure de catégorie de modèles du théorème
\ref{structurecatmodeleabelienne} sur~$\prefab{A}$. Son adjoint à
droite $L^*$ est donc également une équivalence de Quillen, et la
counité de l'adjonction dérivée 
\[
L_! : \Hotab_A \to \Hotab \mdvirg L^* : \Hotab \to \Hotab_A 
\]
est un isomorphisme dans $\Hotab$. On peut alors appliquer le lemme
\ref{CounitePasBesoin} pour conclure que $L$ est un intégrateur test
faible.
\end{proof}

\begin{prop}\label{prop:testWhiteheadTestHomologiqueFaible}
Toute catégorie test de Whitehead est une catégorie test homologique faible.
\end{prop}
\begin{proof}
Cela découle immédiatement de la proposition précédente ainsi que du
théorème \ref{thm:TestWhiteheadPseudoTestHomologique}.
\end{proof}

Faute d'une compréhension plus profonde des catégories de Whitehead,
le résultat suivant constitue pour l'instant notre outil principal
pour exhiber des catégories test homologiques faibles.

\begin{coro}\label{coro:CatTestHomFaibledeWhiteheadFoncteurAspherique}
Soient $A_0$ une catégorie de Whitehead et $A$ une catégorie test. On
suppose qu'il existe un foncteur asphérique $A_0 \to \pref{A}$ au sens
du paragraphe \ref{defFoncteurAspheriqueÂ}. Alors $A$ est une
catégorie de Whitehead test homologique faible.

En particulier, si il existe un morphisme asphérique $A_0 \to A$ de
$\Cat$, alors~$A$ est une catégorie de Whitehead test homologique
faible.
\end{coro}
\begin{proof}
Il suffit d'appliquer le corollaire \ref{coro:ResumeCatPseudoTest}
affirmant que sous ces conditions, $A$ est une catégorie pseudo-test
homologique de Whitehead. La proposition
\ref{prop:WhiteheadPseudoTestSsiTestFaible} permet alors de conclure.
\end{proof}

\begin{remark}
Par opposition à la proposition
\ref{prop:intCoinduitAsphEqTestHomFaible}, on souligne qu'on a bien
besoin que le morphisme $u$ soit asphérique, et pas seulement
asphérique en homologie.
\end{remark}

On va maintenant s'intéresser à la stabilité par produits de la classe
des catégories test homologiques faibles. 

\medskip 

\emph{
Jusqu'au corollaire \ref{coro:produitTestHomologiqueFaible}, on fixe deux petites catégories $A$ et $B$, et deux intégrateurs $L_A$
sur $A$ et~$L_B$ sur $B$. 
}

\paragr \label{notationsProduitIntegrateurs}
On rappelle (voir \ref{integrateurProduit}) qu'on
peut alors définir un intégrateur noté $L_{A\times B}$ sur le
produit~$A\times B$ en posant
\begin{align*}
L_{A\times B} = L_{A}\boxtimes L_B : A\times B &\to \Ch(\Ab) \\
(a,b) &\mapsto L_A(a)\otimes L_B(b) \pbox{.}
\end{align*}
Notons alors 
\[
\widetilde{L_A}_! : \prefab{A\times B} \to \Hom({B}^{\op},\Ch(\Ab))
\]
le foncteur associant à un préfaisceau en groupes abéliens $X$ sur
$A\times B$
le préfaisceau en complexes de chaînes 
\[
\widetilde{{L_A}_!} X : b \mapsto {L_A}_!X(-,b) \mdvirg
\]
c'est-à-dire le foncteur défini par le diagramme commutatif 
\[
\xymatrix{
\prefab{A\times B} \ar[r]^{\simeq} \ar[rd]_-{\widetilde{L_A}_!}  & \Hom({B}^{\op}, \prefab{A})
\ar[d]^{- \circ {L_A}_!} \\
& \Hom({B}^{\op}, \Ch(\Ab)) \pbox{.}
}
\]

On rappelle également qu'on note 
\[
\underline{L_B}_! : \Hom({B}^{\op},\Ch(\Ab)) \to \Ch(\Ab) 
\]
le foncteur obtenu en appliquant terme à terme le foncteur ${L_B}_!$ à
un préfaisceau en complexes de chaînes sur $B$, puis en prenant le
complexe total du complexe obtenu, défini au paragraphe
\ref{def:LSouligne}.
\begin{lemme}
Le diagramme suivant est commutatif à isomorphisme près : 
\[
\xymatrix{
{\prefab{A\times B}} \ar[rr]^{\widetilde{L_A}_!} \ar[rd]_{{L_{A\times
B}}_!} &&
{\Hom({B}^{\op},\Ch(\Ab))} \ar[ld]^{\underline{L_B}_!} \\
& {\Ch(\Ab)} & \pbox{.}
} 
\]
\end{lemme}
\begin{proof}
Puisque tous les foncteurs en jeu commutent aux limites inductives, il suffit
de vérifier que pour tout couple $(a,b)$ d'objets de $A\times B$, on a bien 
\[
\underline{L_B}_!\widetilde{L_A}_!\Wh{A\times B}{a,b} = L_A(a)\otimes
L_B(b) \pbox{.} 
\]
Or, le préfaisceau en complexe de chaînes $\widetilde{L_A}_!\Wh{A\times
B}{a,b}$ est le préfaisceau 
\[
b' \mapsto {L_A}_!\Wh{A\times B}{a,b}(-,b') = L_A(a) \otimes
\Wh{B}{b}(b') \pbox{.}
\]
Appliquer le foncteur $\underline{L_B}_!$ revient alors à prendre le
complexe total du complexe double obtenu en appliquant terme à terme
le foncteur ${L_B}_!$. Ainsi, l'image de~$\Wh{A\times B}{a,b}$ par le
foncteur $\underline{L_B}_!\widetilde{L_A}_!$ est le complexe total du
complexe double
\[
{L_B}_!(L_A(a)_0 \otimes \Wh{B}{b}) \leftarrow
{L_B}_!(L_A(a)_1 \otimes \Wh{B}{b}) \leftarrow
\cdots
\]
et on retrouve alors bien l'expression du produit tensoriel
$L_A(a)\otimes L_B(b)$.
\end{proof}

\paragr On note à présent
\[
  \widetilde{L_A^*} : \Hom({B}^{\op},\Ch(\Ab)) \to \prefab{A\times B}
\]
l'adjoint à droite du foncteur $\widetilde{L_A}_!$, obtenu également en appliquant terme à terme le
foncteur $L_A^*$, ainsi que 
\[
\widetilde{\epsilon} : \widetilde{L_A}_!\widetilde{L_A^*} \to \id 
\]
le morphisme d'adjonction correspondant. On rappelle (voir également
au paragraphe \ref{def:LSouligne}) que le foncteur
$\underline{L_B}_!$ admet pour adjoint à droite le foncteur 
\begin{align*}
\underline{L_B^*} : \Ch(\Ab) &\to \Hom({B}^{\op},\Ch(\Ab)) \\
C &\mapsto (b \mapsto \Homi_{\Ch(\Ab)}(L_B(b),C)) \mdvirg
\end{align*}
et on note alors
\[
  \underline{\epsilon} : \underline{L_B}_!\underline{L_B^*} \to \id
\]
la counité associée à cette adjonction. Par composition des
adjonctions, le diagramme 
\[
\xymatrix{
{\prefab{A\times B}} && 
\Hom({B}^{\op},\Ch(\Ab))
\ar[ll]_{\widetilde{L_A^*}}
\\
& \Ch(\Ab)
\ar[lu]^{L_{A\times B}^*} \ar[ur]_{{\underline{L_B^*}}}
} 
\]
est commutatif à isomorphisme près. De plus, pour tout
objet $C$ de $\Ch(\Ab)$, le diagramme 
\begin{equation}\label{diag:couniteProduit}
\xymatrix{
{L_{A\times B}}_!{L_{A\times B}}^*(C) \ar@{=}[r] \ar[rd]_{\epsilon_C} 
& \underline{L_B}_!\widetilde{L_A}_!\widetilde{L_A^*}\underline{L_B^*}(C)
\ar[r]^-{\lambda_C} 
& \underline{L_B}_!\underline{L_B^*}(C)
\ar[ld]^{\underline{\epsilon}_C} 
\\ & \Ch(\Ab) & \pbox{,}
} 
\end{equation}
est aussi commutatif à isomorphisme près, où on a noté $\lambda_C$ l'image par le foncteur $\underline{L_B}_!$
du morphisme 
\[
\widetilde{\epsilon}_{\underline{L_B^*}(C)} :
\widetilde{L_A}_!\widetilde{L_A^*}\underline{L_B^*}(C)  \to
\underline{L_B^*}(C) \pbox{.}
\]

\begin{prop}
Soient $A$ une petite catégorie et $B$ une catégorie asphérique en
homologie. Si $L_A$ est un intégrateur test faible sur $A$ et $L_B$
est un intégrateur sur $B$, alors l'intégrateur $L_A\boxtimes L_B$ est
un intégrateur test faible sur $A\times B$.
\end{prop}
\begin{proof}
On va montrer que pour tout
complexe de chaînes $C$, les morphismes~$\lambda_C$ et $\underline{\epsilon}_C$
du diagramme \ref{diag:couniteProduit} sont
des quasi-isomorphismes.

D'une part, on vérifie que pour tout complexe de groupes abéliens $C$,
le préfaisceaux en complexes de chaînes de groupes abéliens 
\[
\widetilde{L_A}_!\widetilde{L_A^*}\underline{L_B}^*(C) : {B}^{\op} \to \Ch(\Ab)
\]
est le préfaisceau 
\[
b \mapsto {L_A}_!{L_A}^*\left(\underline{L_B}^*(C)(b)\right) \pbox{.}
\]
Puisque $L_A$ est un intégrateur test faible, la counité
$\epsilon_A : {L_A}_!{L_A}^* \to \id$ induit un quasi-isomorphisme
argument par argument 
\[
\widetilde{L_A}_!\widetilde{L_A}^*\underline{L_B}^*(C) \to
\underline{L_B}^*(C) \pbox{.}
\]
Puisque $L_B$ est un intégrateur sur $B$, en vertu de la proposition
\ref{colimiteHomotopiqueDegreparDegre}, le foncteur
$\underline{L_B}_!$
calcule la limite inductive homotopique des préfaisceaux en complexes
de chaînes, et envoie donc les quasi-isomorphismes argument par
argument sur des quasi-isomorphismes. Ainsi, le morphisme $\lambda_C$
est un quasi-isomorphisme,
pour tout complexe $C$.

Pour le morphisme $\underline{\epsilon}_C$, puisque $B$ est une
catégorie asphérique en homologie, alors pour tout complexe $C$, le morphisme
$\underline{\epsilon}_C$ est un quasi-isomorphisme : on renvoie pour
cela à l'annexe \ref{annexe:categoriesAspheriquesIntegrateurs}.
\end{proof}

\begin{coro}\label{coro:produitTestHomologiqueFaible}
Si $A$ est une catégorie test homologique faible et $B$ est une catégorie
asphérique en homologie, alors $A\times B$ est une catégorie test
homologique faible.
\end{coro}

\section{Catégories test homologiques}

\paragr \label{integrateurTranches}
Soit $A$ une petite catégorie. Pour tout préfaisceau $F$ sur $A$, on
dispose d'un foncteur d'oubli
\begin{align*}
\alpha : \tranche{A}{F} &\to A \\
(x, x \xrightarrow{u} F) &\mapsto x 
\end{align*}
qui nous permet, étant donné un intégrateur $L_A$ sur $A$, de définir
un foncteur~$L_{\tranche{A}{F}}$ par le diagramme commutatif 
\[
\xymatrix{
\tranche{A}{F} \ar[rd]_{L_{\tranche{A}{F}}} \ar[r]^{\alpha} & A \ar[d]^{L_A} \\
& \Ch(\Ab) \pbox{.}
} 
\]

\begin{prop}\label{prop:integrateurTranches}
Pour tout intégrateur $L_A$ sur $A$, le foncteur $L_{\tranche{A}{F}}$
défini ci-dessus est un intégrateur sur $\tranche{A}{F}$. 
\end{prop}
\begin{proof}
Voir le corollaire \ref{coro:inductionTranche}. 
\end{proof}

\begin{example}\label{exemple:homologieTrancheDelta}
On peut ainsi calculer l'homologie d'un préfaisceau abélien sur la tranche
$\tranche{\Delta}{\Delta_n}$ en utilisant l'intégrateur libre $L_\Delta$
associant à tout objet de $\Delta$ son complexe non normalisé. On
obtient alors, pour tout préfaisceau abélien
$X$ sur $\tranche{\Delta}{\Delta_n}$, un isomorphisme
\[
\H{\tranche{\Delta}{\Delta_n}}{X} \simeq \cdots \leftarrow
\bigoplus_{\mathclap{\Delta_i \xrightarrow{u}{\Delta_n}}} X(\Delta_i,u)
\leftarrow 
\bigoplus_{\mathclap{\Delta_{i+1} \xrightarrow{u} \Delta_n}} X(\Delta_{i+1},u)
\leftarrow \cdots
\]
dans $\Hotab$, et on va voir en \ref{DeltaTestHomologiqueStricte} que
cet intégrateur est bien un intégrateur test faible.
\end{example}

\begin{defin}\label{defIntegrateurTestLocal}
On dit qu'un intégrateur $L_A$ sur $A$ est un
\ndef[intégrateur!test!local]{intégrateur test local} si pour tout
objet $A$ de $a$, l'intégrateur $L_{\tranche{A}{a}}$ induit par $L$
est un intégrateur test faible. 
On dit que $L_A$ est un \ndef[intégrateur!test]{intégrateur test} si
c'est un intégrateur test faible et test local. 

On dit qu'une catégorie $A$ est une \ndef[catégorie!test
homologique!locale]{catégorie test homologique
locale} (resp. \ndef[catégorie!test homologique]{test homologique}) si
il existe un intégrateur test local (resp. test) sur 
$A$.
\end{defin}

\paragr On rappelle (voir au paragraphe
\ref{def:PrefaisceauAspherique}) qu'un préfaisceau $F$ sur $A$ est un
préfaisceau localement $\W_\infty^\ab$-asphérique si le morphisme 
\[
\alpha : \tranche{A}{F} \to A 
\]
est un morphisme asphérique en homologie.

\begin{prop}\label{prop:tranchesLocAspheriquesTestHomologique}
Soient $A$ une petite catégorie, $L_A$ un intégrateur test local
sur $A$ et $F$ un préfaisceau sur $A$. Alors l'intégrateur
$L_{\tranche{A}{F}}$ est un intégrateur test local. De plus, si $F$
est un préfaisceau localement $\W_\infty^\ab$-asphérique,
alors~$L_{\tranche{A}{F}}$ est un intégrateur test si et seulement si
$L_A$ est un intégrateur test.
\end{prop}
\begin{proof}
La première assertion vient de l'observation qu'on a 
\[
\tranche{(\tranche{A}{F})}{(a,x)} = \tranche{A}{a}
\]
pour tout objet $(a,a\xrightarrow{x} F)$ de $\tranche{A}{F}$. Pour la
deuxième assertion, si $F$ est localement asphérique, alors le
morphisme $\alpha : \tranche{A}{F}\to A$ est un morphisme asphérique,
et on peut appliquer la proposition
\ref{prop:intCoinduitAsphEqTestHomFaible} pour conclure.
\end{proof}

\begin{remark}
Un élément de frustration se doit d'être explicité : on peut montrer
qu'une catégorie $A$ est une catégorie test \emph{si et seulement si}
$A$ est une catégorie test locale asphérique. Nous n'avons pas
d'analogue de cette propriété pour l'instant. La preuve (voir
\cite[remarque 1.5.4]{maltsiniotis2005}) utilise la caractérisation
des catégories test faibles dont on a déjà dit qu'elle nous manquait
dans la remarque \ref{etsionavaitlabonnecondition}. 

Si l'analogue abélien de cette caractérisation tient, c'est-à-dire si
toute catégorie test homologique locale qui est aussi asphérique en homologie est
une catégorie test homologique, on pourra en particulier montrer
que si $A$ est une catégorie test homologique locale, alors
$\tranche{A}{F}$ est une catégorie test homologique pour tout
préfaisceau $\W_\infty^\ab$-asphérique $F$ (c'est l'analogue de la
proposition \ref{prop:testLocalTranche} pour les catégories test).
\end{remark}

\begin{prop}
Soit $A$ une catégorie test homologique locale (resp. test
homologique) et $B$ une catégorie arbitraire (resp. asphérique en
homologie). Alors la catégorie $A\times B$ est une catégorie test
homologique locale (resp. test homologique).
\end{prop}
\begin{proof}
Il suffit de remarquer qu'on a 
\[
\tranche{(A\times B)}{(a,b)}= \tranche{A}{a}\times \tranche{B}{b} 
\]
pour tout couple $(a,b)$ d'objets de $A\times B$, avant d'appliquer le
corollaire~\ref{coro:produitTestHomologiqueFaible}, et de remarquer
que
l'intégrateur induit par $L_A\boxtimes L_B$ sur la
catégorie $\tranche{(A\times B)}{(a,b)}$ est bien l'intégrateur
$L_{\tranche{A}{a}}\boxtimes L_{\tranche{B}{b}}$ \pbox{.}
\end{proof}

\begin{remark}\label{blablaWhiteheadLocal}
Pour utiliser la proposition
\ref{prop:testWhiteheadTestHomologiqueFaible} affirmant que toute
catégorie test de Whitehead est une catégorie test homologique, il
serait évidemment intéressant de caractériser les petites catégories
\emph{localement de Whitehead}. Nous ne disposons pas d'une telle
caractérisation pour l'instant. En revanche, on peut montrer qu'il n'est pas
raisonnable de demander que pour tout préfaisceau $F$ sur $A$, la
catégorie $\tranche{A}{F}$ soit une catégorie de Whitehead, puisque
même la catégorie $\Delta$ ne vérifie pas cette condition.
En effet, on a vu au paragraphe~\ref{oubliDelta/NAaspherique} que pour
toute petite catégorie $A$, le foncteur 
\[
\sigma : \tranche{\Delta}{\nerf A} \to A \mdvirg (\Delta_n, \Delta_n
\xrightarrow{u} A) \mapsto u(n)
\]
est un foncteur asphérique. La proposition
\ref{aspheriqueSourceWhitehead} impliquerait alors que
toute catégorie est une catégorie de Whitehead. Or, il existe des
petites catégories qui ne sont pas de Whitehead (on va monter par exemple que c'est le cas de la catégorie $\Grefl$ des globes réflexifs
dans la proposition \ref{globespaswhitehead}).
La stratégie consisterait donc à s'intéresser seulement aux catégories $A$
telles que pour tout objet $a$ de $A$, la catégorie $\tranche{A}{a}$
est de Whitehead. 
\end{remark}

\begin{prop}
Soit $A$ une catégorie test. Si~$\tranche{\Delta}{\nerf A}$ est une
catégorie de Whitehead, alors $A$ est une catégorie test homologique
faible de Whitehead.
\end{prop}
\begin{proof}
Comme expliqué au paragraphe \ref{oubliDelta/NAaspherique}, le
morphisme $\sigma : \tranche{\Delta}{\nerf A}\to
A$ est alors asphérique, et $A$ est donc une catégorie test de
Whitehead. La proposition
\ref{prop:testWhiteheadTestHomologiqueFaible} permet alors de
conclure.
\end{proof}

\section{Catégories test homologiques strictes}
\label{secCatTestHomStrictes}

La notion de catégorie totalement asphérique que l'on a introduit au
paragraphe \ref{def:totAspherique} peut être définie pour tout
localisateur fondamental. 

\begin{prop}[Grothendieck]\label{TotAspheriqueEqDef}
Soient $A$ une petite catégorie et~$\W$ un localisateur
fondamental (\ref{defLocFondamental}). Les conditions suivantes sont
équivalentes : 
\begin{enumerate}
\item pour tous objets $a$ et $b$ de $A$, le préfaisceau $a\times b$
est $\W$-asphérique;
\item le produit de deux préfaisceaux $\W$-asphériques est $\W$-asphérique;
\item pour tous préfaisceaux $F$ et $G$ sur $A$, le foncteur
canonique 
\[
\tranche{A}{(F\times G)} \to \tranche{A}{F}\times \tranche{A}{G} 
\]
est dans $\W$;
\item pour tous préfaisceaux $F$ et $G$ sur $A$, le foncteur
canonique 
\[
\tranche{A}{(F\times G)} \to \tranche{A}{F}\times \tranche{A}{G} 
\]
est $\W$-asphérique;
\item tout préfaisceau représentable de $\pref{A}$ est localement
$\W$-asphérique;
\item tout préfaisceau asphérique de $\pref{A}$ est localement
$\W$-asphérique;
\item le foncteur diagonal $A \to A\times A$ est $\W$-asphérique.
\end{enumerate}
De plus, si $A$ est non vide, alors toutes ces conditions impliquent
que $A$ est une catégorie $\W$-asphérique.
\end{prop}
\begin{proof}
Voir \cite[1.6.1]{maltsiniotis2005}.
\end{proof}

On peut en particulier appliquer cette proposition au localisateur
fondamental $\W_{\infty}^{\ab}$ défini au paragraphe
\ref{localisateurFondamentalAbelien}. 

\begin{defin}
Une catégorie $A$ est dite \ndef[catégorie!totalement asphérique!en homologie]{totalement asphérique en homologie} si elle est
asphérique en homologie et si elle vérifie les conditions
équivalentes de la proposition \ref{TotAspheriqueEqDef}

Les catégories totalement asphériques coïncident donc exactement avec
les catégories non vides vérifiant les conditions équivalentes de la
proposition \ref{TotAspheriqueEqDef}.
\end{defin}

\begin{remark}
On souligne que si $A$ est totalement asphérique, alors elle est
totalement asphérique en homologie. 
\end{remark}

\begin{defin}
On appelle \ndef[catégorie!test homologique!stricte]{catégorie test
homologique stricte} une catégorie test homologique totalement
asphérique en homologie.
\end{defin}

\begin{prop}\label{foncteurAspheriqueSourceTHS}
Si $u : A \to B$ est un foncteur asphérique en homologie, et si $A$
est totalement asphérique en homologie, alors il en est de même pour
$B$.
\end{prop}
\begin{proof}
Ce résultat reste également vrai pour tout localisateur fondamental
faible, on trouvera une preuve dans \cite[1.6.5]{maltsiniotis2005}.
\end{proof}

\begin{prop}\label{propTestHomologiqueTotAspherique}
Soit $A$ une catégorie totalement asphérique en homologie. Alors les conditions
suivantes sont équivalentes : 
\begin{enumerate}
  \item $A$ est une catégorie test homologique faible;
  \item $A$ est une catégorie test homologique locale;
  \item $A$ est une catégorie test homologique;
  \item $A$ est une catégorie test homologique stricte,
\end{enumerate}
Si, de plus, $A$ est une catégorie de Whitehead test locale, alors les conditions
précédentes sont encore équivalentes à la condition 
\begin{enumerate}[resume]
\item $A$ est une catégorie pseudo-test homologique.
\end{enumerate}
\end{prop}
\begin{proof}
L'équivalence $(5) \iff (1)$ a déjà été montrée dans la proposition
\ref{prop:WhiteheadPseudoTestSsiTestFaible}, et il suffit donc de
montrer l'équivalence $(1)\iff(2)$. Si~$A$ est totalement asphérique
en homologie, alors tout préfaisceau représentable est localement
$\W_\infty^\ab$-asphérique, ce qui signifie que pour tout objet $a$ de
$A$, le foncteur~$\alpha: \tranche{A}{a}\to A$ est asphérique
en homologie. En fixant un intégrateur~$L_A$ sur $A$, la proposition
\ref{prop:intCoinduitAsphEqTestHomFaible} permet de conclure que
$L_A$ est un
intégrateur test faible si et seulement si l'intégrateur
$L_{\tranche{A}{a}}$ défini au paragraphe~\ref{defIntegrateurTranche}
est un intégrateur test faible. En d'autres termes, $L_A$ est un
intégrateur test faible si et seulement si c'est un intégrateur test
local. 
\end{proof}

\begin{prop}
Toute catégorie test stricte de Whitehead est une catégorie test
homologique stricte.
\end{prop}
\begin{proof}
C'est un corollaire immédiat de la proposition
\ref{prop:testWhiteheadTestHomologiqueFaible} et du fait que toute
catégorie totalement asphérique est totalement asphérique en
homologie.
\end{proof}

\begin{coro}\label{DeltaTestHomologiqueStricte}
La catégorie $\Delta$ est une catégorie test homologique stricte. 
\end{coro}
\begin{proof}
Puisque $\Delta$ est une catégorie test stricte (voir
\cite[proposition 1.6.5]{maltsiniotis2005}) de Whitehead (voir la
proposition \ref{DeltaWhitehead}), la proposition précédente permet de
conclure.
\end{proof}

\begin{prop}
Soit $A$ une catégorie test homologique stricte. Alors pour tout
préfaisceau $\W_\infty^\ab$-asphérique $F$ sur $A$, la catégorie
$\tranche{A}{F}$ est une catégorie test homologique.
\end{prop}
\begin{proof}
Puisque $A$ est totalement asphérique en homologie, tout préfaisceau
$\W_\infty^\ab$-asphérique est localement $\W_\infty^\ab$-asphérique,
et on peut appliquer la proposition
\ref{prop:tranchesLocAspheriquesTestHomologique}.
\end{proof}

\begin{remark}
La proposition précédente permet de montrer qu'il
existe des catégories test homologiques qui ne sont pas des catégories
de Whitehead : comme expliqué dans l'exemple
\ref{blablaWhiteheadLocal}, la catégorie $\tranche{\Delta}{\nerf
\Grefl}$ n'est pas une catégorie de Whitehead. Mais c'est une
catégorie test homologique, puisqu'on verra que $\Grefl$ est
asphérique, ce qui implique que $\nerf \Grefl$ est un préfaisceau
localement asphérique sur $\Delta$. On note de plus que
$\tranche{\Delta}{\nerf \Grefl}$ est également une catégorie test.
\end{remark}

Finalement, les résultats principaux que nous utiliserons afin
d'exhiber des correspondances de Dold-Kan homotopiques sont les
suivants.

\begin{prop}\label{propMorphismeAsphSourceTHS}
Soient $A_0$ une catégorie de Whitehead et $A$ une catégorie test. On
suppose que $A$ est totalement asphérique en homologie, et qu'il
existe un foncteur asphérique $A_0 \to \pref{A}$ au sens du paragraphe
\ref{defFoncteurAspheriqueÂ}. Alors~$A$ est une catégorie test
homologique stricte de Whitehead.
\end{prop}
\begin{proof}
Sous ces conditions, le corollaire
\ref{coro:CatTestHomFaibledeWhiteheadFoncteurAspherique} implique que
$A$ est une catégorie test homologique faible de Whitehead. Puisque
$A$ est une catégorie totalement asphérique en homologie, c'est donc
une test homologique stricte en vertu de la proposition
\ref{propTestHomologiqueTotAspherique}.
\end{proof}

\begin{coro}
Soient $A_0$ une catégorie de Whitehead totalement
asphérique et $A$ une catégorie test locale. On suppose
qu'il existe un morphisme asphérique $A_0 \to A$ de $\Cat$. Alors $A$
est une catégorie test homologique stricte de Whitehead.
\end{coro}
\begin{proof}
Sous ces conditions, $A$ est une catégorie totalement asphérique en
vertu de la proposition \ref{totAspheriqueMorphismeAspherique}, et est
donc une catégorie test stricte. La proposition précédente permet donc
de conclure.
\end{proof}

En particulier, on appliquera souvent les deux résultats précédents aux
foncteurs asphériques de source $\Delta$.

\chapter{Théorèmes de Dold-Kan homotopiques}\label{chapExemples}

\section{Orientations sur les petites catégories}
\label{secOrientation}

\paragr Si $A$ est une catégorie, on introduit une relation
d'équivalence sur les monomorphismes de but fixé de la manière
suivante. Si $a$ est un objet de~$A$ et $\varphi' : a'\hookrightarrow
a$ et $\varphi'' : a'' \hookrightarrow a$ sont deux monomorphismes, on pose
$\varphi' \sim \varphi''$ s'il existe un isomorphisme $f : a \to a'$
faisant commuter le triangle
\[
\xymatrix{
a' \ar@{^(->}[rd]_{\varphi'} \ar[rr]^{f}_{\simeq} && a''
\ar@{_(->}[ld]^{\varphi''} \\ & a & \pbox{.} 
}  
\]

On appelle \ndef{sous-objet} de $a$ une classe d'équivalence de
monomorphismes de but~$a$ pour cette relation d'équivalence. On dit
qu'un sous-objet de $a$ est trivial s'il est équivalent à l'identité
de $a$, c'est-à-dire s'il est dans la classe d'équivalence des
automorphismes de $a$.

On effectuera l'abus de notation consistant à noter de la même
manière les monomorphismes de but $a$ et les sous-objets de $a$. En
particulier, un monomorphisme est trivial si c'est un isomorphisme.
 
\paragr \label{defDimensionSSObjets}
Si $C = a_0 \hookrightarrow a_1\hookrightarrow \cdots 
\hookrightarrow a_{n-1} \hookrightarrow a_n$ est une chaîne d'inclusions
de sous-objets de \emph{source} $a_0$ et de \emph{but} $a_n$, on note
$\ell(C)$ l'entier $n$ et on dit que~$\ell(C)$ est la \emph{longueur} de
la chaîne $C$. 

\paragr Étant donnés deux objets $a'$ et $a$ de $A$, on dit qu'une chaîne
d'inclusions de sous-objets de source $a'$ et de but $a$ est \emph{maximale}
si on ne peut pas la raffiner en une chaîne de sous-objets de source
$a'$ et de but $a$ de longueur strictement supérieure,
c'est-à-dire si on ne peut décomposer aucun des monomorphismes qui la
composent comme composées de deux monomorphismes non triviaux.

On dit
qu'une chaîne d'inclusion de sous-objets de but $a$ est \emph{totalement
maximale} si on ne peut pas la raffiner en une chaîne de longueur
strictement supérieure, et si la source de cette chaîne n'est le but d'aucun
monomorphisme non trivial.

\paragr\label{def:catenaire} On dit que $A$ est une
\ndef[catégorie!caténaire]{catégorie caténaire} si elle
vérifie les conditions suivantes :
\begin{enumerate}
\item tout objet $a$ est le but d'une chaîne totalement maximale;
\item toutes les chaînes totalement maximales de but $a$ ont la même
longueur, que l'on note $\dim(a)$;
\item pour tout objet $a$ de $A$ et tout chaîne $C$ d'inclusions de
sous-objets de but $a$, on a $\ell(C) \leq \dim(a)$;
\item pour tous objets $a$ et $a'$ de $A$, l'ensemble $\Hom_A(a,a')$
est fini.
\end{enumerate}

\begin{lemme}\label{lemme:catenaire}
Soient $A$ une catégorie caténaire et $a' \hookrightarrow a$ un
sous-objet de $a$ dans $A$. Alors il existe une chaîne maximale $C$ de
$a'$ vers $a$, et, pour une telle chaîne, on a
\[
  \dim a = \dim a' + \ell(C) \pbox{.}
\]
\end{lemme}
\begin{proof}
On fixe une chaîne totalement maximale
\[
  C'=a'_0 \hookrightarrow \cdots \hookrightarrow a'_n = a'
  \]
de but $a'$. On a donc $\dim a' = \ell(C')$. La condition $3$
ci-dessus implique qu'il existe une chaîne maximale de $a'$ vers $a$.
Soit alors $C$ une telle chaîne. La chaîne obtenue en concaténant $C'$
et~$C$ 
  \[
a'_0 \hookrightarrow \cdots \hookrightarrow a' \hookrightarrow
\cdots \hookrightarrow a
  \]
est totalement maximale, et est donc de longueur $\dim(a)$.
\end{proof}

\paragr En particulier, si $A$ est une catégorie caténaire et si $a'
\hookrightarrow a$ est un sous-objet de $a$ dans $A$, toutes les
chaînes maximales de $a'$ vers $a$ ont la même longueur, que l'on note
alors~$\codim(a',a)$. Le lemme précédent implique alors que pour tout
sous-objet~$a' \hookrightarrow a$ dans $A$,
on a 
\[
\codim(a',a) = \dim a - \dim a' \pbox{.}
\]
Si $\varphi : a' \hookrightarrow a$ est un sous-objet de $a$ dans $A$,
on notera aussi
\[
  \codim(\varphi) = \codim(a',a) \pbox{.}
\]
En particulier, $\varphi$ est un monomorphisme de codimension $1$ si
et seulement si on ne peut pas le décomposer comme une composition de
deux monomorphismes non triviaux.

\begin{example}\label{ex:orientationDelta}
On peut vérifier que $\Delta$ est une catégorie caténaire, et que pour
tout entier $n\geq 0$, l'objet $\Delta_n$ est de dimension $n$. Les
monomorphismes de codimension $1$ sont alors les cofaces (voir
\ref{defDelta})
\[
  \delta_i : \Delta_{n} \to \Delta_{n+1}
\]
pour $n\geq 0$ et $0 \leq i \leq n+1$.
\end{example}

\bigskip 
\emph{Jusqu'à la fin de cette section, on fixe une petite catégorie
caténaire $A$}.
\medskip

\paragr Pour tout entier $n\geq0$, on définit un préfaisceau en groupes
abélien sur~${A}^{\op}$ en posant 
\[
L_n = \bigoplus_{\dim a = n}{\Wh{{A}^{\op}}{a}} \pbox{.}
\]
Concrètement, pour tout objet $x$ de $A$ et pour tout entier $n\geq0$,
$L_n(x)$ est le groupe abélien libre engendré par les morphismes de but
$x$ et de source un objet de dimension~$n$
\[
L_n(x) = \bigoplus_{\dim a = n} \Z^{(a \to x)} \pbox{.}
\]

\paragr Le lemme \ref{lemme:catenaire} implique que si $\varphi : a'
\hookrightarrow a$ est un monomorphisme de codimension $1$ et $a$
est de dimension $n >0$, alors $a'$ est de dimension $n-1$. La
précomposition par $\varphi$ induit donc un morphisme dans
$\prefab{{A}^{\op}}$
\[
\varphi^* : L_n \to L_{n-1}
\]
défini pour tout objet $x$ de $A$ et pour tout objet $a$ de dimension
$n$ de $A$ par 
\[
\varphi^*_x \langle a \xrightarrow{u}
x \rangle = \langle a' \xrightarrow{\varphi} a \xhookrightarrow{u} x \rangle \pbox{.}
\]

\paragr On appelle \emph{préorientation} sur $A$ la
donnée, notée $\mathcal{O}$, d'un signe \[
\sg(\varphi) \in \lbrace \pm 1 \rbrace 
\]
pour tout sous-objet $\varphi : a' \hookrightarrow a$ de codimension
$1$.
Étant donnée une telle préorientation, on définit un morphisme 
\[
d_\bullet : L_\bullet \to L_{\bullet-1} 
\]
de préfaisceaux abéliens sur ${A}^{\op}$ en posant, pour tout entier
$n>0$,
\[
d_n = \,\sum_{\mathclap{\codim\varphi=1}}\,\sg(\varphi)\varphi^* : L_n \to L_{n-1} \pbox{.}
\]

\paragr\label{def:orientation} On dit qu'une préorientation $\mathcal{O}$ sur $A$ est une
\ndef{orientation} si pour tout objet $a$ de $A$, le morphisme $d$
défini ci-dessus est une
différentielle, c'est-à-dire si on a~$d\circ d=0$. Dans ce cas, on
note 
\[
L = L_0 \xleftarrow{d_1} L_1 \xleftarrow{d_2} \cdots 
\]
le complexe dans $\Ch(\Addinf({A}^{\op}))$ obtenu (voir
\ref{defAddinf}). Concrètement, si $x$ est un objet de $A$, on obtient
le complexe de chaînes de groupes abéliens 
\[
L(x) = \cdots \xleftarrow{d_x}  \bigoplus_{\dim a = n-1} \Z^{(a\to x)}
\xleftarrow{d_x} \bigoplus_{\dim a = n}\Z^{(a\to x)} \xleftarrow{d_x}
\cdots
\]
où pour tout objet $a$ de dimension $n>0$ et tout morphisme $u : a\to x$, on a 
\[
d_x \langle u \rangle = \sum_{\codim \varphi = 1} \sg(\varphi) \langle u \circ
\varphi \rangle \pbox{.}
\]

\paragr\label{def:orientationAspherique} On dit qu'une orientation sur $A$ est
\ndef[orientation!asphérique]{asphérique}
si le complexe $L$ a l'homologie du point. En d'autres termes,
$\mathcal{O}$ est une orientation asphérique si le complexe
$L_\bullet$ est un intégrateur libre (voir \ref{defIntegrateurs}) sur $A$.

\begin{prop}
Si $\mathcal{O}$ est une orientation asphérique sur $A$, alors pour tout
préfaisceau abélien $X$ sur $A$, on a un isomorphisme naturel 
\[
\H{A}{X} \simeq \cdots \leftarrow \bigoplus_{\dim a = n-1} Xa 
\leftarrow \bigoplus_{\dim a = n}Xa \leftarrow \cdots
\]
dans $\Hotab$, où la différentielle est donnée par $\Addinf(X)(d)$.
\end{prop}
\begin{proof}
C'est un cas particulier de la proposition
\ref{propExpressionHomologieDerivateurLibre}.
\end{proof}

\begin{example}
En reprenant l'exemple de la catégorie $\Delta$,
l'orientation~$\mathcal{O}$ donnée pour toute coface $\delta_i :
\Delta_n \to \Delta_{n+1}$ avec $0 \leq i \leq n+1$ par 
\[
\sg (\delta_i) = (-1)^i
\]
est bien une orientation asphérique. En effet, on vérifie que le
complexe obtenu est l'intégrateur libre standard $L_\Delta$ sur
$\Delta$ (voir \ref{ex:IntegrateurLibreDelta}).
\end{example}

\section{Ensembles multisimpliciaux}
Le cas des ensembles multisimpliciaux est entièrement contenu dans le
théorème de Dold-Kan original. On va ici simplement regrouper les
énoncés, et les réinterpréter légèrement avec le vocabulaire que nous
avons introduit dans cette thèse.

\paragr On peut prouver, à partir du théorème de Dold-Kan
(\ref{doldkan}), que l'on a une équivalence de catégories 
\[
\prefab{\Delta\times\Delta} \simeq \Ch(\Ch(\Ab)) 
\]
entre la catégorie des groupes abéliens bisimpliciaux et celle des
complexes doubles de groupes abéliens du premier quadrant,
c'est-à-dire les complexes de chaînes doubles $C_{\bullet,\bullet}$
tels que $C_{i,j} = 0$ dès que $i<0$ ou $j<0$. 

En effet, cette correspondance est obtenue en appliquant terme à terme le
foncteur 
\[
\dk : \prefab{\Delta}\to\Ch(\Ab) \mdvirg
\]
qui, on le rappelle (\ref{doldkan}), est une équivalence de catégories. Précisément,
on obtient la chaîne d'équivalences de catégories 
\begin{align*}
\prefab{\Delta\times\Delta} &\simeq \Hom({\Delta}^{\op},\prefab{\Delta})
\\
&\simeq\Hom({\Delta}^{\op},\Ch(\Ab)) 
\\
&\simeq \Ch(\prefab{\Delta}) 
\\
&\simeq \Ch(\Ch(\Ab)) 
\end{align*}
où, dans l'avant dernière ligne, on a simplement utilisé
l'isomorphisme de catégories $\Ch(\prefab{A}) \simeq
\Hom({A}^{\op},\Ch(\Ab))$ appliqué à la catégorie $A=\Delta$.
On alors peut itérer cette construction pour obtenir une équivalence de
catégories
\[
  \prefab{\Delta^n} \simeq \underbrace{\Ch(\Ch(\dots(\Ch}_{n
          \text{ fois}}(\Ab)))) \pbox{.}
  \]

\paragr On peut introduire un intégrateur sur la catégorie
$\Delta\times\Delta$ de la manière suivante. On fixe un
intégrateur $L$ sur $\Delta$, par exemple l'intégrateur normalisé $c : \Delta
\to \Ch(\Ab)$ (\ref{integrateurNormalise}) ou l'intégrateur non
normalisé $L_\Delta$ (\ref{ex:IntegrateurLibreDelta}).
On peut alors, comme expliqué au paragraphe \ref{integrateurProduit},
définir un intégrateur sur la catégorie~$\Delta\times\Delta$ en posant 
\begin{align*}
L_{\Delta\times\Delta} :=L \boxtimes L : \Delta\times\Delta &\to \Ch(\Ab) \\
\Delta_{n,m} &\mapsto
L(\Delta_n)\otimes_{\Ch(\Ab)} L(\Delta_m)
\end{align*}
où on note
$\Delta_{n,m}$ le couple $(\Delta_n, \Delta_m)$ pour $n,m\geq 0$.

\begin{prop}\label{homologieBiSimpliciaux}
Si $X$ est un groupe abélien bisimplicial et $L$ est un intégrateur
sur $\Delta$, on a un isomorphisme naturel 
\[
\H{\Delta\times\Delta}{X} \simeq \Tot({L_!}_\bullet X) 
\]
dans $\Hotab$, où on a noté ${L_!}_\bullet X$ le complexe double
obtenu en appliquant terme à terme le foncteur ${L_!} :
\prefab{\Delta}\to\Ch(\Ab)$ à $X$.
\end{prop}
\begin{proof}
Il suffit de remarquer que les foncteurs ${L_!}_\bullet$ et $\Tot$
commutent aux limites inductives, et que leurs restrictions aux groupes
abéliens bisimpliciaux représentables coïncident avec le foncteur
$L\boxtimes L$ défini ci-dessus, et on peut alors conclure grâce à la
proposition \ref{propIntegrateurCocontinuExtKan}.
\end{proof}

\begin{remark}
On peut évidemment itérer cette construction pour la catégorie
$\Delta^n$ pour tout entier $n>0$.
\end{remark}

\paragr En choisissant l'intégrateur non normalisé $L_\Delta$ sur
$\Delta$ (\ref{ex:IntegrateurLibreDelta}), qui est un
intégrateur libre, cela revient à introduire l'orientation suivante
sur la catégorie $\Delta\times\Delta$, au sens de la section
\ref{secOrientation}. 

L'objet~$\Delta_{n,m}$ de $\Delta\times\Delta$ a pour
dimension l'entier $n+m$. Les monomorphismes de codimension $1$ sont
les morphismes de la forme 
\begin{align*}
(\delta_i,\id) : \Delta_{n-1,m}\to\Delta_{n,m} &\mdvirg n\geq i\geq 0
\mdvirg\\
(\id,\delta_j) : \Delta_{n,m-1}\to\Delta_{n,m} &\mdvirg m \geq j \geq
0 \mdvirg
\end{align*}
et on pose alors
\begin{align*}
\begin{cases}
\sg(\delta_i,\id_{\Delta_m}) &= (-1)^i \\
\sg(\id_{\Delta_n},\delta_j) &= (-1)^n(-1)^j \pbox{.}
\end{cases}
\end{align*}

On vérifie alors que le complexe $L$ associé à cette orientation
(\ref{def:orientation}), défini
pour tout entier $n\geq 0$ par 
\[
(L_{\Delta\times\Delta})_n(\Delta_{i,j}) = \bigoplus_{k+l=n}
\Z^{\left(\Delta_{k,l}\to\Delta_{i,j}\right)}
\]
et dont la différentielle est donnée par la somme des
compositions par les monomorphismes de codimension $1$ comptés avec
leur signes, est bien un complexe, et qu'il correspond au complexe
produit $L_\Delta \boxtimes L_\Delta$. Autrement dit, l'orientation
définie ci-dessus est bien asphérique
(\ref{def:orientationAspherique}).

\paragr On obtient par récurrence, une orientation sur la
catégorie~$\Delta^n$ pour tout entier $n>0$. Concrètement, si $S$ est
un objet de $\Delta^{n-1}$ et~$k\geq 0$, on a 
\[
\dim (\Delta_k, S) = k + \dim S \mdvirg
\]
et on
définit les signes sur les monomorphismes $\varphi$ de codimension $1$
de $\Delta^n$ de la manière suivante :
\begin{itemize}
\item soit $\varphi$ est de la forme $\varphi=(\delta^k_i, \id_S) :
(\Delta_{k-1},S) \to (\Delta_{k},S)$ où $k$ est un entier strictement
positif, $0 \leq i \leq k$ et
$S\in\Ob(\Delta^{n-1})$, et on pose
\[
\sg (\delta_i^k, \id_S) = (-1)^i \pbox{;}
\]
\item soit $\varphi$ est de la forme $(\id_{\Delta_k}, \delta) : (\Delta_k, S') \to
(\Delta_k, S)$ où $k\geq 0$ et $\delta : S \to S'$ est un monomorphisme de
codimension $1$ de $\Delta^{n-1}$, et on pose 
\[
\sg (\id_{\Delta_n}, \delta) = (-1)^n\sg\delta \pbox{.} 
\]
\end{itemize}

On peut alors vérifier que cette orientation est bien asphérique, et
que l'intégrateur associé correspond à l'intégrateur produit
\[
  \underbrace{L_\Delta
    \boxtimes \cdots \boxtimes L_\Delta}_{n \text{ fois}} : \Delta^n
    \to \Ch(\Ab) \pbox{.}
\]

\paragr Puisque la catégorie $\Delta$ est totalement asphérique
(\ref{exTotAspherique}), pour tout entier $n>0$, le morphisme 
\[
\diag : \Delta\to\Delta^n 
\]
est asphérique, et il est donc asphérique en homologie par la
proposition \ref{propAspheriqueImpliqueAsphHomologie}. Cela signifie
qu'on a un diagramme commutatif à isomorphisme naturel près
\[
\xymatrix{
\Hotab_{\Delta^n} \ar[rr]^{\diag^*} \ar[rd]_{\Hf{\Delta^n}} 
&& \Hotab_{\Delta} \ar[ld]^{\Hf{\Delta}}
\\ & \Hotab & \pbox{.}
} 
\]
Dans le cas de la catégorie $\Delta\times\Delta$, on retrouve donc,
comme expliqué dans l'exemple
\ref{exTotalementAspheriqueEilenbergZilber}, l'énoncé du théorème
d'Eilenberg-Zilber\footnote{Le théorème d'Eilenberg-Zilber montre en
fait qu'il existe une équivalence d'homotopie entre le complexe total
du complexe normalisé d'un groupe abélien bisimplicial et le complexe
normalisé de sa diagonale, ce qui est un résultat plus fort.}
(voir \cite{eilenbergZilber1953products} et \cite{doldpuppe1961}): 
\begin{prop}[Eilenberg-Zilber]
Pour tout groupe abélien bisimplicial $X$, on a une chaîne
d'isomorphismes
\[
\Tot\dk_\bullet(X) \simeq \H{\Delta\times\Delta}{X} \simeq
\H{\Delta}{\diag^*(X)} \simeq \dk\diag^*(X)
\]
dans $\Hotab$. 
\end{prop}
\begin{proof}
Le premier isomorphisme provient de la proposition
\ref{homologieBiSimpliciaux} appliquée à l'intégrateur normalisé. Le
deuxième provient de l'asphéricité du foncteur diagonal et de la proposition
\ref{morphismesAsphEnHomologieEquivalences}. On rappelle que le
troisième isomorphisme est un résultat classique de Bousfield et Kan
\cite{bousfield1972homotopy}, que l'on a également prouvé au
paragraphe \ref{integrateurNormalise}.
\end{proof}
\begin{remark}
On peut aussi remplacer le foncteur complexe normalisé dans la
proposition précédente pas le foncteur complexe non
normalisé~$L_\Delta$.
\end{remark}

L'existence de ce foncteur asphérique implique également la
proposition suivante : 
\begin{prop}
Pour tout entier $n\geq1$, la catégorie $\Delta^n$ est une catégorie
test homologique stricte de Whitehead.
\end{prop}
\begin{proof}
C'est un corolaire immédiat de l'existence d'un foncteur asphérique
$\Delta \to \Delta^n$, grâce à la proposition
\ref{propMorphismeAsphSourceTHS}. 
\end{proof}

\section{Ensembles semi-simpliciaux}
Une propriété importante de l'homologie des ensembles simpliciaux est
que le processus de normalisation n'est pas nécessaire, contrairement
par exemple aux ensembles cubiques. On montre ici que cela implique
que la sous-catégorie des monomorphismes de $\Delta$ est aussi une
catégorie test homologique faible.  

\paragr On note $\Delta'$\notindex{$\Delta'$} la sous-catégorie des monomorphismes de
$\Delta$. Les préfaisceaux (resp. les préfaisceaux en groupes
abéliens) sur $\Delta'$ sont appelés des
\ndef[ensemble!semi-simplicial]{ensembles semi-simpliciaux} (resp.
\ndef[groupe abélien!semi-simplicial]{groupes abéliens
semi-simpliciaux}). Pour éviter toute confusion, on notera
$\Delta'_n$ les objets de
$\Delta'$ pour $n\geq 0$. On note
\[
  j : \Delta' \hookrightarrow \Delta
  \]
le foncteur d'inclusion évident.

\begin{prop}
Le foncteur
${L_{\Delta'}}:\Delta'\xrightarrow{j}\Delta\xrightarrow{\dk}\Ch(\Ab)$
est un intégrateur sur $\Delta'$.
\end{prop}
\begin{proof}
Par construction, pour tout objet $\Delta'_n$ de $\Delta'$, le
complexe~$L_{\Delta'}(\Delta'_n)$ a bien l'homologie du point. De
plus, on remarque que, pour tout
entier~$k$, le groupe abélien ${L_{\Delta'}}(\Delta'_n)_k$ est le groupe abélien
libre engendré par les monomorphismes de $\Delta_i$ dans $\Delta_n$,
c'est-à-dire 
\[
{L_{\Delta'}}(\Delta'_n)_k = \Z^{(\Hom_{\Delta'}(\Delta'_k,\Delta'_n))} \pbox{.}
\]
et $L$ est donc bien un intégrateur libre sur $\Delta'$. 
\end{proof}
\begin{remark}
L'orientation standard sur $\Delta$ (voir \ref{ex:orientationDelta})
induit une orientation sur $\Delta'$, et on peut vérifier que
l'intégrateur $L_{\Delta'}$ introduit ci-dessus correspond à
l'intégrateur associé à cette orientation, qui est donc asphérique.
\end{remark}
\begin{coro}
Si $X$ est un groupe abélien semi-simplicial, on a
un isomorphisme  
\[
\H{\Delta'}{X} \simeq X_0 \leftarrow X_1
\leftarrow X_2 \leftarrow \cdots
\]
dans $\Hotab$, où la différentielle est donnée par la somme
des faces comptées avec leur signe.
\end{coro}
\begin{proof}
C'est un cas particulier de la proposition
\ref{propExpressionHomologieDerivateurLibre}.
\end{proof}

\paragr On remarque que le diagramme suivant est commutatif à
isomorphisme naturel près 
\[
\xymatrix{
\prefab{\Delta'} \ar[rd]_{{L_{\Delta'}}_!} && \prefab{\Delta}
\ar[ll]_{j^*} \ar[ld]^{{L_\Delta}_!} \\
& \Ch(\Ab) & \pbox{.}
} 
\]
Cela implique en particulier que le morphisme $j:\Delta'\to\Delta$ est
un morphisme asphérique en homologie, en vertu de la proposition
\ref{morphismesAsphEnHomologieEquivalences}. Le corollaire suivant en
découle.

\begin{coro}
L'intégrateur $L_\Delta'$ est un intégrateur test faible, et la
catégorie $\Delta'$ est une catégorie test homologique faible.
\end{coro}
\begin{proof}
Puisque le foncteur
\[
  L_{\Delta'} = \dk j : \Delta' \to \Ch(\Ab)
\]
est un intégrateur et que $j$ est
asphérique en homologie, on peut appliquer la proposition
\ref{prop:intCoinduitAsphEqTestHomFaible}. Puisque $L_\Delta$ est un
intégrateur test faible, il en est donc de même pour $L_{\Delta'}$.
\end{proof}
\begin{remark}
Le morphisme $j : \Delta' \to \Delta$ est bien asphérique, et pas
seulement asphérique en homologie. On trouvera une preuve dans
\cite[proposition 1.7.25]{maltsiniotis2005}.
\end{remark}
 \begin{remark}
Un argument semblable devrait pouvoir se reproduire dans un cadre
assez général où $A$ est une petite catégorie test homologique, $A'$
désigne la sous-catégorie des monomorphismes de $A$, et où
les épimorphismes de $A$ se comportent bien, au sens où on peut
définir un complexe normalisé de la même manière que dans le cas où
$A$ est la catégorie $\Delta$, induisant ainsi un intégrateur libre
sur $A'$. On pense notamment au cadre étudié par Gagna dans sa thèse
\cite[partie 2.2]{gagna2019homotopie3cat}, qui inclurait notamment le
cas où $A$ est la catégorie $\Theta$. L'espoir
est surtout d'obtenir des conditions pour que les complexes normalisés
et non normalisés soient quasi-isomorphes.
 \end{remark}

\begin{prop}
La catégorie $\Delta'$ n'est pas une catégorie test homologique
locale.
\end{prop}
\begin{proof}
Soit $\Delta'_n$ un objet de $\Delta'$. Comme détaillé en
\ref{integrateurTranches}, on dispose d'un intégrateur sur la
catégorie $\tranche{\Delta'}{\Delta'_n}$ défini comme la composée 
\[
\tranche{\Delta'}{\Delta'_n} \xrightarrow{\alpha} \Delta'
\xrightarrow{L_{\Delta'}} \Ch(\Ab) 
\]
où on note $\alpha$ le morphisme d'oubli. On obtient alors un
isomorphisme  
\[
\H{\tranche{\Delta'}{\Delta'_n}}{X} \simeq \cdots \leftarrow
\bigoplus_{\Delta'_i \to \Delta'_n} X_i \leftarrow
\bigoplus_{\Delta'_{i+1} \to \Delta'_n} X_i \leftarrow \cdots
\]
dans $\Hotab$. Pour tout préfaisceau abélien $X$ sur
$\tranche{\Delta'}{\Delta'_n}$ et pour $i>n$, on a donc
$\mathsf{H}_i(\Delta',X)=0$, ce qui montre que
$\Delta'$ n'est pas une catégorie pseudo-test homologique.
 \end{proof}

\begin{remark}
Il est intéressant de noter la commutativité du diagramme suivant,
reliant la correspondance de Dold-Kan, et le processus d'ajout libre
de dégénérescences à un ensemble semi-simplicial : 
\[
\xymatrix{
\prefab{\Delta'} \ar[rd]_{{L_{\Delta'}}_!} \ar[rr]^{j_!^\ab} 
&&
\prefab{\Delta} 
\\
& \Ch(\Ab) \ar[ru]_{\dk^*} & \pbox{.}
} 
\]
\end{remark}

\section{Cubes}\label{secCubes}
Dans \cite{brown2003cubical}, Brown et Higgins prouvent que la
catégorie des groupes abéliens cubiques avec connexions est équivalente
à la catégorie des complexes de chaînes de groupes abéliens. On montre
dans cette section que cette correspondance calcule bien l'homologie
de la catégorie cubique avec connexions, et donne donc lieu à un
théorème de Dold-Kan homotopique. On montre aussi que la catégorie des
cubes, sans connexion, est une catégorie test homologique faible.

\paragr Pour tout entier $n\geq0$, on note $\cub_n$\notindex{$\cub_n$}
l'ensemble ordonné produit~$\lbrace 0 < 1 \rbrace^n$. Pour $1 \leq i
\leq n$ et $\epsilon \in \lbrace 0,1\rbrace$, on note 
\[
\delta^n_{i,\epsilon} : \cub_{n-1} \to \cub_n 
\]
le morphisme défini par 
\[
\delta^n_{i,\epsilon} (x_1,\dots,x_{n-1}) = (x_1, \dots,
x_{i-1},\epsilon,x_{i},\dots,x_{n-1})
\]
et, pour $n\geq 0$, et $1\leq i \leq n+1$, on note 
\[
\sigma_i^n : \cub_{n+1} \to \cub_n 
\]
le morphisme défini par 
\[
\sigma_i^n(x_1,\dots,x_{n+1}) = (x_1, \dots, x_{i-1}, x_{i+1},
\dots, x_{n+1}) \pbox{.}
\]

\paragr La \emph{catégorie des cubes} $\cub$\notindex{$\cub$} est la
sous-catégorie de la catégorie des ensembles ordonnés dont les objets
sont les ensembles $\cub_n$ pour $n\geq 0$, et les morphismes sont
engendrés par les applications croissantes 
\begin{align*}
\delta_{i,\epsilon}^n &\mdvirg n \geq i \geq 1 \mdvirg \epsilon=0,1 \\
\sigma_i^n &\mdvirg n\geq0 \mdvirg n+1\geq i\geq 1 \pbox{.}
\end{align*}

\paragr Pour $n\geq 1$ et $1\leq i \leq n$, on note 
\[
\gamma_i^n : \cub_{n+1} \to \cub_n 
\]
le morphisme défini par 
\[
\gamma_i^n(x_1,\dots,x_{n+1}) =
(x_1,\dots,x_{i-1},\sup(x_i,x_{i+1}),x_{i+2},\dots,x_{n+1}) 
\]
appelé morphisme de \emph{connexion}. La \emph{catégories des cubes
avec connexions} $\cubc$\notindex{$\cubc$} est la sous-catégorie de la catégorie des
ensembles ordonnés dont les objets sont les ensembles $\cub_n$ pour
$n\geq 0$ et dont les morphismes sont engendrés par les
applications croissantes 
\begin{align*}
\delta_{i,\epsilon}^n &\mdvirg 1 \leq i \leq n \mdvirg \epsilon=0,1 \pbox{,} \\
\sigma_i^n &\mdvirg n\geq0 \mdvirg n+1\geq i\geq 1 \pbox{,}\\
\gamma_i^n & \mdvirg n\geq 1 \mdvirg n \geq i \geq 1 \pbox{.}
\end{align*}

\paragr Un \ndef[ensemble!cubique]{ensemble cubique} (resp. un
\ndef[groupe abélien!cubique]{groupe abélien
cubique}) est un préfaisceau (resp. un préfaisceau en groupes
abéliens) sur la catégorie $\cub$. Un
\ndef[ensemble!cubique!avec connexions]{ensemble cubique avec
connexions} (resp. un \ndef[groupe abélien!cubique!avec connexions]{groupe abélien cubique avec connexions})
est un préfaisceau (resp. un préfaisceau en groupes abéliens) sur la
catégorie $\cub^c$. Étant donné un tel préfaisceau $X$, on notera $X_n
= X(\cub_n)$, et 
\[
d_{i,\epsilon}^n = X(\delta_{i,\epsilon}^n) : X_n\to X_{n-1}
\]
les applications \emph{faces},
\[
s_i^n = X(\sigma_i^n) : X_n \to X_{n+1}
\]
les applications \emph{dégénérescences} et 
\[
\Gamma_i^n = X(\gamma_i^n) : X_n \to X_{n+1} 
\]
les morphismes de \emph{connexion}.

\begin{prop}[Cisinski]
La catégorie $\cub$ est une catégorie test.
\end{prop}
\begin{proof}
Voir \cite[corollaire 8.4.13]{cisinskipref}.
\end{proof}

\begin{prop}[Maltsiniotis]\label{prop:cubesTestStricte}
La catégorie $\cub$ n'est pas une catégorie test stricte. En revanche,
la catégorie $\cubc$ est une catégorie test stricte. 
\end{prop}
\begin{proof}
Voir \cite{maltsiniotis2009cubique}.
\end{proof}

\paragr On introduit l'orientation
suivante sur les catégories $\cub$ et $\cubc$. L'ajout des morphismes
de connexion n'ajoute pas de monomorphismes, et ne modifie donc pas la
dimension au sens de la section \ref{secOrientation} : la
dimension de
l'objet $\cub_n$ vu dans $\cub$ ou dans $\cubc$ est l'entier $n\geq
0$. On pose alors 
\[
\sg(\delta^n_{i,\epsilon}) = (-1)^{i+\epsilon} \mdvirg 1\leq i \leq n
\mdvirg \epsilon=0,1\pbox{.}
\]
On peut alors montrer qu'on a bien défini une orientation sur les
catégories~$\cub$ et~$\cubc$. Le complexe de chaînes associé
(\ref{def:orientation}) se décrit de la manière suivante :
pour tout entier $n\geq 0$ et pour tout entier $k\geq0$, on rappelle
qu'on a
\[
L_\cub(\cub_n)_k = \Z^{(\Hom_\cub(\cub_k,\cub_n))}  \mdvirg
L_{\cubc}(\cub_n)_k = \Z^{(\Hom_{\cubc}(\cub_k,\cub_n))} 
\]
\notindex{$L_\cub$}%
\notindex{$L_\cubc$}%
et on obtient, pour tout
entier~$k>0$ et pour tout morphisme $u :\cub_k \to \cub_n$ de~$\cub$
(ou de $\cubc$)
\[
d \langle u \rangle =
\sum_{i=1}^k(-1)^i \left(\langle d^k_{i,0}u\rangle - \langle
d^k_{i,1}u\rangle\right) \pbox{.}
\]

En revanche, cette orientation n'est pas asphérique
au sens du paragraphe~\ref{def:orientationAspherique}. On peut s'en convaincre
rapidement : puisque $\cub_0$ est un objet final dans $\cub$ et dans
$\cubc$, les préfaisceaux $\Wh{\cub}{\cub_0}$ et $\Wh{\cubc}{\cub_0}$
coïncident avec le préfaisceau constant de valeur $\Z$ respectivement sur
$\cub$ et sur $\cubc$. Or, puisque la différentielle définie ci-dessus
est une somme paire en tout degré, on a des isomorphismes de complexes
\[
(L_\cub)_!(\Z) \simeq \Z \xleftarrow{0} \Z \xleftarrow{0} \Z
\xleftarrow{0} \cdots
\]
et 
\[
(L_{\cubc})_!(\Z) \simeq \Z \xleftarrow{0} \Z \xleftarrow{0} \Z
\xleftarrow{0} \cdots
\]
et les groupes d'homologie valent $\Z$ en tout degré.
Cela implique que ni $L_\cub$, ni $L_{\cubc}$ ne sont
des intégrateurs.

\paragr Pour obtenir un complexe calculant la bonne homologie, il faut
donc normaliser. Pour les groupes abéliens cubiques avec connexions, on
dispose alors d'un choix : on peut quotienter par les dégénérescences,
ou quotienter par les dégénérescences et les connexions.

Pour tout groupe abélien cubique avec connexions $X$, on note alors~$D_\sigma(X)_n$ le sous-groupe engendré par les éléments de $X_n$ qui
sont dans l'image d'une application de dégénérescence $s_i : X_{n-1}
\to X_n$ pour $1 \leq i \leq n$, et~$D_\gamma(X)_n$ le sous-groupe
engendré par les éléments qui sont dans l'image d'une
connexion~$\Gamma_i : X_{n-1} \to X_n$ pour $1 \leq i \leq n$.

\begin{prop}[Barcelo-Greene-Jarrah-Welker]\label{normalisationsCubesConnexionsOuPas}
Soit $X$ un groupe abélien cubique avec connexions. Alors : 
\begin{enumerate}
\item les sous-groupes $D_\sigma(X)_\bullet$ et $D_\gamma(X)_\bullet$
de $X_\bullet$ forment des sous-complexes de ${L_{\cub}}_!(X)$;
\item il existe un isomorphisme naturel 
\[
{L_{\cub}}_!(X) / {D_\sigma(X)} \simeq 
{L_{\cub}}_!(X) / ({D_\sigma(X)} \oplus D_\gamma(X) )  
\]
dans $\Hotab$.
\end{enumerate}
\end{prop}
\begin{proof}
Voir \cite[corollaire 3.10]{barcelo2018homologygroupscubicalsets}.
\end{proof}

\begin{remark}
Si $X$ est un groupe abélien cubique sans connexion, on peut aussi
montrer que les sous-groupes $D_\sigma(X)_\bullet$ forment bien un
sous-complexe du complexe~${L_{\cub}}_!(X)$.
\end{remark}

\paragr Si $X$ est un groupe abélien cubique, on note
\termindex{complexe normalisé!d'un groupe abélien cubique}
\begin{align*}
\dkc : \prefab{\cub}&\to\Ch(\Ab) \\
X &\mapsto {L_{\cub}}_!(X) / D_\sigma(X)
\end{align*}
\notindex{$\dkc$}%
\notindex{$\dkcc$}%
et si $X$ est un groupe abélien cubique avec connexions, on note 
\begin{align*}
\dkcc : \prefab{\cubc}&\to\Ch(\Ab) \\
X &\mapsto {L_{\cubc}}_!(X) / (D_\sigma(X) \oplus D_\gamma(X) ) \pbox{.}
\end{align*}

\medskip
Pour montrer que ces foncteurs calculent la bonne homologie,
on doit montrer que leurs restrictions aux préfaisceaux représentables
sont des intégrateurs. Dans le cas des groupes abéliens cubiques avec
connexions, on va introduire un autre complexe, qui est l'objet du
théorème de Dold-Kan cubique prouvé par Brown et Higgins dans
\cite{brown2003cubical}. 

\paragr\label{defCubesintersectionNoyaux} Pour tout groupe abélien avec connexions $X$, on pose~$\normcc
(X)_0 = X_0$ et, pour tout entier~$n>0$, 
\[
\normcc(X)_n = \bigcap_{\substack{1 \leq i \leq n\\
\epsilon\in\lbrace0,1\rbrace\\\mathclap{(i,\epsilon)\neq(n,0)}}} \ker d^n_{i,\epsilon}
\subset X_n\pbox{.}
\]
On peut alors vérifier qu'on obtient bien un sous-complexe de
${L_{\cubc}}_!(X)$ que l'on note $\normcc X$.

\begin{prop}[Carranza-Kapulkin-Tonks]\label{normaliseCubesIntersectionNoyaux}
Pour tout groupe abélien cubique avec connexions~$X$, on a un
isomorphisme naturel de complexes de chaînes
\[
{L_{\cubc}}_! X \simeq \normcc X \oplus (D_\sigma X \oplus D_\gamma X) \pbox{.}
\]
En particulier, on a un isomorphisme naturel de complexes de chaînes
\[
\dkcc X \simeq \normcc X \pbox{.}
\]
\end{prop}
\begin{proof}
Voir \cite[introduction, théorème 4.8 et remarque
4.9]{carranza2023hurewicz}.
\end{proof}

\paragr On note $c_{\cubc}$ la restriction du foncteur $\dkcc$ aux
préfaisceaux représentables, c'est-à-dire la composée
\[ 
\cubc \xrightarrow{\Whf{\cubc}} \prefab{\cubc}
\xrightarrow{\dkcc} \Ch(\Ab) \pbox{.}
\]

\begin{coro}
Le foncteur $c_\cubc$, vu comme un complexe de préfaisceaux sur
$(\cubc)^{\op}$, est un complexe de préfaisceaux projectifs.
\end{coro}
\begin{proof}
Par la proposition précédente, $c_{\cubc}$ est un rétracte de $L_{\cubc}$,
qui est un complexe de préfaisceaux libres, et on peut conclure grâce
à la proposition \ref{projectifPrefab}.
\end{proof}

\paragr Pour montrer que $\dkcc$ calcule bien l'homologie au sens de la
section \ref{sectypehomologieassocie}, il reste à montrer que
$c_\cubc$ est une résolution du préfaisceau constant de valeur $\Z$.
Pour cela, on peut invoquer l'existence d'un foncteur de réalisation
géométrique~$|-| : \pref{\cubc}\to\Top$, et un résultat classique
(voir par exemple~\cite[théorème
3.17]{barcelo2018homologygroupscubicalsets})
affirme que pour tout ensemble cubique avec connexions $X$, l'homologie
du complexe $\dkcc X$ coïncide avec l'homologie de l'espace
topologique~$|X|$. Ainsi, le complexe $c_\cubc(\cub_n)$ a le type
d'homologie du $n$-cube standard, c'est-à-dire le type d'homologie du
point. 

\begin{prop}\label{homologieCubiquesConnexions}
Le foncteur 
\[
c_{\cubc} : \cubc \xrightarrow{\Whf{\cubc}} \prefab{\cubc} \xrightarrow{\dkcc}
\Ch(Ab) 
\]
\notindex{$c_{\cubc}$}%
est un intégrateur sur $\cubc$. Par conséquent, pour tout groupe abélien cubique
avec connexions $X$, on a un isomorphisme naturel 
\[
\H{\cubc}{X}\simeq \dkcc X
\]
dans $\Hotab$.
\end{prop}
\begin{proof}
Pour conclure à partir du dernier paragraphe, il suffit d'appliquer
la proposition \ref{propIntegrateurCocontinuExtKan} affirmant que
puisque le foncteur $\dkcc$ commute aux limites inductives, il
coïncide avec l'extension de Kan à gauche de sa restriction.
\end{proof}

\paragr En combinant la proposition ci-dessus et la proposition
\ref{normaliseCubesIntersectionNoyaux}, on obtient finalement, pour
tout groupe abélien cubique avec connexions, des isomorphismes naturels
\[
\H{\cubc}{X}\simeq \dkcc X \simeq \normcc X
\]
dans $\Hotab$. En vertu de la proposition
\ref{normalisationsCubesConnexionsOuPas}, on peut aussi simplement
quotienter par le sous-complexe $D_\sigma$ des dégénérescences sans
quotienter par les connexions.

\begin{theorem}[Brown-Higgins]
Le foncteur
\[
\normcc : \prefab{\cubc}\to\Ch(\Ab) 
\]
est une équivalence de catégories.
\end{theorem}
\begin{proof}
On renvoie à \cite{brown2003cubical} pour la preuve originale. Le
foncteur $\normcc$ y est obtenu en deux étapes : on associe d'abord à
un groupe abélien cubique avec connexions un \emph{complexe croisé} en
groupes abéliens, puis un complexe de chaînes en groupes abéliens.
Chacune de ces étapes constitue une équivalence de catégories. 
\end{proof}

\begin{remark}
La formule que donnent Brown et Higgins dans \cite{brown2003cubical}
et~\cite{brown2011nonabelian} est en fait légèrement différente de
celle donné au paragraphe \ref{defCubesintersectionNoyaux} : en degré
$n\geq1$, Brown et Higgins considèrent l'intersection des noyaux de
toutes les faces $d_{i,\epsilon}$ pour $1\leq i \leq n$ et
$\epsilon\in\lbrace0,1\rbrace$ sauf pour $(i, \epsilon) = (1,0)$. Pour passer d'une
version à l'autre, on peut utiliser l'automorphisme de la
catégorie~$\cubc$ induit par le renversement des composantes des
cubes.
\end{remark}

\begin{coro}
La catégorie $\cubc$ est une catégorie test homologique stricte.
\end{coro}
\begin{proof}
Puisque le foncteur $\normcc$ est une équivalence de catégories, il en
est de même pour le foncteur $\dkcc$, et
on a montré que ce dernier calcule l'homologie des groupes abéliens
avec connexions (\ref{homologieCubiquesConnexions}). Le foncteur induit
entre les catégories localisées est donc encore une équivalence de
catégories. La catégorie $\cubc$ est donc une catégorie test
homologique faible. De plus, $\cubc$ est totalement asphérique (voir
\ref{prop:cubesTestStricte}), et elle est donc totalement asphérique
en homologie. On peut donc conclure que la catégorie $\cubc$ est une
catégorie test homologique stricte grâce à la
proposition~\ref{propTestHomologiqueTotAspherique}.
\end{proof}

Le but de ce qui suit est de montrer que la catégorie $\cubc$
est une catégorie de Whitehead.

\paragr On note 
\[
j : \cubc \to \Top 
\]
le foncteur envoyant, pour tout entier $n \geq 0$, l'objet $\cub_n$
sur l'espace topologique produit $\lbrack0,1\rbrack^n$. On obtient
alors un foncteur
\[
  j_! : \pref{\cubc} \to \Top
  \]
appelé foncteur de \emph{réalisation géométrique}, qui admet pour
adjoint à droite le foncteur \emph{complexe singulier cubique avec connexions}
\[
j^* : \Top \to \pref{\cubc} \pbox{.}
\]

Dans \cite{carranza2023hurewicz}, Carranza et Kapulkin donnent
plusieurs descriptions des groupes d'homotopie des ensembles cubiques
avec connexions, c'est-à-dire des groupes d'homotopie de leur
réalisation géométrique. Dans le cas des groupes abéliens cubiques
avec connexions, ces groupes d'homotopie admettent une description
combinatoire particulièrement simple, à l'image des groupes
d'homotopie des groupes abéliens simpliciaux (voir
\ref{homologieHomotopieNormalise}).

\begin{prop}[Carranza-Kapulkin-Tonks]\label{homologieHomotopieNormaliseCub}
Si $X$ est un groupe abélien cubique avec connexions, alors
\begin{enumerate}
\item il existe une bijection naturelle $\pi_0(\U X) \simeq
H_0(\normcc X)$;
\item pour tout entier $n>0$, il existe un isomorphisme naturel 
\[
\pi_n(\U X,0) \simeq H_n(\normcc X) \pbox{.}
\]
\end{enumerate}
\end{prop}
\begin{proof}
Voir \cite[théorème 4.1]{carranza2023hurewicz}.
\end{proof}

On va utiliser cette proposition pour montrer que $\cubc$ est une
catégorie de Whitehead. Pour cela, on a besoin de montrer que les
groupes d'homotopie des ensembles cubiques avec connexions coïncident
avec les groupes d'homotopie de leur catégorie des éléments. Ce
résultat est peut-être bien connu des spécialistes, mais nous n'avons
pas trouvé de preuve dans la littérature. 
\begin{lemme}
On a l'égalité 
\[
\W_{\cubc} = j_!^{-1}(\W_\Top).
\]
En d'autres termes, les équivalences test (\ref{def:eqTest})
de $\pref{\cubc}$ coïncident avec les morphismes d'ensembles cubiques
avec connexions
dont l'image par le foncteur de réalisation géométrique est une
équivalence faible topologique.
\end{lemme}
\begin{proof}
On peut montrer que le foncteur $j_! : \pref{\cubc}\to\Top$ envoie les
équivalences test sur des équivalences faibles topologiques. En effet,
\hbox{Cisinski} donne dans \cite[théorème 1.7]{cisinski2014univalent} des
générateurs des cofibrations et
des cofibrations triviales de la structure de catégorie de
modèles de Cisinski sur les ensembles cubiques avec connexions (voir théorème
\ref{structureCisinskiTestLocal}). On peut alors vérifier que $j_!$
est un foncteur de Quillen à gauche pour cette structure de catégorie
de modèles, et donc qu'il préserve les équivalences faibles, puisque
tous les ensembles cubiques avec connexions sont cofibrants pour cette
structure. Cela implique qu'on a~$\W_{\cubc} \subset
j_!^{-1}(\W_{\Top})$, et que $j_!$ induit donc un foncteur
\[
\overline{j_!} : \Hot_\cubc \to \Hot_\Top \pbox{.}
\]
Pour conclure, on veut montrer que ce foncteur est une équivalence de
catégories.
On va en fait
montrer que le foncteur $j^* : \Top \to \pref{\cubc}$
envoie les équivalences faibles topologiques sur des équivalences
test, et qu'il induit une équivalence de catégories entre les
catégories localisées. 

Pour cela, on montre d'abord que le foncteur $j^*$ envoie les espaces
topologiques contractiles sur des préfaisceaux
asphérique (au sens du paragraphe~\ref{def:PrefaisceauAspherique}).
Puisque~$j(\cub_1)=\lbrack0,1\rbrack$, on dispose d'un
morphisme d'adjonction 
\[
\eta_{1} : \cub_1 \to j^*j_!(\cub_1) = j^*(\lbrack0,1\rbrack) \pbox{.}
\]
Ainsi, étant donné une homotopie $h$ entre deux morphismes $f,g : X \to Y$ dans~$\Top$
\[
\xymatrix{
X \ar[d]_{i_0} \ar[rd]^{f} \\
\lbrack0,1\rbrack \times X \ar[r]^-{h} & Y \\
X \ar[u]^{i_1} \ar[ru]_{g} & \pbox{,}
}
\]
on obtient (en utilisant le fait que le foncteur $j^*$ commute aux produits) un diagramme commutatif dans $\pref{\cubc}$
\[
\xymatrix@C=4em{
&j^*(X) \ar[ld]_{\delta^1_{1,0}\times\id} \ar[rd]^{f}
\ar[d]^{j^*(i_0)}
\\
\cub_1 \times j^*(X) \ar[r]^-{\eta_1 \times \id} 
& j^*(\lbrack0,1\rbrack\times X) \ar[r]^{j^*(h)}
& j^*(Y) 
\\
&j^*(X) \ar[lu]^{\delta^1_{1,1}\times\id} \ar[ru]_{g} 
\ar[u]_{j^*(i_1)} & \pbox{.}
}
\]
Puisque~$\cub_1$ est un préfaisceau représentable et que $\cubc$ est
une catégorie totalement asphérique, alors $\cub_1$ est un préfaisceau
localement asphérique (voir proposition
\ref{prop:totalementAspheriqueEquivalences}). Cela
implique que les flèches obliques de gauche dans le diagramme
ci-dessus sont dans $\W_{\cubc}$ (voir \cite[paragraphes 1.2.5 et
1.2.6]{maltsiniotis2005}), et donc, par deux-sur-trois, que $f$ est une équivalence test si
et seulement si il en est de même pour $g$. Ainsi, si~$X$ est un
espace topologique contractile, l'identité de $X$ est homotope à un
morphisme constant, et on peut déduire par faible saturation
(\ref{defFaibleSaturation}) que le morphisme $j^*(X) \to e_{\pref{\cubc}}$ est
donc une équivalence test, ce qui signifie que le préfaisceau~$j^*(X)$ est
asphérique.

Dans \cite{ara2022compnerfs}, Ara et Maltsiniotis établisent un cadre
général pour comparer des foncteurs nerfs, que l'on peut
utiliser dans ce contexte de la manière suivante. Puisque les cubes et les simplexes topologiques sont contractiles,
l'argument précédent permet d'appliquer la
proposition \cite[proposition 4.13]{ara2022compnerfs} aux foncteurs
\[
\Delta \xrightarrow{|-|} \Top \xleftarrow{j} \cubc 
\]
(en utilisant également le fait que le foncteur complexe singulier
préserve les équivalences faibles), qui garantit alors la
commutativité à isomorphisme naturel près du diagramme
\[
\xymatrix{
& \Top \ar[rd]^{j^*} \ar[ld]_{\mathsf{Sing}} \\
\pref{\Delta} \ar[r]_{i_\Delta} & \Cat \ar[d]& \pref{\cubc} \ar[l]^{i_\cubc} \\
& \Hot & \pbox{.}
} 
\]
En particulier, le foncteur $j^*$ envoie les équivalences faibles
topologiques sur des équivalences test, et induit donc un foncteur
s'insérant dans le diagramme commutatif à isomorphisme naturel près
\[
\xymatrix{
& \Hot_\Top \ar[rd]^{\overline{j^*}} \ar[ld]_{\overline{\mathsf{Sing}}}^\simeq \\
\Hot_{\Delta} \ar[rd]_{\overline{i_\Delta}}^\simeq && \Hot_{\cubc}
\ar[ld]^{\overline{i_\cubc}}_\simeq \\
& \Hot & \pbox{.}
} 
\]
Le foncteur $\overline{j^*}$ est donc une équivalence de catégories, et il en est
de même pour son adjoint à gauche, qui est le foncteur 
$\overline{j_!} : \Hot_{\cubc} \to \Hot_\Top$. Ce dernier reflète donc les isomorphismes. 

Pour conclure,
on utilise le fait que, puisque $\cubc$ est une catégorie test, un
morphisme d'ensembles cubiques avec connexions est dans $\W_{\cubc}$
si et seulement si son image dans $\Hot_{\cubc}$ est un isomorphisme
(voir
\cite[corollaire 1.4.7 et proposition 4.2.3]{cisinskipref}), ce qui
permet de montrer que $j_!^{-1}(\W_\Top) \subset \W_{\cubc}$. 
\end{proof}

\begin{coro}
La catégorie $\cubc$ est une catégorie de Whitehead.
\end{coro}
\begin{proof}
La proposition \ref{homologieHomotopieNormaliseCub} et le
lemme précédent permettent d'affirmer
que pour tout morphisme $f$ de groupes abéliens cubiques avec
connexions, le morphisme $\normcc(f)$ est un quasi-isomorphisme si et
seulement si le morphisme d'ensembles cubiques avec connexions $\U f$
est un élément de~$\W_{\cubc}$.
Puisqu'on a montré que le complexe $\normcc X$ calcule bien l'homologie
des groupes abéliens cubiques au sens du paragraphe \ref{defH_A}, on peut conclure de la même manière que
pour $\Delta$ (voir la proposition \ref{DeltaWhitehead}).
\end{proof}

\medskip

On va maintenant traiter le cas des groupes abéliens cubiques sans
connexion. À titre de comparaison, on énonce le résultat suivant : 

\begin{theorem}[Lack-Street]
Il existe une équivalence de catégories 
\[
\prefab{\cub} \simeq \prefab{\Delta'} \mdvirg
\]
où $\Delta'$ désigne la sous-catégorie des monomorphismes de $\Delta$.
\end{theorem}
\begin{proof}
Voir \cite[exemple 7.3]{lack2014combinatorial}.
\end{proof}

\paragr\label{defc_cub} On note $c_\cub$ la restriction du foncteur $\dkc$ aux
\notindex{$c_\cub$}%
préfaisceaux représentables, c'est-à-dire la composée
\[
c_\cub : \cub \xrightarrow{\Whf{\cub}} \prefab{\cub}
\xrightarrow{\dkc} {\Ch(\Ab)} \pbox{.}
\]

Un résultat de Husainov \cite{husainov2018homologycubical}
affirme déjà que $c_\cub$ est un intégrateur sur~$\cub$, et donc que
le foncteur $\dkc$ calcule l'homologie au sens introduit dans la
section~\ref{sectypehomologieassocie}. En particulier, on trouvera une
preuve combinatoire du fait que pour tout entier $n\geq0$, le
complexe $c_\cub(\cub_n)$ a bien l'homologie du point.

\begin{prop}[Husainov]
Le foncteur $c_\cub : \cub \to \Ch(\Ab)$ est un intégrateur sur
$\cub$. En particulier, pour tout groupe abélien cubique $X$, on a un
isomorphisme naturel 
\[
\H{\cub}{X} \simeq \dkc X 
\]
dans $\Hotab$.
\end{prop}
\begin{proof}
Voir \cite{husainov2018homologycubical}.
\end{proof}

\paragr On va maintenant montrer que $\cub$ est une catégorie test
homologique faible. Notons 
\[
c_{\cubc,\gamma} : \cubc \to \Ch(\Ab) 
\]
le foncteur associant à tout objet $\cub_n$ de $\cub$ le complexe
$L(\cub_n)/D_\sigma(\cub_n)$, c'est-à-dire le complexe obtenu en
quotientant par les dégénérescences mais pas par les connexions. En
vertu de la proposition \ref{normaliseCubesIntersectionNoyaux} et de
la proposition \ref{normalisationsCubesConnexionsOuPas},
$c_{\cubc,\gamma}$ est encore un intégrateur sur $\cubc$. De plus, en
notant $j : \cub \to \cubc$ le foncteur évident, on voit que le
diagramme 
\begin{equation}\label{diag:cubtesthomfaible}
\xymatrix{
{\cub} \ar[rr]^{j} \ar[rd]_{c_\cub} && {\cubc}
\ar[ld]^{c_{\cubc,\gamma}} \\
& {\Ch(\Ab)}
} 
\end{equation}
est commutatif.
\begin{prop}
L'intégrateur $c_{\cub}$ défini au paragraphe \ref{defc_cub} est un
intégrateur test faible, et la catégorie $\cub$ est une catégorie test
homologique faible.
\end{prop}
\begin{proof}
Le foncteur $j$ induit à son tour un foncteur
\[
  j^* : \prefab{\cubc}\to\prefab{\cub} 
\]
et on remarque qu'on obtient un diagramme commutatif
\[
\xymatrix{
{\prefab{\cubc}} \ar[rr]^{j^*} \ar[rd]_{\dkccg} && {\prefab{\cub}}
\ar[ld]^{\dkc} \\
& {\Ch(\Ab)} & \pbox{,}
} 
\] 
où on a noté $\dkccg$ l'extension de Kan à gauche de
$c_{\cubc,\gamma}$, c'est-à-dire le foncteur associant à tout groupe
abélien cubique avec
connexions $X$ le complexe obtenu en quotientant le complexe non
normalisé par les dégénérescences, et pas par les
connexions. Comme expliqué au paragraphe
\ref{normalisationsCubesConnexionsOuPas}, le foncteur $\dkccg$ calcule
également l'homologie des groupes abéliens avec connexions. Ainsi,
par la proposition \ref{morphismesAsphEnHomologieEquivalences}, le
foncteur $j$ est asphérique en homologie. Or, l'intégrateur
$c_{\cubc,\gamma}$ est un intégrateur test faible sur $\cubc$, et il
en est donc de même pour l'intégrateur $c_\cub$ grâce à la
commutativité du diagramme \ref{diag:cubtesthomfaible} et à la
proposition \ref{prop:intCoinduitAsphEqTestHomFaible}.
\end{proof}

\begin{remark}
On peut montrer que la catégorie $\cubc^+$ obtenue en remplaçant, dans
la définition des connexions, les $\sup$ par des $\inf$, c'est-à-dire
en rajoutant à la catégorie $\cub$ des \emph{connexions positives},
définies pour $n\geq 0$ et~$1 \leq i \leq n$ par 
\begin{align*}
\gamma_{i,+} : \cub_{n+1} &\to \cub_n \\
(x_1,\dots,x_{n+1}) &\mapsto (x_1,\dots,\inf(x_i,x_{i+1}),\dots,x_{n+1})
\end{align*}
est également une catégorie test stricte, ainsi qu'une catégorie test
homologique stricte de Whitehead. En effet, les objets $\cub_n$ de
$\cubc$, vus comme des ensembles ordonnés et donc comme des
catégories, sont isomorphes à leur catégories opposées.
On peut alors vérifier que l'automorphisme de $\Cat$
envoyant une catégorie sur sa catégorie opposée induit
alors un isomorphisme de la catégorie $\cubc$ vers la catégorie
$\cubc^+$.

De même, on peut considérer la catégorie
$\cubc^{\pm}$, obtenue en considérant les connexions positives et les
connexions que nous avons utilisées ici, appelées alors connexions
négatives. C'est avec cette catégorie que travaillent Carranza,
Kapulkin et Tonks dans \cite{carranza2023hurewicz}, et on y trouve un
preuve de la proposition \ref{homologieHomotopieNormaliseCub} pour la
catégorie $\cubc^\pm$, qui implique que cette dernière est une
catégorie de Whitehead. De plus, d'après Maltsiniotis \cite[commentaire
2, p 392]{maltsiniotisPS}, cette catégorie est également une catégorie
test stricte, et le théorème
\ref{thm:TestWhiteheadPseudoTestHomologique} implique donc que la
catégorie $\cubc^\pm$ est une catégorie test homologique stricte de
Whitehead. 
\end{remark}

\section{\texorpdfstring{Globes et $\omega$-catégories }%
                               {Globes et ω-catégories }}
                               \label{secGlobes}
La catégorie $\Grefl$ des \emph{globes réflexifs} fait partie de ce
que Grothendieck appelle la \og trinité des catégories test \fg{}
(avec $\Delta$ et $\cub$) dans \emph{Pursuing Stacks}. Toutefois,
Maltsiniotis prouve dans \cite[commentaire 1]{maltsiniotisPS} que
celle-ci n'est pas une catégorie pseudo-test. Cependant, le théorème
de Bourn (\ref{th:Bourn}) implique que la catégorie des préfaisceaux
abéliens sur $\Grefl$ est équivalente à la catégorie des complexes de
chaînes de groupes abéliens. On va voir que $\Grefl$ est bien une
catégorie test homologique, mais qu'elle n'est pas
une catégorie de Whitehead. On en profite aussi pour définir la notion
de $\omega$-catégorie stricte, qui nous servira pour définir les exemples
suivants.

\paragr On rappelle que la catégorie des \ndef{globes} $\G$ est la petite catégorie 
\notindex{$\G$}%
engendrée par le graphe 
\[
\xymatrix{
\D_0 \ar@/^/[r]^{\sigma_1} \ar@/_/[r]_{\tau_1} & 
\D_1 \ar@/^/[r]^{\sigma_2} \ar@/_/[r]_{\tau_2} & 
\D_2 \ar@/^/[r]^{\sigma_3} \ar@/_/[r]_{\tau_3} & 
\cdots
\ar@/^/[r]^{\sigma_n} \ar@/_/[r]_{\tau_n} & 
\D_n 
\ar@/^/[r]^{\sigma_n} \ar@/_/[r]_{\tau_n} & 
\cdots
}
\]
\notindex{$\D_n$}%
soumise aux relations \emph{coglobulaires}
\[
\sigma_{n+1}\circ \sigma_{n} = \tau_{n+1} \circ \sigma_{n} \mdvirg
\sigma_{n+1}\circ \tau_n = \tau_{n+1} \circ \tau_n \mdvirg n \geq 1
\pbox{.}
\]
On appelle
\ndef[$\omega$-graphe]{$\omega$-graphes} ou \ndef[ensemble!globulaire]{ensembles
globulaires} les préfaisceaux sur $\G$. Si $X$ est un ensemble
globulaire, on note $X_n$ l'ensemble $X(\D_n)$ pour $n\geq0$,
et, pour tout entier $n\geq 1$, on note
\[
  t_n=X(\tau_n) : X_n \to X_{n-1}
\]
  et
\[
  s_n=X(\sigma_n) : X_n \to X_{n-1}
\]
les applications \emph{but} et \emph{source}.

\paragr On a déjà montré (voir l'exemple \ref{exGlobesHomologie}) que 
si $X$ est un \ndef[$\omega$-graphe!en groupes abéliens]{$\omega$-graphe
abélien},\termindex{groupe abélien!globulaire}%
c'est-à-dire un préfaisceau en groupes abéliens sur $\G$,
alors son homologie est celle du
complexe
\[
\H{\G}{X} = X_0 \xleftarrow{t-s} X_1 \xleftarrow{t-s} X_2 \leftarrow
\cdots 
\]
ce qui implique en particulier que la catégorie $\G$ a pour groupes
d'homologie \[
\mathsf{H}_i(\G,\Z) \simeq \Z 
\]
pour tout entier $i \geq 0$.
Par conséquent, en vertu de la proposition \ref{prop:testHomFaibleImpliqueAspherique}, $\G$
n'est pas une catégorie test homologique faible.
On va s'intéresser dans cette section à la catégorie des globes
\emph{réflexifs}, utilisée pour définir les $\omega$-catégories
strictes. Pour un exposé détaillé, on pourra se reporter à 
\cite[chapitre 1]{ara2010infini} sur lequel nous nous sommes largement
basés.

\paragr La catégorie $\Grefl$ des \ndef[globes!réflexifs]{globes réflexifs} est la catégorie
\notindex{$\Grefl$}%
obtenue en ajoutant à $\G$ des morphismes 
\[
\kappa_n : \D_{n+1} \to \D_n \mdvirg n \geq 0
\]
soumis aux relations suivantes : 
\[
\kappa_n \circ \sigma_{n+1} = \id \mdvirg \kappa_n \circ \tau_{n+1} =
\id \mdvirg n \geq 0\pbox{.} 
\]
En d'autres termes, il s'agit de la catégorie 
engendrée par le graphe 
\[
\xymatrix@C=3em{
\D_0 \ar@/^1pc/[r]^{\sigma_1} \ar@/_1pc/[r]_{\tau_1} & 
\D_1 \ar@/^1pc/[r]^{\sigma_2} \ar@/_1pc/[r]_{\tau_2} \ar[l]|{\kappa_0}  & 
\D_2 \ar@/^1pc/[r]^{\sigma_3} \ar@/_1pc/[r]_{\tau_3} \ar[l]|{\kappa_1} & 
\cdots
\ar@/^1pc/[r]^{\sigma_n} \ar@/_1pc/[r]_{\tau_n} \ar[l]|{\kappa_2}& 
\D_n \ar@/^1pc/[r]^{\sigma_{n+1}} \ar@/_1pc/[r]_{\tau_{n+1}}
\ar[l]|{\kappa_{n-1}}& \cdots \ar[l]|{\kappa_n}
}
\]
soumise aux relations \emph{coglobulaires réflexives}
\begin{align*}
\sigma_{n+1}\circ \sigma_{n} = \tau_{n+1} \circ \sigma_{n} &\mdvirg
\sigma_{n+1}\circ \tau_n = \tau_{n+1} \circ \tau_n \mdvirg n\geq 1 \\
\kappa_n \circ \sigma_{n+1} = \id &\mdvirg \kappa_n \circ \tau_{n+1} =
\id \mdvirg n \geq 0\pbox{.} 
\end{align*}
Si $n \leq m$ sont deux entiers, on notera 
\begin{align*}
\sigma_n^m &= \sigma_m\sigma_{m-1}\cdots\sigma_{n+1} : \D_n \to \D_m
\mdvirg\\ 
\tau_n^m &= \tau_m\tau_{m-1}\cdots\tau_{n+1} : \D_n \to \D_m \mdvirg\\
\kappa^n_m &= \kappa_n\kappa_{n+1}\cdots\kappa_{m-1} : \D_m \to \D_n
\pbox{.}
\end{align*}

\paragr Un préfaisceau d'ensembles (resp. de groupes abéliens) sur
$\Grefl$ est appelé un \ndef[ensemble!globulaire!réflexif]{ensemble
globulaire} (resp. un \ndef[groupe abélien!globulaire!réflexif]{groupe
abélien globulaire) réflexif}. Si~$X$ est un ensemble ou un groupe abélien
globulaire réflexif, on note, pour tout entier $n\geq0$, $X_n =
X(\D_n)$ l'ensemble (resp. le groupe abélien) des \emph{$n$-cellules}
de $X$. On note également, pour
tout entier $n\geq1$,
\begin{align*}
s_n&=X\sigma_n : X_n \to X_{n-1} \\
t_n&=X\tau_n : X_n \to X_{n-1} 
\end{align*}
les applications \emph{source} et \emph{but} associés, ainsi que 
\[
k_n = X\kappa_n : X_n \to X_{n+1} 
\]
l'application \emph{identité}. Si $m$ et $n$ sont deux entiers tels
que $m \geq n \geq 0$, on note
\begin{align*}
s_n^m &= X\sigma_n^m : X_m \to X_n \\
t_n^m &= X\tau_n^m : X_m \to X_n \\
k^n_m &= X\kappa^m_n : X_n \to X_m
\end{align*}
les morphismes source, but et identité itérés. La donnée d'un
tel préfaisceau correspond donc à la donnée d'un diagramme d'ensembles
(ou de groupes abéliens) 
\[
\xymatrix{
X_0 \ar[r]|{k_0} 
& X_1 \ar@/^1pc/[l]^{t_1} \ar@/_1pc/[l]_{s_1} \ar[r]|{k_1} 
& X_2 \ar@/^1pc/[l]^{t_2} \ar@/_1pc/[l]_{s_2} \ar[r]|{k_2}
& \cdots \ar@/^1pc/[l]^{t_3} \ar@/_1pc/[l]_{s_3} 
}
\]
vérifiant les relations \emph{globulaires réflexives} 
\begin{align*}
s_n \circ s_{n+1} = s_n \circ t_{n+1} &\mdvirg
t_n \circ t_{n+1} = t_n \circ s_{n+1} \mdvirg n \geq 1 \\
t_{n+1}\circ k_n = \id &\mdvirg s_{n+1}\circ k_n=\id \mdvirg n\geq 0
\pbox{.}
\end{align*}

\begin{defin}
Une \ndef{$\omega$-catégorie} (sous-entendu, \emph{stricte}) est un ensemble globulaire
réflexif $X$ muni d'applications de composition
\[
*_n^m : (X_m, s_n^m) \times_{X_n} (t_n^m, X_m) \to X_m \mdvirg m > n
\geq 0
\]
soumis aux axiomes suivants : 
\begin{enumerate}
\item 
pour tout $(x,y)$ dans $(X_m,s_n^m)\times_{X_n}(t_n^m,X_m)$ avec $m>n\geq
0$, on a 
\begin{enumerate}
\item $s_m(x*_n^m y) = \begin{cases}
s_m(y) & \mdvirg n=m-1 \\
s_m(x) *_n^{m-1} s_n(y) & \mdvirg n<m-1
\end{cases}$ ;
\item $t_m(x*_n^m y) = \begin{cases}
t_m(x) & \mdvirg n=m-1 \\
t_m(x) *_n^{m-1} t_n(y) & \mdvirg n<m-1
\end{cases} \pbox{;}$ 
\end{enumerate}
\item pour tout triplet $(x,y,z)$ dans 
$(X_m, s^m_n)\times_{X_n}(t^m_n,X_m,s^m_n)\times_{X_n}(t^m_n,X_m)$, on
a 
\[
(x *^m_n y) *^m_n z = x *^m_n (y*^m_n z) ;
\]
\item pour tout quadruplet $(x,x',y,y')$ dans 
\[
(X_m,s^m_n) \times_{X_n} (t^m_n, X_m, s^m_l) 
\times_{X_l} (t^m_l , X_m, s^m_n) \times_{X_n} (t^m_n, X_m)
\]
avec $ m > n > l \geq 0 $, on a 
\[
(x *^m_n x') *^m_l (y *^m_n y') = (x *^m_l y) *^m_n (x' *^m_l y') ;
\]
\item pour tout $x$ dans $X_m$ et tout $n$ tel que $m>n \geq 0$, on a
\[
x *^m_n k^n_m s^m_n x = x = k^n_m t^m_n x *^m_n x ;
\]
\item pour tout $(x,y)$ dans $(X_m,
s^m_n)\times_{X_n}(t^m_n,X_m)$ avec $m > n\geq 0$, on a 
\[
k_m (x *^m_n y) = k_m(x) *^{m+1}_n k_m(y). 
\]
\end{enumerate}
Un \emph{$\omega$-foncteur} est un morphisme d'ensembles globulaires
réflexifs commutant aux opérations $*^m_n$ pour $m > n \geq 0$. On note $\wcat$ la catégorie des $\omega$-catégories. 
\end{defin}

\paragr Si $X$ est une $\omega$-catégorie, et si $m$ et $n$ sont deux
entiers tels que~$m>n\geq 0$, on dit qu'une $m$-cellule $x$ admet un
$*^m_n$-inverse s'il existe une $m$-cellule $y$ de $X$ vérifiant les
conditions suivantes : 
\begin{align*}
s^m_n(y)=t^m_n(x) &\mdvirg t^m_n(y)=s^m_n(x) \mdvirg\\
x *^m_n y = k^n_m t^m_n(x) &\mdvirg y *^m_n x = k^n_m s^m_n(x) \pbox{.}
\end{align*}
On vérifie qu'un tel inverse, s'il existe, est unique.

Un~\ndef{$\omega$-groupoïde} (sous-entendu, \emph{strict}) est une
$\omega$-catégorie telle que pour tous entiers $m$, $n$ tels que $m
>n\geq 0$, toute $m$-cellule admet un~$*^m_n$-inverse. On note $\wgpd$
la catégorie des $\omega$-groupoïdes.

\paragr\label{isowCatAbwGpdAbAbReflexif} Une
\ndef[$\omega$-catégorie!en groupes abéliens]{$\omega$-catégorie
en groupes abéliens} est une $\omega$-catégorie 
interne à~$\Ab$, ou, ce qui revient au même, un groupe abélien interne
à $\wcat$. Un~\ndef[$\omega$-groupoïde!en groupes
abéliens]{$\omega$\nobreakdash-groupoïde en groupes abéliens} est un
$\omega$-groupoïde interne à $\Ab$ ou, ce qui revient au même, une
$\omega$-catégorie en groupes abéliens dont la~$\omega$\nobreakdash-catégorie
sous-jacente est un $\omega$-groupoïde. On note $\wcat(\Ab)$ la
catégorie des $\omega$-catégories en groupes abéliens,
et $\wgpd(\Ab)$ celle des $\omega$-groupoïdes en groupes abéliens. 

\paragr Il s'avère que ces deux dernières notions sont en fait
dégénérées : si~$X$ est un groupe abélien globulaire réflexif,
on vérifie qu'en posant, pour tout couple
$(x,y)\in(X_m,s^m_n)\times_{X_n}(t^m_n,X_m)$, 
\[
x *^m_n y = x+y-k^n_m z \mdvirg 
\]
où on a noté
\[
  z = s^m_n(x) = t^m_n(y)\in X_n \mdvirg
  \]
alors on obtient bien de cette
manière une~$\omega$\nobreakdash-catégorie en groupes abéliens, et que
toute~$\omega$\nobreakdash-catégorie en groupes abéliens peut être
obtenue de
cette manière (voir par exemple \cite{milburn2012abelian}). De plus,
si~$C$ est une $\omega$-catégorie en groupes abéliens, alors pour~$m > n \geq 0$, et pour
toute~$m$-cellule $x$ de $X$, la $m$-cellule 
\[
y = k^n_m(s^m_n(x)+t^m_n(x))-x 
\]
est un $*^m_n$-inverse de $x$. En résumé, on a des isomorphismes de
catégories 
\[
\prefab{\Grefl} \simeq \wcat(\Ab) \simeq \wgpd(\Ab) \pbox{.}
\]

\medskip

On va à présent décrire une construction, due à Bourn
\cite{bourn1990denormalization}, établissant une équivalence de
catégories entre la catégorie des $\omega$-groupoïdes en groupes
abéliens et la catégorie $\Ch(\Ab)$, que nous interpréterons donc ici
comme un théorème de Dold-Kan pour la catégorie $\Grefl$ des globes
réflexifs. Nous montrerons par la suite que le complexe associé à un
groupe abélien globulaire réflexif par cette construction calcule la limite
inductive homotopique, établissant ainsi que la catégorie $\Grefl$ est
une catégorie pseudo-test homologique.

\paragr Si $X$ est un $\omega$-groupoïde en groupes abéliens, on lui
associe un complexe de chaînes de groupes abéliens $\bourne{X}$ en posant
$\bourne{X}_0 = X_0$, et, pour $n\geq 1$,
\notindex{$\bournef$}%
\begin{align*}
\bourne{X}_n &= X_n/k(X_{n-1}) \mdvirg \\
d_n &= t_n-s_n : \bourne{X}_n \to \bourne{X}_{n-1} \pbox{.}
\end{align*}
On peut vérifier grâce aux relations globulaires qu'on obtient bien de
cette manière un complexe de chaînes. 

Ainsi, si $X$ est un groupe
abélien globulaire réflexif, on peut lui associer le complexe obtenu
en appliquant le foncteur $\bournef$ au
$\omega$-groupoïde en groupes abéliens associé à $X$, décrit au
paragraphe précédent. On obtient ainsi un
foncteur encore noté
\[
  \bournef : \prefab{\Grefl}\xrightarrow{\simeq}\wgpd(\Ab)\to\Ch(\Ab) \pbox{.}
\]

\paragr Inversement, étant donné un complexe de chaînes de groupes abéliens
$C$, on peut lui associer un groupe abélien globulaire réflexif $X$ de
la manière suivante :  
\begin{itemize}
\item pour tout $n\geq0$, on pose $X_n = C_n \oplus \cdots \oplus C_1
\oplus C_0$; 
\item pour $n\geq 1$, les morphismes 
\[
s_n, t_n : C_n \oplus \cdots \oplus C_0 \to C_{n-1} \oplus \cdots
\oplus C_0  
\]
sont donnés respectivement par la projection, et par la somme de la
projection et de la différentielle de $C$;
\item pour $n\geq 0$, le morphisme
\[
  k_n : C_n \oplus \cdots \oplus C_0 \to C_{n+1} \oplus \cdots \oplus C_0 
  \]
est l'inclusion canonique.
\end{itemize}
On peut vérifier qu'on a bien défini de cette manière un préfaisceau
en groupes abéliens sur $\Grefl$. En utilisant l'isomorphisme décrit
au paragraphe \ref{isowCatAbwGpdAbAbReflexif}, on obtient donc un
$\omega$-groupoïde en groupes abéliens associé, dont la composition
est induite par la structure abélienne de $X$, et on a ainsi défini
un foncteur
\[
  \bournef^{-1} : \Ch(\Ab) \to \prefab{\Grefl} \xrightarrow{\simeq}
  \wgpd(\Ab) \pbox{.}
\]

\begin{theorem}[Bourn]\label{th:Bourn}
Le foncteur $\bournef: \wgpd(\Ab)\to\Ch(\Ab)$ est une équivalence
de catégories, dont un quasi-inverse est donné par le
foncteur~$\bournef^{-1} : \Ch(\Ab) \to \wgpd(\Ab)$.
\end{theorem}
\begin{proof}
Voir \cite[théorème 3.3]{bourn1990denormalization} pour la preuve
originale, ou encore~\cite[1.4.4]{ara2010infini} pour une preuve plus
directe.
\end{proof}
En particulier, le foncteur $\bournef$ induit une équivalence de
catégories entre la catégorie des groupes abéliens globulaires
réflexifs et celle des complexes de chaînes. On va maintenant voir que ce
dernier calcule bien la limite inductive homotopique, établissant
ainsi que la catégorie $\Grefl$ est une catégorie pseudo-test
homologique.
\begin{prop}
Le foncteur $c_{\Grefl} : \Grefl \to \Ch(\Ab)$, obtenu par restriction du foncteur
$\bournef$ aux préfaisceaux représentables, est un intégrateur sur~$\Grefl$.
\end{prop}
\begin{proof}
On peut facilement vérifier que pour tout groupe abélien globulaire
réflexif $X$, le complexe $\bournef(X)$ est isomorphe au sous-complexe
\[
X_0 \xleftarrow{t} \ker s_1 \xleftarrow{t} \ker{s_2} \leftarrow \cdots
\]
du complexe 
\[
X_0 \xleftarrow{t-s} X_1 \xleftarrow{t-s} X_2 \leftarrow \cdots \pbox{.}
\]
En particulier,
pour tout entier~$n\geq0$, $c_{\Grefl}(-)_n$ est un
facteur direct de $\Wh{\Grefl^{\op}}{\D_n}$, et est donc projectif.

Il reste donc à vérifier que, pour tout entier positif $n$, le
complexe~$c(\D_n)$ a le type d'homologie du point. Si $i$
est un entier tel que~$i<n$, alors $\D_n$ possède deux
$i$-cellules non-dégénérées (n'étant pas dans l'image du foncteur
$k_{i-1}$) que l'on note $x^0_i$ et $x^1_i$. À celles-ci s'ajoutent
les~$i$-cellules dégénérées de la forme $k_l^ix^\epsilon_l$ pour
$l<i$ et $\epsilon \in \lbrace 0,1 \rbrace$. Pour~$i=n$,~$\D_n$
possède une unique $n$-cellule non dégénérée $x^0_n=x^1_n$,
qui donne alors $2n+1$ cellules de dimension $n$ au total. Enfin, pour $i>n$,
les $i$-cellules sont toutes de la forme $k_n^ik_j^nx_j^\epsilon$.
Ainsi, le groupe abélien
des~$i$-cellules de $c_{\Grefl}(\D_n)$ est le quotient du groupe
abélien
\[
\Z^{(\D_i \to \D_n)} \simeq \begin{cases}
\Z^{(2i+2)} & \mdvirg i<n \\
\Z^{(2i+1)} & \mdvirg i\geq n \\
\end{cases} \pbox{.}
\]
par l'image du foncteur $k_{i-1}$. On obtient alors le complexe 
\[
{c_{\Grefl}}(\D_n) = 
\Z^{2} 
\xleftarrow{\begin{bsmallmatrix*}[r]
1 & 1 \\ -1 & -1
\end{bsmallmatrix*}} 
\Z^{2} 
\xleftarrow{\begin{bsmallmatrix*}[r]
1 & 1 \\ -1 & -1
\end{bsmallmatrix*}} 
\cdots \xleftarrow{}
\Z^{2} 
\xleftarrow{\begin{bsmallmatrix*}[r]
1 \\ -1
\end{bsmallmatrix*}} 
\Z 
\leftarrow 0 \leftarrow \cdots 
\]
dont l'homologie est bien celle du point.
\end{proof}

\begin{coro}
Si $X$ est un groupe abélien globulaire réflexif, on a un isomorphisme naturel
\[
\H{\Grefl}{X} \simeq \bourne{X} 
\]
dans $\Hotab$. En particulier, la catégorie $\Grefl$ est une catégorie
pseudo-test homologique faible.
\end{coro}
\begin{proof}
Comme énoncé dans la proposition \ref{propIntegrateurCocontinuExtKan},
il suffit de remarquer que puisque le foncteur $\bournef$ est
cocontinu, il coïncide avec l'extension de Kan de sa restriction à
$\Grefl$, dont on a montré qu'elle était un intégrateur. Le foncteur
$\bournef$ calcule donc bien la limite inductive homotopique. De plus,
on sait que $\bournef$ est une équivalence de catégories, et,
puisqu'elle préserve et reflète les équivalences faibles, elle induit
une équivalence de catégorie entre les catégories localisées.
\end{proof}

\begin{remark}
On pourrait être tentés d'essayer de calculer l'homologie des groupes
abéliens globulaires réflexifs de la même manière que dans le cas non
réflexif, c'est-à-dire en calculant l'homologie du complexe 
\[
X_0 \xleftarrow{t-s} X_1 \xleftarrow{t-s} X_2 \leftarrow \cdots 
\]
et s'attendre, comme dans le cas de $\Delta$, à ce que ce complexe \og
non normalisé \fg{} soit quasi-isomorphe au complexe $\bourne{X}$. On
peut rapidement s'apercevoir qu'il n'en est rien : par exemple,
puisque $\Grefl$ possède un objet
final, c'est une catégorie contractile, et elle a donc le type
d'homologie du point. Mais si $X$ est le préfaisceau
constant de valeur $\Z$, le complexe ci-dessus a le même type
d'homologie que $\H{\G}{\Z}$, dont tous les groupes d'homologie sont
isomorphes à $\Z$. Cela implique que le complexe ci-dessus ne calcule
pas l'homologie des groupes abéliens globulaires réflexifs.
\end{remark}

\paragr En vertu de la proposition \ref{propAdjIntegrateurs}, le
foncteur
\[
  \bournef : \prefab{\Grefl}\to\Ch(\Ab)
\]
admet pour adjoint à droite le foncteur  
\begin{align*}
c_{\Grefl}^* : \Ch(\Ab)\to\prefab{\Grefl}
\end{align*}
défini pour $C\in\Ch(\Ab)$ par
\[
c_{\Grefl}^*(C) : \D_i \mapsto \Hom(c_{\Grefl}(\D_i),C) \mdvirg
\]
où on rappelle qu'on note $c_{\Grefl}$ la restriction de $\bournef$ aux
préfaisceaux représentables. Le foncteur $c_{\Grefl}^*$ coïncide donc, à
isomorphisme naturel près, avec le foncteur $\bournef^{-1}$. On peut
alors conclure que $c_{\Grefl}^*$ préserve les équivalences faibles, et
on obtient une équivalence de catégories 
\[
\bournef : \Hotab_{\Grefl} \xrightarrow{\simeq} \Hotab \mdvirg
c_{\Grefl}^* : \Hotab
\xrightarrow{\simeq} \Hotab_{\Grefl} \mdvirg
\]
ce qui signifie qu'on a prouvé la proposition suivante : 
\begin{prop}
L'intégrateur $c_{\Grefl}$ est un intégrateur test faible, et la catégorie $\Grefl$ est test homologique faible.
\end{prop}

\paragr\label{debutPreuveGlobesPasWhitehead} On va maintenant montrer que $\Grefl$ n'est pas une catégorie de
Whitehead. Pour cela, considérons le complexe de chaînes 
\[
C = \Z \xleftarrow{\id} \Z \leftarrow 0 \leftarrow \cdots
\]
et notons $X=\bournef^{-1}(C)$ le groupe abélien globulaire réflexif
associé 
\[
X = \xymatrix{
\Z \ar@{^{(}->}[r]^{}  
& \Z\oplus\Z \ar@/_1pc/[l]_{\pr} \ar@/^1pc/[l]^{(\pr,\id)} \ar[r]^{\id} 
& \Z\oplus\Z \ar@/_1pc/[l]_{\id} \ar@/^1pc/[l]^{\id} \cdots
} \pbox{.}
\]
Puisque $C$ est un complexe exact, on sait que le morphisme $X\to0$
est dans~$\Wab_{\Grefl}$, et on va montrer que l'ensemble globulaire
$\U X$ a un type d'homotopie non
trivial. Pour cela on utilise le résultat suivant : 
\begin{prop}[Maltsiniotis]
Soit $X$ un ensemble globulaire réflexif, et $R$ la relation
d'équivalence sur $X_1$ engendrée par la relation 
\[
x \sim y \iff \exists u \in X_2 \mdvirg \text{$x = s(u)$ et $y = t(u)$} \pbox{.}
\]
Considérons alors le graphe $\tilde{X}$ défini par
\[
\tilde{X}=\xymatrix{
X_0 & \tilde{X_1} \ar@<0.5ex>[l]^-{s} \ar@<-0.5ex>[l]_-{t} 
}
\]
où $\tilde{X_1}$ désigne l'ensemble quotient de $X_1 - k(X_0)$
par la relation $R$. Alors on a un isomorphisme 
\[
\pi_1(\tranche{\Grefl}{X}) \simeq \pi_1(\tranche{\G_1}{\tilde{X}})
\]
et, en conséquence, le $1$-type d'homotopie de $X$ est celui du groupoïde
libre engendré par le graphe $\tilde{X}$.
\end{prop}
\begin{proof}
Voir \cite[commentaire 1, p 383]{maltsiniotisPS}. 
\end{proof}

\begin{prop}\label{globespaswhitehead} La catégorie $\Grefl$
n'est pas une catégorie de Whitehead. 
\end{prop}
\begin{proof}
On reprend l'exemple du groupe abélien globulaire réflexif~$X$
introduit au paragraphe \ref{debutPreuveGlobesPasWhitehead}. Toutes les
$2$-cellules de~$\U X$ étant des identités, il n'y a rien à identifier
dans le quotient définissant le graphe~$\U\tilde{X}$ de la
proposition précédente, qui se décrit donc de la manière suivante :
les sommets sont les entiers relatifs, et les flèches de $m$ vers $n$
sont les entiers relatifs $k$ tels que $m+k = n$.

Le groupoïde libre sur ce graphe possède un groupe d'automorphismes
non trivial, et la proposition précédente implique donc que
$\U X$ n'est pas simplement connexe. Par conséquent, le morphisme $X
\to 0$ de $\prefab{\Grefl}$ est dans $\W_{\Grefl}^\ab$, mais pas dans
$\U^{-1} \W_{\Grefl}$.
\end{proof}

\section{\texorpdfstring{La catégorie $\Theta$ de Joyal}%
                               {La catégorie ϴ de Joyal}}
                               \label{secTheta}
La catégorie $\Theta$, introduite par Joyal dans
\cite{joyal1997disks}, joue le même rôle relativement à la catégorie
$\wcat$ que $\Delta$ pour la catégorie $\Cat$ : on dispose d'un
foncteur nerf pleinement fidèle $\wcat \to \pref{\Theta}$ permettant
d'étudier la théorie de l'homotopie des $\omega$-catégories strictes,
et Cisinski et Maltsiniotis \cite{cisinski2011theta} ont montré que
$\Theta$ est une catégorie test stricte. De plus, la catégorie
$\Theta$ est intimement liée à diverses notions de $\omega$-catégories
faibles, notamment via les travaux de Rezk~\cite{rezk2010weak} et Ara~\cite{ara2014quasi}.
On montre dans cette section que $\Theta$ est une catégorie
test homologique stricte.

\paragr Soit $m$ un entier positif. Un \ndef{tableau de dimensions} de
\emph{largeur} $m$ est un tableau de la forme
\[
\left(
\begin{matrix}
i_1 && i_2 && i_3 && \cdots && i_m \\
& i'_1 && i'_2 && \cdots && i'_{m-1}
\end{matrix}
\right)
\]
dont les entrées sont des entiers positifs vérifiant 
\[
i_k > i'_k \mdvirg i_{k+1} > i'_k \mdvirg 1 \leq k < m \pbox{.}
\]
La \emph{hauteur} d'un tableau de dimensions est la plus grande de
ses entrées. 

\paragr On renvoie à la section \ref{secGlobes} pour la définition de
la catégorie $\G$.
Si
\[
T = \left(
\begin{matrix}
i_1 && i_2 && i_3 && \cdots && i_m \\
& i'_1 && i'_2 && \cdots && i'_{m-1}
\end{matrix}
\right)
\]
est un tableau de dimensions, la \ndef{somme globulaire} associée
à~$T$ est la somme amalgamée itérée
\[
(\D_{i_1},\sigma_{i'_1}^{i_1}) \amalg_{\D_{i'_1}} 
(\tau_{i'_1}^{i_2},\D_{i_2},\sigma_{i'_1}^{i_2}) \amalg_{\D_{i'_2}} 
\cdots 
\amalg_{\D_{i'_{m-1}}} (\tau_{i'_{m-1}}^{i_m}, \D_{i_m}) 
\]
dans $\pref{\G}$, c'est-à-dire la limite inductive du diagramme
d'ensembles globulaires 
\[
\xymatrix@R=.2pc@C=1pc{
\D_{i_1} && \D_{i_2} && \D_{i_3} && \D_{i_{m-1}} && \D_{i_m} \\
&&&&& \cdots \\
& 
\D_{i'_1} \ar[luu]^{\sigma^{i_1}_{i'_1}} \ar[ruu]_{\tau^{i_2}_{i'_1}}
&&
\D_{i'_2} \ar[luu]^{\sigma^{i_2}_{i'_2}} \ar[ruu]_{\tau^{i_3}_{i'_2}}
&& &&
\D_{i'_{m-1}}\ar[luu]^{\sigma^{i_{m-1}}_{i'_{m-1}}}
\ar[ruu]_{\tau^{i_{m}}_{i'_{m-1}}} & \pbox{.}
} 
\]
Par la suite, on notera cette somme globulaire simplement
par 
\[
\D_{i_1} \amalg_{\D_{i'_1}} \D_{i_2} \amalg_{\D_{i_2}} \cdots 
\amalg_{\D_{i_{m-1}}} \D_{i_m} \pbox{.}
\]

\begin{example}\label{ex3graphe}
Les objets de $\G$ correspondent aux tableaux de dimensions à une
seule entrée. Le $2$-graphe 
\[
\xymatrix{
\bullet \ar@/^1pc/[r]_{}="a" \ar@/_1pc/[r]^{}="b" & \bullet 
\ar@{=>}"a";"b"
\ar@/^1pc/[r]_{}="c" \ar@/_1pc/[r]^{}="d" & \bullet 
\ar@{=>}"c";"d"
}
\]
est la somme globulaire engendrée par le tableau de dimensions
\[
\left( \begin{matrix}
2 && 2 \\
& 0
\end{matrix} 
\right)
\]
tandis que pour le tableau 
\[
T = \left(
\begin{matrix}
2 && 2 && 1 && 3 && 2 \\
& 1 && 0  && 0 && 1 
\end{matrix}
\right) \mdvirg
\]
le $3$-graphe associé est 
\[
\xymatrix@C=4pc@H=5pc{
\bullet 
\ar@/^2.5pc/[r]_{}="f"
\ar[r]^{}="g"_{}="h"
\ar@/_2.5pc/[r]^{}="i"
\ar@2"f";"g"
\ar@/_1pc/@2"h";"i"_{}="0"
\ar@/^1pc/@2"h";"i"_{}="1"
\ar@3"0";"1"
&
\bullet 
\ar[r]
&
\bullet
\ar@/^1.5pc/[r]_{}="u"
\ar[r]^{}="v"_{}="w"
\ar@/_1.5pc/[r]^{}="x"
\ar@2"u";"v"
\ar@2"w";"x"
&
\bullet \pbox{.}
} 
\]
En revanche, les graphes
\[
  \xymatrix{
    \bullet & \bullet \ar[l]^{} \ar[r]^{} & \bullet
  }
  \hspace{2pc}\text{et} \hspace{2pc}
  \xymatrix@C=4em{
  \bullet \ar@/^1.5pc/[r]_{}="a" \ar[r]^{}="b"_{}="c" \ar@/_1.5pc/[r]^{}="d"
  & \bullet
  \ar@{=>}"b";"a"
  \ar@{=>}"c";"d"
  }
\]
ne sont pas des sommes globulaires.
\end{example}

\paragr Le foncteur d'oubli 
\[
\U : \wcat \to \pref{\G} \mdvirg
\]
associant à une $\omega$-catégorie stricte son ensemble globulaire sous-jacent,
admet un adjoint à gauche noté
\[
L : \pref{\G} \to \wcat 
\]
associant à un ensemble globulaire $X$ la \emph{$\omega$-catégorie
librement engendrée} par~$X$. 

\paragr La catégorie $\Theta$\notindex{$\Theta$} est définie comme la
sous-catégorie pleine de $\wcat$ dont les objets sont les
$\omega$-catégories librement engendrées par les sommes globulaires. 

\paragr Pour tout entier $n \geq 0$, on note $\Theta_n$ la
sous-catégorie pleine de~$\Theta$ dont les objets sont les
$\omega$-catégories libres sur les sommes globulaires de hauteur
inférieure ou égale à $n$. Les objets de $\Theta_n$ sont donc des
\emph{$n$-catégories}, c'est-à-dire des $\omega$-catégories dont
toutes les flèches en dimension strictement supérieure à $n$ sont des
identités. 

La catégorie $\Theta_0$ est la catégorie
terminale, et $\Theta_1$ est canoniquement isomorphe à~$\Delta$, comme
on peut le voir en représentant, pour tout entier $i\geq0$, l'objet
$\Delta_i$ de $\Delta$ par le tableau de dimension de longueur~$i$
\[
\left( 
\begin{matrix}
1 && 1 && \cdots && 1 && 1\\
& 0 && 0 && \cdots && 0
\end{matrix}
\right) \pbox{.}
\]

\begin{remark}
La catégorie $\Theta$ est isomorphe à la catégorie introduite par
Joyal dans \cite{joyal1997disks} comme catégorie opposée à celle des
disques combinatoires finis. L'équivalence de la définition donnée ici
et de celle de Joyal a été conjecturée par Batanin et Street dans
\cite{batanin2000multitude}, puis prouvée indépendamment par Makkai et
Zawadowski dans \cite{makkai2001duality} et par Berger dans
\cite{berger2002cellular}. Une définition alternative, que l'on va
exposer ci-dessous, est donnée par
Berger dans \cite{berger2007iterated}, que Cisinski et Maltsiniotis
utilisent dans \cite{cisinski2011theta} pour prouver que la
catégorie~$\Theta$ est une catégorie test stricte. Une preuve
alternative de ce résultat, utilisant une définition plus proche de
celle donnée ci-dessus, est dans la thèse de Ara~\cite{ara2010infini}.
\end{remark}

\paragr La combinatoire de $\Theta$ est assez sauvage. Pour illustrer
cela, on peut calculer numériquement que l'objet de $\Theta_2$
représentant la loi d'échange, c'est-à-dire la $2$-catégorie librement
engendrée par le $2$-graphe 
\[
\xymatrix@C=4pc@R=5pc{
\bullet \ar@/^1.5pc/[r]_{}="a" \ar[r]^{}="b"_{}="c" \ar@/_1.5pc/[r]^{}="d" & 
\bullet \ar@/^1.5pc/[r]^{}="u" \ar[r]^{}="v"_{}="w" \ar@/_1.5pc/[r]^{}="x" & 
\bullet
\ar@{=>}"a";"b"
\ar@{=>}"c";"d"
\ar@{=>}"u";"v"
\ar@{=>}"w";"x"
} 
\]
possède $169$ sous-objets (voir \cite{dimitriHtheta} pour une
bibliothèque Haskell permettant de jouer avec la combinatoire de
$\Theta$). Aussi, on n'introduira pas une orientation sur $\Theta_n$
pour $n > 1$. En revanche, on va voir dans la suite qu'on peut
calculer l'homologie des groupes abéliens cellulaires, c'est-à-dire
des préfaisceaux en groupes abéliens sur $\Theta$, grâce au produit en
couronne introduit par Berger dans
\cite{berger2007iterated}.

\paragr Si $A$ est une petite catégorie, on
définit une nouvelle catégorie notée~$\Delta \wr A$, appelée
\notindex{$\Delta\wr A$}%
\ndef{produit en couronne} de $\Delta$ et $A$, de la manière suivante.
Les objets de $\Delta\wr A$ sont les couples $[\Delta_n;(a_1, \cdots
a_n)]$ où les $a_i$ sont des objets de $A$. Un morphisme
\[
  [\Delta_n;(a_1,\cdots,a_n)] \to [\Delta_m;(a'_1,\cdots,a'_m)]
\]
est un couple $[\varphi, \mathbf{f}]$ où $\varphi
: \Delta_n \to \Delta_m$ est un morphisme de $\Delta$ et 
\[
    \mathbf{f}=(f_{ji}:a_i \to a'_j)
    _{
      1 \leq i \leq n,\,\varphi(i-1)<j\leq\varphi(i)
    }
\]
est un famille de morphismes de $A$. Étant données deux flèches
composables \[
  [\Delta_n;( a_1, \cdots,a_i)] \xrightarrow{[\varphi, \mathbf{f}]}
  [\Delta_{n'};( a'_1, \cdots,a'_j)] \xrightarrow{[\varphi', \mathbf{f'}]}
  [\Delta_{n''};( a''_1, \cdots,a''_k)] \mdvirg
\]
leur composée est définie par le couple
\[
(\varphi'', \mathbf{f''}) : [\Delta_n;(a_1,\cdots,a_n)] \to
[\Delta_{n''};(a''_1,\cdots,a''_{n''})]  \mdvirg
\]
avec $\varphi''=\varphi'\varphi$ et 
\[
\mathbf{f''}=(f''_{i''i})_{
    1 \leq i \leq n~,~
    \varphi''(i-1) < i'' \leq \varphi''(i) } \mdvirg
    f''_{i''i} = f'_{i''i'}f_{i'i}
\]
où $i'$ est l'unique entier tel que $\varphi(i-1) < i' \leq
\varphi(i)$ et $ \varphi'{(i'-1)} < i'' \leq \varphi'(i')$. 

\begin{remark}
Si $\varphi : \Delta_n \to \Delta_m$ est un morphisme constant de
$\Delta$, alors pour toutes familles $(a_i)_{1 \leq i \leq n}$ et
$(a'_j)_{1\leq j \leq m}$ d'objets de $A$, on dispose d'un
morphisme 
\[
\lbrack \varphi; ()\rbrack : \lbrack \Delta_n;(a_1,\dots,a_n)\rbrack \to
\lbrack \Delta_m;(a'_1,\dots,a'_m)\rbrack
\]
qui ne dépend pas du choix des familles $(a_i)$ et $(a'_j)$.
\end{remark}

\paragr Pour se représenter plus facilement les morphismes et leur
composition, on peut
représenter les objets de $\Delta\wr A$ comme
des graphes linéaires dont les arrêtes sont décorées par des
objets de $A$. L'objet $[\Delta_n;(a_1,\cdots,a_n)]$ est alors
représenté par le graphe
\[
0 \xrightarrow{a_1}
1 \xrightarrow{a_2}
2 \xrightarrow{a_3}
\cdots \xrightarrow{a_n}
n \mdvirg
\]
et les morphismes peuvent alors être représentés de manière graphique.
Par exemple, si $f_{12} : a_2 \to a'_1$ est un morphisme de $A$, le
morphisme
\[
  [\Delta_2;(a_1,a_2)]\xrightarrow{[s_0,\mathbf{f}]}
  [\Delta_{1},(a'_1)]
  \mdvirg \mathbf{f}=(f_{12})
\]
est représenté par
\[
\xymatrix{
0 \ar[r]^{a_1} \ar[d]^{} & 1 \ar[ld] \ar[r]^{a_2}_{}="a" & 2 \ar[ld]_{} 
\\
0 \ar[r]_{a'_1}^{}="b" & 1 
\ar@{-->}"a";"b"|{f_{12}} 
} 
\]
tandis qu'un morphisme 
\[
(\Delta_1, a_1) \xrightarrow{(d_1,\mathbf{f})} (\Delta_2, a'_1,a'_2) \mdvirg
\mathbf{f}=(f_{11},f_{21})
\]
est représenté par
\[
\xymatrix{
0 \ar[r]^{a_1}_{}="a" \ar[d]^{} & 1 \ar[rd]^{} 
\\
0 \ar[r]_{a'_1}^{}="b" & 1 \ar[r]_{a'_2}^{}="c" & 2
\ar@{-->}|{f_{11}}"a";"b"
\ar@{-->}|{f_{21}}"a";"c"
\pbox{.}} 
\]

\paragr\label{def:mu_A} La construction du produit en couronne est fonctorielle :
étant donné un foncteur $F : A\to B$ entre petites catégories, on
obtient un foncteur 
\[
\Delta\wr F : \Delta\wr A \to \Delta\wr B  
\]
défini de la manière suivante. Si $[\Delta_n;(a_1,\cdots,a_n)]$ est un
objet de $\Delta\wr A$, son image par le foncteur $\Delta\wr F$ est
l'objet 
\[
[\Delta_n;(F(a_1),\cdots,F(a_n)] 
\]
de $\Delta\wr B$. Si $[\varphi, \mathbf{f}=(f_{ji})_{i,j}]$ est un morphisme de
$\Delta\wr A$, on pose 
\[
(\Delta\wr F)[\varphi,\mathbf{f}] = [\varphi, F(\mathbf{f}) =
\left(F(f_{ji})\right)_{i,j}] 
\]
et on peut alors vérifier qu'on a bien défini de cette manière un
foncteur.

\paragr Si $F : A \to B$ est un foncteur fidèle (resp. pleinement fidèle), alors il en
est de même pour le foncteur $\Delta \wr F$. 
En particulier, si $a$ est un objet de~$A$, alors on dispose d'un
foncteur fidèle 
\[
I_a : \Delta \wr e \simeq \Delta \xrightarrow{\Delta \wr a} \Delta \wr A
\]
en notant $e$ la catégorie ponctuelle. Ce foncteur est pleinement fidèle si $a$
n'admet pas d'endomorphisme non trivial. Explicitement, pour tout
objet $\Delta_n$ de $\Delta$, on a 
\[
I_a(\Delta_n) = [\Delta_n;(\underbrace{a,a,\cdots,a}_{n\text{ fois}})] 
\]
et si $\varphi : \Delta_n \to \Delta_m$ est un morphisme de $\Delta$,
on a 
\[
I_a(\varphi) = [\Delta_n;(\underbrace{a,a,\cdots,a}_{n\text{ fois}})] \xrightarrow{[\varphi,
\mathbf{i}]} [\Delta_m;(\underbrace{a,a,\cdots,a}_{m\text{ fois}})] 
\]
où toutes les composantes de $\mathbf{i}$ sont l'identité de $a$.

Cette construction est également fonctorielle en $a$, et on obtient en
fait un foncteur qu'on note
\notindex{$\mu_A$}%
\[
\mu_A : \Delta \times A \to \Delta \wr A \mdvirg (\Delta_n, a) \mapsto
I_a(\Delta_n) \pbox{.}
\]

\paragr\label{notationCouronne} On utilisera les notations suivantes par
la suite. On note simplement~$\Delta_0$ l'objet $\lbrack
\Delta_0;()\rbrack$ de $\Delta\wr A$ et, si $a$ est un objet de~$A$, on
notera, pour tout entier $n > 0$,
\[
\Delta_n \wr a = \mu_A(\Delta_n, a) =
[\Delta_n;(\underbrace{a,a,\cdots,a}_{n\text{ fois}})] \pbox{.}
\]
\notindex{$\Delta_n\wr a$}%
Si $\varphi : \Delta_n \to \Delta_m$ est un morphisme constant de
$\Delta$, alors pour tous objets $a$ et~$a'$ de $A$, on note
\[
  \varphi \wr - : \Delta_n\wr a \to
\Delta_m\wr a'
  \]
le morphisme $\lbrack \varphi;()\rbrack$, qui ne dépend donc pas du
choix des objets $a$ et $a'$.
Enfin, si $\varphi : \Delta_n \to \Delta_m$ est un morphisme non constant de $\Delta$
et $f : a\to a'$ est un morphisme de~$A$, on notera
\[
\varphi \wr f = \mu_A(\varphi,f) : \Delta_n \wr a \to \Delta_m \wr a'
\pbox{.}
\]
\notindex{$\varphi\wr f$}%
\notindex{$\varphi\wr -$}%
On insiste sur le fait que, si $\varphi : \Delta_n \to \Delta_m$ est
un morphisme constant de $\Delta$, alors pour tous morphismes $f : a
\to b$ et $f' : a' \to b'$ de $A$, le diagramme 
\[
\xymatrix{
\Delta_n\wr a \ar[d]_{\id\wr f} \ar[r]^{\varphi \wr -} & \Delta_m \wr
a' \ar[d]^{\id \wr f'} \\
\Delta_n\wr b \ar[r]_{\varphi\wr -} & \Delta_m \wr b'
} 
\]
dans $\Delta\wr A$ est commutatif.

\paragr \label{defThetaCouronne} On peut maintenant donner la définition de $\Theta$ utilisant
le produit en couronne : puisque $\Delta_0$ est l'objet final de
$\Delta$, le foncteur $I_{\Delta_0}:\Delta \to \Delta\wr\Delta$ est pleinement
fidèle, et on le traitera comme une inclusion. On dispose alors d'une
séquence d'inclusions de petites catégories 
\[
\Delta \xhookrightarrow{I_{\Delta_0}} \Delta \wr \Delta 
\xhookrightarrow{\Delta\wr I_{\Delta_0}} \Delta \wr ( \Delta \wr \Delta )
\hookrightarrow \cdots
\]
et on peut alors redéfinir la catégorie $\Theta$ de la manière suivante
: 
\begin{itemize}
\item $\Theta_0 = e$ est la catégorie ponctuelle; 
\item $\Theta_{n+1} = \Delta \wr \Theta_n$ pour $n\geq 0$.
\end{itemize}
La catégorie $\Theta$ est alors obtenue comme la limite inductive du
diagramme de petites catégories 
\[
\Theta_0 \hookrightarrow \Theta_1 \hookrightarrow \Theta_2
\hookrightarrow \cdots \pbox{.}
\]

\medskip
On trouvera dans \cite{berger2007iterated} une preuve du fait que
cette définition de $\Theta$ est équivalente à celle que nous
avons exposée précédemment. On va montrer dans la suite comment cette
correspondance opère sur les objets.

\paragr Si $A$ est une petite catégorie, on peut représenter les
objets de~$\Delta\wr A$ par des \emph{arbres} planaires de hauteur au plus $1$ dont les
\emph{feuilles} sont décorées par des objets de $A$.
Si~$[\Delta_n;(a_1,\cdots,a_n)]$ est un objet de $\Delta\wr A$, on le
représente par l'arbre à $n$ \emph{branches}

\shorthandoff{;}
\[
\begin{xy}
(0,0)*=0{\bullet\annot{\Delta_n}}="0";
(-20,10)*=0{\bullet\annot{a_1}}="1";
(-10,10)*=0{\bullet\annot{a_2}}="2";
 (-0,10)*=0{\bullet\annot{a_3}}="3";
 (10,10)*=0{\cdots};
 (20,10)*=0{\bullet\annot{a_n}}="4";
\ar@{-}"0";"1"
\ar@{-}"0";"2"
\ar@{-}"0";"3"
\ar@{-}"0";"4"
\end{xy} \mdvirg
\]
\shorthandon{;}%
et en itérant cette construction à partir de $\Delta$, c'est-à-dire si
les $a_i$ sont des objets de $\Delta$, puis de $\Delta\wr\Delta$ et
ainsi de suite, on augmente progressivement la hauteur de l'arbre en
les concaténant.

\begin{example}\label{exArbreTheta}
Notons $T$ l'objet $[\Delta_3;(\Delta_1,\Delta_0,\Delta_3)]$ de
$\Delta\wr\Delta$. Alors $T$ est représenté par l'arbre
\shorthandoff{;}
\[
\begin{xy}
(0,0)*=0{\bullet}="12";
(-10,10)*=0{\bullet}="21";
(0,10)*=0{\bullet}="22";
(10,10)*=0{\bullet}="23";
(-20,20)*=0{\bullet}="31";
(0,20)*=0{\bullet}="32";
(10,20)*=0{\bullet}="33";
(20,20)*=0{\bullet}="34";
\ar@{-}"12";"21"
\ar@{-}"12";"22"
\ar@{-}"12";"23"
\ar@{-}"21";"31"
\ar@{-}"23";"32"
\ar@{-}"23";"33"
\ar@{-}"23";"34"
\end{xy} 
\]
\shorthandon{;}%
et l'objet $[\Delta_2;(\Delta_0,T)]$ de $\Delta\wr\Delta\wr\Delta$ est
alors représenté par l'arbre

\shorthandoff{;}
\[
\begin{xy}
(-10,-10)*=0{\bullet}="01";
(25,-10)*=0{.}="";
(-20,0)*=0{\bullet}="11";
(0,0)*=0{\bullet}="12";
(-10,10)*=0{\bullet}="21";
(0,10)*=0{\bullet}="22";
(10,10)*=0{\bullet}="23";
(-20,20)*=0{\bullet}="31";
(0,20)*=0{\bullet}="32";
(10,20)*=0{\bullet}="33";
(20,20)*=0{\bullet}="34";
\ar@{-}"12";"21"
\ar@{-}"12";"22"
\ar@{-}"12";"23"
\ar@{-}"21";"31"
\ar@{-}"23";"32"
\ar@{-}"23";"33"
\ar@{-}"23";"34"
\ar@{-}"01";"11"
\ar@{-}"01";"12" 
\end{xy}
\]
\shorthandon{;}%
où on a noté $\Delta_1$ l'objet $[\Delta_1;(\Delta_0)]$ de
$\Delta\wr\Delta$.
\end{example}

\paragr Si $T$ est un arbre planaire, on peut lui associer un tableau
de dimensions de la manière
suivante. On remarque que l'orientation du plan induit un ordre total
sur l'ensemble des feuilles de l'arbre 
(c'est-à-dire, sur les nœuds qui ne sont reliés à aucun nœud de hauteur
supérieure). On considère alors les feuilles 
\[
F_0 < F_1 < \cdots < F_m 
\]
et on note $i_k$ la hauteur de la feuille $F_k$ pour $0 \leq k \leq
m$. Si $F_k$ et $F_{k+1}$ sont deux feuilles, on note $i'_k$ la
hauteur du plus haut nœud en dessous de~$F_k$ et $F_{k+1}$. Autrement
dit, $i'_k$ est la hauteur du point d'intersection des deux chemins
reliant $F_k$ et $F_{k+1}$ à la racine.
On associe alors à $T$ le tableau de dimensions 
\[
\left(
\begin{matrix}
i_0 && i_1 && i_2 && \cdots && i_m \\
& i'_0 && i'_1 && \cdots && i'_{m-1}
\end{matrix}
\right) \pbox{.}
\]

\begin{example}
Le tableau de dimensions associé à l'arbre

\shorthandoff{;}
\[
\begin{xy}
(-10,-10)*=0{\bullet}="01";
(-20,0)*=0{\bullet}="11";
(0,0)*=0{\bullet}="12";
(-10,10)*=0{\bullet}="21";
(0,10)*=0{\bullet}="22";
(10,10)*=0{\bullet}="23";
(-20,20)*=0{\bullet}="31";
(0,20)*=0{\bullet}="32";
(10,20)*=0{\bullet}="33";
(20,20)*=0{\bullet}="34";
\ar@{-}"12";"21"
\ar@{-}"12";"22"
\ar@{-}"12";"23"
\ar@{-}"21";"31"
\ar@{-}"23";"32"
\ar@{-}"23";"33"
\ar@{-}"23";"34"
\ar@{-}"01";"11"
\ar@{-}"01";"12" 
\end{xy}
\]
\shorthandon{;}%
de l'exemple
\ref{exArbreTheta} est 
\[
\left(
\begin{matrix}
1 && 3 && 2 && 3 && 3 && 3 \\
& 0 && 1 && 1 && 2 && 2
\end{matrix}
\right)
\]
et l'arbre représentant le $3$-graphe 
\[
\xymatrix@C=4pc@H=5pc{
\bullet 
\ar@/^2.5pc/[r]_{}="f"
\ar[r]^{}="g"_{}="h"
\ar@/_2.5pc/[r]^{}="i"
\ar@2"f";"g"
\ar@/_1pc/@2"h";"i"_{}="0"
\ar@/^1pc/@2"h";"i"_{}="1"
\ar@3"0";"1"
&
\bullet 
\ar[r]
&
\bullet
\ar@/^1.5pc/[r]_{}="u"
\ar[r]^{}="v"_{}="w"
\ar@/_1.5pc/[r]^{}="x"
\ar@2"u";"v"
\ar@2"w";"x"
&
\bullet \pbox{.}
} 
\]
de l'exemple \ref{ex3graphe} est le
suivant :
\shorthandoff{;}
\[
\begin{xy}
(0,0)*=0{\bullet}="01";
(-10,10)*=0{\bullet}="11";
(0,10)*=0{\bullet}="12";
(10,10)*=0{\bullet}="13";
(-15,15)*=0{\bullet}="21";
(-5,15)*=0{\bullet}="22";
(5,15)*=0{\bullet}="23";
(15,15)*=0{\bullet}="24";
(0,20)*=0{\bullet}="31";
\ar@{-}"01";"11"
\ar@{-}"01";"12"
\ar@{-}"01";"13"
\ar@{-}"11";"21"
\ar@{-}"11";"22"
\ar@{-}"13";"23"
\ar@{-}"13";"24"
\ar@{-}"23";"31"
\end{xy} \pbox{.}
\]
\shorthandon{;}%
\end{example}

Maintenant que l'on a donné les deux définitions de $\Theta$ qui nous
serviront dans la suite, on va s'intéresser aux préfaisceaux sur $\Theta$.

\paragr On appelle \ndef[ensemble!cellulaire]{ensembles cellulaires}
(resp. \ndef[groupe abélien!cellulaire]{groupes
abéliens cellulaires}) les préfaisceaux d'ensembles (resp. de groupes
abéliens) sur $\Theta$. 

\begin{prop}[Cisinski-Maltsiniotis]\label{ThetaTestStricte}
Pour tout entier $n\geq1$, la catégorie $\Theta_n$ est une catégorie
test stricte. De plus, la catégorie $\Theta$ est aussi une catégorie
test stricte.
\end{prop}
\begin{proof}
Voir \cite{cisinski2011theta}.
\end{proof}

Dans ce qui suit, on va montrer que les catégories $\Theta_n$, pour
$n\geq1$, sont des catégories test homologiques strictes.

\begin{prop}[Ara-Maltsiniotis]
Soit $A$ une catégorie totalement asphérique (voir
\ref{TotAspheriqueEqDef}). Alors le foncteur
\[
\mu_A : \Delta \times A \to  \Delta \wr A
\]
est un foncteur asphérique.
\end{prop}
\begin{proof}
Voir \cite[proposition 5.4]{ara2022compnerfs}.
\end{proof}

\paragr\label{def:m_n} On définit, pour tout entier $n$, un foncteur
\begin{align*}
m_n : \Delta^n &\to \Theta_n 
\end{align*}
par récurrence de la manière suivante : 
\begin{itemize}
\item $m_0$ est l'identité de $e=\Delta^0=\Theta_0$;

\item $m_{n+1}$ est la composée
\[
\Delta^{n+1} = \Delta \times \Delta^n \xrightarrow{\Delta\times m_n}
\Delta\times\Theta_n \xrightarrow{\mu_{\Theta_n}} \Theta_{n+1} \pbox{.}
\]
\end{itemize}
En utilisant la notation introduite en \ref{notationCouronne}, cela
signifie qu'on a 
\[
m_n(\Delta_{i_0},\Delta_{i_1},\cdots,\Delta_{i_n}) 
= \Delta_{i_0}\wr\Delta_{i_1}\wr\cdots\wr\Delta_{i_n}
\]
où l'opération est parenthésée à droite. On fera toutefois attention à la
subtilité suivante : si $i_k=0$ pour un certain $k\leq n$, alors on a
\[
m_n(\Delta_{i_0},\Delta_{i_1},\cdots,\Delta_{i_n}) = 
m_k(\Delta_{i_0},\Delta_{i_1},\cdots,\Delta_{i_k}) = 
\Delta_{i_0}\wr\Delta_{i_1}\wr\cdots\wr\Delta_{i_{k}} \pbox{.}
\]
Par exemple, on a 
\shorthandoff{;}
\begin{align*}
m_2(\Delta_3,\Delta_2) = \Delta_3 \wr \Delta_2&=
\xymatrix{
\bullet 
\ar@/^1.5pc/[r]|{\vphantom{X}}="a"
\ar[r]|{\vphantom{X}}="b"
\ar@/_1.5pc/[r]|{\vphantom{X}}="c"
\ar@2"a";"b"
\ar@2"b";"c"
&
\bullet
\ar@/^1.5pc/[r]|{\vphantom{X}}="a1"
\ar[r]|{\vphantom{X}}="b1"
\ar@/_1.5pc/[r]|{\vphantom{X}}="c1"
\ar@2"a1";"b1"
\ar@2"b1";"c1"
&
\bullet
\ar@/^1.5pc/[r]|{\vphantom{X}}="a2"
\ar[r]|{\vphantom{X}}="b2"
\ar@/_1.5pc/[r]|{\vphantom{X}}="c2"
\ar@2"a2";"b2"
\ar@2"b2";"c2"
&
\bullet
}
\\
m_3(\Delta_2,\Delta_1,\Delta_1) = \Delta_2\wr\Delta_1\wr\Delta_1&=
\xymatrix@=4em{
\bullet 
\ar@/^1.5pc/[r]_{}="a"
\ar@/_1.5pc/[r]^{}="b"
\ar@/^1pc/@2"a";"b"_{}="u"
\ar@/_1pc/@2"a";"b"^{}="v"
\ar@3"v";"u"
&
\bullet
\ar@/^1.5pc/[r]_{}="a1"
\ar@/_1.5pc/[r]^{}="b1"
\ar@/^1pc/@2"a1";"b1"_{}="u1"
\ar@/_1pc/@2"a1";"b1"^{}="v1"
\ar@3"v1";"u1"
&
\bullet
}
\\
m_3(\Delta_2,\Delta_0,\Delta_i) = \Delta_2 &= 
\xymatrix@C=3em{
\bullet \ar[r]^{} & \bullet \ar[r]^{} & \bullet \mdvirg {i \in \N \pbox{.}}
} 
\end{align*}
\shorthandon{;}%

\begin{coro}[Ara-Maltsiniotis]
Pour tout entier $n\geq 0$, le foncteur 
\[
m_n : \Delta^n \to \Theta_n 
\]
est asphérique.
\end{coro}
\begin{proof}
Voir \cite[corollaire 5.5]{ara2022compnerfs}. 
\end{proof}

Puisque la catégorie $\Delta$ est totalement asphérique, le foncteur
diagonal~$\Delta\to\Delta^n$ est asphérique pour tout entier $n\geq
1$, et on obtient finalement un foncteur asphérique
\begin{align*}
\Delta &\xrightarrow{\diag} \Delta^n \xrightarrow{m_n} \Theta_n \\
\Delta_i &\mapsto \Delta_i \wr \Delta_i \wr \cdots \wr \Delta_i \pbox{.}
\end{align*}
\begin{coro}
Pour tout entier $n\geq 1$, la catégorie $\Theta_n$ est une catégorie test
homologique stricte. De plus, $\Theta_n$ est également une catégorie de
Whitehead.
\end{coro}
\begin{proof}
C'est un cas particulier de la proposition
\ref{propMorphismeAsphSourceTHS}, puisque la catégorie $\Theta_n$ est une catégorie
test et que $m_n : \Delta \to \Theta_n$ est un foncteur asphérique.
\end{proof}

\paragr En vertu de la proposition \ref{coro:integrateurInduit}, on
peut définir un intégrateur sur $\Theta_n$ pour $n\geq 1$ comme la
composée 
\[
L_{\Theta_n} : \Theta_n \xrightarrow{\mathsf{Wh}_{{\Theta}_n}} \prefab{\Theta_n}
\xrightarrow{m_n^*} \prefab{\Delta} \xrightarrow{{L_{\Delta}}_!}
\Ch(\Ab) \mdvirg
\]
\notindex{$L_{\Theta_n}$}%
en notant toujours ${L_\Delta}_!$ le foncteur complexe de chaînes non
normalisé (on peut aussi utiliser le complexe normalisé). 
Si $X$ est un préfaisceau abélien sur~${\Theta_n}$, alors le
complexe ${L_{\Theta_n}}_!X$ est le complexe non normalisé associé
au groupe abélien simplicial
\[
m_n^*\diag^*(X) : {\Delta}^{\op}\to\Ab \mdvirg \Delta_i \mapsto
X(m_n\diag(\Delta_i)) \pbox{.}
\]
Explicitement, on a un isomorphisme naturel
\[
{L_{\Theta_n}}_!X \simeq
X(\Delta_0) 
\leftarrow
X(\Delta_1\wr\cdots\wr\Delta_1) 
\leftarrow
X(\Delta_2\wr\cdots\wr\Delta_2) 
\leftarrow
\cdots
\]
dans $\Ch(\Ab)$,
où la différentielle est donnée au rang $k>1$ par
\[
d_k = \sum_{i=0}^k (-1)^i
X(\delta_i\wr\delta_i\wr\cdots\wr\delta_i) 
\]
en notant $\delta_i = \delta^k_i: \Delta_{k-1} \hookrightarrow \Delta_k$ la
$i$-ème coface de $\Delta_k$ dans
$\Delta$. Au rang $1$, on trouve
\[
d_1 = X(\delta_0) - X(\delta_1)
\]
en notant $\delta_i$ le morphisme $\Delta_0 \to \Delta_1 \wr \Delta_1 \wr \cdots
\wr \Delta_1$ de $\Theta_n$ défini par $\delta_i : \Delta_0 \to
\Delta_1$ pour~$i=0,1$.
\begin{coro}
Pour tout entier $n\geq 1$ et pour tout préfaisceau en groupes abélien $X$ sur
$\Theta_n$, on a un isomorphisme 
\[
\H{\Theta_n}{X} \simeq {L_{\Theta_n}}_! X
\]
dans $\Hotab$.
\end{coro}

\medskip
En revanche, on ne dispose pas d'un foncteur
asphérique de $\Delta$ dans $\Theta$. Pour montrer que $\Theta$ est
une catégorie test homologique stricte, on va utiliser la notion de
foncteur asphérique introduite au paragraphe
\ref{defFoncteurAspheriqueÂ}.

\paragr\label{paragrCompNerfs} L'inclusion $\Theta \hookrightarrow \wcat$
permet de définir un foncteur pleinement fidèle appelé \emph{nerf cellulaire} 
\[
\nerf_{\Theta} : \wcat \to \pref{\Theta}
\]
permettant d'étudier la théorie de l'homotopie des $\omega$-catégories
strictes. Par ailleurs, Street introduit dans \cite{street1987algebra}
un foncteur 
\[
\orient : \Delta \to \wcat \mdvirg \Delta_n \mapsto \orient_n 
\]
appelé \emph{objet cosimplicial des orientaux}. On renvoie à l'article
de Street pour la définition des $\omega$-catégories $\orient_n$,
ainsi qu'à \cite{ara2018thmAsimp} pour une description utilisant la
théorie de Steiner~\cite{steiner2004omegacat}. Ce foncteur permet à
son tour de définir un foncteur appelé \ndef[nerf!de Street]{nerf de Street} 
\[
\nerf : \wcat \to \pref{\Delta} 
\]
par le même procédé. Explicitement, si $C$ est une $\omega$-catégorie,
alors $\nerf_{\Theta}(C)$ est l'ensemble cellulaire 
\[
{\Theta}^{\op} \to \Ens \mdvirg T \mapsto \Hom_{\wcat}(T,C) 
\]
et $\nerf(C)$ est l'ensemble simplicial 
\[
{\Delta}^{\op} \to \Ens \mdvirg \Delta_n \mapsto
\Hom_{\wcat}(\orient_n, C) \pbox{.}
\]

Dans \cite{ara2022compnerfs}, Ara et Maltsiniotis prouvent que ces
deux foncteurs nerfs sont équivalents, au sens où le carré 
\[
\xymatrix{
\wcat \ar[r]^{\nerf_\Theta} \ar[d]_{\nerf} & \pref{\Theta}
\ar[d]^{i_{\Theta}} \\
\pref{\Delta} \ar[r]_{i_{\Delta}} & \Hot
} 
\]
est commutatif à isomorphisme près. Une des étapes de
leur démonstration est la suivante : 

\begin{prop}[Ara-Maltsiniotis]
Le foncteur 
\[
j : \Delta \xrightarrow{\orient} \wcat \xrightarrow{N_\Theta} \pref{\Theta} 
\]
est un foncteur asphérique au sens du paragraphe \ref{defFoncteurAspheriqueÂ}.
\end{prop}
\begin{proof}
Voir \cite[théorème~5.8]{ara2022compnerfs}. On rappelle que cela
signifie que \begin{itemize}
\item pour tout objet $T$ de $\Theta$, le préfaisceau $j^*(T)$ est asphérique;
\item pour tout objet $\Delta_n$ de $\Delta$, les préfaisceaux
$j(\Delta_n)$ et $j^*j(\Delta_n)$ est asphérique. \qedhere
\end{itemize}
\end{proof}

\begin{coro}
La catégorie $\Theta$ est une catégorie test homologique stricte de Whitehead.
\end{coro}
\begin{proof}
Puisque $\Theta$ est une catégorie test stricte (voir \ref{ThetaTestStricte}),
il suffit donc d'utiliser la proposition
\ref{propMorphismeAsphSourceTHS}, appliquée au foncteur asphérique $j
: \Delta \to \pref{\Theta}$ ci-dessus.
\end{proof}

\section{\texorpdfstring{La catégorie $\Xi$ }%
                               {La catégorie Ξ}}
\label{secXi}
Dans cette section, on introduit une sous-catégorie notée $\Xi_n$ de $\Theta_n$
définie comme l'image de $\Delta^n$ par le produit en couronne. Cette
catégorie est isomorphe à la catégorie introduite par Simpson (qu'il
note $\Theta_n$) dans
\cite{simpson1997closed}. On va
montrer que $\Xi_n$ est une catégorie test stricte, ainsi qu'une catégorie
test homologique stricte. De plus, on donne une orientation sur
$\Xi_n$ et on montre qu'elle est asphérique au sens du paragraphe
\ref{def:orientationAspherique}, donnant alors une formule
explicite pour calculer l'homologie des préfaisceaux abéliens sur
$\Xi_n$ de manière intrinsèque.

\paragr 
\notindex{$\Delta\xiwr A$}%
Pour toute petite catégorie $A$, on note $\Delta\xiwr A$ la
sous-catégorie de~$\Delta\wr A$
engendrée par l'image du foncteur $\mu_A$ introduit au paragraphe
\ref{def:mu_A}.

\begin{prop}
Pour toute petite catégorie $A$, la catégorie $\Delta\xiwr A$ peut
être décrite de la manière suivante : les objets de $\Delta\xiwr A$
sont $\Delta_0$ et $\Delta_n\wr a$ pour tout $n> 0$ et tout objet
$a$ de $A$. La classe des morphismes de $\mu_A(\Delta_i,a)$
vers~$\mu_A(\Delta_j,a')$ est l'union de :
\begin{itemize}
\item la classe des morphismes de la forme $\varphi \wr \psi$ où
$\varphi : \Delta_i \to \Delta_j$ est un morphisme non constant et
$\psi : a \to a'$ est un morphisme de $A$;
\item la classe des morphismes de la forme $\varphi\wr -$ où $\varphi :
\Delta_i \to \Delta_j$ est un morphisme constant de $\Delta$.
\end{itemize}
\end{prop}
\begin{proof}
On vérifie facilement que ces classes de morphismes sont stables par
composition.
De plus, si $C$ est une sous-catégorie de $\Delta\wr A$ contenant l'image du
foncteur $\mu_a$, on montre que $C$ contient les morphismes décris
ci-dessus : en
effet, si~$(\varphi,\psi) : (\Delta_n,a) \to (\Delta_m,a')$ est un
morphisme de $\Delta\times A$ tel que $\varphi$ n'est pas un morphisme
constant de $\Delta$, alors $C$ contient le
morphisme~$\varphi\wr\psi$. De plus, pour tout morphisme constant
$\varphi : \Delta_n \xrightarrow{r} \Delta_0 \xrightarrow{s} \Delta_m$
de $\Delta$, on obtient, pour tous objets $a$ et $a'$ de $A$, une
chaîne de morphismes composables dans l'image de $\Delta\wr A$
\[
\mu_A(\Delta_n, a) \xrightarrow{\mu_A(r,\id_a)}
\mu_A(\Delta_0,a)=\mu_A(\Delta_0,a') \xrightarrow{\mu_A(s,\id_{a'})}
\mu_A(\Delta_m, a')
\]
dont la composée $\varphi \wr - : \Delta_n\wr a \to \Delta_m \wr a'$
est donc un morphisme de $C$. 
\end{proof}

\begin{remark}\label{remark:imagemu_Astable}
Si $\varphi : \Delta_n
\to \Delta_m$ est un morphisme constant de $\Delta$ et si~$u : a \to
a'$ est un morphisme de $A$, alors on a $\varphi \wr
- = \mu_A(\varphi,u)$. En particulier, si $A$ est une catégorie telle que pour tous
objets $a$ et $a'$, l'ensemble~$\Hom_A(a,a')$ est non vide, alors
$\Delta\xiwr A$ est l'image du foncteur $\mu_A$, qui est donc stable par
composition. C'est le cas par exemple si $A=\Delta$, ou si
$A=\Delta\wr A'$ pour toute petite catégorie $A'$.
\end{remark}

\paragr\label{secmu'} Pour toute petite catégorie $A$, on note
\begin{align*}
\ximuf_A : \Delta\times A &\to \Delta\xiwr A 
\\(\Delta_n,a) &\mapsto [\Delta_n;(a,\cdots,a)] 
\end{align*}
\notindex{$\ximuf_A$}%
le foncteur induit par le foncteur
$\mu$. On vérifie alors que si $F : A \to B$ est un foncteur fidèle
(resp. pleinement fidèle), alors le foncteur 
\[
\Delta \xiwr F : \Delta \xiwr A \to \Delta \xiwr B 
\]
est également un foncteur fidèle (resp. pleinement fidèle).

\paragr Pour tout entier $n \geq 0$, on définit alors par récurrence
une catégorie notée $\Xi_n$ de la manière suivante. Puisque $\Delta_0$
est l'objet final de $\Delta$, le foncteur~$\Delta\xiwr\Delta_0 : \Delta \to
\Delta\xiwr\Delta$ (en notant $\Delta_0 : e \to \Delta$ le foncteur
constant de valeur $\Delta_0$) est pleinement fidèle. On dispose alors d'une
chaîne d'inclusions 
\[
\Delta \xhookrightarrow{\Delta\xiwr\Delta_0} \Delta\xiwr\Delta
\xhookrightarrow{\Delta\xiwr(\Delta\xiwr\Delta_0)} \Delta\xiwr(\Delta\xiwr\Delta)
\hookrightarrow \cdots 
\]
\notindex{$\Xi_n$}%
et on définit alors la catégorie $\Xi_n$ en posant :
\begin{itemize}
\item $\Xi_0 = e$ est la catégorie ponctuelle;
\item $\Xi_{n+1} = \Delta \xiwr \Xi_n$ pour $n\geq 0$.
\end{itemize}

En d'autres termes, $\Xi_n$ est la sous-catégorie de $\Theta_n$
obtenue comme l'image (voir la remarque \ref{remark:imagemu_Astable})
du foncteur $m_n : \Delta^n \to \Theta_n$ introduit au
paragraphe~\ref{def:m_n}.

\paragr\label{defxim} Pour tout entier $n\geq 0$, on note alors 
\[
\xim_n : \Delta^n \to \Xi_n 
\]
le foncteur défini par récurrence à partir de $\ximuf_\Delta$,
c'est-à-dire qu'on pose : \begin{itemize}
\item $\xim_0 : e \to e$ est le foncteur identité;
\item $\xim_{n+1} : \Delta^{n+1}\to\Xi_{n+1}$ est la composée 
\[
\Delta^{n+1} = \Delta \times \Delta^{n} \xrightarrow{\Delta\times \xim_n}
\Delta\times\Xi_n \xrightarrow{\ximuf_{\Xi_n}} \Xi_{n+1} \pbox{.}
\]
\end{itemize}

\paragr La catégorie $\Xi_n$ admet la description explicite suivante.
On rappelle que si~$(\Delta_{i_1},\dots,\Delta_{i_n})$ est un objet de
$\Delta^n$ et si $1 \leq k \leq n-1$ est un entier tel que~$i_{k+1}=0$,
alors on a 
\[
\xim_n(\Delta_{i_1},\dots,\Delta_{i_n}) =
\xim_k(\Delta_{i_1},\dots,\Delta_{i_{k}}) \pbox{.}
\]
Ainsi, chaque objet $S \neq \Delta_0$ de $\Xi_n$ peut s'écrire
\emph{de manière unique}, en utilisant la notation introduite en
\ref{notationCouronne}, comme un produit en couronnes de la forme 
\[
\Delta_{i_1}\wr\Delta_{i_2}\wr\cdots\wr\Delta_{i_k} 
\]
où $1 \leq k\leq n$ et $i_l$ est un entier non nul pour $1 \leq l \leq
k$. On appellera \emph{forme normale} des objets de $\Xi_n$ une
telle écriture en produit en couronnes. On précise que la forme
normale de $\Delta_0$ vu comme objet de $\Xi_n$ est $\Delta_0$.
Si
  \[
S = \xim_n(\Delta_{i_1},\Delta_{i_2},\dots,\Delta_{i_n})
  \]
et~
  \[
T = \xim_n(\Delta_{j_1},\Delta_{j_2},\dots,\Delta_{j_n})
  \]
sont deux objets de $\Xi_n$, un morphisme $S\to T$ de $\Xi_n$ est un
morphisme de la forme 
\[
u = \xim_n(u_1,u_2,\dots,u_n) 
\]
où $u_l : \Delta_{i_l}\to\Delta_{j_l}$ est un morphisme de $\Delta$
pour $1 \leq l \leq n$. On rappelle à nouveau que si~$u_{k} :
\Delta_{i_{k}} \to \Delta_{j_{k}}$ est un morphisme constant pour
un certain~$1 \leq k \leq n$, alors on a 
\[
\xim_n(u_1,u_2,\dots, u_n) = \xim_n(u_1,u_2,\dots,u_{k},v_{k+1},\dots,v_n) 
\]
pour tous morphismes $v_{k+1},\dots,v_n$ de $\Delta$. 
Les morphismes $S \to T$ de $\Xi_n$ correspondent donc exactement aux
familles de morphismes~$(u_l : \Delta_{i_l} \to \Delta_{j_l})_{1\leq
l\leq k}$ de morphismes de $\Delta$ avec $k\leq n$, tels que~$u_l$ ne
se factorise pas par $\Delta_0$ pour $1\leq l < k$. On pourra
alors écrire $u$ sous la forme normale
\[
\Delta_{i_1} \wr \Delta_{i_2} \wr \cdots \wr \Delta_{i_S}
\xrightarrow{u_1\wr u_2\wr\cdots\wr u_k}
\Delta_{j_1} \wr \Delta_{j_2} \wr \cdots \wr \Delta_{j_T} 
\]
où $i_S$ et $j_T$ sont des entiers inférieurs à $n$.

\paragr 
On note $\Xi$ la limite inductive du diagramme de petite catégories 
\[
\Xi_0 \hookrightarrow 
\Xi_1 \hookrightarrow 
\Xi_2 \hookrightarrow 
\Xi_3 \hookrightarrow 
\cdots
\]
\notindex{$\Xi$}%
et, pour tout entier $n\geq 0$, on note $j_n$ le morphisme canonique
$\Xi_n \to \Theta_n$. On note également $j : \Xi
\to \Theta$ le morphisme induit par passage à la limite inductive.

\begin{example}
Les objets de $\Xi_2$ sont de la forme 
\shorthandoff{;}
\[
\resizebox{8cm}{!}{
\xymatrix@=5em{
\bullet 
\ar@/^3pc/[r]|{\vphantom{X}}="a"
\ar@/^2pc/[r]|{\vphantom{X}}="b"
\ar@/^1pc/[r]|{\vphantom{X}}="c"
\ar@/_1pc/[r]|{\vphantom{X}}="d"
\ar@/_2pc/[r]|{\vphantom{X}}="e"
\ar@/_3pc/[r]|{\vphantom{X}}="f"
\ar@{}"c";"d"|<<<{\vdots}
\ar@2"a";"b"
\ar@2"b";"c"
\ar@2"d";"e"
\ar@2"e";"f"
&
\bullet
\ar@/^3pc/[r]|{\vphantom{X}}="a1"
\ar@/^2pc/[r]|{\vphantom{X}}="b1"
\ar@/^1pc/[r]|{\vphantom{X}}="c1"
\ar@/_1pc/[r]|{\vphantom{X}}="d1"
\ar@/_2pc/[r]|{\vphantom{X}}="e1"
\ar@/_3pc/[r]|{\vphantom{X}}="f1"
\ar@{}"c1";"d1"|<<<{\vdots}
\ar@2"a1";"b1"
\ar@2"b1";"c1"
\ar@2"d1";"e1"
\ar@2"e1";"f1"
&
\cdots
\ar@/^3pc/[r]|{\vphantom{X}}="a2"
\ar@/^2pc/[r]|{\vphantom{X}}="b2"
\ar@/^1pc/[r]|{\vphantom{X}}="c2"
\ar@/_1pc/[r]|{\vphantom{X}}="d2"
\ar@/_2pc/[r]|{\vphantom{X}}="e2"
\ar@/_3pc/[r]|{\vphantom{X}}="f2"
\ar@{}"c2";"d2"|<<<{\vdots}
\ar@2"a2";"b2"
\ar@2"b2";"c2"
\ar@2"d2";"e2"
\ar@2"e2";"f2"
&
\bullet 
}
} \pbox{,}
\]
\shorthandon{;}%
et voici l'objet $\Delta_3\wr\Delta_2\wr\Delta_1$ de $\Xi_3$ :
\shorthandoff{;}%
\[
T = \xymatrix@=5em{
  \bullet 
    \ar@/^2.5pc/[r]_{}="a"
    \ar[r]^{}="b"_{}="c"
    \ar@/_2.5pc/[r]^{}="d"
      \ar@/^1pc/@2"a";"b"_{}="u"
      \ar@/_1pc/@2"a";"b"^{}="v"
      \ar@/_1pc/@2"c";"d"^{}="w"
      \ar@/^1pc/@2"c";"d"_{}="x"
        \ar@3"v";"u"
        \ar@3"w";"x"
    &
  \bullet 
    \ar@/^2.5pc/[r]_{}="a1"
    \ar[r]^{}="b1"_{}="c1"
    \ar@/_2.5pc/[r]^{}="d1"
      \ar@/^1pc/@2"a1";"b1"_{}="u1"
      \ar@/_1pc/@2"a1";"b1"^{}="v1"
      \ar@/_1pc/@2"c1";"d1"^{}="w1"
      \ar@/^1pc/@2"c1";"d1"_{}="x1"
        \ar@3"v1";"u1"
        \ar@3"w1";"x1"
&
    \bullet
    \ar@/^2.5pc/[r]_{}="a2"
    \ar[r]^{}="b2"_{}="c2"
    \ar@/_2.5pc/[r]^{}="d2"
      \ar@/^1pc/@2"a2";"b2"_{}="u2"
      \ar@/_1pc/@2"a2";"b2"^{}="v2"
      \ar@/_1pc/@2"c2";"d2"^{}="w2"
      \ar@/^1pc/@2"c2";"d2"_{}="x2"
        \ar@3"v2";"u2"
        \ar@3"w2";"x2"
&
    \bullet
} \pbox{.}
\]
\shorthandon{;}%
Si $S = \Delta_{i_1}\wr\Delta_{i_2}\wr\Delta_{i_3}$ est un objet de
$\Xi_3$ écrit sous forme normale, alors les morphismes $S \to T$ dans
$\Xi_3$ sont de trois formes :
\begin{itemize}
\item les morphismes de la forme 
\[
{\varphi_1\wr\varphi_2\wr\varphi_3} : S \to T
\]
où les $\varphi_i$ sont des morphismes de $\Delta$ non constants pour
$i=1$ et $2$;
\item les morphismes de la forme 
\[
 \varphi_1\wr\varphi_2\wr- : S \to T
\]
où $\varphi_1$ est un morphisme non constant de $\Delta$ et
$\varphi_2$ est un morphisme constant de $\Delta$;
\item les morphismes de la forme 
\[
\varphi\wr -  : S \to T 
\]
où $\varphi$ est un morphisme constant de $\Delta$. 
\end{itemize}

On souligne que les objets suivants de $\Theta$ ne sont pas des objets de
$\Xi$ :
\[
\xymatrix@=4em{
\bullet
\ar[r]^{} 
&
\bullet 
\ar@/^1.5pc/[r]|{\vphantom{X}}="a"
\ar[r]|{\vphantom{X}}="b"
\ar@/_1.5pc/[r]|{\vphantom{X}}="c"
\ar@2"a";"b"
\ar@2"b";"c"
&
\bullet
} 
\mdvirg
\xymatrix@=5em{
\bullet 
\ar@/^1.5pc/[r]_{}="a"
\ar@/_1.5pc/[r]^{}="b"
\ar@/^1pc/@2"a";"b"_{}="u"
\ar@/_1pc/@2"a";"b"^{}="v"
\ar@3"v";"u"
&
\bullet
\ar@/^1pc/[r]|{\vphantom{X}}="a"
\ar@/_1pc/[r]|{\vphantom{X}}="b"
\ar@2"a";"b"
&
\bullet
} \pbox{,}
\]
et que le monomorphisme représenté par 
\[
\xymatrix@=5em{
  \bullet 
    \ar@/^2.5pc/[r]_{}="a"
    \ar[r]^{}="b"_{}="c"
      \ar@/^1pc/@2"a";"b"_{}="u"
      \ar@/_1pc/@2"a";"b"^{}="v"
        \ar@3"v";"u"
&
  \bullet 
    \ar[r]^{}="b1"_{}="c1"
    \ar@/_2.5pc/[r]^{}="d1"
      \ar@/_1pc/@2"c1";"d1"^{}="w1"
      \ar@/^1pc/@2"c1";"d1"_{}="x1"
        \ar@3"w1";"x1"
&
    \bullet
    \ar@/^2.5pc/[r]_{}="a2"
    \ar[r]^{}="b2"_{}="c2"
      \ar@/^1pc/@2"a2";"b2"_{}="u2"
      \ar@/_1pc/@2"a2";"b2"^{}="v2"
        \ar@3"v2";"u2"
&
    \bullet
} 
\]
entre les images de $\Delta_3\wr\Delta_1\wr\Delta_1$ et de $T = \Delta_3
\wr \Delta_2 \wr \Delta_1$ dans $\Theta$ n'est pas un morphisme de~$\Xi$.
\end{example}
\medskip

Dans ce qui suit, on va montrer que, pour tout entier $n\geq
0$, le foncteur $\xim_n$ est un foncteur asphérique. On va en fait
essentiellement adapter la preuve de Ara et Maltsiniotis
\cite[corollaire 5.5]{ara2022compnerfs} à cette variation du produit
en couronne.

\paragr Si $A$ est une petite catégorie, on dispose d'un triangle commutatif{
\begin{equation}\label{lemme:triangleXiFibGrothendieck}
\xymatrix{
\Delta\times A \ar[rr]^{\ximuf} \ar[rd]_{p_1} 
&& \Delta\xiwr  A \ar[ld]^{\pi} \\ & \Delta
}
\end{equation}
où $p_1$ est la première projection, et le foncteur $\pi$ est défini
par 
\begin{align*}
\ximu{\Delta_i}{a} &\mapsto \Delta_i \\
\ximu{\varphi}{f} &\mapsto \varphi \pbox{.}
\end{align*}
Pour tout objet $\ximu{\Delta_i}{ a}$ de $\Delta\xiwr  A$, on en déduit un
foncteur 
\[
P : \tranche{(\Delta\times A)}{\ximu{\Delta_i}{ a}}\to
\tranche{\Delta}{\Delta_i} 
\]
qui envoie un objet 
\[
((\Delta_{i'},a'),(\ximu{\Delta_{i'}}{a'} \xrightarrow{\ximu{\varphi}{\psi}}
\ximu{\Delta_i}{ a)}
\]
sur l'objet $(\Delta_{i'}, \Delta_{i'}\xrightarrow{\varphi}\Delta_i)$, et
un morphisme $(\psi, g)$ sur le morphisme $\psi$.

\begin{lemme}
Pour tout objet $\ximu{\Delta_i}{ a}$ de $\Delta\xiwr  A$, le foncteur
\[
P : \tranche{(\Delta\times A)}{\ximu{\Delta_i}{a}}\to
\tranche{\Delta}{\Delta_i} 
\]
est une fibration de Grothendieck.
\end{lemme}
\begin{proof}
C'est un cas particulier du lemme prouvé par Ara et Maltsiniotis
\cite[lemme 5.3]{ara2022compnerfs} appliqué au triangle
\ref{lemme:triangleXiFibGrothendieck}.
\end{proof}
On renvoie à \cite[paragraphe 1.1.15]{maltsiniotis2005} pour la
définition des fibrations de Grothendieck. L'élément important ici est
que les tranches d'une fibration de Grothendieck coïncident, à
équivalence de Thomason près, avec ses
fibres (voir~\cite[proposition 1.1.17]{maltsiniotis2005}).

\begin{prop}
Pour toute petite catégorie asphérique $A$, le foncteur 
\[
\ximuf_A : \Delta\times A \to \Delta\xiwr  A 
\]
est asphérique.
\end{prop}
\begin{proof}
Considérons un objet $\ximu{\Delta_i,a}$ de $\Delta\xiwr A$. Il s'agit de
montrer que la catégorie $\tranche{(\Delta\times A)}{\ximu{\Delta_i}{a}}$ est
asphérique. Pour cela, on va montrer que le foncteur
\[
P : \tranche{(\Delta\times A)}{\ximu{\Delta_i}{a}}\to \tranche{\Delta}{\Delta_i} 
\]
du lemme précédent est
une équivalence faible. En vertu de~\cite[corollaire~1.1.24]{maltsiniotis2005}, puisque $P$ est une fibration de
Grothendieck, il suffit de montrer que ses fibres sont asphériques. On
vérifie alors que la fibre en~$(\Delta_k,
\Delta_k\xrightarrow{u}\Delta_i)$ de $P$
est la catégorie $\tranche{A}{a}$ si $u$ est un morphisme non constant, et
$A$ si $u$ est constant. Dans tous les cas, la fibre de $P$ en $(\Delta_k,
u)$ est une catégorie asphérique, ce qui montre que~$P$ est une
équivalence faible, et donc que $\tranche{(\Delta\times
A)}{\ximu{\Delta_i}{a}}$ est une catégorie asphérique. 
\end{proof}

\begin{remark}
Contrairement au résultat analogue pour le produit en couronne usuel, on
n'a pas besoin de supposer que la catégorie $A$ est totalement
asphérique, puisqu'ici aucun produit de préfaisceaux représentables
n'intervient dans le calcul de la fibre du foncteur $P$.
\end{remark}

\begin{coro}
Pour tout entier $n\geq 1$, le foncteur 
\[
m_n' : \Delta^n \to \Xi_n 
\]
est un foncteur asphérique.
\end{coro}
\begin{proof}
Ce foncteur est obtenu comme composée des foncteurs asphériques
\[
\Delta^n 
= \Delta \times \Delta^{n-1} \xrightarrow{\id\times \xim_{n-1}}
\Delta\times\Xi_{n-1} \xrightarrow{\ximuf_{\Xi_{n-1}}}
\Delta\xiwr \Xi_{n-1} = \Xi_n \mdvirg
\]
et une récurrence montre que $\xim_n$ est bien asphérique.
\end{proof}

\paragr On dispose d'un triangle commutatif, pour tout entier $n\geq 1$, 
\[
\xymatrix{
{\Delta^n} \ar[rr]^{\xim_n} \ar[rd]_{m_n} && {\Xi_n} \ar[ld]^{j_n} \\
& {\Theta_n}
} 
\]
où les foncteurs $m_n$ et $\xim_n$ sont asphériques, et il en est donc de
même pour le foncteur $j_n:\Xi_n \to \Theta_n$ en vertu de
\cite[proposition 1.1.8]{maltsiniotis2005}. 
\begin{prop}
Pour tout entier $n\geq 1$, la catégorie $\Xi_n$ est une catégorie test
stricte.
\end{prop}
\begin{proof}
Puisque le foncteur $\xim_n$ est asphérique et que $\Delta^n$ est une
catégorie totalement asphérique, la proposition
\ref{totAspheriqueMorphismeAspherique} implique qu'il en est de même
pour $\Xi_n$. Puisque le foncteur $j_n$
est asphérique et que $\Theta_n$ est une catégorie test, la totale
asphéricité de $\Xi_n$ implique, grâce à la
proposition~\ref{propFonctAspheriqueSourceTotAsphTest}, que $\Xi_n$ est une
catégorie test stricte.
\end{proof}

\begin{coro}
Pour tout entier $n \geq 1$, la catégorie $\Xi_n$ est une catégorie test
homologique stricte de Whitehead.
\end{coro}
\begin{proof}
Puisque $\Xi_n$ est une catégorie test stricte et que le foncteur 
\[
\Delta \xrightarrow{\diag} \Delta^n \xrightarrow{\xim_n} \Xi_n 
\]
est asphérique, cela résulte de la proposition
\ref{propMorphismeAsphSourceTHS}.
\end{proof}

\medskip
\emph{Dans la suite, on fixe une petite catégorie $A$, et on suppose
que $A$ admet une orientation asphérique au sens de la section
\ref{secOrientation}.}
\medskip

Le but de ce qui suit est de donner une orientation sur
la catégorie~$\Xi_n$ pour tout entier $n\geq 1$, et de montrer que
cette orientation est asphérique.

\paragr On peut vérifier que $\Delta \xiwr A$ est une catégorie caténaire au
sens du paragraphe \ref{def:catenaire}. Si $T=\ximu{\Delta_i}{a}$, on a
\[
\dim T = \begin{cases}
i + \dim a \mdvirg & i \neq 0 \pbox{;} \\
0 \mdvirg & i = 0 \pbox{.}
\end{cases}
\]

On peut décrire les monomorphismes de codimension $1$ de $\Delta\xiwr
A$ de la manière suivante.  
\begin{itemize}
\item Si $T$ est un objet de $\Delta$, un monomorphisme
$S\hookrightarrow T$ de codimension $1$ correspond
exactement à un monomorphisme de codimension $1$ de but~$T$ dans
$\Delta$.

\item Si $T=\ximu{\Delta_1}{a}$ avec $\dim a\neq 0$, les monomorphismes
de codimension $1$ sont les morphismes 
\[
\ximu{\id}{ \psi} : \ximu{\Delta_1}{ a'} \to \ximu{\Delta_1}{ a}
\]
où $\psi : a' \to a$ est un monomorphisme de
codimension $1$ dans $A$. Néanmoins, pour $\delta
: \Delta_0 \to \Delta_1$, le morphisme $\ximu{\delta}{\id_a}$ n'est pas
de codimension~$1$ puisque $\ximu{\Delta_0}{a}$ est de dimension $0$.

\item Dans tous les autres cas, les monomorphismes de codimension $1$
correspondent aux monomorphismes 
\[
\ximu{\id_{\Delta_i}}{ \psi} : \ximu{\Delta_i}{ a'} \to \ximu{\Delta_i}{a}
\]
pour $\psi : a' \hookrightarrow a$ de codimension $1$, et aux
morphismes
\[
\ximu{\delta_k}{\id_{a}} : \ximu{\Delta_{i-1}}{a} \to
\ximu{\Delta_i}{a} \mdvirg 0 \leq k \leq i \pbox{.}
\]
\end{itemize}

\paragr \label{signesXi} On définit une préorientation $\mathcal{O}$ sur $\Delta \xiwr
A$ en utilisant les règles de Koszul : si $i$ est un entier non nul et
\[
  \delta_k : \Delta_{i-1} \to \Delta_{i}
\]
est une coface de $\Delta$, on pose, pour tout objet $a$ de $A$,
\[
\sg \ximu{\delta_k}{\id_{a}) = \sg(\delta_k) = (-1}^k \mdvirg 
0 \leq k \leq i \mdvirg
\]
et si $\varphi : a' \hookrightarrow a$ est un monomorphisme de
codimension $1$ de $A$, on pose 
\[
\sg\ximu{\id_{\Delta_i}}{\varphi) = (-1)^{i}\sg(\varphi} \pbox{.}
\]

\begin{example}
L'objet $\Delta_3 \wr \Delta_1$ de $\Xi_2$
\[
\xymatrix{
\bullet \ar@/^1pc/[r]_{}="a" \ar@/_1pc/[r]^{}="b" &
\bullet \ar@/^1pc/[r]_{}="c" \ar@/_1pc/[r]^{}="d" &
\bullet \ar@/^1pc/[r]_{}="e" \ar@/_1pc/[r]^{}="f" 
& \bullet
\ar@2"a";"b"
\ar@2"c";"d"
\ar@2"e";"f"
}
\]
est de dimension $4$ et possède six sous-objets de codimension $1$,
correspondant aux quatre morphismes 
\[
\Delta_2 \wr \Delta_1 
\xrightarrow{\delta_i\wr\id} 
\Delta_3\wr\Delta_1
\mdvirg 0 \leq i \leq 3 \mdvirg \]
de signes $(-1)^i$, et aux deux morphismes
\[
\Delta_3 \wr \Delta_0 
\xrightarrow{\id\wr \delta_j} \Delta_3\wr\Delta_1
\mdvirg 0 \leq j \leq 1 
\]
de signes $(-1)^3(-1)^j=(-1)^{j+1}$.
On rappelle que, par exemple, les sous-objets suivants de l'image de
$T$ dans $\Theta_2$ \emph{ne sont pas} des sous-objets de $T$ :
\[
\xymatrix{
  \bullet \ar@/^1pc/[r] &
  \bullet \ar@/_1pc/[r] &
  \bullet \ar@/^1pc/[r] &
  \bullet
  } \mdvirg
\xymatrix{
\bullet \ar[r] &
\bullet \ar@/^1pc/[r]_{}="c" \ar@/_1pc/[r]^{}="d" &
\bullet \ar@/^1pc/[r]_{}="e" \ar@/_1pc/[r]^{}="f" 
& \bullet
\ar@2"c";"d"
\ar@2"e";"f"
} \pbox{.}
\]
\end{example}

\paragr On va maintenant montrer que, sous des conditions raisonnables
sur $A$, la préorientation définie ci-dessus est bien une
orientation. On rappelle (voir \ref{def:orientation}) que cela
signifie qu'en notant, pour tout objet $T$ de $\Delta\xiwr A$, 
\[
L(T)_k = \bigoplus_{\dim S = k}\Z^{(S\to T)}
\]
et, pour tout morphisme $u : S\to T$ avec $\dim S=k$, 
\[
{d} \langle u \rangle =
\sum
    _{{\substack{S'\xhookrightarrow{\delta}S
    \\\mathclap{\codim\delta = 1}}}} 
    \sg(\varphi)\langle u\circ\delta\rangle 
\in L(T)_{k-1}\mdvirg
\]
on obtient alors un complexe de chaînes. Pour simplifier les notations, on va
étendre par linéarité le produit en couronnes en un foncteur 
\[
\ximuf : \Add(\Delta)\times\Add(A) \to \Add(\Delta\xiwr A) \pbox{.}
\]
De cette manière, si 
\[
u = \mu(\varphi,\psi) : \mu(\Delta_i,a') \to \mu(\Delta_j,a)
\]
est un morphisme de $\Delta\xiwr A$, alors on peut noter\footnote{On se permet également de sommer
directement les morphismes au lieu de sommer les générateurs du groupe
$L(T)_k$ , c'est-à-dire que si $u$ et $v$ sont deux morphismes de but
$T$ et de source des objets de dimension $k$, on note $u+v$
pour l'élément $\langle u \rangle + \langle v \rangle$ du groupe
$L(T)_k$.}
\[
\ximu{d\varphi}{\psi} =
\sum_{k=0}^i(-1)^k \langle\ximu{\varphi\circ\delta_k}{\psi} \rangle
\]
et 
\[
\ximu{\varphi}{d\psi}=
\sum
    _{\substack{a'' \xhookrightarrow{\delta} a'
    \\\mathclap{\codim \delta = 1}}}
    \sg(\delta)\left\langle\ximu{\varphi}{\psi\circ\delta}\right\rangle
\pbox{.}\]
Avec ces notations, et la description des monomorphismes de
codimension $1$ et de leur signe faite précédemment, on trouve
l'expression suivante pour la différentielle : 
\[
d\ximu{\varphi}{\psi} = 
\begin{cases}
\ximu{d\varphi}{\psi}                           \mdvirg & i\geq1 \text{
et } \dim a' = 0 \pbox{;}\\
- \ximu{\varphi}{d\psi}                         \mdvirg & i=1 \text{ et }
\dim a' \neq 0\pbox{;}\\
\ximu{d\varphi}{\psi) +(-1}^i \ximu{\varphi}{d\psi}\mdvirg & i>1 \text{ et }
\dim a'\neq 0 \pbox{.}
\end{cases} 
\]
\begin{lemme}\label{lemme:OrientationDeltaWrA}
On suppose que $A$ admet un objet final $e_A$ qui est son seul objet de
dimension $0$, et que tout objet de dimension $1$ de $A$ possède
exactement une face positive et une face négative. Alors la
préorientation sur~$\Delta \xiwr A$ définie au
paragraphe~\ref{signesXi} est une orientation.
\end{lemme}

\begin{proof}
Si $T$ est de la forme $\ximu{\Delta_i}{e_A}$, alors la différentielle coïncide avec
la différentielle du complexe non normalisé~$L_\Delta(\Delta_i)$. Sinon, $T$
est de la forme 
\[
T = \ximu{\Delta_j}{a} 
\]
avec $j>0$ et $\dim a > 0$. On fixe alors un générateur $u$ de $L(T)_k$
pour $k\geq 2$. Cela signifie qu'on a un morphisme
\[
u = \ximu{\varphi}{\psi} : \ximu{\Delta_i}{a'}\to\ximu{\Delta_j}{a}
\]
où $a'$ est un objet de $ A$, $i$ est un entier non nul, et $i
+ \dim a' \geq 2$. On distingue alors plusieurs cas
: \begin{itemize}
\item Si $\dim a' = 0$, alors on a 
\[
d^2\ximu{\varphi}{\psi} = \ximu{d^2\varphi}{\psi} =0 \pbox{.} 
\]
\item Si $i=1$ et $\dim a' =1$, alors on a 
\begin{align*}
d^2\ximu{\varphi}{\psi} 
&= - d\ximu{\varphi}{d\psi} \\
&= - \ximu{d\varphi}{d\psi} \\
&= -\Big( \ximu{\varphi\circ\delta_1}{ d\psi} -
\ximu{\varphi\circ\delta_0}{d\psi}\Big)
\end{align*}
et puisque $\varphi\circ\delta_k$ est un morphisme constant, on trouve
\[
d^2\ximu{\varphi}{\psi} = \Big(\sum_{e_A\xhookrightarrow{\delta}a'}
\sg(\delta)\Big)
(\varphi\circ \delta_0 - \varphi\circ\delta_1) \pbox{.}
\]
L'hypothèse sur $A$ implique que la somme des signes des
monomorphismes $e_A \to a'$ est nulle, et permet donc bien de conclure que
$d^2\mu'(\varphi,\psi)=0$.

\item Si $i=1$ et $\dim a' \geq 2$, on trouve 
\[
d^2\mu'(\varphi,\psi)= \mu'(\varphi,d^2\psi) = 0 \pbox{.}
\]
\item Si $i =2$ et $\dim a' = 1$, on a 
\begin{align*}
d^2\mu'(\varphi,\psi) &= d\mu'(d\varphi,\psi)+d\mu'(\varphi,d\psi)\\
&=-\mu'(d\varphi,d\psi) + \mu'(d\varphi,d\psi) = 0 \pbox{.}
\end{align*}

\item Si $i = 2$ et $\dim a' \geq 2$, on a 
\begin{align*}
d^2\mu'(\varphi,\psi) &= d\mu'(d\varphi,\psi)+d\mu'(\varphi,d\psi)\\
&= -\mu'(d\varphi,d\psi) + \mu'(d\varphi,d\psi) - \mu'(\varphi,d^2\psi)
=0 \pbox{.}
\end{align*}

\item Si $i>2$ et $\dim a' = 1$, on a 
\begin{align*}
d^2\mu'(\varphi,\psi) &= d\mu'(d\varphi,\psi) + (-1)^i
d\mu'(\varphi,d\psi) \\
&= 
\left(\mu'(d^2\varphi,\psi) +(-1)^{i-1}\mu'(d\varphi,d\psi)\right)
\\& \quad+
(-1)^i\mu'(d\varphi,d\psi) 
\\&=0 \pbox{.}
\end{align*}

\item Si $i>2$ et $\dim a' \geq 2$, on a 
\begin{align*}
d^2\mu'(\varphi,\psi) &= d\mu'(d\varphi,\psi) + (-1)^i
d\mu'(\varphi,d\psi) \\
&= 
\left(\mu'(d^2\varphi,\psi) +(-1)^{i-1}\mu'(d\varphi,d\psi)\right)
\\& \quad+
(-1)^i\left(\mu'(d\varphi,d\psi) + (-1)^{i}\mu'(\varphi,d^2\psi)\right)
\\&=0 \pbox{.}
\end{align*}
\end{itemize}
On a bien montré que $d$ définit une
différentielle, et les signes introduits précédemment définissent donc bien
une orientation de $\Delta\xiwr A$.
\end{proof}
\begin{remark}
On aurait pu supposer seulement que tout objet de dimension $1$ de $A$
possède autant de faces positives que de faces négatives, mais il
semble raisonnable de demander que tout objet de dimension $1$ se
comporte comme un segment.
\end{remark}

\emph{
Dans la suite de cette section, on suppose que la catégorie~$A$ admet
une orientation asphérique au sens du paragraphe
\ref{def:orientationAspherique}, et on note~$L_A : A\to \Ch(\Ab)$
l'intégrateur associé. On suppose de plus que $A$ possède un objet
final $e_A$ qui est le seul objet de dimension $0$, et que tout objet
de dimension $1$ possède exactement une face positive et une face
négative.
}

\medskip

\paragr Pour tout objet $T$ de $\Delta\xiwr  A$, on dispose d'un
morphisme canonique 
\[
r : L(T) \to \Z[0] 
\]
dans $\Ch(\Ab)$, envoyant les générateurs de dimension $0$ sur
$1\in\Z$. On construit alors une section de $r$ de la manière suivante :
si $\Delta_i$ est un objet de~$\Delta$, on note~$\sigma^i_0 :
\Delta_0\to\Delta_i$ le~$0$-simplexe minimal de $\Delta_i$,
c'est-à-dire l'injection~$\Delta_0 \to \Delta_i$ prenant la valeur
$0$. 
Si $T=\ximu{\Delta_{j}}{a}$ est un objet de $\Delta\xiwr A$, alors on note 
\[
\sigma_0^T : \Delta_0 \to T 
\]
le $0$-simplexe minimal de $T$, défini par le morphisme $\sigma_0^j :
\Delta_0 \to \Delta_j$ de $\Delta$. On définit de cette manière un
morphisme 
\[
s = s_T : \Z[0] \to L(T) 
\]
dans $\Ch(\Ab)$, envoyant le générateur de degré $0$ sur $\sigma_0^T$.
On a bien $rs=\id_{\Z[0]}$, et on va produire une homotopie du
morphisme $sr$ vers l'identité de $L(T)$, c'est-à-dire une famille de
morphismes de groupes abéliens
\[
h_k : L(T)_k \to L(T)_{k+1} \mdvirg k\geq 0
\]
vérifiant les relations 
\[
\begin{cases}
d_{k+1}h_k + h_{k-1}d_k = \id_{L(T)_k} \mdvirg k > 1 \\
d_1h_0 = \id_{L(T)_0} - sr \pbox{.}
\end{cases} 
\]

Si $T$ est un objet de $\Delta$, on a déjà vu l'existence d'une telle
rétraction par déformation dans l'exemple
\ref{ex:IntegrateurLibreDelta}. On rappelle que si $\varphi :
\Delta_i \to \Delta_j$ est un~$i$\nobreakdash-simplexe de $\Delta_j$ représenté
par 
\[
\varphi = \left\langle x_0 \to x_1 \to  \cdots \to x_i \right\rangle
\]
avec $0 \leq x_0 \leq x_1 \leq \cdots \leq x_i \leq j$, alors on a
posé
\[
h(\varphi) =  \left\langle 0 \to x_0 \to x_1 \to \cdots \to x_i
\right\rangle .
\]

\paragr\label{notationA} Puisque $L_A$ est un intégrateur sur $A$, on dispose d'un
quasi-isomorphisme $r_A : L_A \to \Z_{{A}^{\op}}$. Pour tout objet $a$ de
$A$, puisque $L_A(a)$ et
$\Z$ sont des complexes de groupes abéliens projectifs, le
morphisme $r_{A,a}$ est
une équivalence d'homotopie. On fixe alors un inverse à homotopie
près $s_a : \Z[0]\to L_A(a) $ de $r_{A,a}$ pour tout objet $a$ de $A$
(on précise que $s_a$ n'est pas naturel en $a$). Puisque~$L_A$ est un
intégrateur
libre, le morphisme
\[
  s_{a} : \Z[0] \to {L_A}(a)
\]
est déterminé par une somme de faces de dimension $0$ de $a$, que l'on
note
\[
x_a = (s_a)_0(1) \in {L_A}(a)_0 \pbox{.}
\]
On fixe aussi, pour tout objet $a$ de $A$, une homotopie $h_a :
L_A(a)_\bullet \to L_A(a)_{\bullet+1}$ de
$s_ar_{A,a}$ vers $\id_{L_A(a)}$.

\paragr Si $T=\ximu{\Delta_j}{a}$ est un objet de $\Delta\xiwr  A$ avec $a\neq
e_A$ et si $k$ est un entier positif, on définit alors le morphisme 
\[
h_k : L(T)_k \to L(T)_{k+1} 
\]
de la manière suivante. Pour tout morphisme 
\[
  \ximu{\varphi}{\psi} : \ximu{\Delta_i}{a'}\to\ximu{\Delta_j}{a}
\]
de $\Delta\xiwr  A$, on pose
\[
h\mu'(\varphi,\psi) = \begin{cases}
\mu'(h\varphi, x_a) \mdvirg
& i= 0 \pbox{;} \\
(-1)^i\mu'(\varphi,h\psi) \mdvirg
& i \neq 0 \text{ et } a'\neq e_A \pbox{;}\\
\mu'(h\varphi, x_a) + (-1)^i\mu'(\varphi, h\psi) \mdvirg
& i\neq 0 \text{ et } a' =  e_A \pbox{;}
\end{cases} 
\]
où on a noté $h(\varphi)$ l'image de $\varphi : \Delta_i \to \Delta_j$
par l'opérateur d'homotopie introduit dans l'exemple
\ref{ex:IntegrateurLibreDelta} (et que l'on vient de
rappeler plus haut), ainsi que $h(\psi)$ l'image de $\psi : a' \to
a$ par l'opérateur d'homotopie $h_a$ que l'on a fixé. On insiste sur
le fait que $h(\psi)$ est un élément du groupe abélien~$\bigoplus_{\dim a''=\dim a'+1}\Z^{(a'' \to a)}$.

Le morphisme $h$ est bien défini, puisque si $\varphi$ est un morphisme
constant, alors on a $\ximu{\varphi}{h\psi} = \ximu{\varphi}{h\psi'}$ pour
tous morphismes $\psi$ et $\psi'$ de $ A$.
\begin{example}
Dans le cas où $A=\Delta$ et $a = \Delta_k$, et en notant 
\[
T = \Delta_j\wr\Delta_k =
\xymatrix@=5em{
0
\ar@/^3pc/[r]^{f_1^0}|{\vphantom{X}}="a"
\ar@/^1pc/[r]|{\vphantom{X}}="b"
\ar@/_1pc/[r]|{\vphantom{X}}="d"
\ar@/_3pc/[r]_{f_1^k}|{\vphantom{X}}="e"
\ar@{}"b";"d"|<<<{\vdots}
\ar@2"a";"b"^{\alpha_1^1}
\ar@2"d";"e"^{\alpha_1^k}
&
1
\ar@/^3pc/[r]^{f_2^0}|{\vphantom{X}}="a"
\ar@/^1pc/[r]|{\vphantom{X}}="b"
\ar@/_1pc/[r]|{\vphantom{X}}="d"
\ar@/_3pc/[r]_{f_2^k}|{\vphantom{X}}="e"
\ar@{}"b";"d"|<<<{\vdots}
\ar@2"a";"b"^{\alpha_2^1}
\ar@2"d";"e"^{\alpha_2^k}
&
\cdots
\ar@/^3pc/[r]^{f_j^0}|{\vphantom{X}}="a"
\ar@/^1pc/[r]|{\vphantom{X}}="b"
\ar@/_1pc/[r]|{\vphantom{X}}="d"
\ar@/_3pc/[r]_{f_j^k}|{\vphantom{X}}="e"
\ar@{}"b";"d"|<<<{\vdots}
\ar@2"a";"b"^{\alpha_j^1}
\ar@2"d";"e"^{\alpha_j^k}
&
j
}
 \pbox{,}
\]
l'action de l'opérateur d'homotopie $h$ se représente de cette
manière. On va faire une rétraction par déformation de $T$ sur la
section
\[
  s_T : \Z[0]\to L(T)
\]
définie par le $0$-simplexe $\langle 0 \rangle$ de $T$. 

On doit d'abord choisit un opérateur d'homotopie sur la deuxième
coordonnée, et on utilise encore l'opérateur d'homotopie introduit dans
l'exemple~\ref{ex:IntegrateurLibreDelta}. Cela revient à faire un
rétracte par
déformation du complexe $L_\Delta(\Delta_k)$ sur la section définie
par le $0$-simplexe $\sigma_0^k$ de $\Delta_k$.
En d'autres termes, la $0$-face notée~$x_{\Delta_k}$ en
\ref{notationA} est la $0$-cellule minimale $\sigma_0^k : \Delta_0 \to
\Delta_k$ de $\Delta_k$. Les morphismes 
\[
  \ximu{\varphi}{\delta_0^k} : \Delta_1\xiwr\Delta_0 \to \Delta_j \xiwr
  \Delta_k
\]
correspondent alors aux flèches $f^0_k$ de $T$ pour~$1 \leq k \leq j$.

On peut alors décrire l'action de $h$ de la manière suivante.

Pour toute $0$-face~$x : \Delta_0 \to T$ de $T$, c'est-à-dire pour
tout entier $0 \leq x \leq j$, on a
\begin{align*}
h(x) = \ximu{x}{\sigma_0^k} = \left\langle \xymatrix@C=3em{
0 \ar@/^1pc/[r]^{f^0_xf^0_{x-1}\cdots f^0_1} & x
}\right\rangle
\end{align*}
si $x\neq0$, et $h(x)$ est l'identité de $0$ sinon. 

Pour tout morphisme $u : \Delta_i\wr\Delta_l \to T$ avec $l\neq
0$ et $i\neq 0$, représenté par 
\[
u = \left\langle
\xymatrix@=5em{
x_0
\ar@/^3pc/[r]^{f_{x_1}^{m_0}}|{\vphantom{X}}="a"
\ar@/^1pc/[r]|{\vphantom{X}}="b"
\ar@/_1pc/[r]|{\vphantom{X}}="d"
\ar@/_3pc/[r]_{f_{x_1}^{m_l}}|{\vphantom{X}}="e"
\ar@{}"b";"d"|<<<{\vdots}
\ar@2"a";"b"^{}
\ar@2"d";"e"^{}
&
x_1
\ar@/^3pc/[r]^{f_{x_2}^{m_0}}|{\vphantom{X}}="a"
\ar@/^1pc/[r]|{\vphantom{X}}="b"
\ar@/_1pc/[r]|{\vphantom{X}}="d"
\ar@/_3pc/[r]_{f_{x_2}^{m_l}}|{\vphantom{X}}="e"
\ar@{}"b";"d"|<<<{\vdots}
\ar@2"a";"b"^{}
\ar@2"d";"e"^{}
&
\cdots
\ar@/^3pc/[r]^{f_{x_i}^{m_0}}|{\vphantom{X}}="a"
\ar@/^1pc/[r]|{\vphantom{X}}="b"
\ar@/_1pc/[r]|{\vphantom{X}}="d"
\ar@/_3pc/[r]_{f_{x_i}^{m_l}}|{\vphantom{X}}="e"
\ar@{}"b";"d"|<<<{\vdots}
\ar@2"a";"b"^{}
\ar@2"d";"e"^{}
&
x_i
}\right\rangle 
\]
où $x_0 \leq x_1 \leq \cdots \leq x_i$ sont des éléments de
$\lbrace0,\dots,j\rbrace$ et $m_0\leq m_1 \leq \cdots \leq m_l$ sont
des flèches de $T$ (potentiellement des identités), on a 
\[ h(u) = 
(-1)^i\left\langle\xymatrix@=6em{
x_0
\ar@/^5pc/[r]^{f_{x_1}^{0}}|{\vphantom{X}}="0"
\ar@/^3pc/[r]_{}|{f_{x_1}^{m_0}}="a"
\ar@/^1pc/[r]|{\vphantom{X}}="b"
\ar@/_1pc/[r]|{\vphantom{X}}="d"
\ar@/_3pc/[r]_{f_{x_1}^{m_l}}|{\vphantom{X}}="e"
\ar@{}"b";"d"|<<<{\vdots}
\ar@2"0";"a"^{c_1^0}
\ar@2"a";"b"^{}
\ar@2"d";"e"^{}
&
x_1
\ar@/^5pc/[r]^{f_{x_2}^0}|{\vphantom{X}}="0"
\ar@/^3pc/[r]^{}|{f_{x_2}^{m_0}}="a"
\ar@2"0";"a"^{c_2^0}
\ar@/^1pc/[r]|{\vphantom{X}}="b"
\ar@/_1pc/[r]|{\vphantom{X}}="d"
\ar@/_3pc/[r]_{f_{x_2}^{m_l}}|{\vphantom{X}}="e"
\ar@{}"b";"d"|<<<{\vdots}
\ar@2"a";"b"^{}
\ar@2"d";"e"^{}
&
\cdots
\ar@/^5pc/[r]^{f_{x_i}^0}|{\vphantom{X}}="0"
\ar@/^3pc/[r]^{}|{f_{x_i}^{m_0}}="a"
\ar@2"0";"a"^{c_i^0}
\ar@/^1pc/[r]|{\vphantom{X}}="b"
\ar@/_1pc/[r]|{\vphantom{X}}="d"
\ar@/_3pc/[r]_{f_{x_i}^{m_l}}|{\vphantom{X}}="e"
\ar@{}"b";"d"|<<<{\vdots}
\ar@2"a";"b"^{}
\ar@2"d";"e"^{}
&
x_i
}\right\rangle \pbox{,}
\]
où on a noté $c^0_p$ la composée
$\alpha_{x_p}^{m_0}\alpha_{x_p}^{m_0-1}\cdots\alpha_{x_p}^{1}$ pour $
1\leq p\leq i$ avec la convention que si $m_0=0$, alors $c_p^0$ est la
$2$-cellule identité de $f^0_p$. 

Enfin, pour un morphisme $u :
\Delta_i\to T$ représenté par 
\[
u = \left\langle \xymatrix{
x_0 \ar[r]^{f_{x_1}^m} & \cdots \ar[r]^{f_{x_i}^m} & x_i
} \right\rangle
\]
avec $x_0 \leq x_1 \leq \cdots \leq x_i \in \lbrace0,\dots,j\rbrace$
et $0 \leq m \leq k$, on a 
\[
h(u)=
\left\langle\xymatrix@C=1em{
0 \ar@/^1em/[r]^{f} & x_0 \ar@/^1em/[r]^{f_{x_1}^0} & \cdots \ar@/^1em/[r]^{f_{x_i}^0} & x_i
}\right\rangle
+ (-1)^i
\left\langle\xymatrix{
x_0 \ar@/^1pc/[r]^{f_{x_1}^0}_{}="a" \ar@/_1pc/[r]_{f_{x_1}^m}^{}="a1" & 
\cdots 
    \ar@/^1pc/[r]^{f_{x_i}^0}_{}="c" \ar@/_1pc/[r]_{f_{x_i}^m}^{}="c1" & x_i
\ar@{=>}"a";"a1"^{c_1}
\ar@{=>}"c";"c1"^{c_i}
}\right\rangle
\]
où $f=f_{x_0}^0f_{x_{0}-1}^0\cdots f_{1}^0$ (ou l'identité si $x_0=0)$,
et $c_p$ est la composée~$\alpha_{x_p}^m\alpha_{x_p}^{m-1}\cdots\alpha_{x_p}^{1}$ pour $1 \leq
p \leq i$, ou l'identité de $f^0_{x_p}$
si $m=0$.
\end{example}

\begin{prop}
On suppose que la catégorie $A$ vérifie les hypothèses du lemme
\ref{lemme:OrientationDeltaWrA}. 
Alors $h$ est une homotopie du morphisme~$s_Tr$ vers l'identité
de $L(T)$.
\end{prop}
\begin{proof}
On fixe un morphisme
\[
u : \mu(\Delta_i,a') \xrightarrow{\mu(\varphi,\psi)} \mu(\Delta_j,a) 
\]
où $a$ et $a'$ sont des objets de $A$. On vérifie alors que
$h$ définit une homotopie de l'identité de $L(T)$ vers le morphisme
$sr$ : 
\begin{itemize}
\item si $i=0$, alors $\mu(\Delta_0,a') = \Delta_0$ et on a alors,
comme dans le cas de $\Delta$,
\begin{align*}
{d} h\ximu{\varphi}{\psi} &= d\ximu{h\varphi}{x_a} \\
&= \ximu{dh\varphi}{x_a}\\
&= \ximu{\varphi - sr(\varphi)}{x_a} \\
&= \ximu{\varphi - sr(\varphi)}{\psi} \\
&= (\id - sr)(u) \mdvirg
\end{align*}
où on a utilisé que puisque $\varphi$ et $sr(\varphi)$ sont deux
morphismes constants, on a
$\ximu{\varphi-sr(\varphi)}{f}=\ximu{\varphi-sr(\varphi)}{g}$ pour tous
morphismes $f$ et $g$ de~$A$.
\item si $i>1$ et $a' \neq e_A$, on a 
\begin{align*}
{d} h \ximu{\varphi}{\psi}
&= (-1)^i {d} \ximu{\varphi}{ h \psi}
\\
&= (-1)^i \ximu{ {d}\varphi}{ h\psi } 
+  \ximu{\varphi}{ {d} h\psi}\pbox{.}
\end{align*}
Par ailleurs, on a
\begin{align*}
h{d}\ximu{\varphi}{\psi} &= 
h\left(\ximu{d\varphi}{\psi} 
+ (-1)^i\ximu{\varphi}{d\psi)\right} \\
&= (-1)^{i-1}\ximu{{d}\varphi}{ h\psi} 
+ \ximu{\varphi}{h{d}\psi}
\end{align*}
tout en faisant attention au second terme dans le cas où $a'$ est de
dimension $1$ : on a dans ce cas
\begin{align*}
h\ximu{\varphi}{d\psi} 
&= h\ximu{\varphi}{ t\psi - s\psi}\\
&=\ximu{h\varphi}{ x_a -  x_a} +
(-1)^i\ximu{\varphi}{hd\psi} \\
&= (-1)^i\ximu{\varphi}{hd\psi} 
\end{align*}
en notant $t\psi$ (resp. $s\psi$) l'unique $0$-face positive (resp.
négative) de $\psi$. Puisqu'on a 
\[
hd \psi + dh \psi = \psi \mdvirg
\]
on trouve donc bien que $(dh+hd)\ximu{\varphi}{\psi} = \ximu{\varphi}{\psi}$.

\item Si $i=1$ et $a' \neq e_A$, on a cette fois 
\begin{align*}
{d} h\ximu{\varphi}{\psi} 
&= - {d} \ximu{\varphi}{ h\psi } \\
&= \ximu{\varphi}{dh\psi}
\end{align*}
et on a bien 
\begin{align*}
h{d}\ximu{\varphi}{\psi} 
&= -h\ximu{ \varphi}{ {d}\psi} \\
&= \ximu{\varphi}{ h{d}\psi}
\end{align*}
même si $a'$ est de dimension $1$, puisque dans ce cas on a 
\begin{align*}
h\ximu{ \varphi}{{d}\psi} 
&= h\ximu{\varphi}{t\psi - s\psi}\\
&= \ximu{h\varphi}{  x_a -  x_a} 
+ \ximu{\varphi}{ h{d} \psi} \pbox{.}
\end{align*}
On peut donc conclure, toujours en utilisant que
\[
  dh\psi + hd\psi = \psi \mdvirg
\]
que $(dh + hd)\ximu{\varphi}{\psi} = \ximu{\varphi}{\psi}$.

\item Si $i>1$ et $a'=e_A$, on obtient cette fois 
\begin{align*}
{d} h \ximu{ \varphi}{\psi}
&= {d} \Big(  \ximu{h\varphi}{ x_a}
+ (-1)^i\ximu{\varphi}{ h\psi) \Big}
\\
&= \ximu{dh\varphi}{ x_a}
+(-1)^i\ximu{d\varphi}{ h\psi}
+\ximu{\varphi}{ {d} h\psi}
\end{align*}
où on a utilisé que $s_{A,a}r_{A,a}(\psi)=x_a$. Par ailleurs, on a
  \begin{align*}
hd\ximu{\varphi}{\psi}
&=h\ximu{d\varphi}{\psi}\\
&= \ximu{hd\varphi}{ x_a} 
+ (-1)^{i-1}\ximu{d\varphi}{h\psi}
  \end{align*}
ce qui donne bien
\begin{align*}
  (dh+hd)\ximu{\varphi}{\psi}
  &= \ximu{dh\varphi+hd\varphi}{x_a}
  + \ximu{\varphi}{dh\psi} \\
  &= \ximu{\varphi}{x_a} + \ximu{\varphi}{\psi-x_a} \\
  &= \ximu{\varphi}{\psi} \pbox{.}
  \end{align*}

\item Enfin, si $i=1$ et $a'=e_A$, on trouve 
\begin{align*}
{d} h \ximu{\varphi }{ \psi} &=
{d} \ximu{h\varphi}{  x_a} 
- d\ximu{\varphi}{ h\psi} \\
&= \ximu{{d} h\varphi}{  x_a} 
+ \ximu{\varphi}{ {d} h \psi}
\\
&= \ximu{\varphi - h{d}\varphi}{  x_a} 
+ \ximu{\varphi}{ \psi - x_a} \mdvirg
\end{align*}
et on a
\begin{align*}
  h{d}\ximu{\varphi}{ \psi}
  &= h\ximu{d\varphi}{\psi} \\
  &= \ximu{h{d}\varphi}{ x_a}
  \end{align*}
ce qui donne bien à nouveau.
  \begin{align*}
(dh+hd)(\ximu{\varphi}{\psi)} 
&=\ximu{\varphi}{x_a} + \ximu{\varphi}{\psi-x_a}\\
&=\ximu{\varphi}{\psi} \pbox{.}
  \end{align*}
\end{itemize} 

On a bien montré que $h$ est une homotopie du morphisme $sr : L(T) \to
L(T)$ vers l'identité.
\end{proof}

\begin{prop}
Sous les hypothèses du lemme \ref{lemme:OrientationDeltaWrA},
l'orientation introduite au paragraphe \ref{signesXi} sur la
catégorie~$\Delta\xiwr  A$ est une orientation asphérique au sens du
paragraphe \ref{def:orientationAspherique}. 
\end{prop}
\begin{proof}
On a montré que pour tout objet $T$ de $\Delta\xiwr  A$, le morphisme
$r : L(T) \to \Z[0]$ admet un inverse à homotopie près, en notant
$L(T)$ le complexe défini par l'orientation de $\Delta \xiwr  A$.
\end{proof}

\begin{coro}
Pour tout entier $n\geq 1$, l'orientation sur $\Xi_n$ définie au
paragraphe \ref{signesXi} est une orientation asphérique.
\end{coro}

\begin{theorem}\label{integrateurSousObjetsXi_n}
Si $X$ est un préfaisceau en groupes abéliens sur $\Xi_n$, alors on a
un isomorphisme canonique naturel
\[
\H{\Xi_n}{X} \simeq \bigoplus_{\dim T = 0} X(T)
\leftarrow
\bigoplus_{\dim T=1} X(T) \leftarrow 
\bigoplus_{\dim T=2} X(T) \leftarrow \cdots
\]
\notindex{$L_{\Xi_n}$}%
dans $\Hotab$, où la différentielle est donnée pour $x \in X(T)$ par
\[
dx = \sum_{\substack{S
\xhookrightarrow{\delta}T\\\mathclap{\codim\delta=1}}} \sg(\delta)
X(\delta)(x) \pbox{.}
\]

\end{theorem}
\begin{proof}
On utilise simplement la description des extensions de Kan à gauche
des intégrateurs libres décrite en
\ref{propExpressionHomologieDerivateurLibre}.
\end{proof}

\begin{remark}
On peut aussi utiliser le foncteur $\xim_n : \Delta^n \to \Xi_n$, ainsi
que sa composée avec la diagonale $\xim_n\circ \delta : \Delta\to\Xi_n$
pour calculer l'homologie des préfaisceaux sur $\Xi_n$. Dans $\Xi_2$,
la différence entre les complexes~$L(X)$
et~$L_{\Delta\times\Delta}({\xim_n}^*(X))$ s'exprime de la manière
suivante : par définition, le complexe $L_{\Delta\times\Delta}({\xim_n}^*X)$ est
le complexe total du complexe double
\[
\xymatrix@R=1.5em{
X(\Delta_0) & \ar[l]_{0} X(\Delta_0) & \ar[l]_{\id} X(\Delta_0) &
\ar[l]_{0}
\cdots \\
X(\Delta_1\wr\Delta_0) \ar[u]^{} & \ar[l]^{}
X(\Delta_1\wr\Delta_1) \ar[u]^{0} & \ar[l]^{}
X(\Delta_1\wr\Delta_2) \ar[u]^{0} & \ar[l]^{} 
\cdots \\
X(\Delta_2\wr\Delta_0) \ar[u]^{}  & \ar[l]^{} 
X(\Delta_2\wr\Delta_1) \ar[u]^{} & \ar[l]^{} 
X(\Delta_2\wr\Delta_2) \ar[u]^{} & \ar[l]^{}  
\cdots \\
\vdots \ar[u]^{} & \vdots \ar[u]^{}  & \vdots \ar[u]^{} 
}
\]
alors qu'on peut voir le complexe $L_{\Xi_2}(X)$ induit par
l'orientation asphérique sur $\Xi_2$ comme le complexe
total du complexe double suivant : 
\[
\xymatrix@R=1.5em{
X(\Delta_0) & \ar[l]^{} 0 & \ar[l]^{} 0 & \ar[l]^{}
\cdots \\
X(\Delta_1\wr\Delta_0) \ar[u]^{} & \ar[l]^{}
X(\Delta_1\wr\Delta_1) \ar[u]^{} & \ar[l]^{}
X(\Delta_1\wr\Delta_2) \ar[u]^{} & \ar[l]^{} 
\cdots \\
X(\Delta_2\wr\Delta_0) \ar[u]^{}  & \ar[l]^{} 
X(\Delta_2\wr\Delta_1) \ar[u]^{} & \ar[l]^{} 
X(\Delta_2\wr\Delta_2) \ar[u]^{} & \ar[l]^{}  
\cdots \\
\vdots \ar[u]^{} & \vdots \ar[u]^{}  & \vdots \ar[u]^{} &\pbox{.}
} 
\] 
\end{remark}

\appendix
\let\paragr\paragrA
\let\prop\propA
\let\example\exampleA
\let\lemme\lemmeA
\let\remark\remarkA
\let\defin\definA
\let\coro\coroA
\chapter{Notions de base sur les
dérivateurs}\label{annexeDerivateurs}
Le langage des dérivateurs, introduit pour la première fois par
Grothendieck dans \emph{Pursuing Stacks}
\cite{pursuingstacks,maltsiniotisPS}, a pour but de dégager un cadre
permettant une étude fine des extensions de Kan homotopiques. En
particulier, elle permet de caractériser les morphismes
homotopiquement cofinaux, et la notion de morphisme asphérique en
homologie introduite au chapitre \ref{chapAsphericite} correspond en
fait à la notion de morphisme \emph{asphérique au sens du dérivateur}
$\DerHotab$. Nous allons expliquer ce que cela signifie dans cette
annexe. Pour un exposé détaillé, on pourra consulter
\cite{maltsiniotis2001derivateurs} sur lequel nous nous sommes basés.

Les propositions utilisant explicitement des énoncés de cette annexe
sont les suivantes : la proposition
\ref{morphismesAsphEnHomologieEquivalences} affirmant que les
foncteurs $A\to B$ induisant un isomorphisme $\H{A}{\Z}\to\H{B}{\Z}$
sont homotopiquement cofinaux, la proposition
\ref{localisateurFondamentalAbelien} affirmant que la classe
$\W_\infty^\ab$ forme un localisateur fondamental, ainsi que la
proposition \ref{integrateurInduitMorphismeAspheriqueTF}, explicitant,
pour un morphisme asphérique en homologie $u: A \to B$, un
quasi-isomorphisme entre les complexes calculant l'homologie de $A$
dans $u^*X$ et celle de $B$ dans $X$.

\paragrA Un \ndef{prédérivateur} est la donnée d'un
$2$-foncteur 
\[
\Der : \Cat^{\mathsf{o}} \to \CAT
\]
où on note $\Cat^{\mathsf{o}}$ la $2$-catégorie obtenue à partir de
$\Cat$ en inversant le sens des $1$-flèches et des $2$-flèches, et $\CAT$ la \og $2$-catégorie des catégories \fg{} (modulo
des questions ensemblistes). Cela signifie qu'un prédérivateur
consiste en la donnée
\begin{enumerate}
\item pour toute petite catégorie $A$, d'une catégorie $\Der(A)$;
\item pour tout foncteur $u:A\to B$ de $\Cat$, d'un foncteur
\[
  u^*=u^*_{\Der} : \Der(B) \to \Der(A) \pbox{;}
\]

\item pour toute transformation naturelle 
\[
\xymatrix@=3em{
A \ar@/^0.7pc/[r]^u_{}="a"
\ar@/_0.7pc/[r]_v^{}="b" & B
\ar@{=>}"a";"b"^{\alpha}
} 
\]
dans $\Cat$, d'une transformation naturelle 
\[
\xymatrix@=3em{
\Der(A) & \Der(B) \ar@/^0.8pc/[l]^{v^*}^{}="a" \ar@/_0.8pc/[l]_{u^*}^{}="b"
\ar@{=>}"a";"b"_{\alpha^*}
} \pbox{;}
\]
\end{enumerate}  
vérifiant les conditions de $2$-fonctorialité évidentes sur les
compositions, et les compositions horizontales et verticales des
$2$-flèches.
\begin{example}
Si $\M$ est une petite catégorie, on appelle \emph{prédérivateur représenté} par
$\M$ le prédérivateur défini par
\[
\Der = \Der_\M : A \mapsto \Homi({A}^{\op},\M) 
\]
et où pour tout morphisme $u : A \to B$, le morphisme $u^* : \Der(B) \to
\Der(A)$ est donné par précomposition avec ${u}^{\op} : {A}^{\op} \to
{B}^{\op}$.
\end{example}
\begin{example}
Un \emph{localisateur} est un couple $(\M,\W)$ où $\M$ est une
catégorie et $\W$ est une partie de $\Fl(\M)$. Si $(\M,\W)$ est un
localisateur, on
peut lui associer un prédérivateur de la manière suivante : si $A$ est
une petite catégorie, on note $\W_A$ la classe des morphismes de
foncteurs $\varphi : X\to Y$ tels que pour tout objet $a$ de $A$, le
morphisme $\varphi_a : Xa \to Ya$ soit dans $\W$. On note alors 
\[
\Der_{(\M,\W)} : A \mapsto \W_A^{-1}\Homi({A}^{\op}, \M)
\]
le prédérivateur qui associe à toute petite catégorie $A$ la catégorie
obtenue en localisant la catégorie des foncteurs ${A}^{\op}\to \M$ par
les morphismes de~$\W_A$. La catégorie $\Der(e)$ est alors
canoniquement isomorphe à la 
catégorie homotopique~$\W^{-1}\M$ associée au localisateur $(\M,\W)$.
\end{example}

\paragr\label{extKanptwise} En reprenant l'exemple du dérivateur représenté par une petite
catégorie $\M$, on sait que si $\M$ est une catégorie admettant des
limites inductives et projectives, alors pour tout foncteur $u:A \to B$
dans $\Cat$, le foncteur~$u^* : \Der(B) \to \Der(A)$ admet un adjoint
à droite et un adjoint à gauche 
\[
\xymatrix@=4em{
\Der(A) \ar@/_{1.3pc}/[r]_{u_*}="a" \ar@/^{1.3pc}/[r]^{u_!}="c" 
& \Der(B) \ar[l]|{u^*}="b"
\ar@{}"c";"b"|{\bot}
\ar@{}"b";"a"|{\bot}
} 
\]
qui sont les extensions de Kan à droite et à gauche le long de $u$. On
sait que ces extensions de Kan peuvent alors se calculer point par
point, ce qui signifie qu'on peut
calculer leurs fibres grâce aux limites inductives et projectives
suivantes : pour tout objet $b$ de $B$ et pour tout objet $X$ de
$\Der(A)$, on a des isomorphismes naturels
\[
u_! X (b)= \limind_{b\to ua} X(a) \mdvirg
u_* X (b)= \limproj_{ua\to b} X(a) \pbox{.}
\]
En particulier, en notant $p_A : A \to e$ l'unique morphisme vers la
catégorie ponctuelle,
le foncteur $p_A^* : \Der(e) \to \Der(A)$ correspond au
foncteur qui associe à tout objet $M$ de $\M$ le diagramme constant de
valeur $M$ sur ${A}^{\op}$. Les foncteurs adjoints à gauche et à
droite de $p_A^*$ correspondent donc aux foncteurs de limite inductive
et projectives. Dans le cas du prédérivateur associé à un localisateur,
les foncteurs ${p_A}_!$ et ${p_A}_*$, quand ils existent,
correspondent aux foncteurs de limite inductive et projective
homotopique. 

\begin{defin}
Soit $\Der$ un prédérivateur sur $\Cat$. On dit que $\Der$ est un
dérivateur s'il vérifie les conditions suivantes : 
\begin{enumerate}
\item[\textbf{Der 1}] \begin{enumerate}

\item Si $A$ et $B$ sont deux petites catégories, le morphisme
canonique
\[
\Der(A\textstyle\coprod B) \to \Der(A) \times \Der(B) 
\]
est une équivalence de catégories.

\item $\Der(\emptyset)$ est la catégorie ponctuelle, où $\emptyset$
désigne la catégorie vide.
\end{enumerate} 
\end{enumerate}

Si $A$ est une petite catégorie et $a$ un objet de $A$, on note
$i_{A,a}: e \to A$ le morphisme correspondant au choix de l'objet $a$.
On appellera foncteur \emph{fibre en~$a$} le foncteur $i_{A.a}^*$.
\begin{enumerate}[resume]
\item[\textbf{Der 2}] Si $\varphi$ est un morphisme de $\Der(A)$ tel que pour tout
objet $a$ de $A$, $i_{A,a}^* (\varphi)$ est un isomorphisme de
$\Der(e)$, alors $\varphi$ est un isomorphisme.

\item[\textbf{Der 3}] Pour tout morphisme $u:A\to B$, le morphisme
$u^*:\Der(B)\to\Der(A)$ admet un adjoint à gauche $u_!$ et un adjoint
à droite $u_*$ , appelés respectivement \emph{image directe homologique} et
\emph{image directe cohomologique}.
\end{enumerate}

Un dernier axiome demande que les images directes homologiques et
cohomologiques puissent se calculer point par point, généralisant les
formules données au paragraphe \ref{extKanptwise} dans le cas du
prédérivateur représenté par une catégorie. Dans le cas du
prédérivateur associé à un localisateur, cela revient à dire que les
fibres des extensions de Kan homotopiques sont des limites projectives
et inductives homotopiques indexées par les tranches. Pour un
dérivateur général, on peut formuler cette propriété grâce au principe
des \og mates \fg{} : étant donné un diagramme 
\[
\xymatrix{
A & B \ar[l]_{f} 
\ar@{}[ld]|(.3){}="a"|(.7){}="b" 
\ar@2"a";"b"^{\theta}
\\
C  \ar[u]^{u}  & D \ar[l]^{g} \ar[u]_{u'} 
} 
\]
dans $\Cat$, si $u$ et $u'$ admettent des adjoints à gauche $l$ et $l'$, on obtient
un nouveau diagramme 
\[
\xymatrix{
A 
\ar@{}[rd]|(.3){}="a"|(.7){}="b"
\ar@2"a";"b"^{\bar{\theta}} 
\ar[d]_{l} & B \ar[l]_{f} \ar[d]^{l'} \\
C & D  \ar[l]^{g} 
}
\]
où $\bar{\theta}$ s'exprime de la manière suivante : 
\[
lf  \xrightarrow{lf\star\eta'} lfu'l' \xrightarrow{l\star\theta\star
l'} lugl' \xrightarrow{\epsilon\star gl'} gl' \pbox{.}
\]

Considérons maintenant un morphisme $u:A\to B$ de $\Cat$ et $b$ un
objet de~$B$. On dispose alors d'un diagramme 
\[
\xymatrix{
\cotranche{b}{A} 
\ar[d]_{p_{\cotranche{b}{A}}}  
\ar[r]^{j_{u,b}} & A  \ar[d]^{u} \\
e 
\ar@{}[ru]|(.3){}="a"|(.7){}="b"
\ar@{=>}"a";"b"_{\alpha_{u,b}}
\ar[r]_{i_{B,b}} &  B
}
\]
en notant $j_{u,b}$ le foncteur d'oubli
\[
  \cotranche{b}{A}\to A \mdvirg (a, b \xrightarrow{\varphi} ua) \mapsto
  a
\]
et $i_{B,b}$ le foncteur correspondant au choix de l'objet $b$ de $B$. Le
morphisme de foncteurs $\alpha_{u,b}$ est donné pour tout objet
$(a,\varphi)$ de~$\cotranche{b}{A}$ par 
\[
i_{B,b}\circ p_{\cotranche{b}{A}}(a, \varphi)=b
\xrightarrow{\phantom{iiiii}\varphi\phantom{iiiii}} ua
= u \circ j_{u,b} (a, \varphi) \pbox{.}
\]
On obtient alors par fonctorialité un diagramme 
\[
\xymatrix{
\Der(\cotranche{b}{A}) & \Der(A) 
\ar@{}[ld]|(.3){}="a"|(.7){}="b"
\ar@{=>}"a";"b"_{\alpha_{u,b}^*}\ar[l]_{j_{u,b}^*} 
\\
\Der(e)  \ar[u]^{p_{\cotranche{b}{A}}^*} & \Der(B) \ar[l]^{i_{B,b}^*}
\ar[u]_{u^*} 
} 
\]
auquel on peut appliquer, si le prédérivateur $\Der$ vérifie l'axiome
\textbf{Der 3}, le procédé décrit ci-dessus pour obtenir finalement un
morphisme de foncteurs 
\[
\xymatrix{
\Der(\cotranche{b}{A}) 
\ar[d]_{{p_{\cotranche{b}{A}}}_!} 
\ar@{}[rd]|(.3){}="a"|(.7){}="b"
\ar@2"a";"b"^{c_{u,b}}
& 
\Der(A) \ar[l]_{j_{u,b}^*} \ar[d]^{u_!} 
\\
\Der(e)
&
\Der(B) \ar[l]^{i_{B,b}^*} 
} 
\]
dont le but correspond à la fibre en $b$ de l'image directe
homologique de~$u$ (donc, dans le cas du dérivateur associé à un
localisateur, de l'extension de Kan homotopique à gauche le long de $u$).

\begin{enumerate}
\item[\textbf{Der 4d}] Pour tout morphisme $u:A\to B$ de $\Cat$ et pour tout objet $b$
de $B$, le morphisme $c_{u,b}$ dans le diagramme ci-dessus est un isomorphisme.
\end{enumerate}

La lettre \og d \fg{} pour \og droite \fg{} indique qu'une autre
version de cet axiome, pour les images directes cohomologiques, est
requise.
Une construction duale à celle expliquée ici (voir
\cite[1.4.3]{maltsiniotis2005}) appliquée au carré 
\[
\xymatrix{
\tranche{A}{b} \ar[r]^{j'_{u,b}} \ar[d]_{p_{\tranche{A}{b}}} & A
\ar@{}[ld]|(.3){}="a"|(.7){}="b"
\ar@{=>}"a";"b"^{\alpha'}
\ar[d]^{u} \\
e \ar[r]_{i_{B,b}} & B
} 
\]
fournit un morphisme 
\[
c'_{u,b} : i^*_{B,b} \circ u_* \to {p_{\tranche{A}{b}}}_* \circ
j'^*_{u,b}
\]
faisant donc intervenir la fibre en $b$ de l'image directe cohomologique de~$u$.

\begin{enumerate}
\item[\textbf{Der 4g}] Pour tout morphisme $u :A\to B $ de $\Cat$ et tout objet $b$ de
$B$, le morphisme~$c'_{u,b}$ est un isomorphisme.
\end{enumerate}
\end{defin}

\begin{example}
Cisinski prouve dans \cite{cisinski2003imagesdirectes} que pour tout
localisateur~$(\M,\W)$ admettant une structure de catégorie de
modèles, le prédérivateur associé est un dérivateur. Dans le cas où on
dispose d'un localisateur~$(\M,\W)$ admettant une structure de
catégorie de modèles \emph{combinatoire}, on peut déjà vérifier rapidement
les axiomes \textbf{Der 1} à \textbf{Der 3} : en effet, dans ce cas, les
catégories de foncteurs de but $M$ peuvent être munies de la structure
de catégorie de modèles projectives qui assurent que la localisation a
bien un sens, et que pour tout foncteur $u :A\to B$ le foncteur $u^*$
est un foncteur de Quillen à droite, ce qui assure l'existence du
foncteur $u_!$. On peut aussi utiliser la structure de catégorie de
modèles injective pour laquelle le foncteur $u^*$ est un foncteur de
Quillen à gauche, ce qui cette fois assure l'existence du foncteur~$u_*$.
\end{example}

\begin{example}
Le dérivateur $\DerHot$ peut être représenté par le
localisateur~$(\Top,\W_{\Top})$, qui n'est évidemment ici pas notre
premier choix. On utilisera plutôt le localisateur
$(\Cat,\W_{\infty})$ introduit au paragraphe \ref{def:eqThomason}, ce
qui signifie que si $A$ est une petite
catégorie, on pose 
\[
\DerHot(A) = (\W_A)^{-1}\Homi_{\Cat}({A}^{\op}, \Cat) 
\]
en notant $\W_A$ la classe des morphismes de foncteurs qui sont des
équivalences de Thomason argument par argument. Pour voir qu'on
définit bien un dérivateur de cette manière, on peut utiliser
l'existence de la structure de catégorie de modèles de Thomason
\cite{thomason1980catModel} sur $\Cat$. Ce dernier prouve aussi dans
\cite{thomason1979homotopy} que les limites inductives homotopiques
des diagrammes dans $\Cat$ peuvent être calculées grâce à la
construction de Grothendieck. 
\end{example}

\begin{example}
Si $\mathcal{A}$ est une catégorie abélienne possédant assez d'objets
projectifs, alors le prédérivateur associé au localisateur
$(\Ch(\mathcal{A}), \W_{qis})$, où on note $\W_{qis}$ la classe des
quasi-isomorphismes de complexes, est un dérivateur. À nouveau, on
peut utiliser l'existence de la structure de catégories de modèles
projective sur $\Ch(\mathcal{A})$ introduite par Quillen dans
\cite[remarque 5, p 4.11]{quillenhomotalg}. En particulier, pour
$A=\Ab$, on dispose d'un dérivateur noté
\[
\DerHotab : A \mapsto (\W_{A})^{-1}\Homi({A}^{\op},\Ch(\Ab)) 
\]
où $\W_A$ désigne la classe des morphismes de foncteurs qui sont des
quasi-isomorphismes argument par argument.
\end{example}

\paragr Si $\Der$ est un dérivateur sur $\Cat$, et si $X$ est un objet
de $\Der(A)$, on appelle \emph{homologie} de $A$ à coefficients dans
$X$ relativement au dérivateur $\Der$ l'objet
\[
  \Hder{\Der}{A}{X}:={p_A}_!(X) 
\]
de $\Der(e)$, et \emph{cohomologie} de $A$ à coefficients dans $X$
relativement à $\Der$ l'objet
\[
  \coHder{\Der}{A}{X}:={p_A}_*(X)
\]
de $\Der(e)$, où on rappelle qu'on note $p_A$ l'unique foncteur $A\to e$. Si $M$
est un objet de $\Der(e)$, on notera simplement $\Hder{\Der}{A}{M}$ pour
l'homologie de $A$ à coefficients dans $p_A^*(M)$, c'est-à-dire pour
l'homologie à coefficients dans le diagramme constant de valeur $M$.
En d'autres termes, l'homologie au sens des dérivateurs correspond à
la notion de limite inductive homotopique, et la
cohomologie correspond à la notion de limite projective homotopique.
Dans le cas du dérivateur $\DerHotab$, c'est de cette manière que l'on
a défini l'homologie des préfaisceaux au paragraphe \ref{defH_A} : si
$X$ est un préfaisceau en groupes abéliens sur $A$, on a 
\[
\H{A}{X} = \Hder{\DerHotab}{A}{iX} 
\]
où $i:\Ab \to \Ch(\Ab)$ désigne le morphisme d'inclusion en degré $0$.

\medskip

Nous allons maintenant introduire la notion d'asphéricité au sens des
dérivateurs (ou de cofinalité homotopique), qui est la motivation
principale de cette annexe. On renvoie à la section
\ref{secMorphismesAsphEnHomologie} où on donne des caractérisations
plus précises des morphismes asphériques en homologie.

\paragr \label{defEqDerivateursEtMorphismesAspheriques} On fixe à présent un dérivateur $\Der$ sur $\Cat$. Si $u:A\to
B$ est un morphisme de $\Cat$, on a un triangle commutatif 
\[
\xymatrix{
A \ar[rr]^{u} \ar[rd]_{p_A} && B \ar[ld]^{p_B}  \\ & e
}
\]
qui induit, par fonctorialité, le triangle commutatif 
\[
\xymatrix{
{\Der(A)}  && {\Der(B)} \ar[ll]_{u^*} \\
& {\Der(e)}\ar[ru]_{p_B^*}\ar[lu]^{p_A^*} & \pbox{.}
} 
\]
À son tour, ce dernier induit, par adjonction, les deux triangles
commutatifs à isomorphisme près
\[\label{diag:mateDerivateurs}
\tag{A.10.1}
\xymatrix{
{\Der(A)} \ar[rr]^{u_!} \ar[rd]_{{p_A}_!} && {\Der(B)}
\ar[ld]^{{p_B}_!} \\
& {\Der(e)}
} 
\mdvirg
\xymatrix{
{\Der(A)} \ar[rr]^{u_*} \ar[rd]_{{p_A}_*} && {\Der(B)}
\ar[ld]^{{p_B}_*} \\
& {\Der(e)} & \pbox{,}
} 
\]
et on peut alors définir une transformation naturelle 
\[
{p_A}_!u^* = {p_B}_!u_!u^* \xrightarrow{{p_B}_!\star\epsilon} {p_B}_!
\]
où on a noté $\epsilon : u_!u^*\to\id_{\Der(B)}$ le morphisme d'adjonction. En
particulier, pour tout objet $X$ de $\Der(B)$, on dispose d'un
morphisme 
\[
\Hder{\Der}{A}{u^*(X)} \xrightarrow{{p_B}_!(\epsilon_X)} \Hder{\Der}{B}{X} \pbox{.}
\]
On dira que $u$ est \emph{$\Der$-asphérique} si ce morphisme est
un isomorphisme pour tout objet $X$ de $\Der(B)$.

Dans le cas particulier où $X$ est un coefficient constant,
c'est-à-dire s'il est de la forme $p_B^*(M)$ pour un objet $M$ de
$\Der(e)$, on obtient donc un morphisme
\[
  \Hder{\Der}{A}{M}\to\Hder{\Der}{B}{M}  
\]
défini par la transformation naturelle
\[
{p_A}_!p_A^* = {p_B}_!u_!u^*p_B^*
\xrightarrow{{p_B}_!\star\epsilon\star p_B^*} {p_B}_!{p_B}^* 
\]
et on dit que $u$ est une \ndef[équivalence faible!de
$\Cat$!relativement à un dérivateur]{$\Der$-équivalence} si ce morphisme
est un isomorphisme pour tout objet $m$ de $\Der(e)$. Enfin, on dit
qu'une \emph{petite catégorie} $A$ est
\emph{$\Der$\nobreakdash-asphérique} si le
morphisme $p_A : A\to e$ est une $\Der$-équivalence. Cela signifie
donc que le morphisme d'adjonction 
\[
{p_A}_!p_A^* \to \id_{\Der(e)}
\]
est un isomorphisme argument par argument, autrement dit que pour tout
coefficient constant $M$ dans $\Der(e)$, la counité de l'adjonction
${p_A}_!\dashv {p_A}^*$ induit un isomorphisme $\Hder{\Der}{A}{M} \simeq M$.

\paragr\label{appendix:DerHotab} En particulier, si $A$ est une petite
catégorie, on a, sans le dire, défini au paragraphe \ref{defH_A}
l'homologie de $A$ à coefficients dans un préfaisceau en groupes
abélien $X$ en posant
\[
\H{A}{X} = \Hder{\DerHotab}{A}{iX} 
\]
en notant $i$ le morphisme $\prefab{A}\to\DerHotab(A)$ obtenu en
voyant un préfaisceau en groupes abéliens comme un préfaisceau en
complexes de chaînes concentré en degré $0$. 
Cela implique que si $u : A\to B$ est un morphisme
$\DerHotab$-asphérique, on dispose alors d'un diagramme commutatif à
isomorphisme près
\[
\xymatrix{
\prefab{B} \ar[rr]^{u^*} \ar[d]_{\tilde{i}} && \prefab{A}
\ar[d]^{\tilde{i}} \\
\DerHotab(B) \ar[rr]^{u^*} \ar[rd]_{{p_B}_!} 
&& \DerHotab(A) \ar[ld]^{{p_A}_!} \\ 
& \DerHotab(e) \simeq \Hotab
} 
\]
et donc que pour tout morphisme $f : X \to Y
\in\Fl(\prefab{B})$, on a l'équivalence 
\[
f \in\Wab_B \iff u^*(f)\in\Wab_A \mdvirg
\]
où on rappelle (\ref{defWab_A}) que $\Wab_A$ désigne la classe des
morphismes de $\prefab{A}$ envoyés sur des isomorphismes par le
foncteur $\Hf{A}$.

\paragr\label{annexe:categoriesAspheriquesIntegrateurs} De plus, on a vu au paragraphe \ref{colimiteHomotopiqueDegreparDegre} que
si $L$ est un intégrateur sur une petite catégorie $A$, alors le
foncteur 
\[
\underline{L}_! : \Homi({A}^{\op},\Ch(\Ab)) \to \Ch(\Ab) 
\]
calcule la limite inductive homotopique, c'est-à-dire qu'on a un
diagramme commutatif à isomorphisme près
\[
\xymatrix{
\Homi({A}^{\op},\Ch(\Ab)) \ar[r]^-{\underline{L}_!} \ar[d]^{} 
& \Ch(\Ab) \ar[d]^{} \\
\DerHotab(A) \ar[r]_{{p_A}_!} & \Hotab
}
\]
où les flèches verticales désignent les projections canoniques. Le
foncteur $\underline{L}_!$ admet un foncteur adjoint à droite 
\[
\underline{L}^* :  C \mapsto \big( a \mapsto \Homi_{\Ch(\Ab)}(L(a),C)
\big)
\]
dont le foncteur induit entre les catégories dérivées coïncide donc à
isomorphisme près avec le foncteur dérivé à droite du foncteur
diagramme constant~$p_A^* : \DerHotab(e) \to \DerHotab(A)$.
Ainsi, une petite catégorie $A$ est asphérique pour le localisateur
$\DerHotab$ si et seulement si le morphisme d'adjonction 
\[
\underline{L}_!\underline{L}^*(C) \to C 
\]
est un quasi-isomorphisme pour tout complexe de chaînes $C$.

\paragr De même, dire qu'un foncteur $u : A \to B$ est $\DerHotab$-asphérique
se traduit de la manière suivante. Le foncteur $u$ induit un couple de
foncteurs adjoints 
\begin{align*}
\underline{u}_!^\ab : 
\Homi({A}^{\op},\Ch(\Ab))&\to\Homi({B}^{\op},\Ch(\Ab)) \\
u^* : \Homi({B}^{\op},\Ch(\Ab)) &\to \Homi({A}^{\op},\Ch(\Ab)) 
\end{align*}
où $\underline{u}_!^\ab$ est l'extension de Kan à gauche le long de
${u}^{\op}$, que l'on a étudiée dans le chapitre
\ref{chapHomologiePrefaisceaux}. On sait que le foncteur $u^*$
préserve les équivalences faibles, ce qui n'est pas le cas en général
de son adjoint à gauche. En revanche, ce dernier se dérive à gauche en
un foncteur 
\[
L\underline{u}_!^\ab : \DerHotab(A) \to \DerHotab(B) 
\]
et le diagramme \ref{diag:mateDerivateurs} fournit un triangle commutatif à isomorphisme près
\[
\tag{A.13.1}
\label{diagLu_!}
\xymatrix{
{\DerHotab(A)} \ar[rr]^{L\underline{u}_!^\ab} \ar[rd]_{{p_A}_!} &&
{\DerHotab(B)} \ar[ld]^{{p_B}_!} \\
& {\Hotab} & \pbox{.}
} 
\]
Le foncteur $L\underline{u}_!^\ab$ se décrit de la manière suivante :
si $F_\bullet$ est un préfaisceau en complexes de chaînes sur $A$, on choisit
une résolution $P_\bullet \to F_\bullet$ par des objets projectifs de
$\prefab{A}$. On a alors un isomorphisme canonique 
\[
L\underline{u}_!^\ab (F_\bullet) \simeq \underline{u}_!^\ab(P_\bullet) 
\]
dans $\DerHotab(B)$.

\paragr\label{appendix:morphismeHotabAspheriqueCounite} Ainsi, si $L_B$ est un intégrateur sur $B$, alors le morphisme $u$ est
asphérique en homologie si et seulement si pour tout préfaisceau en
complexes de chaînes $X_\bullet$ sur $B$
et pour toute résolution projective~$ p: P_\bullet \to
u^*(X_\bullet)$ dans~$\Homi({A}^{\op},\Ch(\Ab))$, la
flèche oblique dans le diagramme 
\[
\xymatrix@C=4em{
\underline{L_B}_!\underline{u}_!^\ab P_\bullet
\ar[rd]^{} 
\ar[r]^{{\underline{L_B}}_!\underline{u}_!^\ab(p)} 
& 
{\underline{L}_B}_!\underline{u}_!^\ab u^*(X) \ar[d]^{{L_B}_!(\epsilon)} 
\\
& {\underline{L}_B}_!(X)
}
\]
est un isomorphisme dans $\Hotab$, où on a noté 
\[
\epsilon : \underline{u}_!^\ab u^* \to \id_{\Homi({B}^{\op},\Ch(\Ab))}
\]
le morphisme d'adjonction.

\begin{prop}[Grothendieck]\label{prop:asphericiteDerivateurs}
Pour tout dérivateur $\Der$ et pour tout morphisme $u : A\to B$ de
$\Cat$, les conditions suivantes
sont équivalentes : \begin{enumerate}
\item $u$ est un morphisme $\Der$-asphérique;
\item pour tout objet $b$ de $B$, le morphisme
\[
  \tranche{A}{b} \xrightarrow{\tranche{u}{b}} \tranche{B}{b} \mdvirg
(a, ua \xrightarrow{\varphi} b) \mapsto (ua, \varphi)
\]
est une $\Der$-équivalence; 
\item pour tout objet $b$ de $B$, la catégorie $\tranche{A}{b}$ est
une catégorie $\Der$-asphérique.
\end{enumerate}
\end{prop}
\begin{proof}
Voir \cite[1.5.11]{maltsiniotis2001derivateurs}.
\end{proof}

\paragr \label{colimiteHomotopiqueDiagrammeConstantGenMod}
Pour illustrer ce qu'affirme cette proposition, considérons le
dérivateur $\Der_\M$ représenté par une petite catégorie $\M$ complète
et cocomplète. Un morphisme asphérique $u:A\to B$ est alors un
foncteur \emph{cofinal}, c'est-à-dire un foncteur tel que pour
diagramme $X : {B}^{\op} \to \M $, on a un isomorphisme 
\[
\limind\nolimits_{{B}^{\op}} X \simeq \limind\nolimits_{{A}^{\op}} u^*
(X) \pbox{.}
\]
La proposition ci-dessus affirme alors qu'un foncteur $u:A\to B$ est
cofinal si et seulement si pour tout objet $b$ de $B$, le morphisme
$p_{\tranche{A}{b}} : \tranche{A}{b}\to e$ est une
$\Der_\M$-équivalence, et il reste donc à déterminer la classe des
$\Der_\M$-équivalences. 

Par définition, il s'agit des foncteurs induisant un isomorphisme
entre les limites inductives des diagrammes constants. Puisque $\M$
est une catégorie cocomplète, on dispose d'un \emph{tenseur} 
\[
\otimes : \Ens \times \M \to \M \mdvirg (E, M) \mapsto \coprod_{E}
M \mdvirg
\]
ce qui signifie que pour tout ensemble $E$ et pour tous objets $M_1,
M_2$ de $\M$, on a un isomorphisme naturel 
\[
\Hom_\M(E\otimes M_1, M_2) \simeq \Hom_\Ens(E, \Hom_\M(M_1,M_2) \pbox{.}
\]
On montre alors facilement que, pour tout objet $M$ de $\M =\Der_\M(e)$ et
toute petite catégorie $A$, on a un isomorphisme naturel
\[
{p_A}_!p_A^*(M)
\simeq \pi_0(A)\otimes M
\]
dans $\Der_\M(e)$, c'est-à-dire que la limite inductive du diagramme
constant de valeur $M$ sur ${A}^{\op}$ est le coproduit
$\coprod_{\pi_0(A)}m$ dans $\M$. 

En particulier, cela signifie que les $\Der_\M$-équivalences
coïncident avec les
morphismes~$u: A\to B$ de $\Cat$ induisant un isomorphisme
$\pi_0(A)\simeq \pi_0(B)$, et donc qu'un foncteur $u : A\to B$ est
$\Der_\M$-asphérique (ou cofinal) si et seulement si pour tout objet
$b$ de $B$, la catégorie $\tranche{A}{b}$ est connexe, ce qui est un
résultat bien connu. La proposition suivante correspond alors à la version
homotopique de ce résultat.

\begin{prop}[Grothendieck]
La classe des~$\DerHot$\nobreakdash-équivalences
coïncide avec la classe $\W_{\infty}$ des équivalences de Thomason.
\end{prop}
\begin{proof}
On rappelle (\ref{def:eqThomason}) qu'on voit le dérivateur $\DerHot$ comme induit par le
localisateur $(\Cat, \W_{\infty})$. Si $A$ est
une petite catégorie, alors la limite inductive homotopique des
foncteurs ${A}^{\op}\to\Cat$ peut être calculée grâce à la
construction de Grothendieck (\cite{thomason1979homotopy}, voir aussi
\cite[3.3.18]{cisinskipref}). 

On voit alors très simplement que si $C$ est un objet de
$\Hot=\DerHot(e)$, la catégorie associée à~$p_A^*(C)$ par la
construction de Grothendieck est le produit de catégories $A\times C$. Un
morphisme $u: A\to B$ de $\Cat$ est donc une $\DerHot$\nobreakdash-équivalence si
et seulement si pour tout objet $C$ de $\Hot$, le morphisme~$A\times C \to
B\times C$ est dans~$\W_{\infty}$, ce qui revient à dire que $u$ est
une équivalence de Thomason.
\end{proof}

\begin{remark}
En particulier, on voit que si $u \in\Fl(\Cat)$ est un morphisme de
$\W_{\infty}$, alors c'est une $\Der_\M$-équivalence pour toute
catégorie complète et cocomplète $\M$. De plus, on peut montrer (voir
par exemple \cite[corollaire~1.5.14]{maltsiniotis2001derivateurs}) que
pour tout dérivateur $\Der$, la classe des $\W_{\Der}$-équivalences
forme un localisateur fondamental (\ref{defLocFondamental}). Ainsi, la
minimalité de $\W_\infty$ parmi les localisateurs fondamentaux
(\ref{locFondamentalMinimalThm}) implique que ce phénomène reste vrai
pour tout dérivateur sur~$\Cat$. 
\end{remark}

On va maintenant montrer que la classe des $\DerHotab$-équivalences
coïncide avec la classe $\W_\infty^\ab$ des équivalences faibles
abéliennes de $\Cat$ introduite dans la section 
\ref{secMorphismesAsphEnHomologie}.

\begin{lemme}
Si $C$ est un complexe de chaînes de groupes abéliens, et si~$A$ est une petite catégorie, 
alors on a un isomorphisme naturel
\[
\hocolim_{{A}^{\op}}^{\Hotab} p_A^*(C) \simeq 
\H{A}{\Z} \Lotimes C\pbox{,}
\]
dans $\Hotab$, où $p_A^*(C)$ désigne le diagramme constant sur ${A}^{\op}$ de valeur
$C$.
\end{lemme}
\begin{proof}
On sait (voir \cite[théorème 1.2]{thomason1979homotopy}) que si $A$
est une petite catégorie, on a un isomorphisme naturel en $A$
\[
\hocolim_{{A}^{\op}}^{\Hot_\Delta}p_A^*(\Delta_0) \simeq \nerf A \pbox{.}
\]
De plus, on rappelle qu'on peut munir la catégorie $\pref{\Delta}$ ainsi que la
catégorie $\Ch(\Ab)$ de structures de catégories de modèles pour
lesquelles le foncteur complexe normalisé $\dk : \pref{\Delta}\to\Ch(\Ab)$ est un
foncteur de Quillen à gauche, et le produit tensoriel est un
bifoncteur de Quillen à gauche. Ainsi, si $C$ est un complexe de
chaînes, on peut conclure par la chaîne d'isomorphismes naturels
suivante :  
\begin{align*}
\hocolim_{{A}^{\op}}^{\Hotab} p_A^*(C) &\simeq 
\hocolim_{{A}^{\op}}^{\Hotab} p_A^*(\dk\Delta_0 \Lotimes C) \\
&\simeq \hocolim_{{A}^{\op}}^{\Hotab}(\dk\Delta_0)\Lotimes C\\
& \simeq \dk \nerf A \Lotimes C \pbox{.} 
\end{align*}
On peut alors conclure en rappelant (voir
\ref{homologieZ=HomologieSinguliere}) que pour toute petite catégorie
$A$, on a un isomorphisme naturel $\H{A}{\Z} \simeq \Hsing(A)$.
\end{proof}

\begin{coro}
On a l'égalité $\W_\infty^\ab = \W_{\DerHotab}$. En particulier, un
foncteur $u : A \to B$ est asphérique en homologie
(\ref{defAspheriqueEnHomologie}) si et seulement si $u$ est un
morphisme $\DerHotab$-asphérique.
\end{coro}

\begin{coro}\label{annexeWabLocFond}
La classe $\W_\infty^\ab$ est un localisateur fondamental. 
\end{coro}
\begin{proof}
On utilise le fait que pour tout dérivateur $\Der$ sur $\Cat$, la
classe des $\Der$-équivalences forme un localisateur fondamental (voir
par exemple~\cite[corollaire 1.5.14]{maltsiniotis2001derivateurs}).
\end{proof}

\def\bibname{Références}
\begin{flushleft}
\printindex[not]
\printindex[term]
\end{flushleft}
\bibliographystyle{smfplain}
\bibliography{bibthese}

\end{document}